\documentclass[12pt]{article}
\usepackage{times}
\usepackage[T1]{fontenc}
\usepackage{amsthm}
\usepackage{verbatim}
\usepackage{graphicx}   
\usepackage{amsmath}    
\usepackage{amsfonts}    
\usepackage{mathrsfs}  
\usepackage{pgf,tikz}
\usepackage{amssymb,amsthm,epsfig,epstopdf,titling,url,array}
\usepackage{mathrsfs}
\usepackage[notref, notcite]{showkeys}
\usepackage{subfigure} 
\usepackage{epstopdf}
\usepackage[numbers]{natbib}
\usetikzlibrary{plotmarks}

\graphicspath{{../}}
\usepackage[colorlinks=TRUE,linkcolor=blue,citecolor=blue]{hyperref}
\theoremstyle{plain}

\newtheorem{cla}{Claim}[section]
\theoremstyle{definition}
\newtheorem{defn}{Definition}[section]

\newtheorem{rem}[defn]{Remark}

\newcommand{\um}{\mathcal{U}}

\newcommand{\wm}{\mathcal{W}}
\newcommand{\Hm}{\mathcal{H}}

\providecommand{\testright}[1]{\shortmid\!\xrightarrow{#1}}

\usepackage{graphics,color}
\definecolor{purpv}{rgb}{0.62, 0.0, 1}
\newcommand{\red}[1]{\textcolor{red}{#1}}

\newcommand{\todo}[1]{[*** TO DO: #1 ]}

\newcommand{\sw}{s_{\mathrm{w}}}
\newcommand{\sg}{s_{\mathrm{g}}}
\newcommand{\so}{s_{\mathrm{o}}}

\newcommand{\fw}{f_{\mathrm{w}}}
\newcommand{\fg}{f_{\mathrm{g}}}
\newcommand{\fo}{f_{\mathrm{o}}}
\newcommand{\muw}{\mu_{\mathrm{w}}}
\newcommand{\mug}{\mu_{\mathrm{g}}}
\newcommand{\muo}{\mu_{\mathrm{o}}}
\newcommand{\mut}{\mu_{\mathrm{tot}}}

\newcommand{\lambdas}{\lambda_{\,\mathrm{s}}}
\newcommand{\lambdaf}{\lambda_{\,\mathrm{f}}}

\newcommand{\BO}{$B$\,-\,$O\;$}
\newcommand{\BOc}{$B$\,-\,$O$}

\newcommand{\OBc}{$O$\,-\,$B$}

\newcommand{\WO}{$W$\,-\,$O\;$}
\newcommand{\WOc}{$W$\,-\,$O$}

\newcommand{\WUc}{$W$\,-\,$\um$}
\newcommand{\UW}{$\um$\,-\,$W\;$}
\newcommand{\UWc}{$\um$\,-\,$W$}

\newcommand{\GUc}{$G$\,-\,$\um$}
\newcommand{\OU}{$O$\,-\,$\um\;$}

\newcommand{\BUc}{$B$\,-\,$\um$}
\newcommand{\EU}{$E$\,-\,$\um\;$}
\newcommand{\EUc}{$E$\,-\,$\um$}
\newcommand{\DU}{$D$\,-\,$\um\;$}
\newcommand{\DUc}{$D$\,-\,$\um$}

\newcommand{\GO}{$G$\,-\,$O\;$}
\newcommand{\GOc}{$G$\,-\,$O$}

\newcommand{\GW}{$G$\,-\,$W\;$}
\newcommand{\GWc}{$G$\,-\,$W$}

\newcommand{\EW}{$E$\,-\,$W\;$}
\newcommand{\WE}{$W$\,-\,$E$\;}
\newcommand{\EWc}{$E$\,-\,$W$}
\newcommand{\WEc}{$W$\,-\,$E$}

\newcommand{\GD}{$G$\,-\,$D\;$}
\newcommand{\DGc}{$D$\,-\,$G$}
\newcommand{\GDc}{$G$\,-\,$D$}

\newcommand{\m}{\mathcal}

\providecommand{\testright}[1]{\shortmid\!\xrightarrow{#1}}
\providecommand{\CS}{'\! S}
\title{Displacement of three-phase flow for Heavy Oil: Riemann Solutions}
\author{
Lozano, Luis F.\\
\texttt{luisfer99@gmail.com}
\and
Furtado, Frederico\\
\texttt{furtado@uwyo.edu}
\and
De Souza, Aparecido\\
\texttt{aparecidosouza@ci.ufpb.br}
\and
Marchesin, Dan\\
\texttt{marchesi@impa.br}
}
\date{}

\begin{document}

\maketitle

\begin{abstract}
This work presents the Riemann solution for three-phase flow in porous media under the condition that oil viscosity exceeds that of water and gas. We classify all Riemann solution problems for scenarios where the left states $L$ lie along the edge $G$-$W$, and the right states $R$ span nearly the entire saturation triangle, excluding small regions near the boundaries $G$-$O$ and $W$-$O$. We use the wave curve method to determine the Riemann solution for initial and injection data within the above-mentioned class. This study extends previous analytical solutions, which were limited to right states near the corner $O$ or within the quadrilateral $O$-$E$-$\mathcal{U}$-$D$. Notably, this classification remains valid for all viscosity variations satisfying the inequalities \eqref{eq:classical}, corresponding to viscosity regimes where the umbilic point is close to the vertex $O$. We verify the $L^1_{loc}$-stability of the Riemann solution with respect to variations in the data. While we do not establish the uniqueness of the Riemann solution, extensive numerical experiments confirm its validity. Our findings provide a comprehensive framework for understanding three-phase flow dynamics in porous media under a wide range of conditions.\\

\noindent {{\it keywords}: Riemann solutions; Multiphase flow in porous media; Heavy oil }
\end{abstract}

\tableofcontents
\section{Introduction}
\label{sec:introduction}
This paper is part of our effort to understand the Riemann problem for
a class of conservation laws systems that model the flow of three
immiscible fluids, such as water, gas, and oil, in a porous medium.
Here, we are concerned with Riemann problems modeling the injection of a mixture of water and gas
in oil reservoirs when the oil viscosity is much larger than those of water and gas. The same class of models includes
the problems of $CO_2$ sequestration in mature oil reservoirs, enhanced oil recovery with foam, and underground hydrogen storage \cite{iskandarov2022data,salimi2012influence,davies2024safety,bagchi2025critical,tang2022foam,lozano2024a}.
We want to understand the effect of the viscosity parameters in the solutions for cases of potential engineering interest. 

It is a fact that the system of conservation laws we consider
fails to be strictly hyperbolic at an umbilic point, 
i.e., an isolated point interior to state space 
at which the characteristic speeds coincide. 

For the viscosity parameters we consider, which have one of the fluid viscosities much larger than the other two,
the umbilic point lies in the subregions $\mathcal{A}$-$O$-$\mathcal{B}$ or $\mathcal{A}'$-$O$-$\mathcal{B}'$ of the
saturation triangle $G$-$W$-$O$, shown in Fig.~\ref{fig:ratioviscosities}(b).
So, the umbilic point
corresponds to Case II of the Schaeffer-Shearer classification of 2x2 systems
of non-strictly hyperbolic conservation laws \cite{schaeffer1987classification}.

In \cite{V.2010,Azevedo2014} we established the existence, uniqueness, and $L^1_{loc}$-continuity with respect to data variations for Riemann solutions in green reservoirs, where a combination of water, gas, and oil is injected into an oil-saturated reservoir. Later, in \cite{L.2016}, we presented the Riemann solution for left states along the edge \GW of the saturation triangle and right states representing the displaced mixture in a neighborhood $\mathcal{R}_1$ of vertex $O$, where only classical waves were needed. Subsequently, in \cite{Andrade2018}, we expanded these findings to include right states in two additional regions $\mathcal{R}_2$ and $\mathcal{R}_3$, 
completing a quadrilateral region in the saturation triangle with vertices $O$, $E$, $\um$ and $D$, in which
nonclassical waves (undercompressive and overcompressive shocks \cite{Andrade2018,Matos2016,lozano2024c, L.1990, silva2014riemann}) may appear.
In this paper, we extend our previous work by allowing the region
of right states $R$ to be significantly larger.
However, we have a more stringent restriction on the viscosity parameters
to prevent nonclassical waves from occurring.


In multiphase flow, it is essential to take into account
the effect of capillary pressure between the different fluids.
To take into account such an effect, we deem shocks admissible
if they possess a viscous profile.

All of the results in this work were obtained using quadratic functions to model the permeabilities of the fluid phases, enabling some explicit calculations.
The wave curve method of Liu \cite{Liu1974, Liu1975}, generalized
in \cite{L.1992a, V.2010}, is
combined with numerical techniques \cite{ELI_web} to provide scientific evidence 
for the existence of solutions and $L^1_{loc}$ continuity under change of data. Similar constructions, also based on the wave curve method, appear in references \cite{barros2021analytical, juanes2004analytical, mehrabi2020solution, pires2024approximate}.

The rest of this paper is organized as follows. Section~\ref{sec:math_model} describes the mathematical model and the problem of interest.
In Section~\ref{sec:FundamentalConstructs}, we discuss the rarefaction and
Rankine-Hugoniot curves that underpin the wave curve method.
We also recall some notation and concepts, and introduce the bifurcation loci 
that define the regions of right states with a common Riemann solution structure.
The section ends with an analysis of the dependence of such bifurcation loci
on the viscosity ratios.
Section~\ref{sec:RS}, the main results, details the structure of the Riemann
solution corresponding to the right states in different regions of the
saturation triangle with the left states on the edge \GWc.
Section~\ref{sec:Simulations} presents numerical simulations illustrating our theoretical results. Finally, Section~\ref{sec:Conclusion} contains concluding remarks.

\section{The mathematical model and the problem of interest}
\label{sec:math_model}
In this section, we introduce notation and recall
basic facts about the mathematical model, which is
also studied in \cite{V.2010}, \cite{Azevedo2014},
\cite{DeSouza1992} and \cite{L.1992a}. We consider
the flow of a mixture of three fluid phases (which,
for concreteness, are called water, gas and oil) in
a thin, horizontal cylinder of porous rock.
Let $\sw(x,t)$, $\sg(x,t)$ and $\so(x,t)$ denote the
corresponding saturations at distance $x$ along the
cylinder, at time $t$.
Under our simplifications, the governing equations for three-phase flow become
\begin{equation}\label{eq:system1}
\frac{\partial\sw}{\partial t}+\frac{\partial \fw}{\partial x}=0, \qquad
\frac{\partial \so}{\partial t}+\frac{\partial \fo}{\partial x}=0,
\end{equation}
which is a non-dimensionalized system representing the conservation of water, gas
oil, and Darcy's law. We adopt the following flux functions $\fw$ and $\fo$:
\begin{equation}\label{eq:flowfunct}
\fw = \sw^2/ (\muw \mathscr{D}),~~~~~\fo = \so^2 /(\muo \mathscr{D}),
\mbox{ where } \mathscr{D} =\sw^2  /\muw + \so^2  /\muo + \sg^2  /\mug > 0,
\end{equation}
and the constants $\muw$, $\muo$ and $\mug$ are the phase viscosities.
These represent a simple Corey's model with quadratic relative permeability
functions.
For conciseness, we may denote the flux functions compactly as
$F(S) = (\fw(S), \fg(S))^T$, with $S = (\sw, \sg)^T$.

Because $\sw +\sg +\so = 1$ and $0 \leq\sw , \so , \sg \leq 1$,
following practice in petroleum engineering, we depict the space of states of
the fluid mixture in barycentric coordinates as the saturation triangle
shown in Figs.~\ref{fig:ratioviscosities}(b) and \ref{fig:SecBifurc},
where we use the labels $O$, $W$, and $G$ to denote the states with
saturation $\so = 1$, $\sw = 1$, and $\sg = 1$, respectively;
the state $B$ corresponds to $\so = 0$, $\sw = \muw /(\muw +\mug )$, and
$\sg = \mug /(\muw + \mug )$; the state $D$ to $\so = \muo /(\muw + \muo )$,
$\sw = \muw /(\muw + \muo )$, and $\sg = 0$; the state $E$ to 
$\so = \muo /(\mug + \muo )$, $\sw = 0$, and $\sg = \mug /(\mug + \muo )$.
The state $\um$ is an umbilic point. It has coordinates 
$\so = \muo /\mut, \sw = \muw/\mut$, and $\sg = \mug /\mut$, with
$\mut = \muw + \muo + \mug$.
The state $D_0$ on \GD corresponds to $\sw + \so = \sg =  1/2$
and the state $E_0$ on \WE to $\sg + \so = \sw =  1/2$, see Fig.~\ref{fig:SecBifurc}. 


We are concerned with the injection problem where water and gas are
injected alternately in any proportion to displace the original mixture
consisting of any proportion of each phase. Mathematically this is translated
as a Riemann problem for System~\eqref{eq:system1} with initial data given by
\begin{equation}\label{eq:RP}
S(x, 0) = \begin{cases}
L, \quad \text{if } x < 0,\\
R, \quad \text{if } x > 0,
\end{cases}
\end{equation}
in which the injected mixture is represented by state $L$ and the displaced
one by state $R$. Specifically we consider left states $L$ on the edge \GWc,
and right states $R$ all over the saturation triangle in
Figs.~\ref{fig:ratioviscosities}(b) and \ref{fig:SecBifurc}.

To characterize the classes of Riemann solutions that will be described, we introduce the viscosity ratios: $r_w = \muw/\muo$ and $r_g = \mug/\muo$. 
In \cite{Andrade2018}, the  parameters
\begin{equation}\label{eq:razoes}
\mathcal{R}_{g o}^{\pm}=\left(r_{\mathrm{g}} \pm 1\right) /
r_{\mathrm{w}} \quad \text{and} \quad  \mathcal{R}_{w o}^{\pm}=\left(r_{\mathrm{w}} \pm
1\right) / r_{\mathrm{g}},
\end{equation}
were defined and used to establish that if the inequalities
\begin{equation}\label{eq:classical}
\left(\mathcal{R}_{g o}^{-}\right)^{2} /
\mathcal{R}_{g o}^{+} \ge 8 \quad \text{and}
\left(\mathcal{R}_{w o}^{-}\right)^{2} /
\mathcal{R}_{w o}^{+} \ge 8,
\end{equation}
hold, then the Riemann solutions for left states
$L$ on the edge \GW and right states $R$ in the quadrilateral
$OE\um D$, see Fig.~\ref{fig:SecBifurc}, 
comprise only classical waves.
 That is, the Riemann solutions consist, at most, of a fast wave group, followed by a slow wave group separated by an intermediate constant state.  

In the case of heavy oil reservoirs, we typically have $\mug < \muw <<\muo$, e.g. $r_w = 10^{-4}$ and $ r_g=10^{-5}$ (see for example \cite{miao2024core,zhang2006enhanced} ). Consequently the inequalities \eqref{eq:classical} are satisfied
and only classical
waves occur in the Riemann solution for any right state $R$ along the
saturation triangle and any left state along \GW in
Figs.~\ref{fig:ratioviscosities}(b) and \ref{fig:SecBifurc}.

We notice that in the regime we adopt in this work, the umbilic point $\um$
lies close to edge \WO and vertex $O$.
Thus, to have a better view of the geometric elements involved in the
construction of Riemann problem solutions of 
(\ref{eq:flowfunct}-\ref{eq:RP}) we plot figures with
$\muw = 1$, $\muo = 9.5$ and $\mug = 0.45$.
 With these values we have
\begin{equation}\label{eq:Ourcase}
{(\mathcal{R}_{go}^-)^2} / {\mathcal{R}_{go}^+} \approx 8.23141 \quad \text{and} \quad 
{(\mathcal{R}_{wo}^-)^2} / {\mathcal{R}_{wo}^+} \approx 15.2911,
\end{equation}
which means that the inequalities in \eqref{eq:classical} are verified.

The viscosities are further restricted so that the double contact
locus, see Definition~\ref{def:DoubleContact}, possesses only one branch in the saturation triangle.
This is to reduce the number of right-state regions for which 
the Riemann solutions undergo the same structural bifurcations
as the left state $L$ varies. 
These restrictions on the viscosity ratios limit
the localization of the umbilic point $\um$, whose coordinates are
determined by the viscosities, to the regions $\mathcal{A}\mathcal{B}O$ and $\mathcal{A}'\mathcal{B}'O$
in Fig.~\ref{fig:ratioviscosities}(b). 
Specifically, the boundary $\mathcal{B}'\mathcal{C}\mathcal{B}$ is defined by the equalities in \eqref{eq:classical}, whereas the
boundaries $\mathcal{A}O$ and $\mathcal{A}'O$ are defined by
the additional restriction imposed on the viscosity ratios.




\section{Fundamental Constructs}
\label{sec:FundamentalConstructs}

We begin this section by discussing the rarefaction and Rankine-Hugoniot curves
on which the wave curve method is based. Next, we recall some useful notation
and concepts used in this paper.
Then we introduce the pertinent bifurcation loci that define the boundaries of
regions of right states with a common Riemann solution structure.
We end the section with an analysis of the dependence of the bifurcation loci
on the viscosity ratios.

\subsection{Rarefaction and Rankine-Hugoniot curves}
\label{subsec:rarefactionWaves}

Rarefaction waves are constructed via rarefaction curves defined by
properly oriented eigenvectors of the Jacobian derivative matrix of the flux $F(S)$ for~\eqref{eq:system1}.
Along such a curve, the corresponding eigenvalue ({\it slow-family} or {\it fast-family} characteristic speed) must be monotone, forcing these curves to stop at {\it inflection loci}, where the characteristic speed has an extremum.

\begin{figure}[]
	\centering
	\subfigure[Slow-family integral curves.]
	{\includegraphics[width=0.4\linewidth]{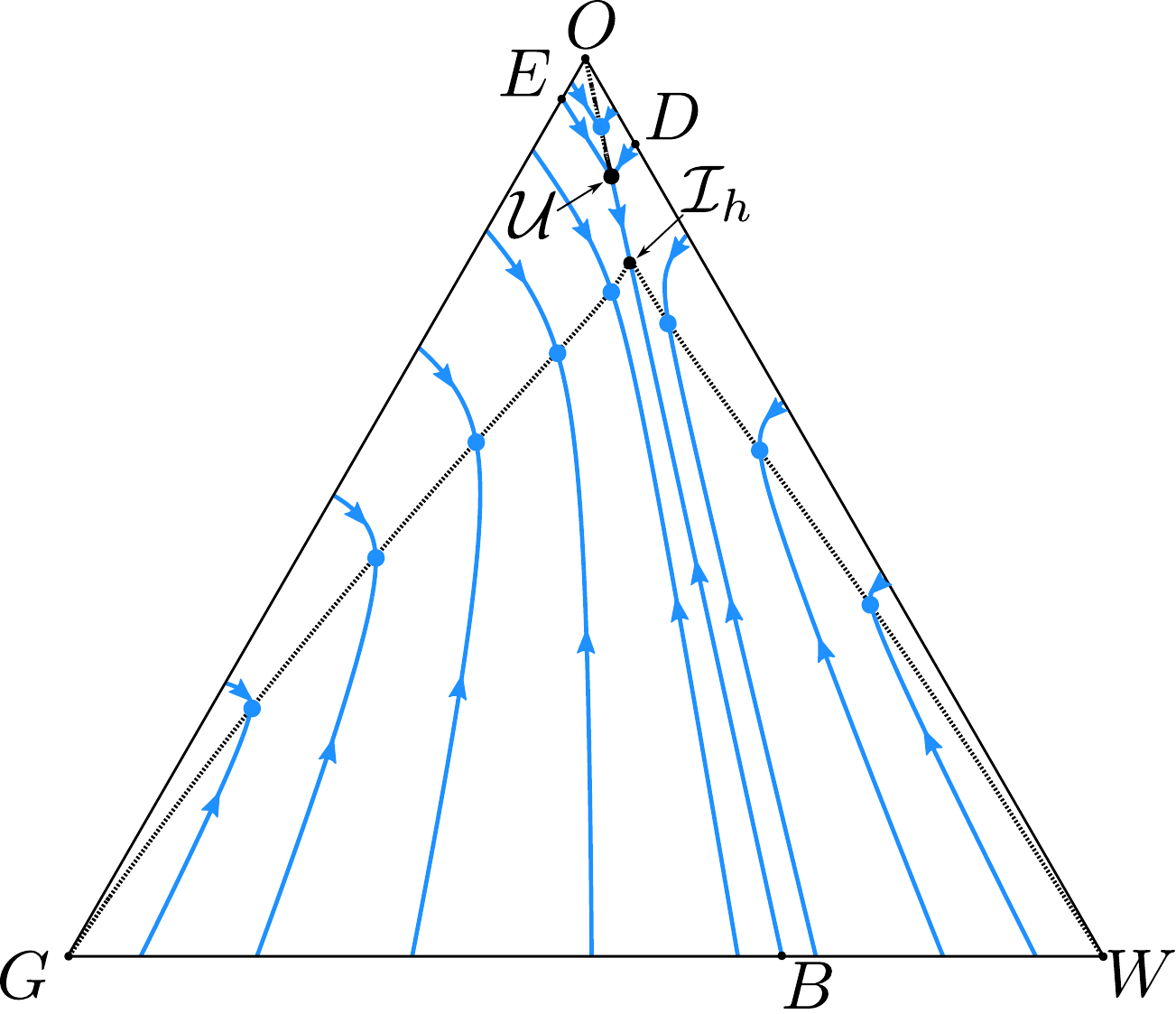}} 
	\hspace{0.8mm}
	\subfigure[Fast-family integral curves.]
	{\includegraphics[width=0.4\linewidth]{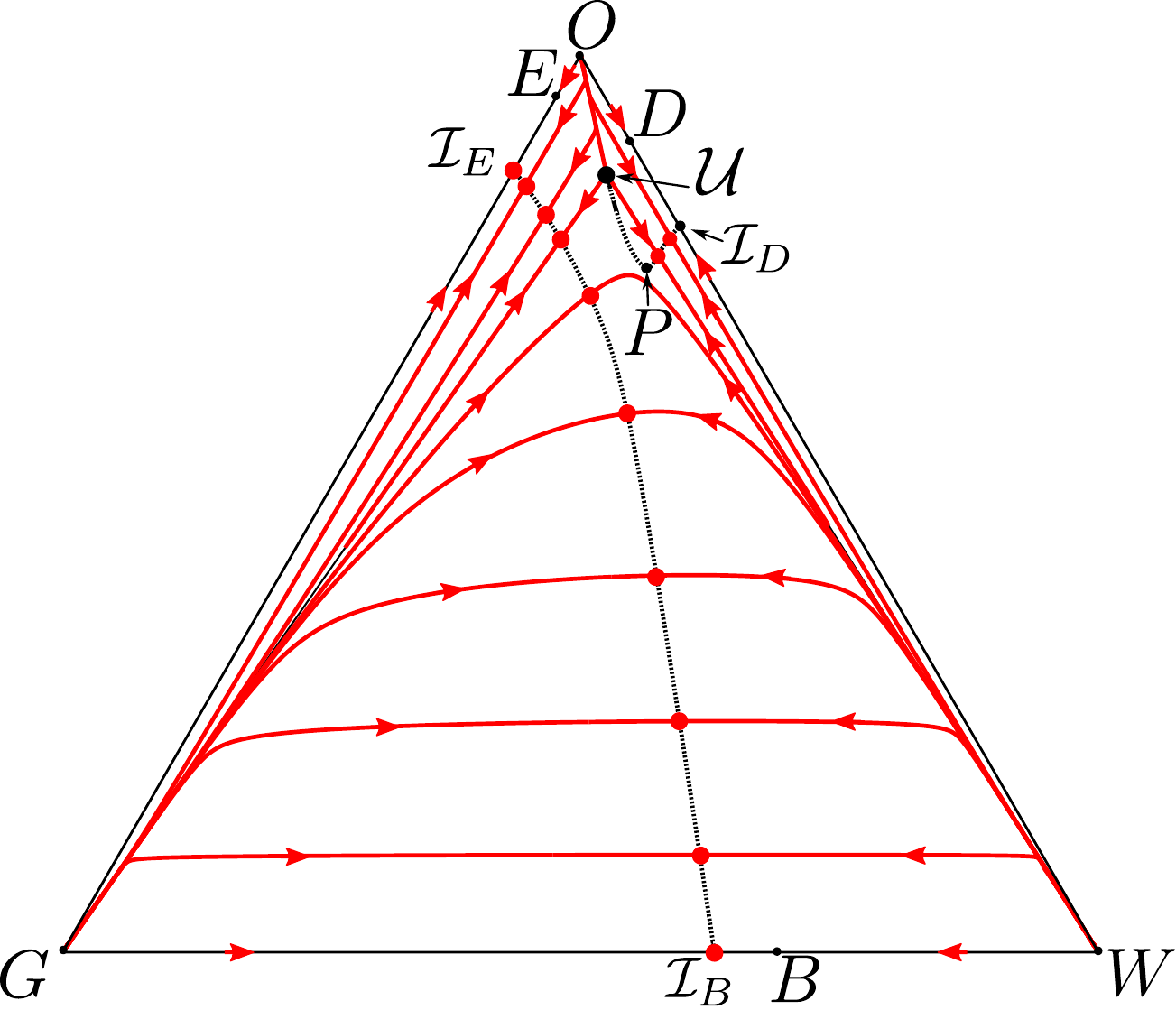}}  
	\caption{ Integral curves and inflection loci within the saturation triangle $GWO$. The arrows indicate the direction of increasing characteristic speed defining rarefaction waves.
	State $\mathcal{I}_h$ in (a) is the ``summit'' of the slow inflection locus. At state $P$ in (b), the fast integral curve is tangent to the fast inflection locus. The fast inflection locus in (b) intersects the edges of the saturation triangle at points $\mathcal{I}_D$,
    $\mathcal{I}_B$ and $\mathcal{I}_E$.
 	}
	\label{fig:IntegralCurves}
\end{figure}

It is proved in \cite{Azevedo2014} that the characteristic speeds are real and positive in the interior of the saturation triangle.
On the edges, the slow-family characteristic speed $\lambdas$ is zero; as there are two fluids, the flow is governed by
the Buckley-Leverett equation. The fast-family characteristic speed $\lambdaf$ is positive, except at the vertices where it also vanishes. At the {\it umbilic point} $\um$ in Fig.~\ref{fig:IntegralCurves}, the
characteristic speeds coincide; see \cite{Azevedo2014} for other properties. On the closed triangle $\lambdas < \lambdaf$ except at
$G$, $W$, $O$, and $\um$.

Figures~\ref{fig:IntegralCurves}(a) and~(b) depict slow- and fast-family integral curves for our model.
They form two foliations with singularities at the umbilic point
and vertices; see \cite{Azevedo2014}.
Each of the foliations maintains topological equivalence, irrespective
of the viscosities; transversality between the two foliations is also preserved. The arrows in Fig.~\ref{fig:IntegralCurves}
indicate the direction of increasing
characteristic speed, which reaches a maximum at the small
dots along the respective inflection locus.

A shock wave consists of two constant states, $M$ and $N$, separated by
a discontinuity traveling with speed $\sigma$;
the Rankine-Hugoniot jump condition relates the states
\begin{equation}\label{eq:RH}
    \sigma(M - N) = F(M) - F(N).
\end{equation}
For a given state $M$, the set of states $N$ satisfying the jump condition
forms the Hugoniot locus based on $M$,
which parameterizes shock waves.
If the base state $M$ is the left state of the discontinuity,
such locus is called the {\it forward}  Hugoniot curve  of $M$, denoted as $ \mathcal{H}^+(M)$,
and we write $\sigma = \sigma(M; N)$.
On the other hand, if $M$ is the right state,
it is called the {\it backward} Hugoniot curve  of $M$, denoted as $ \mathcal{H}(M)$,
and we write $\sigma = \sigma(N; M)$.

We improperly refer to Hugoniot loci as {\it Hugoniot curves} to stress that they are one-dimensional objects,
even though they may have self-intersections and singularities.

To have the uniqueness of Riemann solutions, an admissibility criterion for discontinuities is required.
In multiphase flow, it is essential to take into account the effect of capillary 
pressure between the different fluids. When it is taken into account, instead of
the system~\eqref{eq:system1}, which still appears on the left-hand side below,
we obtain the following parabolic system:
\begin{equation}
\label{eq:parabolic_eq}
\frac{\partial S}{\partial t} + \frac{\partial F(S)}{\partial x} =
\varepsilon\frac{\partial}{\partial x} \left[D(S) \frac{\partial S}{\partial x} \right],
\end{equation}
where the matrix $D(S)$ is positive definite
in the interior of the saturation triangle, \cite{azevedo2002capillary}.

In view of system~\eqref{eq:parabolic_eq}, the natural admissibility criterion
for shocks is the \emph{viscous profile criterion}:
a shock joining a left state $M$ to a right state $N$ must be the limit, as the
positive parameter $\varepsilon$ tends to zero, of traveling wave solutions
$S(x,t) = {\tilde S}\left((x - \sigma t)/\varepsilon\right)$
of \eqref{eq:parabolic_eq}, with the boundary conditions
$S(-\infty) = M$ and $S(+\infty) = N$.
Given that the traveling wave $\tilde{S}$ is a smooth function of $\xi=x - \sigma t)/\varepsilon$, the system can be simplified to a set of ordinary differential equations (ODEs) (omitting hats):
\begin{equation} 
\label{eq:dinamycsyst}
 D(S(\xi)) \dfrac{d S(\xi)}{d \xi} = -\sigma \left(S(\xi) - M\right) + F(S(\xi)) - F(M),   
\end{equation}
where $M$ and $N$ represent equilibrium points of \eqref{eq:dinamycsyst}, meaning that the right-hand side of the system vanishes when evaluated at $M$ or $N$. A traveling wave solution corresponds to a trajectory connecting $M$ and $N$, which is referred to as a viscous profile for the shock wave.

Following the approach established in our previous studies \cite{L.2016, Andrade2018, Azevedo2014, Lozano}, we assume $D(S)$ to be the identity matrix and numerically verify the admissibility of shocks. This simplification ensures that the ODE system \eqref{eq:dinamycsyst} remains invariant along the segments \GDc, \WEc,  and \OBc. Consequently, the admissibility of nonlocal shocks is determined by the relative positions of $M$ and $N$ with respect to the corresponding segments \GDc, \WEc, or \OBc, as discussed in \cite{Lozano}.

The Hugoniot locus is based on a state $M$ along any edge or one of the invariant lines
\EWc, \DGc, and \BO of the saturation triangle is easily obtained;
see \cite{L.2016, Azevedo2014}.
Generically, it consists of such an edge or line together with a hyperbola,
one branch of which contains $M$.
If $M$ is moved away from an edge or from one of the lines \EWc, \DGc, and \BOc,
the Hugoniot locus $ \mathcal{H}(M)$ varies continuously.
In particular, its topological structure with two primary (attached) branches and one
detached branch (which may lie outside the saturation triangle) is preserved;
see figures in Section~\ref{sec:RS} for typical examples.


\begin{figure}[]
	\centering
	 \subfigure[Regions $\mathcal{A}O\mathcal{B}$ and $\mathcal{A}'O\mathcal{B}'$ in the parameter space $r_w \times r_g$ obeying the restriction of the phase viscosities considered in this paper. 
 The point $\um$ represents  the viscosity ratios $r_w = 1/9.5$ and
 $r_g = 0.45/9.5$.]
     {\includegraphics[width=0.35\linewidth]{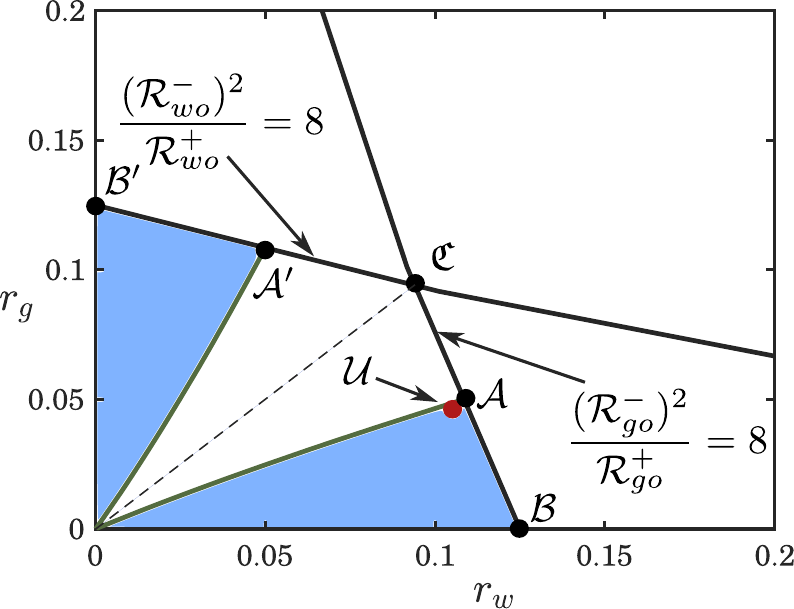}} 
     \hskip 0.8cm
     \subfigure[Subtriangles $I$, $II_G$, $II_W$ and $II_O$ correspond to cases $I$ and $II$ of the Schaeffer-Shearer classification  \cite{schaeffer1987classification} for umbilic points. 
     $\mathcal{A}O\mathcal{B}$
     and $\mathcal{A}'O\mathcal{B}'$ are the possible locations of the umbilic point in the present work.]
          {\includegraphics[width=0.33\linewidth]{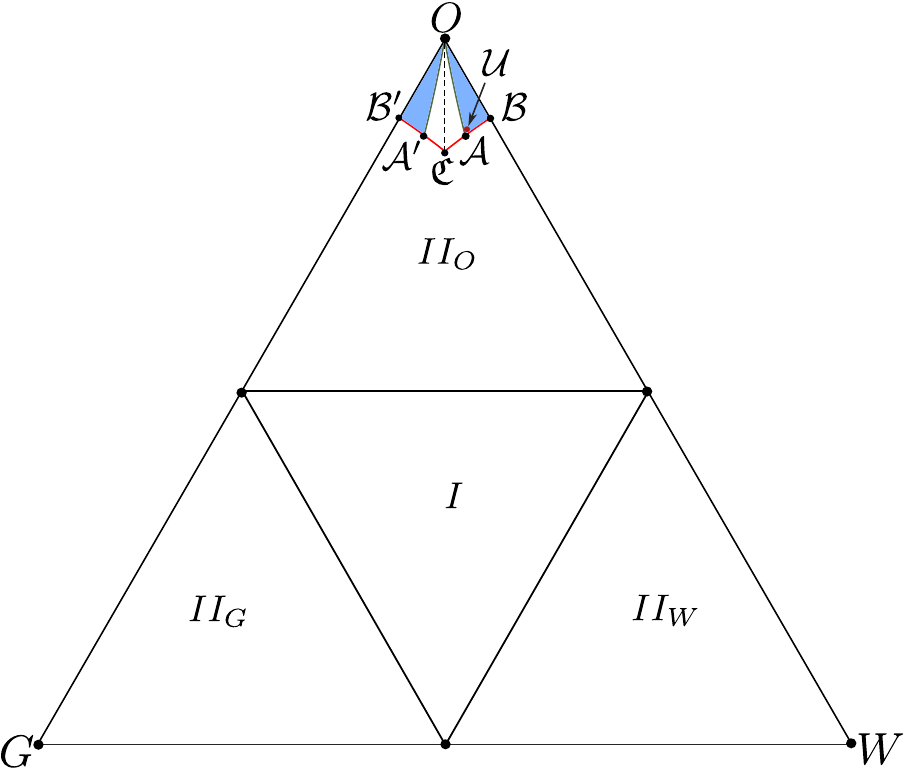}}
          \caption{Geometric representation of the imposed 
          restrictions on the fluid viscosities and the saturation
          triangle subdivision according to the Schaeffer-Shearer
          classification for umbilic points.}
	\label{fig:ratioviscosities}
\end{figure}

\begin{figure}[]
	\centering
     {\includegraphics[width=0.3750\linewidth]{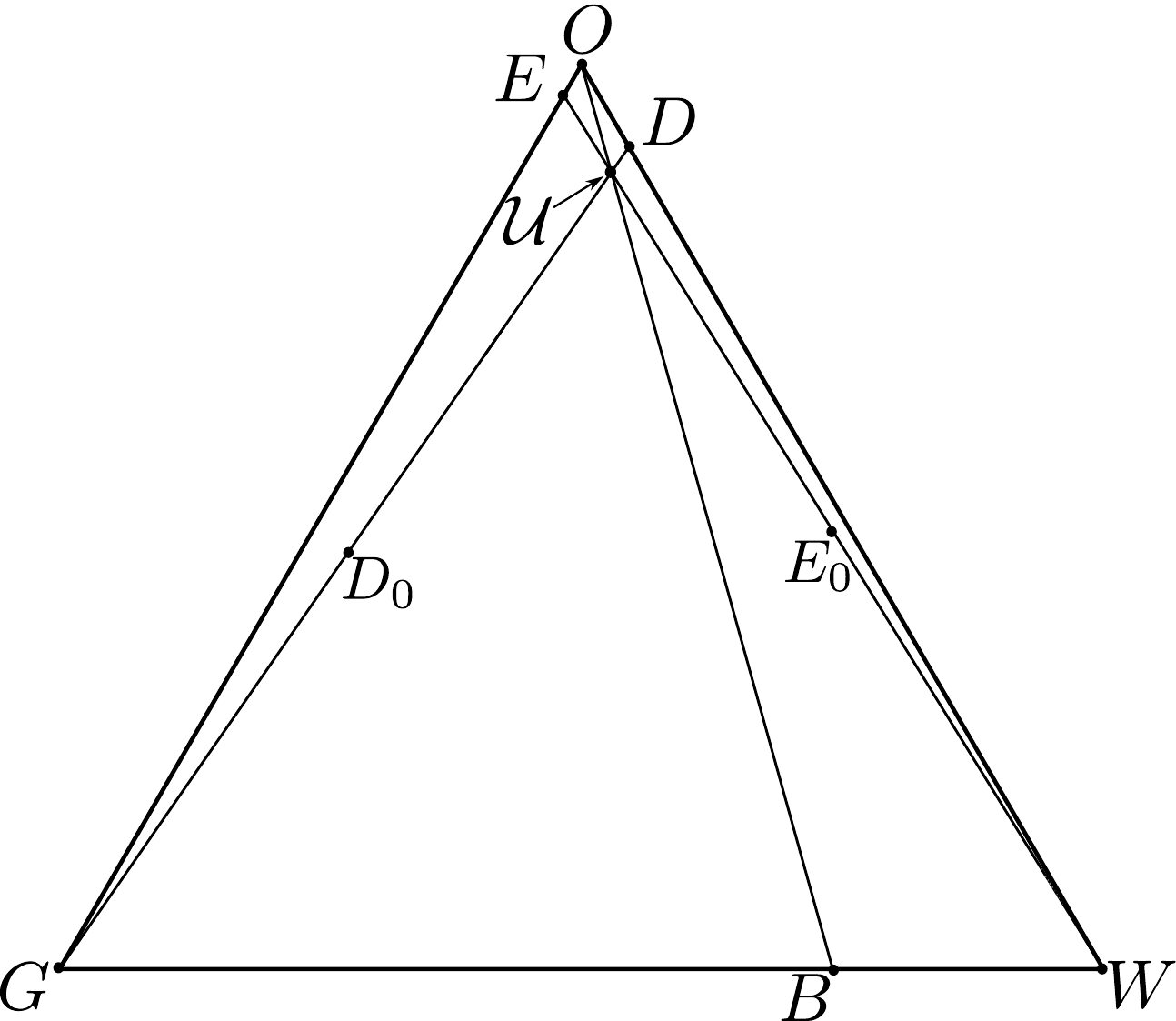}}   
 \caption{Secondary bifurcation loci
 \GDc, \WEc, \OBc.
 Segments \BUc, \DUc, and \EU are invariant under the slow-characteristic vector field.
 Segments \GUc, \WUc,  and \OU are invariant under the fast-characteristic vector field.}
	\label{fig:SecBifurc}
\end{figure}

\begin{figure}[]
	\centering
	\subfigure[$s$-hysteresis locus: $\sigma(\mathcal{I}_i; H_i) = \lambdas(\mathcal{I}_i)$, $\mathcal{I}_i \in \mathcal{I}_s$, $i=1,\dots,6$.]
	{\includegraphics[width=0.4\linewidth]{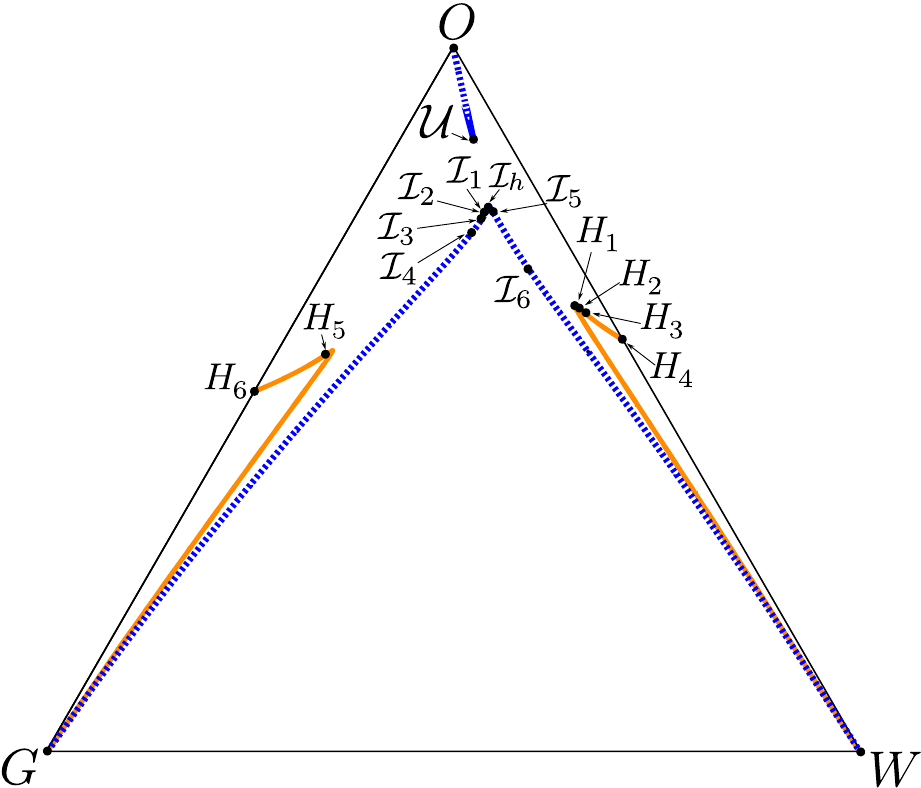}} 
	\hspace{0.8mm}
	\subfigure[$M_E$ is the $f$-extension of $\mathcal{I}_s$
    with $\sigma(\mathcal{I}_i; M_E^i) = \lambdaf(\mathcal{I}_i)$, $\mathcal{I}_i \in \mathcal{I}_s$ $i=a,b,c,d$.]
	{\includegraphics[width=0.4\linewidth]{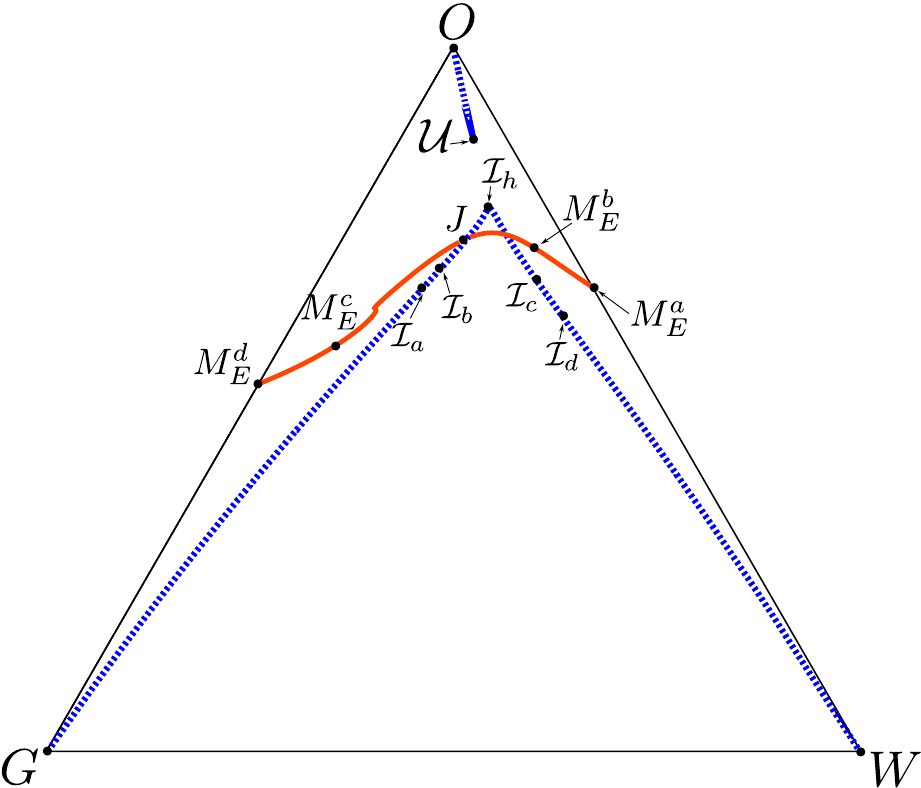}}  
	\caption{(a) State $\mathcal{I}_h$ is the ``summit'' of $\mathcal{I}_s$ in Fig.~\ref{fig:IntegralCurves}(a);
 the state $H_3$ lies in the invariant segment $\um$-$W$;
 states ${H_1}$ and ${H_5}$ are intersection points of the $s$-hysteresis locus
 and the $f$-hysteresis locus in Fig.~\ref{fig:BoundaryExt_MixedDoubleCont}(b); they correspond to $\mathcal{I}_1$ and $\mathcal{I}_5$ in $\mathcal{I}_s$. 
 (b) The state $J$ is the intersection point of $M_E$ with the inflection loci $\mathcal{I}_s$ and $\mathcal{I}_f$, see Fig.~\ref{fig:IntegralCurves}. 
 States $M_E^a$ and $M_E^d$ in $M_E$ lie on the boundary of the saturation triangle; 
 they correspond to states $\mathcal{I}_a$ and $\mathcal{I}_d$ in $\mathcal{I}_s$.
 The state $M_E^b$ is the intersection point of $M_E$ and the segment $P$-$H_1$ of the $f$-hysteresis locus in Fig.~\ref{fig:BoundaryExt_MixedDoubleCont}(b); this state corresponds to $\mathcal{I}_b$ in $\mathcal{I}_s$. The state $M_E^c$ is the intersection point of $M_E$ and the segment $H_5$-$P'$ of the \break
 $f$-hysteresis locus; this state corresponds to
 $\mathcal{I}_c$ in $\mathcal{I}_s$.
  }
\label{fig:IntegralCurves_extensions}
\end{figure}






\subsection{Recalling notation and basic concepts}\label{sec:notation}
We begin by recalling some definitions and notations used in \cite{Matos2016} to describe Riemann solutions. 
We use the following nomenclature for discontinuities
joining a left state $M$  and a right state $N$ with characteristic speeds related to the propagation
speed $\sigma = \sigma(M; N)$ as below:
\begin{itemize}
\item $s$-shock: $\lambda_s(N) < \sigma < \lambda_s(M)$ and $ \sigma < \lambda_f(N)$.
\item $f$-shock: $\lambda_f(N) < \sigma < \lambda_f(M)$ and $\lambda_s(M) < \sigma$.
\item undercompressive or $u$-shock: $\lambda_s(M) < \sigma < \lambda_f(M)$, $\lambda_s(N) < \sigma < \lambda_f(N)$.
\item overcompressive or $o$-shock: $\lambda_f(N) < \sigma < \lambda_s(M)$.
\end{itemize}
Lax, \cite{Lax1957}, used the nomenclature $1$- and $2$-shock for slow and fast shocks, while we
adopted the nomenclature $s$- and $f$-shocks. For more general conservation
laws, characteristics can be tangent to discontinuity. 
In certain cases, we
allow some equalities in the previous Lax configurations, giving way to a
generalized Lax criterion \cite{Liu1974, Liu1975}.
\begin{rem}\label{rem:admissibility}
We emphasize that any shock used in this work
must obey the viscous profile admissibility criterion, with the viscosity matrix being the identity.
\end{rem}

A shock of type $\alpha$, $\alpha \in \{s,f,u,o\}$, connecting a state $M$ to a state $N$ on its right is denoted as $M\xrightarrow{S_\alpha} N$. A {\it prime} on the left or on the right of $S$ indicates a left- or a right-characteristic shock, respectively. For instance, $M\xrightarrow{\CS_f} N$ denotes a $f$-shock joining $M$ to $N$ with propagation speed $\sigma(M;N)=\lambda_f(M)$. On the other hand $M\xrightarrow{R_i} N$ denotes an $i$-rarefaction wave joining $M$ to $N$ associated to the characteristic speed $\lambda_i$.
An $i-${\it wave group}  connecting a state $M$ to a state $N$ on its right, denoted by $M\xrightarrow{\mathcal{G}_i} N$, consists of a sequence of $i-$rarefaction and shock waves with no intermediate constant states filling sectors in the physical space $xt$.

In a wave group $M\xrightarrow{\mathcal{G}_i} N$, we use $v_{ib}$ and $v_{ie}$ to denote the left-most and the right-most wave speeds in the group. This means that if the left-most wave is a shock, then $v_{ib}$ is the shock speed, but if it is 
an $i$-rarefaction $v_{ib} =\lambda_i(M)$; similarly, if the right-most wave is a shock then $v_{ie}$
is the shock speed, but if it is 
an $i$-rarefaction $v_{ie} =\lambda_i(N)$. A wave group sequence
$M\xrightarrow{\mathcal{G}_s} N\xrightarrow{\mathcal{G}_f} P$ is considered \emph{compatible} if and only if
\begin{eqnarray}\label{eq:Cond_Compat}
v_{se} \le v_{fb}.
\end{eqnarray}

The {\it backward} $i$-{\it wave curve} for a fixed state $N$, denoted by $\wm_i(N)$, is a one-dimensional parametrization of the states $M$ which can be connected to $N$ on its left by an $i$-wave group.
The {\it forward} $i$-{\it wave curve} for a fixed state $M$, denoted by $\wm_i^+(M)$, is a one-dimensional parametrization of the states $N$ which can be connected to $M$ on its right by an $i$-wave group.

A Riemann solution
connecting the left state $L$ to the right state $R$ of a Riemann problem consists
of (at most) 2 compatible wave groups generally separated by constant states.
Each wave group is associated with a single wave family.
It is represented, for simplicity, as a sequence of states separated by arrows with the
elementary wave types on top of the arrows. The left-most 
state in a wave group is followed by a right-arrow with
a tail ($\testright{ }$), while states that occur internal to a wave group
are followed by a single right-arrow ($\xrightarrow{ }$). 
For example, a sequence $L\testright{R_s} M \xrightarrow{'S_s} N \testright{S_f}R$ represents a Riemann solution consisting of a $s$-wave group comprising a
rarefaction from $L$ to $M$ adjacent to a left-characteristic shock from $M$ to
$N$, preceded by a $f$-wave group comprising a shock from $N$ to $R$.

A curve segment $A$-$B$ joining $A$ to $B$ can also be
denoted using the interval notation to indicate whether
an endpoint belongs to the segment or not. 
For example, $[A, B)$ denotes the curve segment that joins $A$ to $B$, including $A$ and excluding $B$.

To construct
Riemann solutions, we use the wave curve method \cite{Liu1974} as in \cite{L.2016, Andrade2018, V.2010, Azevedo2014}.
Namely, we employ rarefaction curves,
Hugoniot curves, wave curves, and check the compatibility of the
ending speed $v_{se}$ against the beginning speed $v_{fb}$ of the pertinent slow and fast wave groups.

\subsection{Bifurcation loci}\label{sec:Bifurcation}

In this section, we recall some loci \cite{L.1992a, V.2015} in the state space and define a new one.
These loci play fundamental roles in constructing Riemann solutions.

\begin{defn}\label{def:seconBif}
A state $S$ belongs to the {\it secondary bifurcation locus} for the family $i\in\{s,f\}$ if there exists a state $S'\neq S$ such that
$
S' \in \mathcal{H}(S) \mbox{  with  } \lambda_i(S') = \sigma(S;S') \mbox{  and  } l_i(S')\cdot (S'-S)=0$,
where $l_i(S')$ is the left eigenvector of the Jacobian matrix $DF(S')$ associated with $\lambda_i(S')$.
\end{defn}
In our model, the secondary bifurcation loci comprise the line segments \EWc, \GDc, and \OBc. See Fig.~\ref{fig:SecBifurc}.

\begin{rem}\label{InvariantLines}
The segments \EUc, \BUc, and \DU are invariant under
the slow-characteristic vector field, while \GUc, \WUc, and \OU are invariant under 
the fast-characteristic vector field.
\end{rem}

\begin{defn}\label{def:iinflection}
A state $S$ belongs to the $i$-{\it inflection locus} $\mathcal{I}_i$
if $\nabla \lambda_i(S)\cdot r_i(S) = 0$,
where $r_i(S)$ is the right eigenvector of the Jacobian matrix $DF(S)$ associated with $\lambda_i(S)$.
\end{defn}
The $s$- and $f$-inflection loci for our model are shown in Fig.~\ref{fig:IntegralCurves}(a)-(b).

 \begin{defn}\label{def:DoubleContact}
A state $S$ belongs to the $(i,j)-${\it double contact locus} if there is a state $S'$ in the Hugoniot locus of $S$ with $\lambda_i(S) = \sigma(S;S')=\lambda_j(S')$.
\end{defn}

\begin{figure}[]
	\centering
	\subfigure[]
 	{\includegraphics[width=0.4\linewidth]{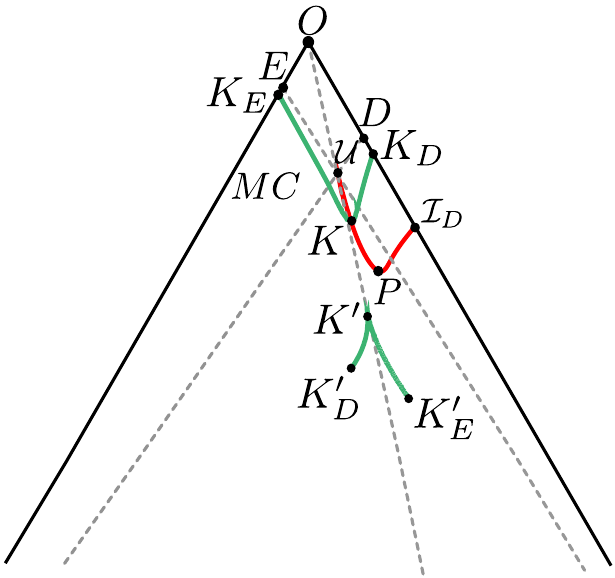}}  
  \hspace{2mm}
	\subfigure[]
 	{\includegraphics[width=0.4\linewidth]{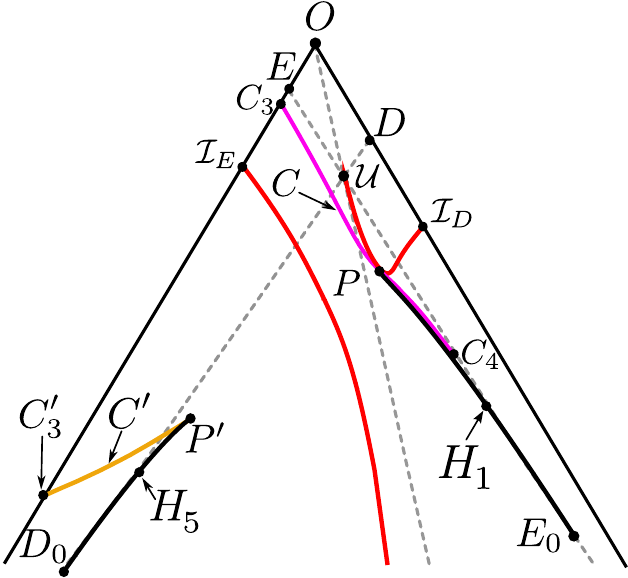}}  
	\caption{The relevant parts of some bifurcation loci. (a) $f$-inflection locus (red) and mixed contact locus (green).
    We have: $\sigma(M'; M) = \lambdas(M') = \lambdaf(M)$; 
    each state $M'$ in $K_E'$-$K_D'$ corresponds univocally to a state $M$ in $K_E$-$K_D$. 
(b) $f$-inflection locus (red), $f$-hysteresis locus (black), and double contact locus (purple).
Regarding the latter, the segment $C_3$-$P$ is an $f$-extension of $P$-$C_4$ and vice-versa.
    The yellow segment $C_3'$-$P'$ is the $f$-extension of  $C_3$-$P$ such that $\sigma(M; M') = \lambdaf(M)$, with $M$ in $C_3$-$P$ and $M'$ in $C_3'$-$P'$.
    The $f$-hysteresis segments $P$-$E_0$ and $P'$-$D_0$ are $f$-extensions
    of the $f$-inflection segment $P$-$\um$, characteristic on 
    $P$-$\um$.
    The state $P$ is defined in Fig.~\ref{fig:IntegralCurves}(b);
the states $H_1$ and $H_5$ are defined in  Fig.~\ref{fig:IntegralCurves_extensions}(a).
}
	\label{fig:BoundaryExt_MixedDoubleCont}
\end{figure}

For brevity, we will refer to the $(f,f)$-double contact locus simply as the {\it double contact} and to the $(s,f)$- (or $(f,s)$)-double contact locus simply as {\it mixed contact}. 
To avoid overloading the figures, we show in Fig.~\ref{fig:BoundaryExt_MixedDoubleCont} only the parts of the mixed and the double contact loci relevant to the present work.
Also, Fig.~\ref{fig:BoundaryExt_MixedDoubleCont} shows the relative position of the relevant parts of the double contact loci and the $f$-inflection locus.

Now, we recall the concept of an extension of a given locus.

\begin{defn}\label{def:Extension-cido}
	Let $\mathcal{L}$ be a locus in state space.
	A state $N$ belongs to the $i$-extension locus of $\mathcal{L}$
 if there exist a scalar $\sigma \in \mathbb{R}$ and a state $M \in \mathcal{L}$ satisfying the Rankine-Hugoniot jump condition \eqref{eq:RH}, with
 $\sigma = \lambda_i(M)$ or $\sigma =\lambda_i(N)$.
 If the state $M$ lying on $\mathcal{L}$ is the right state of the shock, the extension is called a {\it backward extension};
if $M$ is the left state, the extension is called a {\it forward extension}.
If $\sigma = \lambda_i(M)$ as $M$ varies on the base locus $\mathcal{L}$, the extension locus is called characteristic on $\mathcal{L}$;
if $\sigma = \lambda_i(N)$, it is called characteristic on itself.
\end{defn}



Examples of extension loci are provided in the following two definitions.

\begin{defn}\label{def:icomposite-cido}
	A backward (forward) $i$-composite segment is the backward (forward) $i$-extension of an $i$-rarefaction segment, characteristic on the $i$-rarefaction segment.
\end{defn}

\begin{defn}\label{def:ihysteresis-cido}
	The $i$-hysteresis locus $H_i$ is the forward $i$-extension of the $i$-inflection locus, characteristic on the $i$-inflection locus.
\end{defn}

The relevant parts of the $s$-hysteresis locus, and of the $f$-hysteresis locus for the present model
are shown in Fig.~\ref{fig:IntegralCurves_extensions}(a) and in 
Fig.~\ref{fig:BoundaryExt_MixedDoubleCont}(b), respectively.

\begin{rem} \label{rem:M_E}
Another important extension locus is the $f$-extension of the $s$-inflection locus, characteristic at the $s$-inflection that we denote by $M_E$. Relevant parts of this extension locus in the present model are shown in Figs.~\ref{fig:IntegralCurves_extensions}(b), \ref{fig:RRegions Completed}, \ref{fig:RegioesThetas}(b,c), \ref{fig:RRegionsOmegas}(b), and \ref{fig:RRegionsGamma}.
 \end{rem}

Finally, we introduce an exceptional locus defined by the tangency states of Hugoniot curves and a given locus $\mathcal{L}$.

\begin{defn}\label{def:T_Inflection-cido}
Let $\mathcal{L}$ be a locus in state space.
A state $M$ belongs to the tangential extension of $\mathcal{L}$
if there is a state $N$ in $\mathcal{L}$ such that
the tangents of the Rankine-Hugoniot curve based on $M$ and
$\mathcal{L}$ are parallel at $N$.

In our model, the tangential extension of the slow inflection locus
is denoted by $T_I$. The relevant parts of this locus are shown in
Figs.~\ref{fig:RegioesThetas}(b), \ref{fig:RRegionsOmegas} and
\ref{fig:RRegionsGamma}.

\end{defn}



\section{Riemann Solutions}\label{sec:RS}

\begin{figure}
	\centering
	{\includegraphics[scale=0.6]{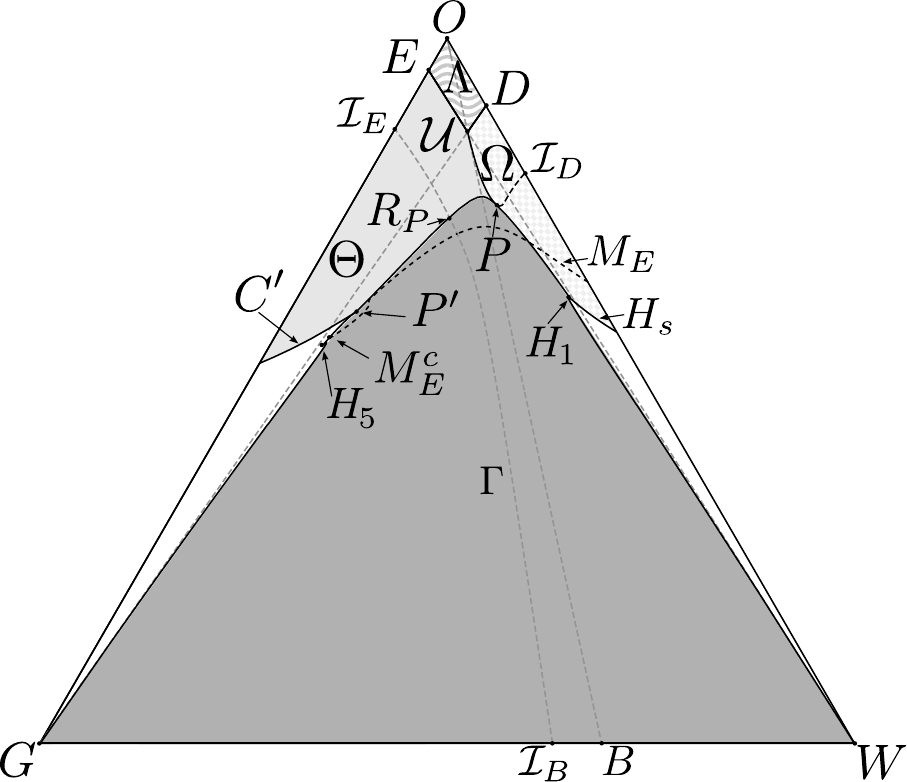}} 
	\caption{Subdivision of $\Delta$ into $R$-regions $\Lambda$, $\Theta$, $\Omega$, and $\Gamma$.}
	\label{fig:RRegions Completed}
\end{figure}

We start by dividing the saturation triangle into four regions, $\Lambda$, $\Theta$, $\Gamma$, and $\Omega$, according to the behavior of the backward fast wave curve of right states $R$; see Fig.~\ref{fig:RRegions Completed}.

\begin{figure}
	\centering
	\subfigure[Subdivision of $\Theta$.
    The boundaries of the subregions are:
the secondary bifurcation segment $E$-$\um$, the $f$-inflection segment $\um$-$P$,
the double contact segment $C_3$-$P$, the $f$-inflection segment $\mathcal{I}_E$-$R_P$,
the $f$-rarefaction segment $P$-$R_P$, the $f$-composite segment $R_P$-$P'$ defined by
$P$-$R_P$, the forward $f$-extension segment $P'$-$C_3'$ of the double contact
segment $P$-$C_3$ characteristic on $P$-$C_3$, and the segment $C_3'$-$E$ of the edge \GOc.
]
	{\includegraphics[scale=0.65]{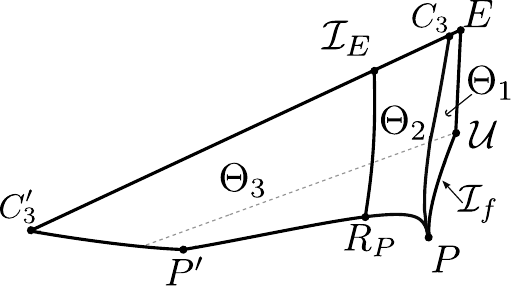}}  
	\hspace{5mm}
	\subfigure[Subdivision of $\Omega$]
	{\includegraphics[scale=0.595]{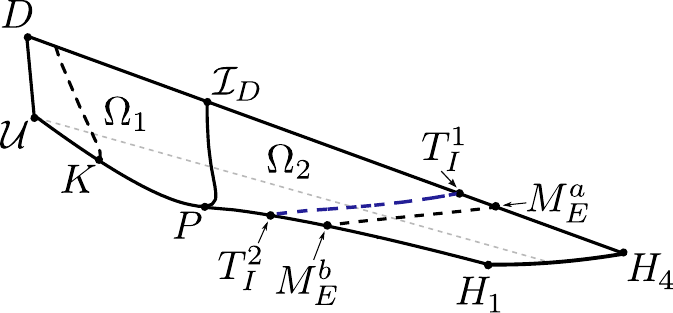}}  
	\hspace{0.75mm}
	\subfigure[Subdivision of $\Gamma$]
	{\includegraphics[scale=0.65]{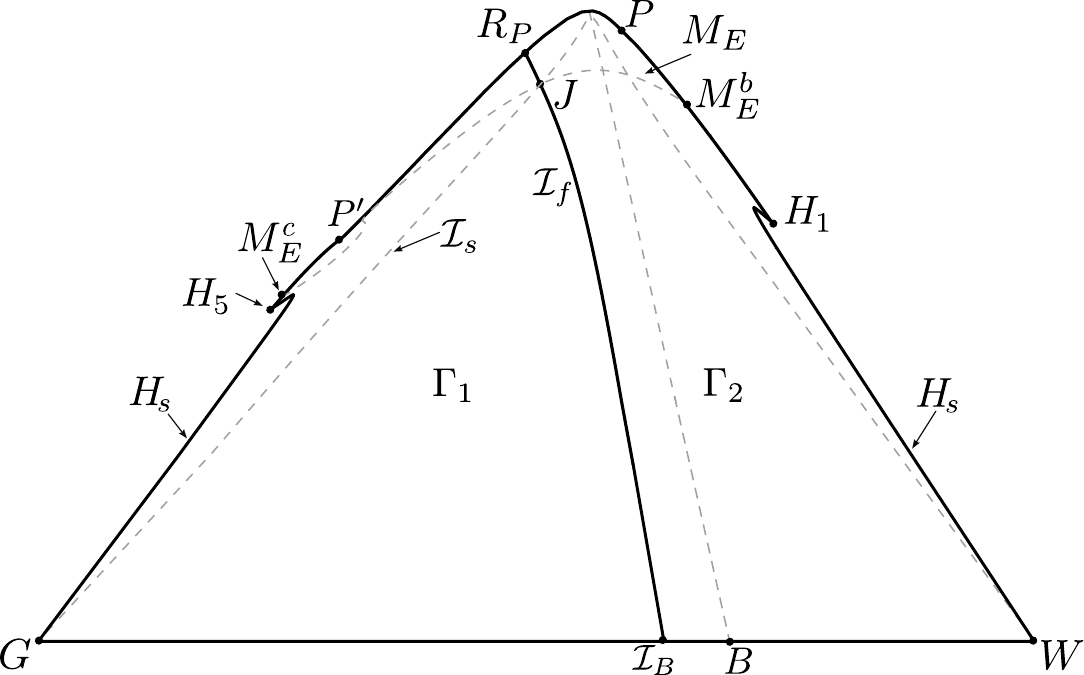}}  
\caption{Subdivisions of regions $\Theta$, $\Omega$ and $\Gamma$. 
}
	\label{fig:RegioesThetas}
\end{figure}

Region $\Lambda$, corresponding to the quadrilateral $OE\um D$, was studied in \cite{L.2016, Andrade2018}, so we do not consider it here. 
Region $\Theta$, see Figs.~\ref{fig:RRegions Completed} and \ref{fig:RegioesThetas}(a), is bounded by the invariant segment \EUc, the segment $\mathcal{U}$-$P$ of the fast inflection locus, the $f$-rarefaction segment $P$-$\mathcal{R}_P$, the forward $f$-composite segment $\mathcal{R}_P$-$P'$ defined by $P$-$\mathcal{R}_P$, the segment $P'$-$C'_3$ of the forward $f$-extension of the fast double contact segment $P$-$C_3$, characteristic on $P$-$C_3$ (shown in Fig.~\ref{fig:BoundaryExt_MixedDoubleCont}(b)), and the segment $C'_3$-$E$ on the edge \GO of the saturation triangle. 
Region $\Omega$, see Figs.~\ref{fig:RRegions Completed} and \ref{fig:RegioesThetas}(b),  is bounded by the invariant segment \DUc; the segment $\mathcal{U}$-$P$ of the fast inflection locus and the segment $P$-$H_1$ of the fast hysteresis locus corresponding to $\mathcal{U}$-$P$; the segment $H_1$-$H_4$ of the slow hysteresis locus $H_s$; and the segment
$H_4$-$D$ on the edge \WO of the saturation triangle.
Region $\Gamma$, see Figs.~\ref{fig:RRegions Completed} and \ref{fig:RegioesThetas}(c), is bounded by the $f$-rarefaction segment $P$-$\mathcal{R}_P$; the forward $f$-composite segment $\mathcal{R}_P$-$P'$ corresponding to $P$-$\mathcal{R}_P$; the segment $P'$-$H_5$ of the fast hysteresis locus corresponding to part of the segment $P$-$\mathcal{U}$ of the fast inflection locus; the segment $H_5$-$G$ of the slow hysteresis locus corresponding to the slow inflection segment $\mathcal{I}_5$-$G$; the edge \GWc; the segment $W$-$H_1$ of the slow hysteresis locus corresponding to the slow inflection segment $W$-$\mathcal{I}_1$;
and the segment $H_1$-$P$ of the fast hysteresis locus corresponding to part of the segment $P$-$\mathcal{U}$ of the fast inflection locus.
Each of the regions $\Theta$, $\Omega$, and  $\Gamma$ is further
subdivided subregions, which will be detailed in the following 
subsections.



 \begin{figure}[]
	\centering
	\subfigure[R-region $\Theta_3$
    subdivided by the secondary bifurcation segment $\mathcal{I}_2^f$-$C_2'$ 
    and the $f$-rarefaction segment $\mathcal{I}_1^f$-$C_1'$.
    ]
	{\includegraphics[scale=0.6]{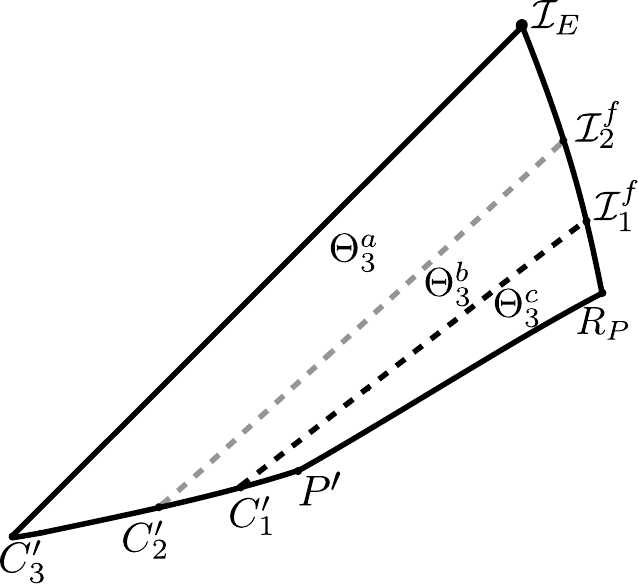}}  
	 \hspace{6mm}
	\subfigure[R-region $\Theta_2$
    subdivided by the secondary bifurcation segment $C_2$-$\mathcal{I}_2^f$ 
    and the $f$-rarefaction segment $C_1$-$\mathcal{I}_1^f$.
    ]
	{\includegraphics[scale=0.6]{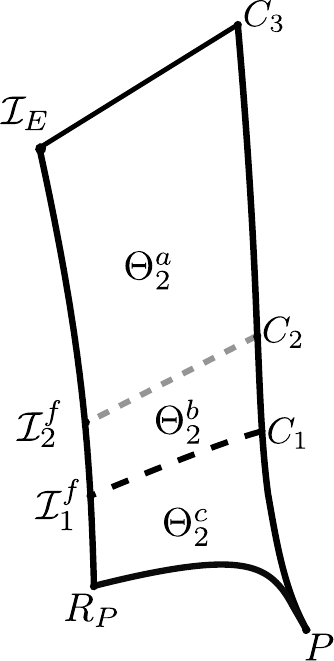}}   \hspace{6mm}
\subfigure[R-region $\Theta_1$
subdivided by the mixed double contact locus segment $K_E$-$K$,
    the secondary bifurcation segment $\um$-$C_2$ and the $f$-rarefaction segment
    $K$-$C_1$.
    ]
	{\includegraphics[scale=0.57]{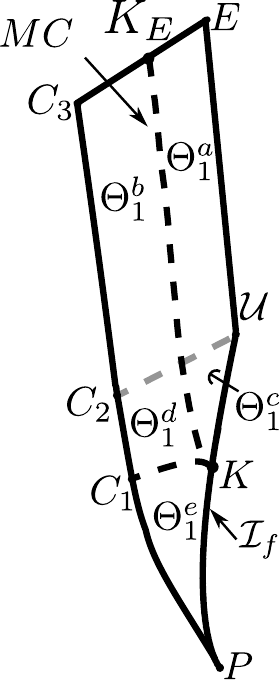}}  
	\caption{$R$-regions $\Theta_1$, $\Theta_2$, $\Theta_3$
 in Fig.~\ref{fig:RegioesThetas}
	displaying in detail their subdivisions.
	} 
	\label{fig:RRegions Completed2}
\end{figure}

\subsection{Right state in region \texorpdfstring{$\Theta$}{O}}

\label{subsec:RinTheta_1}
We divide region $\Theta$ into subregions $\Theta_1,$ $\Theta_2$, and $\Theta_3$; see Fig~\ref{fig:RegioesThetas}(a). These subregions are defined below. 
\subsubsection{Subregion \texorpdfstring{$\Theta_1$}{O}}
The boundary of subregion $\Theta_1$ in Fig~\ref{fig:RegioesThetas}(a)
consists of the following segments: the invariant segment \EUc;
the segment $\um$-$P$ of the $f$-inflection locus;
the segment $P$-$C_3$ of the double contact locus;
and the segment $C_3$-$E$ of the edge \GOc.

The region $\Theta_1$ itself is subdivided into five subregions: $\Theta_1^i$, with $i\in \{a,b,c,d,e\}$, see Fig.~\ref{fig:RRegions Completed2}(c), 
whose boundaries are detailed below.
This subdivision of $\Theta_1$ is related to
changes in the structure of $\wm_f(R)$ that will be described below.

Following the strategy alluded to at the end of Subsection~\ref{sec:notation}, we next describe the backward Hugoniot loci, the backward fast-family wave curves, and Riemann solutions for
representative right states $R$ in regions $\Theta_1$,
$\Theta_2$, and $\Theta_3$.

\medskip

\noindent{\bf The backward Hugoniot curve and the admissible shocks for $R$ in subregion $\Theta_1^a$.}
\medskip

\begin{figure}[]
	\centering
	\subfigure[$R = (0.0402518, 0.913397)$ in $\Theta_1^a$.]
	{\includegraphics[width=0.375\linewidth]{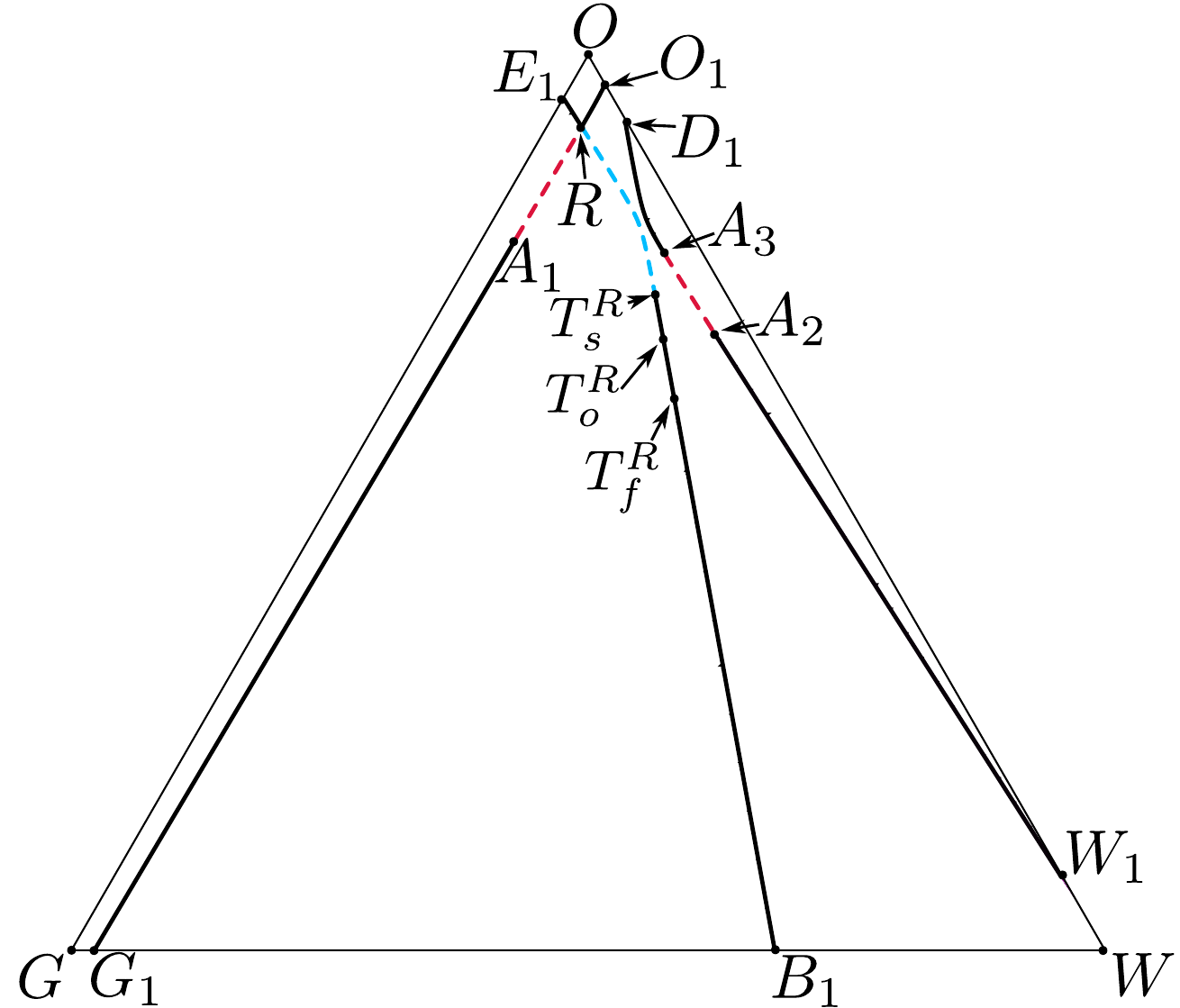}}
	\hspace{1.2mm}
	\subfigure[$R=(0.0414089, 0.90822)$ in $\Theta_1^b$.]
	{\includegraphics[width=0.375\linewidth]{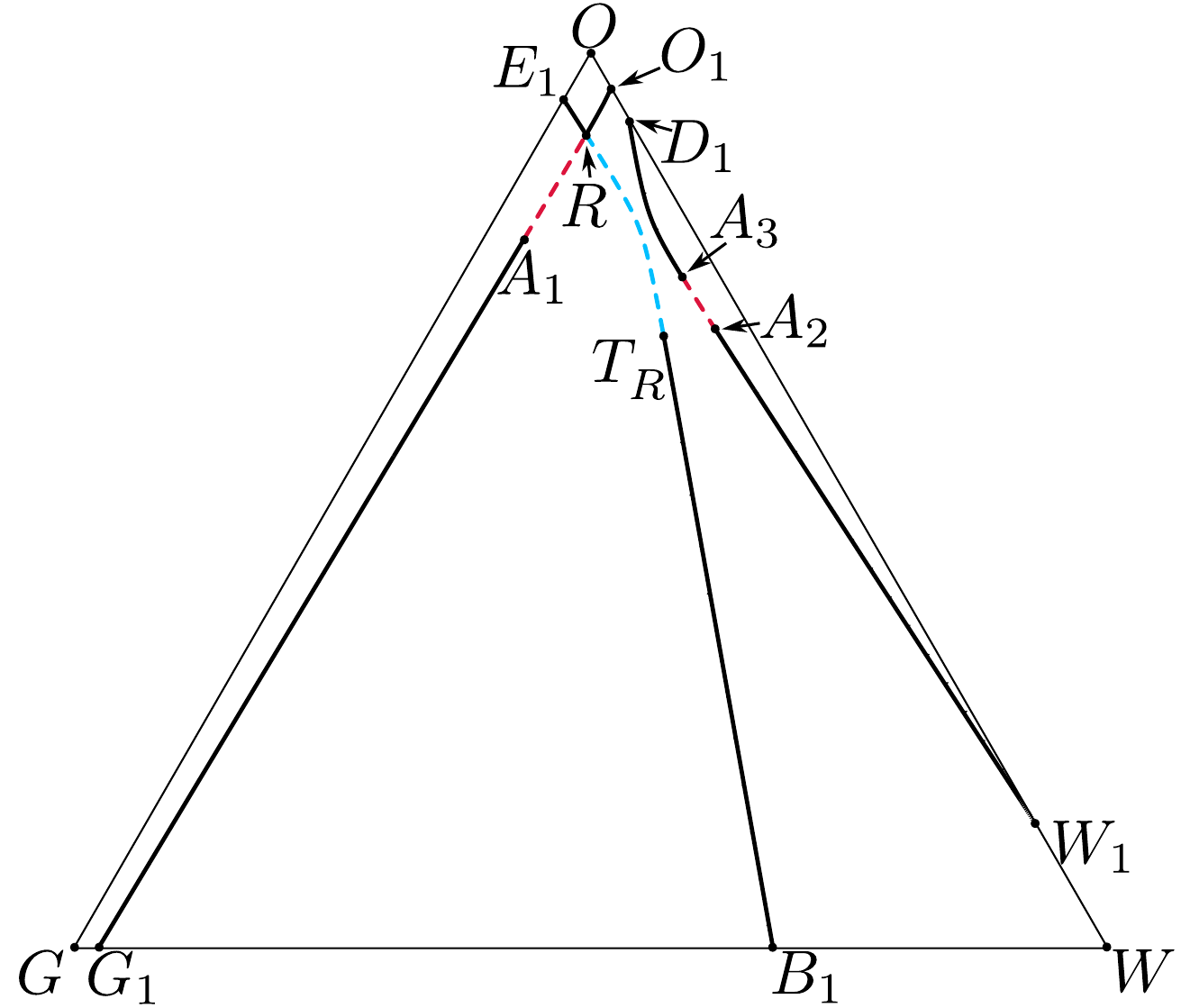}}
	\subfigure[$R=(0.0955103,0.85591)$ in $\Theta_1^c$.]
	{\includegraphics[width=0.375\linewidth]{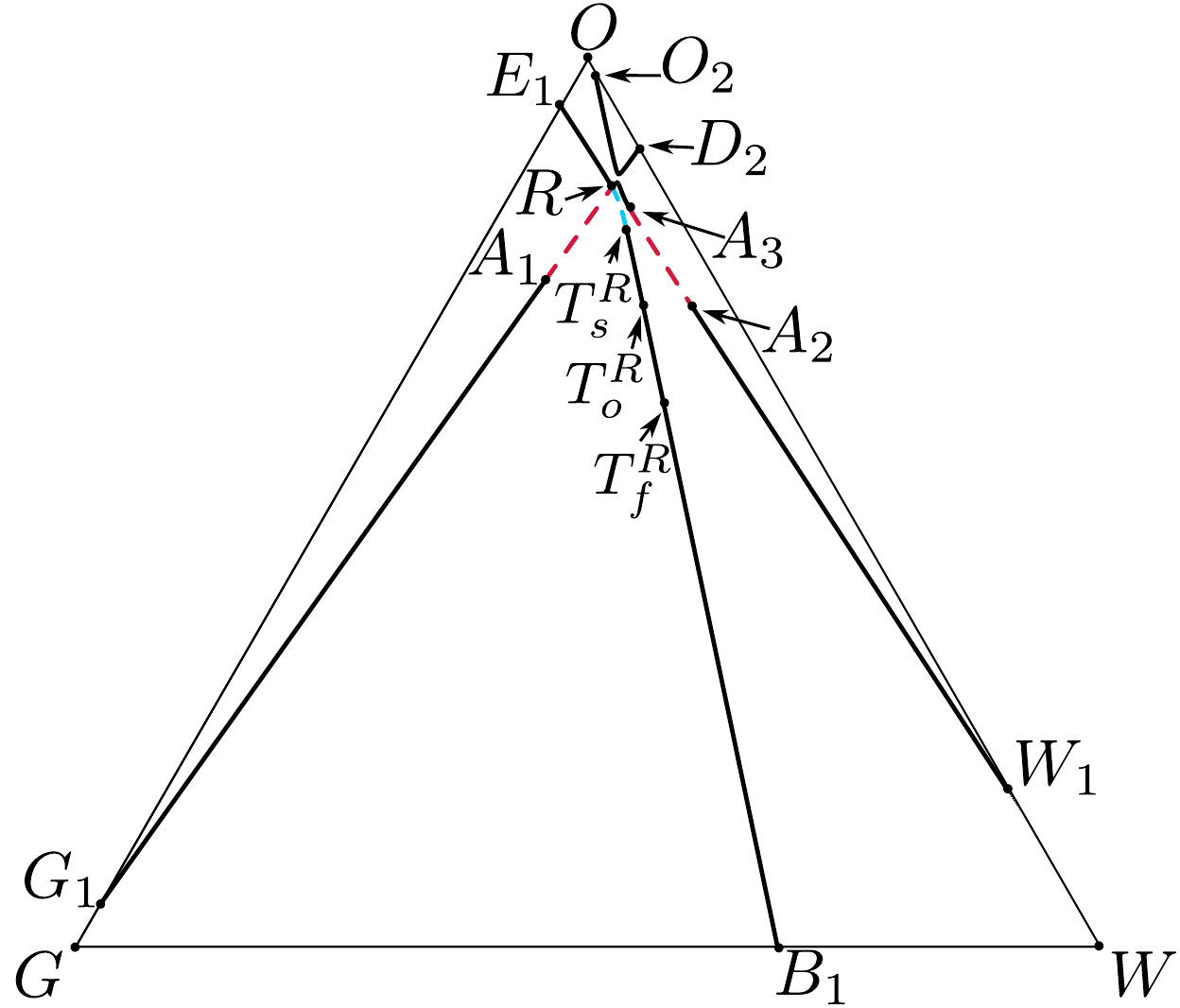}}
	\hspace{1.2mm}
		\subfigure[$R =(0.11209, 0.8299)$ in $\Theta_1^d \cup \Theta_1^e$.]
	{\includegraphics[width=0.375\linewidth]{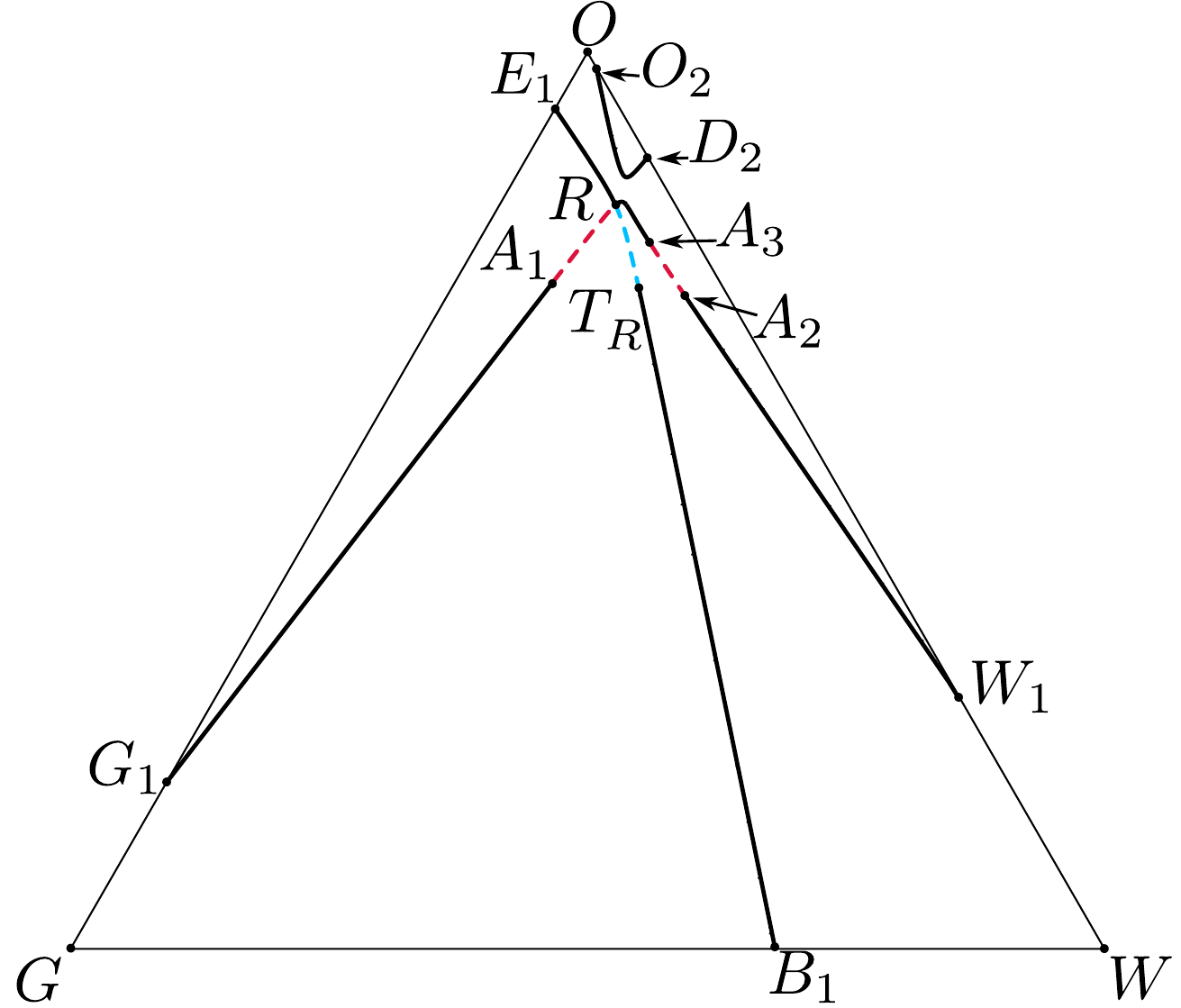}}
	\caption{
	Backward Hugoniot curves for $R$ in subregions of $\Theta_1$, Figs.~\ref{fig:RegioesThetas}(a) and \ref{fig:RRegions Completed2}(a).
	Admissible slow-family shock segments:
	$(R, T_s^R]$ in (a) and (c);  $(R, T_R]$ in (b) and (d).
	Admissible fast-family shock segments:
	$[A_1, R)$ and $[A_2, A_3]$.
    We have: $\sigma(A_1; R) = \lambdaf(A_1)$,
    $\sigma(A_2; R) = \lambdaf(A_2)$,
	$\sigma(A_3; R) = \lambdaf(R)$,
	$\sigma(T_s^R; R) = \lambdaf(R)$,
	$\sigma(T_o^R; R) = \lambdas(T_o^R)$,
	$\sigma(T_f^R; R) = \lambdas(R)$,
	$\sigma(T_R; R) = \lambdaf(R)$.
	The overcompressive shock segment $(T_s^R, T_o^R)$ and 
	the fast-family shock segment $(T_o^R, T_f^R)$ in (a) and (c) are not admissible.
	}
	\label{fig:Hug-R10a}
\end{figure}

A typical backward Hugoniot curve for a generic state $R$ in subregion $\Theta_1^a$ is depicted in Fig.~\ref{fig:Hug-R10a}(a).  
In this case, $\mathcal{H}(R)$ possesses the primary branches (which intersect at $R$) $[G_1, R)$, $(R, O_1]$, $[E_1, R)$,
$(R, B_1]$, and
the detached (nonlocal) branch $[D_1, W_1]$.
For $R$ near \EUc, the Hugoniot curve $\mathcal{H}(R)$ can be
considered as a perturbation of a Hugoniot curve
for $R$ along \EUc, which was obtained explicitly in \cite{Azevedo2014}.
We verify numerically that the $s$-shock segment
$(R, T_s^R]$ in the primary branch $(R, B_1]$, the $f$-shock segment
$[A_1, R)$ in the primary branch $[G_1, R)$, and 
the $f$-shock segment
$[A_2, A_3]$ in the detached branch $[D_1, W_1]$
satisfy the
viscous profile admissibility criterion. The endpoints $A_1$ and $A_2$ of these
segments are Bethe-Wendroff points \cite{Wendroff1972} characterized by the equalities
$\sigma(A_1; R) = \lambdaf(A_1)$ and
$\sigma(T_s^R; R) = \lambdaf(R)$ and $\sigma(A_3;R) = \lambdaf(R)$.
$\sigma(A_2; R) = \lambdaf(A_2)$. The endpoints $A_3$ and $T_s^R$ satisfy the equalities
$\sigma(T_s^R; R) = \lambdaf(R)$ and $\sigma(A_3;R) = \lambdaf(R)$.
\begin{rem}\label{rem:Tstar}
As shown in Fig.~\ref{fig:Hug-R10a}(a),
the backward Hugoniot curve for $R$ in subregion $\Theta_1^a$ also possesses 
the overcompressive shock segment $(T_s^R, T_o^R)$ 
and the $f$-shock segment $(T_o^R, T_f^R)$ along the local branch $(R, B_1]$ such that
$\sigma(T_s^R; R) = \lambdaf(R)$, $\lambdas(T_o^R) = \sigma(T_o^R; R) = \lambdaf(R)$;
however, such segments correspond to non-admissible shock waves.
\end{rem}
\noindent{\bf The backward Hugoniot curve and the admissible shocks for $R$ in subregion
$\Theta_1^b$.} 
\medskip

When the state $R$ lies on the portion of the mixed contact segment $K$-$K_E$ in Figs.~\ref{fig:BoundaryExt_MixedDoubleCont}(a) 
and \ref{fig:RRegions Completed2}(c) that separates $\Theta_1^a$
from $\Theta_1^b$, the states $T_s^R$, $T_o^R$ and $T_f^R$ 
of $\mathcal{H}(R)$ (referred to in Remark~\ref{rem:Tstar}) 
collapse to a single state, denoted by $T_R$. 
Thus, for $R$  in subregion $\Theta_1^b$, see Fig.~\ref{fig:Hug-R10a}(b),
the admissible shock segments are the
$s$-shock segment $(R, T_R]$ and the $f$-shock segments $[A_1, R)$ and 
$[A_2, A_3]$, as for $R$ in $\Theta_1^a$.
\begin{rem}\label{rem:A2A4}
The segment $[A_2, A_3]$ of $\mathcal{H}(R)$ in
Fig.~\ref{fig:Hug-R10a} lies below the invariant line \WUc,
inside the triangle $BW\um$ in Fig.~\ref{fig:SecBifurc}.
\end{rem}
%

\medskip

\noindent{\bf The backward Hugoniot curve and the admissible shocks for $R$ in subregion
$\Theta_1^c$.}
\medskip

The state $R$ crosses the segment $\um$-$C_2$ of the secondary
bifurcation locus as it moves from the subregion $\Theta_1^a$ to the
subregion $\Theta_1^c$. Consequently, the topological structure of the backward
Hugoniot curve $\mathcal{H}(R)$ changes, as seen by comparing
Figs.~\ref{fig:Hug-R10a}(a) and (c). Now, the attached branches are $[G_1, R)$, $(R, W_1]$, $[E_1, R)$, $(R, B_1]$, and the detached branch is $[O_2, D_2]$.

As in the case of $R$ in subregion $\Theta_1^a$, the shock segments
$[A_1, R)$, $(R, T_s^R]$, and $[A_2, A_3]$ are admissible, but
the segments $(T_s^R, T_o^R)$ and $(T_o^R, T_f^R)$ are not.
See Remark~\ref{rem:Tstar}.

\medskip
\noindent{\bf The backward Hugoniot curve and the admissible shocks for $R$ in subregion $\Theta_1^d\cup \Theta_1^e$.}
\medskip

Similar to the case when $R$ crosses from $\Theta_1^a$ to $\Theta_1^b$, when $R$ crosses the remaining portion of the segment $K$-$K_E$ of the mixed contact
locus from $\Theta_1^c$ to $\Theta_1^d \cup \Theta_1^e$, the shock segments $(T_s^R, T_o^R)$ and $(T_o^R, T_f^R)$ in $\mathcal{H}(R)$ 
collapse to the point $T_R$, yielding the admissible $s$-shock segment $(R, T_R)$; see Fig.~\ref{fig:Hug-R10a}(d).



Next, we describe the backward fast wave curves for $R$ in $\Theta_1$.

\medskip

\noindent{\bf The backward fast wave curve $\wm_f(R)$ for $R$ in subregion $\Theta_1^a \cup \Theta_1^b$.}

\medskip

\begin{figure}[ht]
	\centering
	\subfigure[]
	{\includegraphics[width=0.4\linewidth]{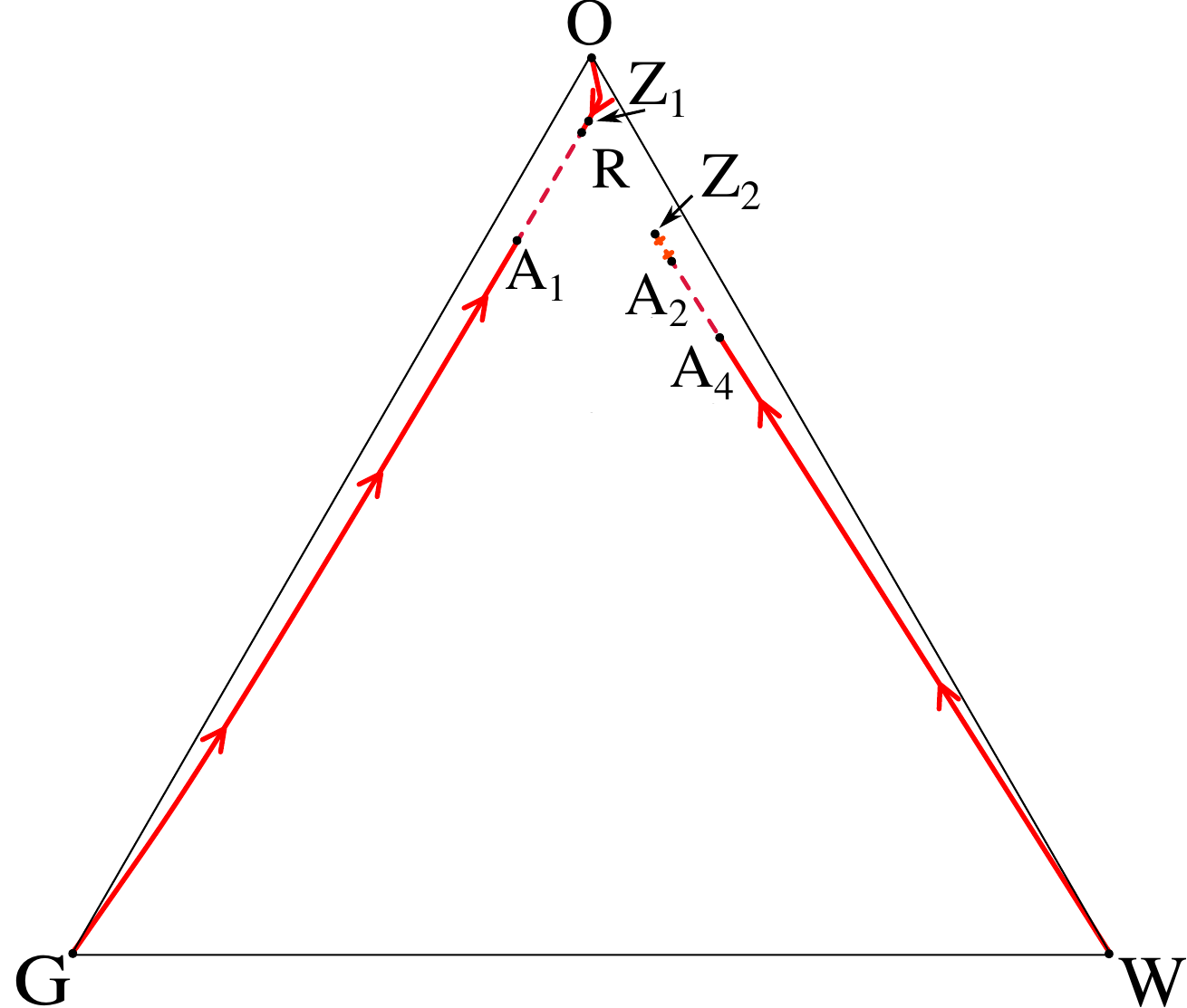}}  
	\hspace{1.75mm}
	\subfigure[]
	{\includegraphics[width=0.4\linewidth]{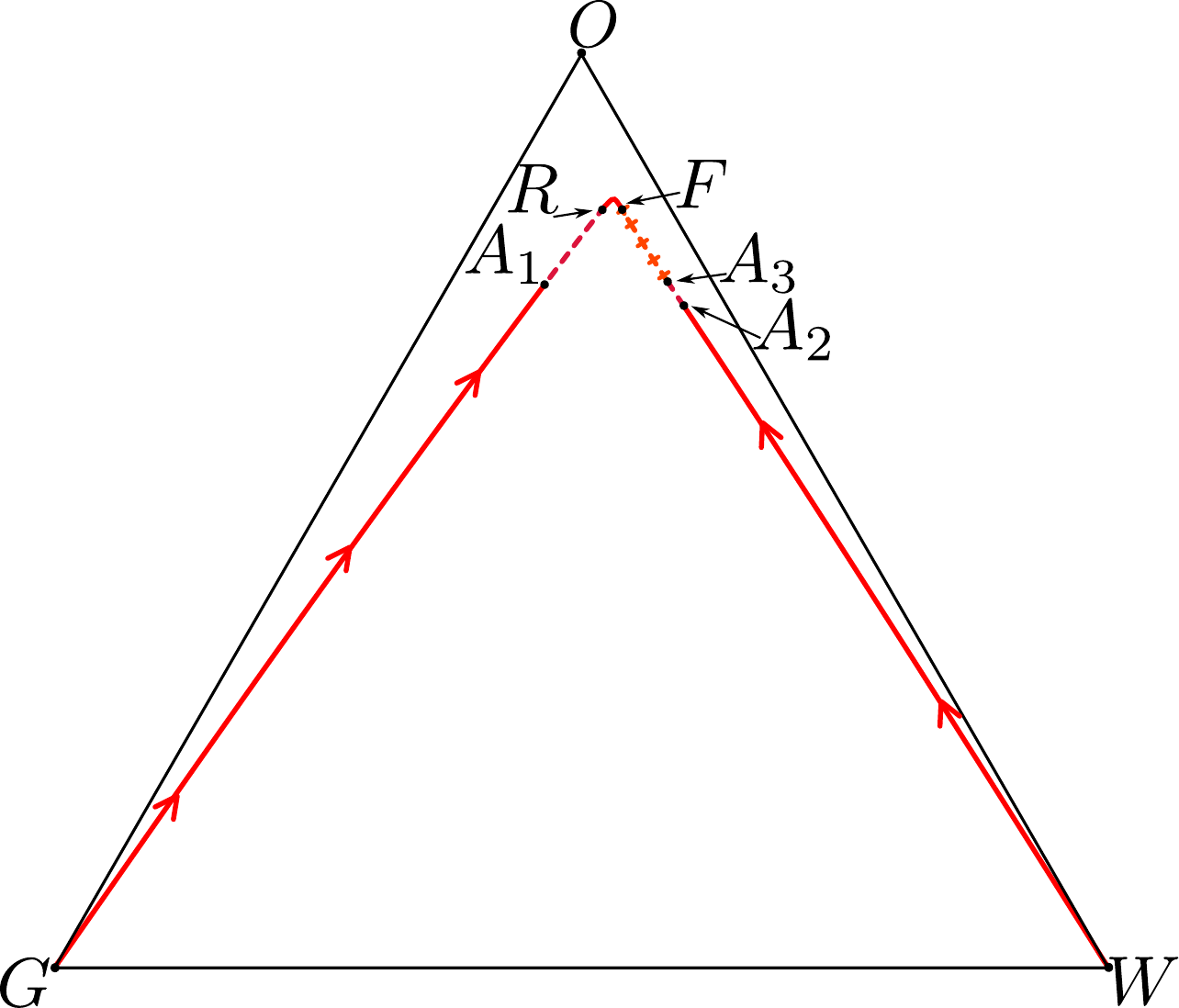}}  
	\caption{
	Backward fast-family wave curve of $R$ for $\Theta_1$.
	(a) $R = (0.0402518, 0.913397)$ in subregion $\Theta_1^a\cup \Theta_1^b$ in Fig.~\ref{fig:RRegions Completed2}(b). The
	rarefaction segment $(R, O]$ crosses the invariant
	segment \EU at a state $Z_1$.
	The state $Z_2$ lies on the invariant segment
	\UWc, with $\sigma(Z_2; Z_1) = \lambdaf(Z_1)$.
	(b) $R=(0.0955103,0.85591)$ in subregion 
	$\Theta_1^c \cup \Theta_1^d$ in Fig.~\ref{fig:RRegions Completed2}(c).
	The state $F$ lies on the $f$-inflection segment $\um$-$P$
 in Fig.~\ref{fig:IntegralCurves}.
	}
	\label{fig:BFWC-R10}
\end{figure}

Refer to Fig.~\ref{fig:BFWC-R10}(a). Notice that 
the $f$-rarefaction curve for a state $R$ in $\Theta_1^a \cup \Theta_1^b$ crosses the invariant
segment \EU at a state $Z_1$ before it reaches the vertex $O$ of
the saturation triangle. 
Corresponding to  $Z_1$ there is a state $Z_2$ on the invariant segment \UW 
such that $\sigma(Z_2; Z_1) = \lambdaf(Z_1)$.
For 
$R$ in $\Theta_1^b$  it also crosses the mixed double
contact segment $MC$
at a state $B_1^*$, which is not shown in Fig.~\ref{fig:BFWC-R10}(a).

\begin{cla}\label{cla:RinR10}
For $R$ in subregion $\Theta_1^a\cup \Theta_1^b$ in Fig.~\ref{fig:RRegions Completed2}(c)
Refer to Fig.~\ref{fig:RRegions Completed} 
$\wm_f(R)$ comprises states $M$ along the
shock segments $[A_1, R)$ and $[A_2, A_3]$ in $\mathcal{H}(R)$,  states $M$ along the $f$-rarefaction
segments $[G, A_1)$, $[O, R)$, $[W, A_2)$,
and states $M'$ along the composite segment
$[Z_2, A_3)$ corresponding to the rarefaction segment $[Z_1, R)$.
\end{cla}
\medskip
\noindent{\bf The backward fast wave curve $\wm_f(R)$ for $R$ in subregion $\Theta_1^c \cup \Theta_1^d  \cup \Theta_1^e$.}

\medskip

Refer to Fig.~\ref{fig:BFWC-R10}(b). The main difference from the previous case is that now
the $f$-rarefaction segment from $R$ does not cross the invariant segment
\EU in Fig.~\ref{fig:RRegions Completed}, but reaches the $f$-inflection segment
$\um$-$P$ at a state denoted by $F$; if $R$ lies in $\Theta_1^c \cup \Theta_1^d$,
	the state $F$ lies in
	$\um $-$ K$, and if $R$ lies in $\Theta_1^e$, $F$ lies in $K $-$ P$. See Fig.~\ref{fig:RRegions Completed2}(c).
Moreover, the $f$-shock segment $[A_2, A_3]$
of $\mathcal{H}(R)$ is such that $A_2$ lies below the  segment $P$-$\mathcal{I}_D$ of the $f$-inflection locus
in Figs.~\ref{fig:IntegralCurves}(a) and \ref{fig:RRegions Completed}.

\begin{cla}\label{cla:BWCR1G}
	For $R$ in subregion $\Theta_1^c \cup \Theta_1^d \cup \Theta_1^e$ in Fig.~\ref{fig:RRegions Completed2}(c)
	$\wm_f(R)$ comprises states
	$M$ along the shock segments $[A_1, R)$ and $[A_2, A_3]$
	in the Hugoniot curve $\mathcal{H}(R)$,
	states $M$ along the $f$-rarefaction segments
	$[G, A_1)$, $[W, A_2)$, and $[F, R)$,
	and states $M'$ along the composite segment
	$(F, A_3)$ corresponding to $(F, A_3)$.
\end{cla}

Concerning $R$ states in region $\Theta_1$, we finally present the Riemann solutions.
We consider three generic right states
and for each one we consider all left states $L$ in the edge \GW of
the saturation triangle, simulating the injection of gas-water mixtures into the reservoir.


\medskip
\noindent{\bf The Riemann solution for $R$ in subregion $\Theta_1^a \cup \Theta_1^c$.}

\medskip


For the sake of clarity we will sometimes denote by $\mathcal{L}_{\text{ext}}$ any pertinent extension of a given locus $\mathcal{L}$.

 Recall the state $T_o^R$ on $\mathcal{H}(R)$ satisfying $\lambdas(T_o^R) = \sigma(T_o^R; R)$;
 see Figs.~\ref{fig:Hug-R10a}(a), (c) and Remark~\ref{rem:Tstar}.
If we consider shock speeds for states $M$ along the segments $[A_1, R)$ and $[A_2, A_3]$
of $\mathcal{H}(R)$,
we find that there is a state $M=A_1^*$ in segment $[A_1, R)$ 
and a state $M=A_2^*$ in segment $[A_2, A_3]$ such that 
$\sigma(A_1^*;R) = \sigma(A_2^*;R) = \sigma(T_o^R; R)$. Thus by the triple shock rule, it follows that
$\sigma(A_2^*; A_1^*)  = 
\sigma(A_1^*;R) = \sigma(A_2^*;R) = \sigma(T_o^R; R)$.
Because of such equalities, the states $A_1^*$, $A_2^*$ play a fundamental role in the
analysis of the speed compatibility between slow and fast wave groups
connecting left states $L$ in the edge \GW to the right state $R$, as we will see next.

\begin{figure}[]
\centering
\subfigure[$R$ in subregion of $\Theta_1^a$.
]
{\includegraphics[width=0.4\linewidth]{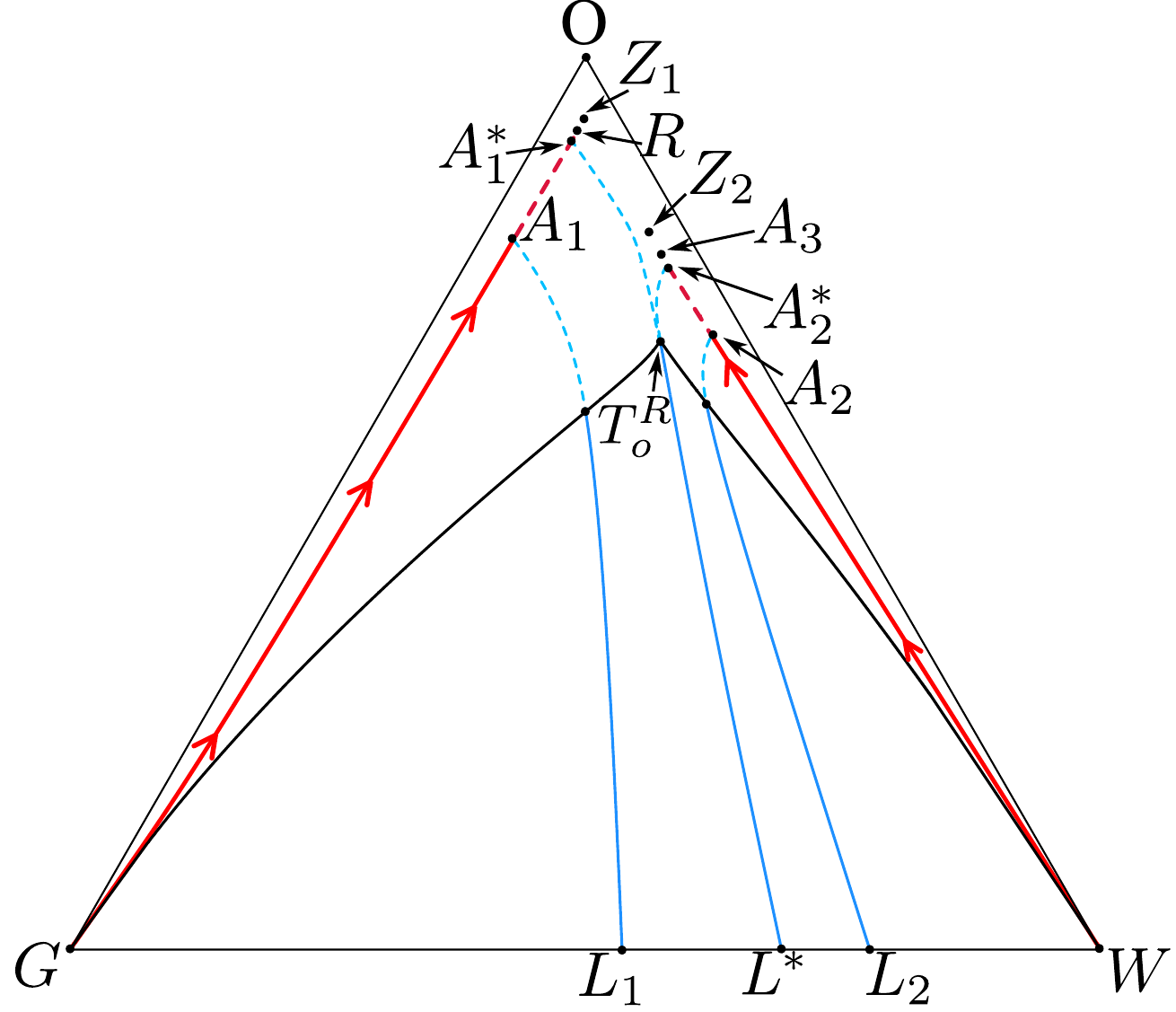}}
\hspace{1.75mm}
\subfigure[$R$ in subregion of $\Theta_1^c$.
]
{\includegraphics[width=0.4\linewidth]{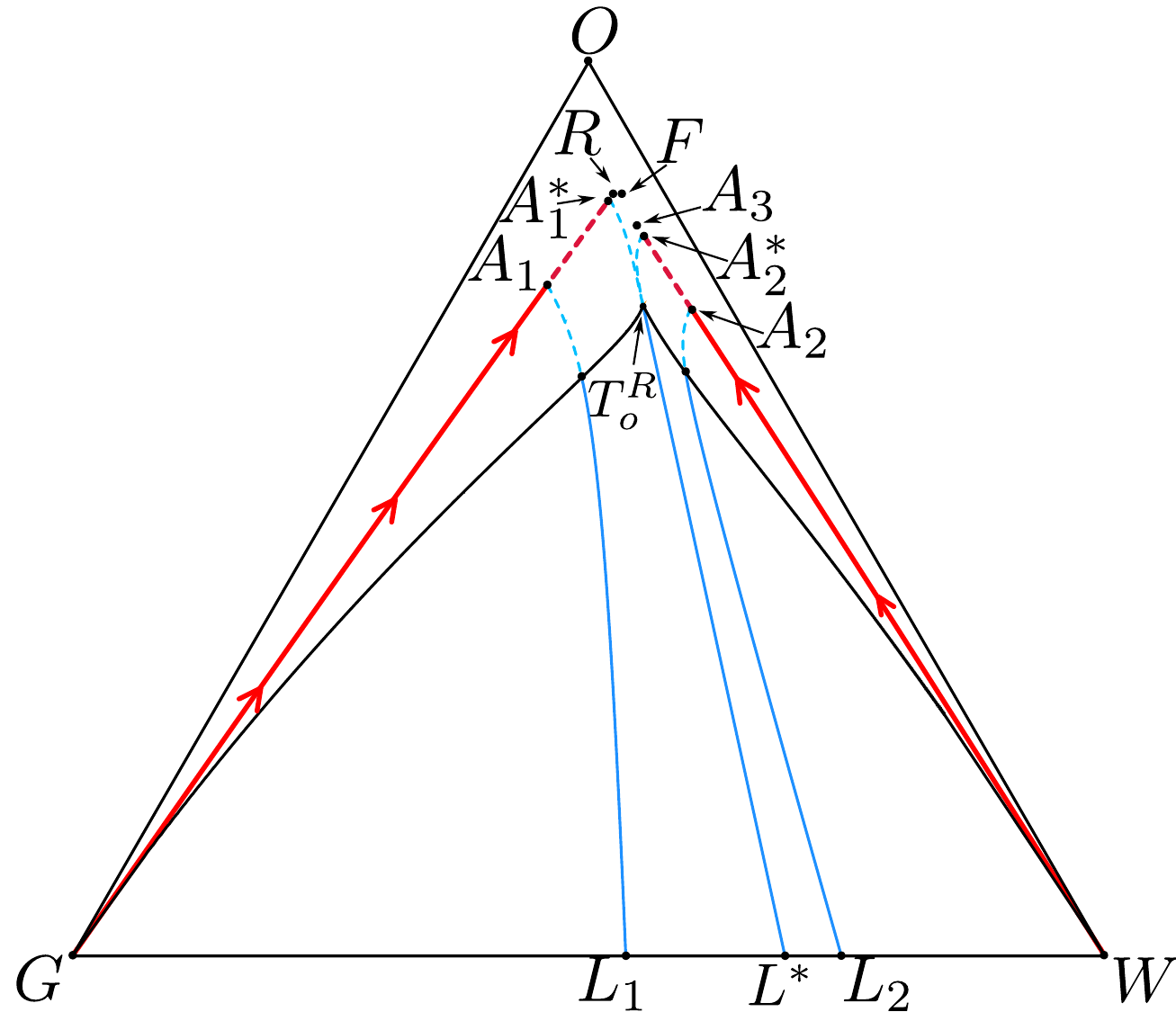}}
\caption{Riemann solution for $R$ in  the subregion $\Theta_1^a\cup \Theta_1^c$ with $L$ along the edge \GWc.
The segments $(G, T_o^R]_{\text{ext}}$ and $(W, T_o^R]_{\text{ext}}$ are backward slow-extensions of segments $(G, A_1^*]$ and $(W, A_2^*]$
on $\wm_f(R)$, with
$\sigma(M'; M) = \lambdas(M')$, for all $M'$ in
$[G, T_R]_{\text{ext}} \cup [W, T_R]_{\text{ext}}$ and all $M$ in $[G, A_1^*] \cup [W, A_2^*]$.
}
\label{fig:RSRinR10}
\end{figure}

For each $M$ along the portions $[O, A_1^*)$ and $[Z_2, A_2^*)$
of $\wm_f(R)$ in Claim~\ref{cla:RinR10},
the beginning speed $v_{fb}$ of the forward fast wave group
$M\xrightarrow{\mathcal{G}_f} R$ 
is smaller than  the ending speed $v_{se}$ of the forward slow wave group $L\xrightarrow{\mathcal{G}_s} M$,
implying the incompatibility of such a sequence of wave groups for constructing a Riemann solution.
Thus, due to the monotonicity of shock speeds along the segments $(A_1, R)$ and $(A_2, A_3)$,
only the portions $[G,A_1^*]$ and $[W, A_2^*]$
of $\wm_f(R)$ and their backward $s$-extension segments $[G, T_o^R]_{\text{ext}}$ and $[T_o^R, W]_{\text{ext}}$, characteristic on these segments,
see Fig.~\ref{fig:RSRinR10},
are used in the construction of
Riemann solutions for $L$ in the edge \GWc.

\begin{cla}\label{cla:RSolution-RinR10a}
Refer to Fig.~\ref{fig:RSRinR10}(a) for $R$
in $\Theta_1^a$ and to 
Fig.~\ref{fig:RSRinR10}(b) for $R$
in $\Theta_1^c$. Let $L$ be a state on the edge
\GW of the saturation triangle and $R$ be a state 
in subregion $\Theta_1^a \cup \Theta_1^c$ in Fig.~\ref{fig:RRegions Completed2}(c).
Let $T_o^R$ on $\mathcal{H}(R)$ and $A_1$,  $A_1^*$, $A_2^*$, $A_2$ on $\mathcal{W}_f(R)$ be the states satisfying $\sigma(A_1;R) = \lambdaf(A_1)$, $\sigma(A_1^*; R) = \sigma(A_2^*; R)
= \sigma(T_o^R; R) = \lambdas(T_o^R)$ and $\sigma(A_2;R) = \lambdaf(A_2)$.
Let $L_1$, $L^*$ and $L_2$ be the intersection points of the backward slow wave curves through $A_1$, $A_1^*$ (or $A_2^*$)
and $A_2$ with the edge \GWc, respectively. Then,

\begin{itemize}
\item[(i)] if $L=G$ or $L = W$, the Riemann solution is
$G\testright{R_f} A_1  \xrightarrow{'S_f} R \quad \mbox{ or }\quad W\testright{R_f} A_2 \xrightarrow{'S_f} R ;
$
\item[(ii)] if $L \in(G, L_1)$ or $L\in (W, L_2)$, the solution is
$
     L\testright{R_s} T_1 \xrightarrow{'S_s} M_1 \testright{R_f} A_1 \xrightarrow{'S_f} R$ or
     $L\testright{R_s} T_2 \xrightarrow{'S_s} M_2 \testright{R_f} A_2 \xrightarrow{'S_f} R,
$
where $T_1\in(G, T_R^o)_{\text{ext}}$, $M_1\in(G, A_1)$,  $T_2\in (W, T_R^o)_{\text{ext}}$ and $M_2\in (W,A_2)$. 
\item[(iii)] if $L \in[L_1, L^*]$ or $L \in [L_2, L^*]$
the Riemann solution is
$
     L\testright{R_s} T_1 \xrightarrow{'S_s} M_1 \testright{S_f} R
     $ or
     $L\testright{R_s} T_2 \xrightarrow{'S_s} M_2 \testright{S_f} R,
$
where $T_1\in(G, T_R^o)_{\text{ext}}$, $M_1\in(A_1, A_1^*)$,  $T_2\in (W, T_R^o)_{\text{ext}}$ and $M_2\in (A_2,A_2^*)$. 
\end{itemize}
\end{cla}
\begin{rem}\label{rem:L*}
If $L = L^*$, see Fig.~\ref{fig:RSRinR10}, there are two paths to reach the right state
$R$, namely 
$L^*\testright{R_s} T_o^R\xrightarrow{'S_s} A_1^* \xrightarrow{S_f} R$ or 
$L^*\testright{R_s} T_o^R \xrightarrow{'S_s} A_2^* \xrightarrow{S_f} R$
but the triple shock rule guarantees they represent the same solution in $xt$-space.
In other words, the solution consists of the single wave group
$L^*\testright{R_s} T_o^R \xrightarrow{'S_o} R$.
\end{rem}

\medskip
\noindent{\bf The Riemann solution for $R$ in subregion $\Theta_1^b \cup \Theta_1^d $.}

\medskip


Recall that for $R$ in the mixed contact segment
$MC$
separating region $\Theta_1^a \cup \Theta_1^c$ from
$\Theta_1^b \cup \Theta_1^d$, see Fig.~\ref{fig:RRegions Completed2}(c),
the states $T_s^R$, $T_o^R$ and $T_f^R$ of $\mathcal{H}(R)$ collapse to a single state, denoted by $T_R$,
with $\lambdas(T_R) = \sigma(T_R; R) = \lambdaf(R)$.
Consequently, relative to the previous case represented in Fig.~\ref{fig:RSRinR10}, the state $A_1^*$ collides with $R$
and state $A_2^*$ with $A_3$.
Thus, the state $T_R$ lies in the segment 
$K_E'$-$K'$ in Fig.~\ref{fig:BoundaryExt_MixedDoubleCont}(b).


When $R$ crosses the mixed contact segment $K_E$-$K$
in Fig.~\ref{fig:RRegions Completed2}(c) from subregion $\Theta_1^a\cup\Theta_1^c$ to the subregion $\Theta_1^b\cup\Theta_1^d$
a state $B_1^*$, defined as the intersection point
of the $f$-rarefaction segment
trough state $R$ 
with $MC$, appears. Thus, there is a correspondent state $T^*$ 
in $MC'$, in Fig.~\ref{fig:BoundaryExt_MixedDoubleCont}(b),
such that $\lambdas(T^*)= \sigma(T^*; B_1^*) = \lambdaf(B_1^*)$.
Moreover, there is a  
state $B_2^*$ in the composite segment $[Z_2, A_3]$
given in Claim~\ref{cla:RinR10} 
or in the composite segment $[F, A_3]$
in Claim~\ref{cla:BWCR1G} such that
$\sigma(B_2^*;B_1^*) = \lambdaf(B_1^*)$. Thus  
the states $T^*$, $B_1^*$ and $B_2^*$ 
satisfy the triple shock rule with
$\lambdas(T^*) = \sigma(T^*; B_1^*) = \lambdaf(B_1^*) = \sigma(T^*; B_2^*)
 = \sigma(B_2^*;B_1^*)$.
The states $B_1^*$ and $B_2^*$ have the same role as the states
$A_1^*$ and $A_2^*$ have for the speed compatibility
of slow and fast wave groups.

\begin{figure}[]
\centering
\subfigure[$R$ in subregion of $\Theta_1$ bounded by $\Theta_1^b$.
]
{\includegraphics[width=0.4\linewidth]{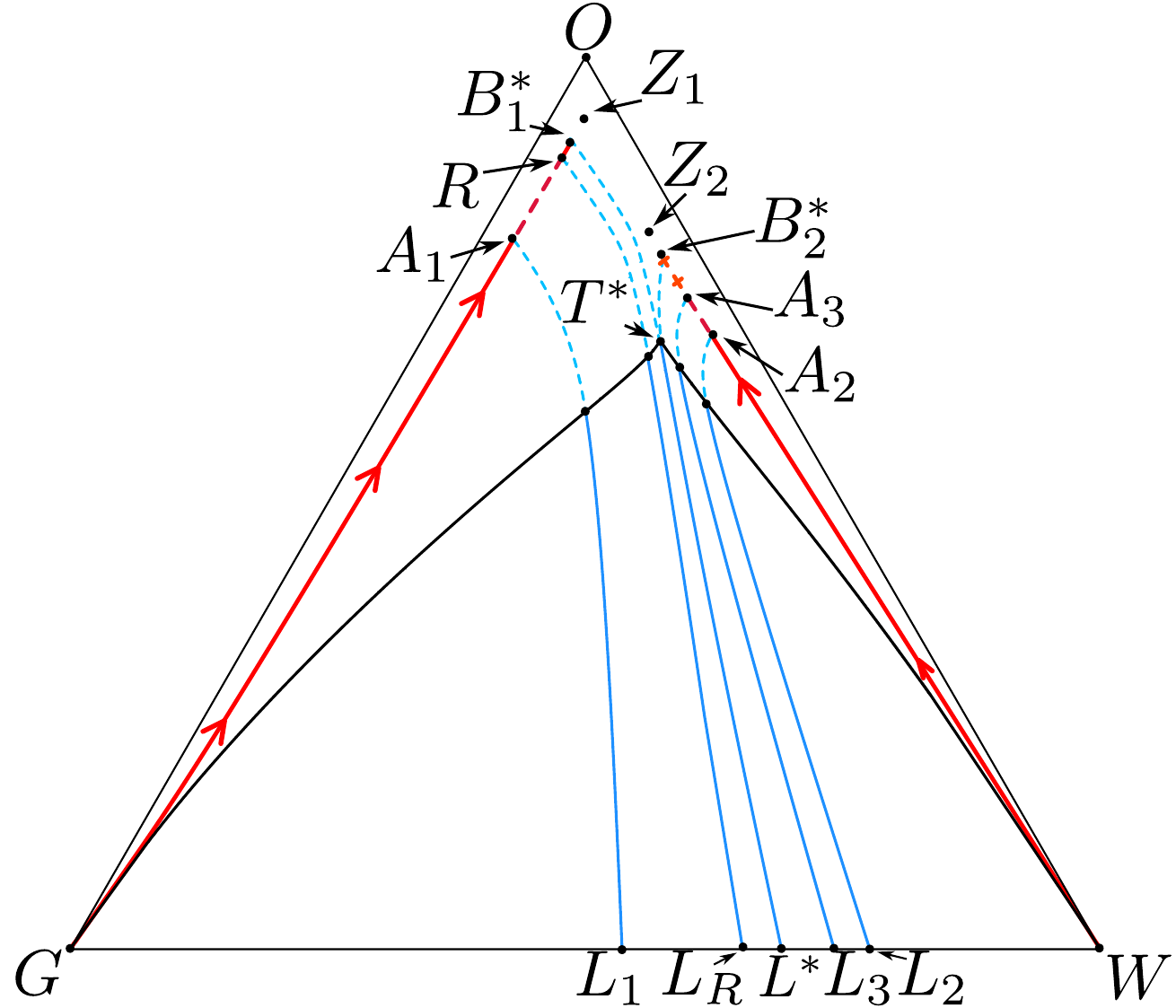}}
\hspace{1.75mm}
\subfigure[$R$ in subregion of $\Theta_1$ bounded by $\Theta_1^d$.
]
{\includegraphics[width=0.4\linewidth]{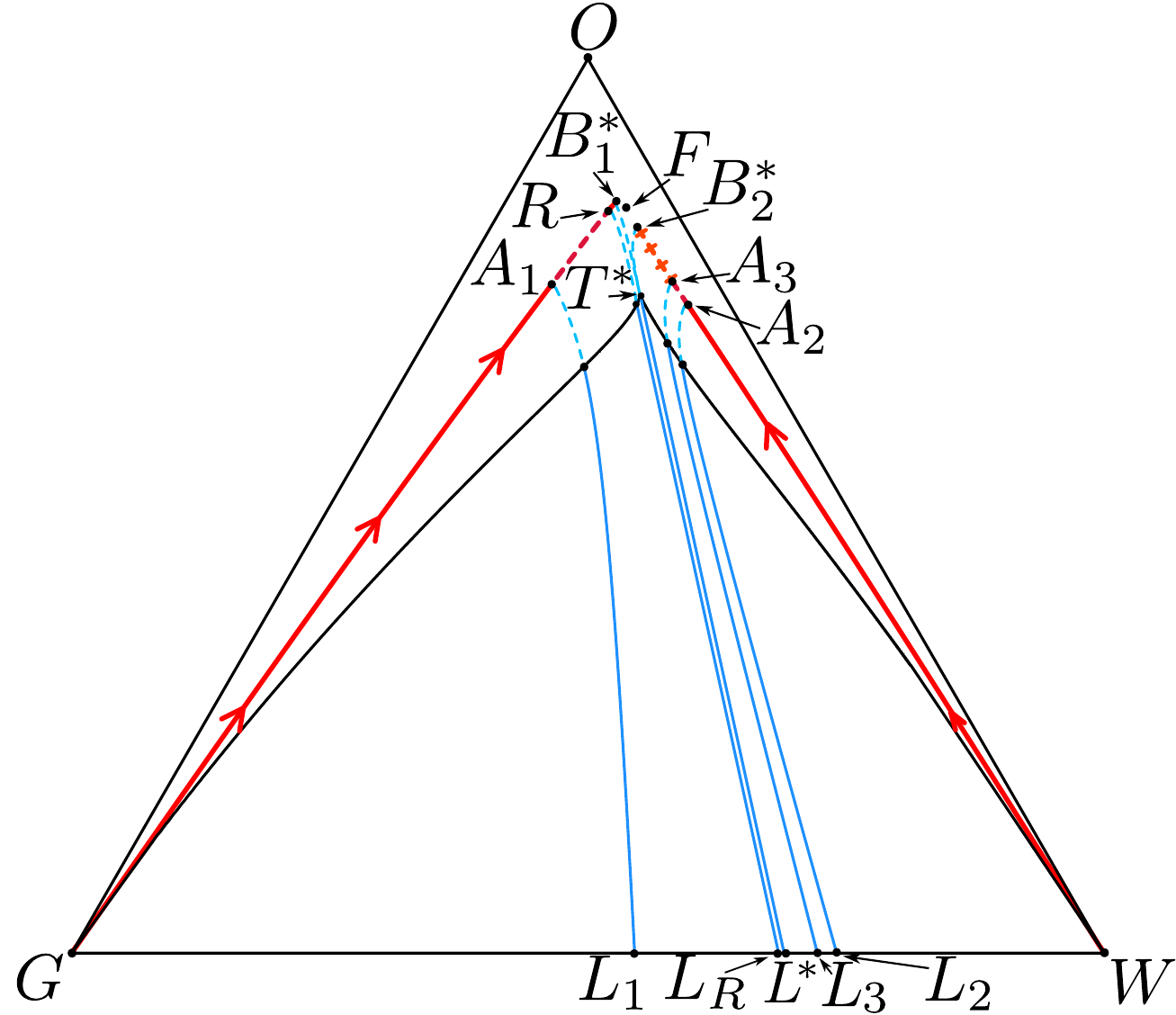}}
\caption{Riemann solutions  for $R$ in  the subregion of $\Theta_1$ with $L$ along the edge \GWc.
}
\label{fig:RSRinR10CandD}
\end{figure}

\begin{cla}\label{cla:RSolution-RinR10b}
Refer to Fig.~\ref{fig:RSRinR10CandD}(a) for $R$
in $\Theta_1^b$ and to 
Fig.~\ref{fig:RSRinR10CandD}(b) for $R$
in $\Theta_1^d$. Let $L$ be a state on the edge
\GW of the saturation triangle and $R$ be a state 
in subregion $\Theta_1^b \cup \Theta_1^d$ in Fig.~\ref{fig:RRegions Completed2}(c). 
Let $A_1$, $B_1^*$, $B_2^*$, $A_2$ and $A_3$ be the states on
$\wm_f(R)$ satisfying $\sigma(A_1;R) = \lambdaf(A_1)$,
$\lambda_s(T^*) = \sigma(T^*; B_1^*) = \lambdaf(B_1^*) = \sigma(T^*; B_2^*)
= \sigma(B_2^*; B_1^*)$, $\sigma(A_2;R) = \lambdaf(A_2)$
and $\sigma(A_3;R) = \lambdaf(R)$.
Let $L_1$, $L_R$, $L^*$, $L_3$ and $L_2$ be the intersection points of the backward slow wave curves through  states $A_1$, $R$, $B_1^*$ (or for $B_2^*$),
$A_3$ and $A_2$ with the edge \GWc, respectively. Then,

\begin{itemize}
\item[(i)]if $L=G$ or $L = W$, the Riemann solution is
$
    G\testright{R_f} A_1  \xrightarrow{'S_f} R$ or $W\testright{R_f} A_2 \xrightarrow{'S_f} R ;
$
\item[(ii)] if $L \in(G, L_1)$ or $L\in (W, L_2)$, the Riemann solution is 
$
     L\testright{R_s} T_1 \xrightarrow{'S_s} M_1 \testright{R_f} A_1 \xrightarrow{'S_f} R$ or  $L\testright{R_s} T_2 \xrightarrow{'S_s} M_2 \testright{R_f} A_2 \xrightarrow{'S_f} R,
$
where $T_1\in(G, T^*)_{\text{ext}}$, $M_1\in(G, A_1)$,  $T_2\in (W, T^*)_{\text{ext}}$ and $M_2\in (W,A_2)$;
\item[(iii)] if $L \in[L_1, L_R)$ or $L \in [L_2, L_3]$
the Riemann solution is
$
     L\testright{R_s} T_1 \xrightarrow{'S_s} M_1 \testright{S_f} R$ or $L\testright{R_s} T_2 \xrightarrow{'S_s} M_2 \testright{S_f} R,
$
where $T_1\in(G, T^*)_{\text{ext}}$, $M_1\in[A_1, R)$,  $T_2\in (W, T^*)_{\text{ext}}$ and $M_2\in [A_2, A_3]$. 
\item[(iv)] if $L \in (L_R, L^*)$, then 
the solution is
$
     L\testright{R_s} T_1 \xrightarrow{'S_s} M_1 \testright{R_f} R,
$
where $T_1\in(G, T^*)_{\text{ext}}$ and $M_1\in(B_1^*, R)$; 
\item[(v)] if $L \in (L_3, L^*)$, then 
the Riemann solution is 
$
     L\testright{R_s} T_1 \xrightarrow{'\!S_s} M_1' \testright{S_f\!\!\!'} M_1 \xrightarrow{R_f} R,
$
where $T_1\in(W, T^*)_{\text{ext}}$, $M_1'\in(B_2^*, A_3)$ and $M_1\in(B_1^*, R)$. 
\end{itemize}
\begin{rem}\label{rem:L*B*}
If $L = L^*$, see Fig.~\ref{fig:RSRinR10CandD}, there are two paths to reach the right state
$R$, namely 
$L^*\testright{R_s} T^*\xrightarrow{'S_o'} B_1^* \xrightarrow{R_f} R$ or 
$L^*\testright{R_s} T^* \xrightarrow{'S_s} B_2^* \xrightarrow{S_f'} B_1^* \xrightarrow{R_f} R$
but the triple shock rule guarantees they represent the same solution in $xt$-space.
In other words, the solution consists of the single wave group $L^*\testright{R_s} T^*\xrightarrow{'S_o'} B_1^* \xrightarrow{R_f} R$.
\end{rem}
\end{cla}

\medskip
\noindent{\bf The Riemann solution for $R$ in subregion $\Theta_1^e$.}

\begin{figure}
	\centering
	\subfigure[]
	{\includegraphics[width=0.4\linewidth]{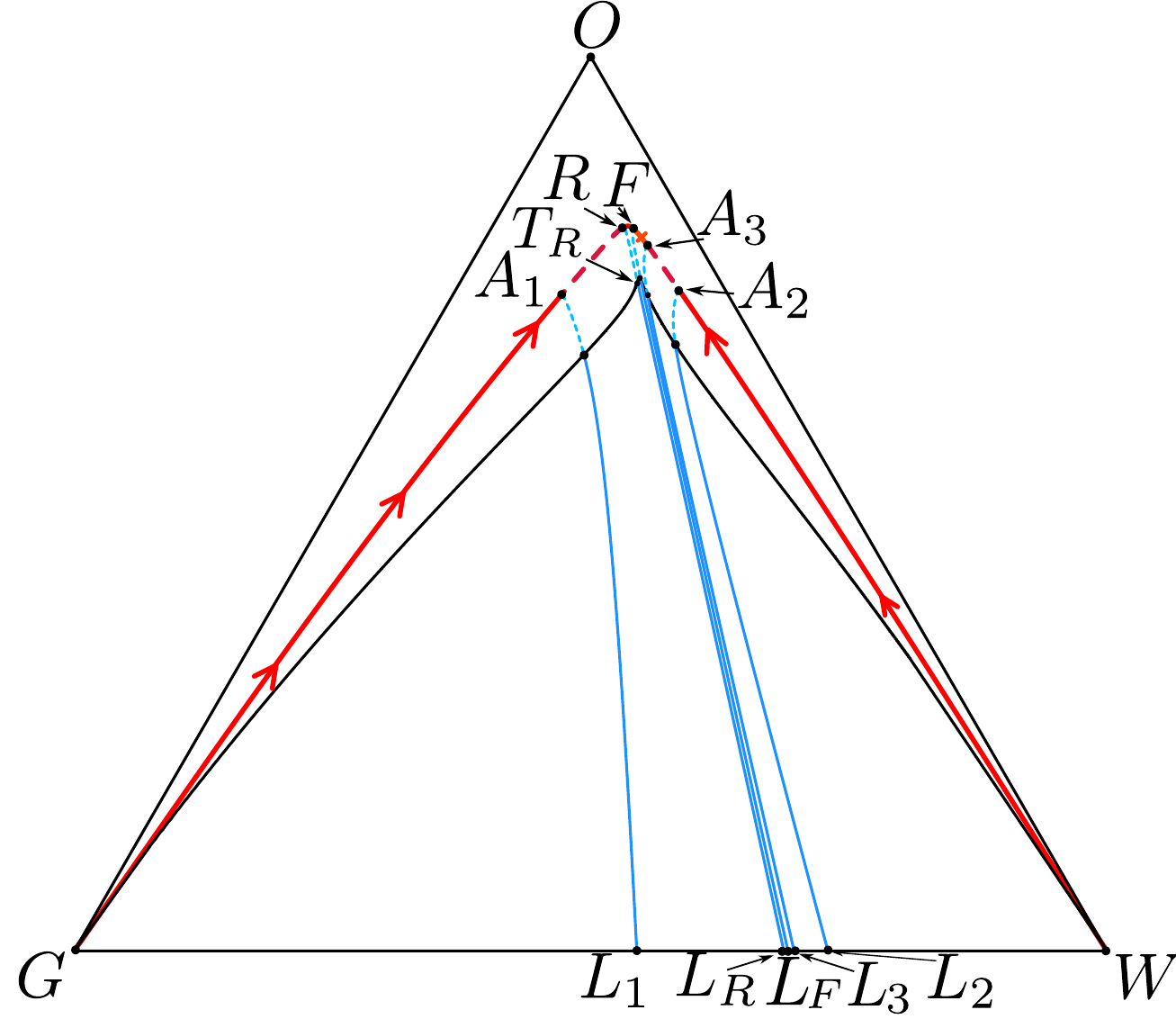}}   \hspace{1.5mm}	
 \subfigure[Zoom]
	{\includegraphics[width=0.4\linewidth]{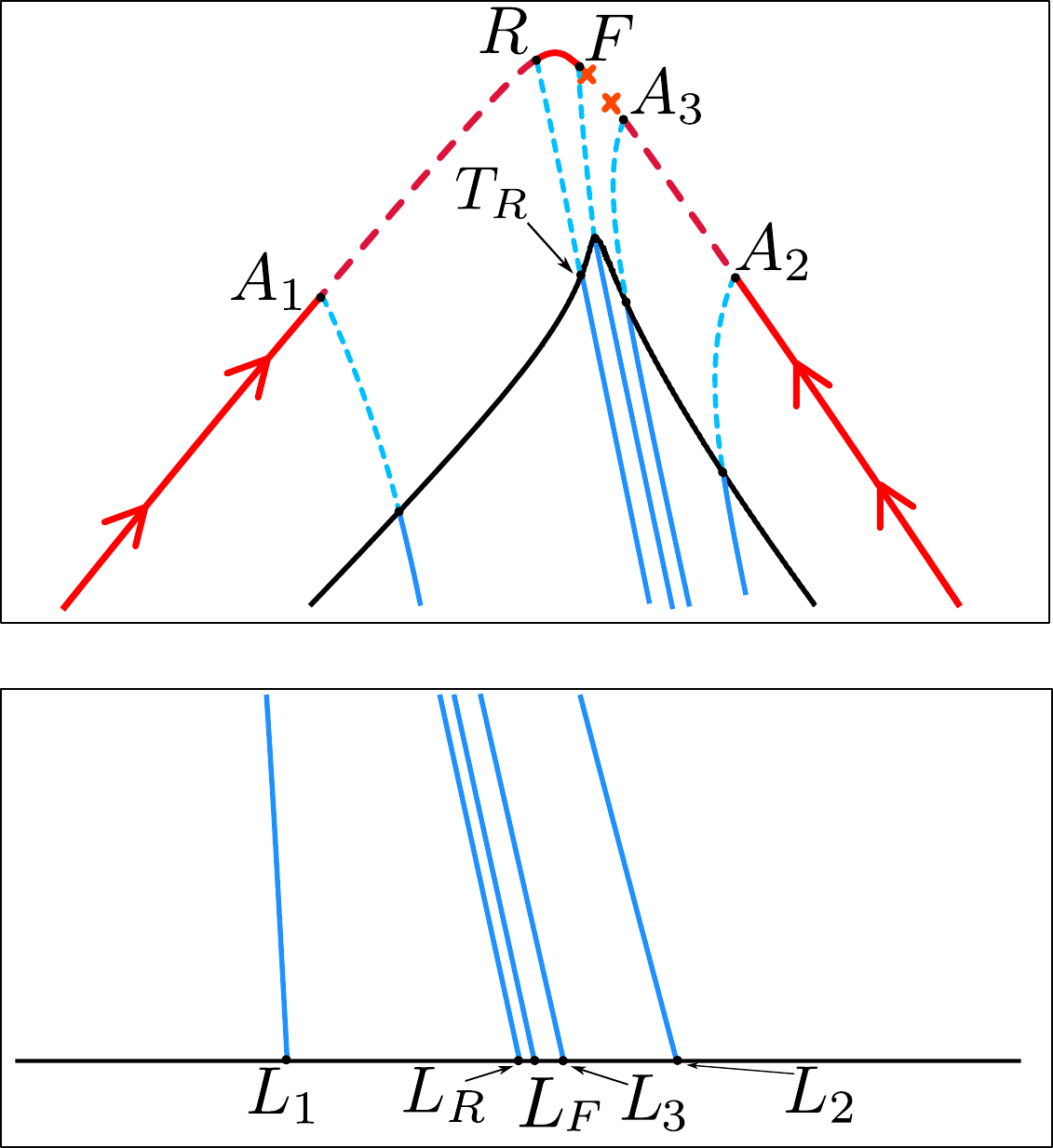}} 
	\caption{Riemann solutions for $R$ in  the subregion $\Theta_1^e$ in Fig.~\ref{fig:RRegions Completed}(b), with $L$ along the edge \GWc.
}
	\label{fig:RSRinR10E}
\end{figure}

When state $R$ is moved from subregion $\Theta_1^d$ onto the $f$-rarefaction segment $K$-$C_1$, the end-state $F$ of
the rarefaction segment $(R, F]$
coincides with $K$. Consequently, the state $B_1^*$,
on the mixed contact segment $K_E$-$K$-$K_D$ (see Fig.~\ref{fig:BoundaryExt_MixedDoubleCont}(a)), the state $B_2^*$ in the composite segment $(F, A_2)$ and state $F$ collide, reducing to $K$.
When state $R$ crosses the $f$-rarefaction segment $K$-$\mathcal{I}_1^f$,
keeping in region $\Theta_1$,
from the subregion $\Theta_1^d$ to $\Theta_1^e$
the state $F$ lies on the $f$-inflection segment $K$-$P$ and now the states $B_1^*$ 
and $B_2^*$ are not defined. Thus, the composite
segment $(F, A_3)$, along with
$\sigma(M';M) = \lambdaf(M)$ for all $M' \in (F, A_3)$
and $M\in (R, F)$, in Fig.~\ref{fig:BFWC-R10}(b), is
entirely used in the construction of the Riemann
solution  for $L$ along \GWc. We have:

\begin{cla}\label{cla:RSolution-RinR10b_casoE}
Refer to Fig.~\ref{fig:RSRinR10E}. 
Let $L$ be a state on the edge
\GW of the saturation triangle and $R$ be a 
state in subregion $\Theta_1^b$ in Fig.~\ref{fig:RRegions Completed2}(c).
Let $A_1$, $F$, $A_2$ and $A_3$ be states on the $\wm_f(R)$
satisfying $\sigma(A_1;R) = \lambdaf(A_1)$, $F$ on the $f$-inflection locus $\um$-$P$,
$\sigma(A_2;R) = \lambdaf(A_2)$
and $\sigma(A_3;R) = \lambdaf(R)$.
Let $L_1$,
$L_R$, $L_F$, $L_3$ and $L_2$ be the intersection points of the backward
slow wave curves through states $A_1$, $R$, $F$,
$A_3$ and $A_2$ with the edge \GWc, respectively. Then,

\begin{itemize}
\item[(i)]if $L=G$ or $L = W$, the Riemann solution is
$
     G\testright{R_f} A_1  \xrightarrow{'S_f} R$ or $W\testright{R_f} A_2 \xrightarrow{'S_f} R ;
$
\item[(ii)] if $L \in(G, L_1)$ or $L\in (L_2, W)$, the Riemann solution is
$
     L\testright{R_s} T_1 \xrightarrow{'S_s} M_1 \testright{R_f} A_1 \xrightarrow{'S_f} R$ or $L\testright{R_s} T_2 \xrightarrow{'S_s} M_2 \testright{R_f} A_2 \xrightarrow{'S_f} R,
$
where $T_1\in(G, T_R)_{\text{ext}}$, $M_1\in(G, A_1)$,  $T_2\in (W, T_R)_{\text{ext}}$ and $M_2\in (W,A_2)$; 
\item[(iii)] if $L \in[L_1, L_R]$ or $L \in [L_3, L_2]$;
the Riemann solution is
$
     L\testright{R_s} T_1 \xrightarrow{'S_s} M_1 \testright{S_f}  R$ or   $L\testright{R_s} T_2 \xrightarrow{'S_s} M_2 \testright{S_f} R,
$
where $T_1\in(G, T_R)_{\text{ext}}$, $M_1\in(A_1,R)$,  $T_2\in (W, T_R)_{\text{ext}}$ and $M_2\in (A_2,A_3)$;
\item[(iv)] if $L \in (L_R, L_F]$, 
the Riemann solution is
$L\testright{R_s} T_1 \xrightarrow{'S_s} M_1 \testright{R_f}  R,
$
where $T_1\in(T_R,W)_{\text{ext}}$ and $M_1\in(R,F)$;
\item[(v)] if $L \in(L_F, L_3)$,
the Riemann solution is
$
     L\testright{R_s} T_1 \xrightarrow{'S_s} M_1' \testright{S_f\,'} M_1 \xrightarrow{R_f} R,
$
where $T_1\in(W, T_R)_{\text{ext}}$, $M_1'\in(F, A_3)$ and $M_1\in(F, R)$. 
\end{itemize}
\end{cla}
%

\subsubsection{Subregion \texorpdfstring{$\Theta_2$}{O}.}
\label{subsec:RinTheta_2}
The boundary of subregion $\Theta_2$
in Fig.~\ref{fig:RegioesThetas}(a) comprises the following parts: segment $C_3$-$P$ of the double contact locus $C$; the $f$-rarefaction segment $P$-$R_P$;   
segment $R_P$-$\mathcal{I}_E$ of the $f$-inflection locus $\mathcal{I}_f$; and segment $\mathcal{I}_E$-$C_3$ of the edge \GOc.

As shown in Fig.~\ref{fig:RRegions Completed2}(b), the segments $C_2$-$\mathcal{I}_2^f$ and $C_1$-$\mathcal{I}_1^f$ split $\Theta_2$ into subregions $\Theta_2^a$, $\Theta_2^b$, $\Theta_2^c$ characterized by certain topological changes in the backward fast wave curve for $R$.
The segment $C_2$-$\mathcal{I}_2^f$ is part of the invariant line $\um$-$D_0$ (see Fig.~\ref{fig:RRegions Completed}); the segment $C_1$-$\mathcal{I}_1^f$ is part of the $f$-rarefaction from $K$, where $K$ 
in Fig.~\ref{fig:RRegions Completed2}(c) is
the intersection of $\mathcal{I}_f$ and $MC$. 

\bigskip

\noindent{\bf The backward Hugoniot curve and the admissible shocks for $R$ in subregion $\Theta_2$.}
\medskip

\begin{figure}
	\centering
	\subfigure[$R=(0.0298543, 0.884736)$ in $\Theta_2^a$]
	{\includegraphics[width=0.4\linewidth]{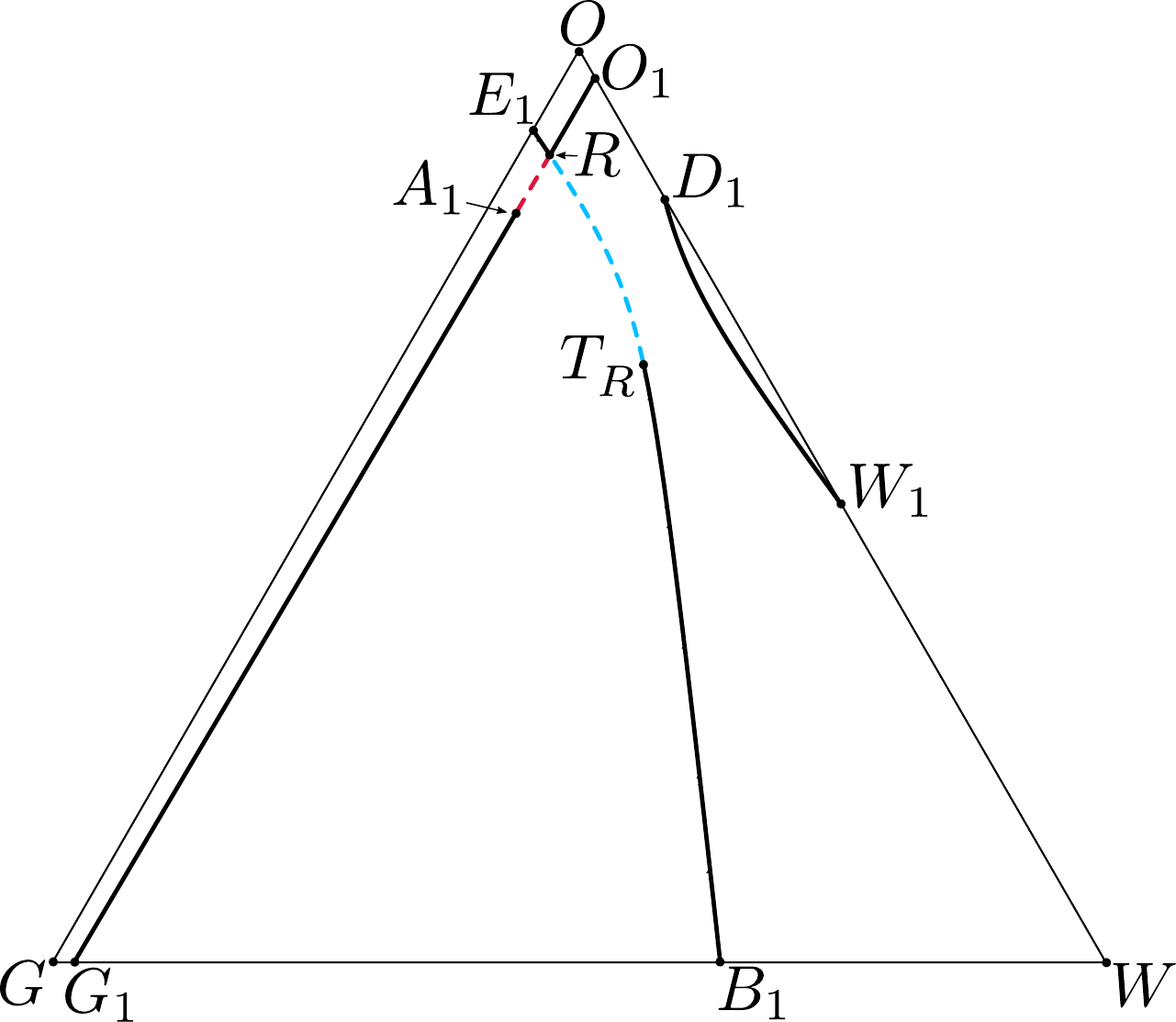}} 
 \hspace{2mm}
		\subfigure[$R=(0.10249, 0.82816)$ in $\Theta_2^b\cup \Theta_2^c$]
	{\includegraphics[width=0.4\linewidth]{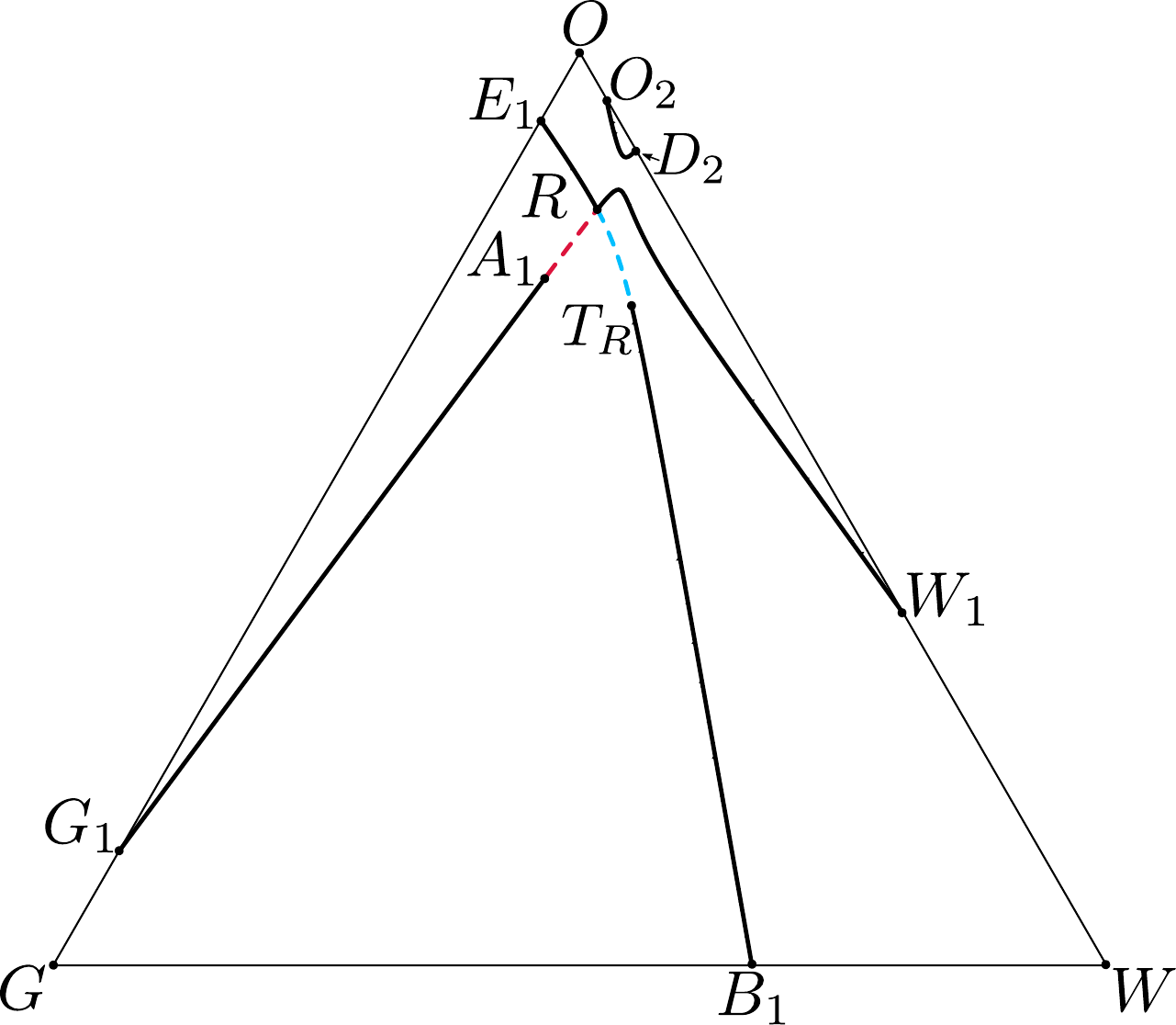}}
	\caption{Hugoniot curve for $R$ in  the subregion of $\Theta_2$. 
		We have $\sigma(A_1;R) = \lambdaf(A_1)$, $\sigma(T_R;R) = \lambdas(T_R)$.
		The segments $(R, A_1]$ and $(R, T_R]$ are admissible.
}
	\label{fig:Hugoniot_R3}
\end{figure}

%
For a right state $R$ lying precisely
on the double contact segment $C$ in
Figs.~\ref{fig:RegioesThetas}, \ref{fig:RRegions Completed2},
the shock segment $[A_2, A_3]$ of $H(R)$ in Figs.~\ref{fig:Hug-R10a}(b), \ref{fig:Hug-R10a}(d)
collapses to a single state $A_2\equiv A_3$ with $\sigma(A_2; R) = \lambdaf(R) = \lambdaf(A_2)$.
When the right state $R$ moves from  $\Theta_1$ to 
$\Theta_2$, the segment $[A_2, A_3]$ disappears and
Figs.~\ref{fig:Hug-R10a}(b), \ref{fig:Hug-R10a}(d) evolve to
Figs.~\ref{fig:Hugoniot_R3}(a), \ref{fig:Hugoniot_R3}(b), respectively.
The fast segment $[A_1, R)$ and the slow segment $(R, T_R]$ of admissible
shocks are preserved as for $R$ in the subregions
$\Theta_1^b$, $\Theta_1^d$ and $\Theta_1^e$.

\begin{rem}\label{rem:theta2}
Notice that for $R$ in 
$\Theta_2$ the local $f$-rarefaction segment through $R$ intersects the double contact segment $C$ at a state $Q_1$ which defines a corresponding state $Q_1'$
in $H(Q_1)$ such that
$\sigma(Q_1';Q_1) = \lambdaf(Q_1) = \lambdaf(Q_1')$, with $Q_1'$
in the
triangle $W\um B$ in Fig.~\ref{fig:RRegions Completed}.
\end{rem}

\medskip

\noindent{\bf The backward fast wave curve $\wm_f(R)$ for $R$ in subregion $\Theta_2^a$.}
\medskip

\begin{figure}
	\centering
	\subfigure[Admissible $\wm_f(R)$ for $R$.]
	{\includegraphics[width=0.4\linewidth]{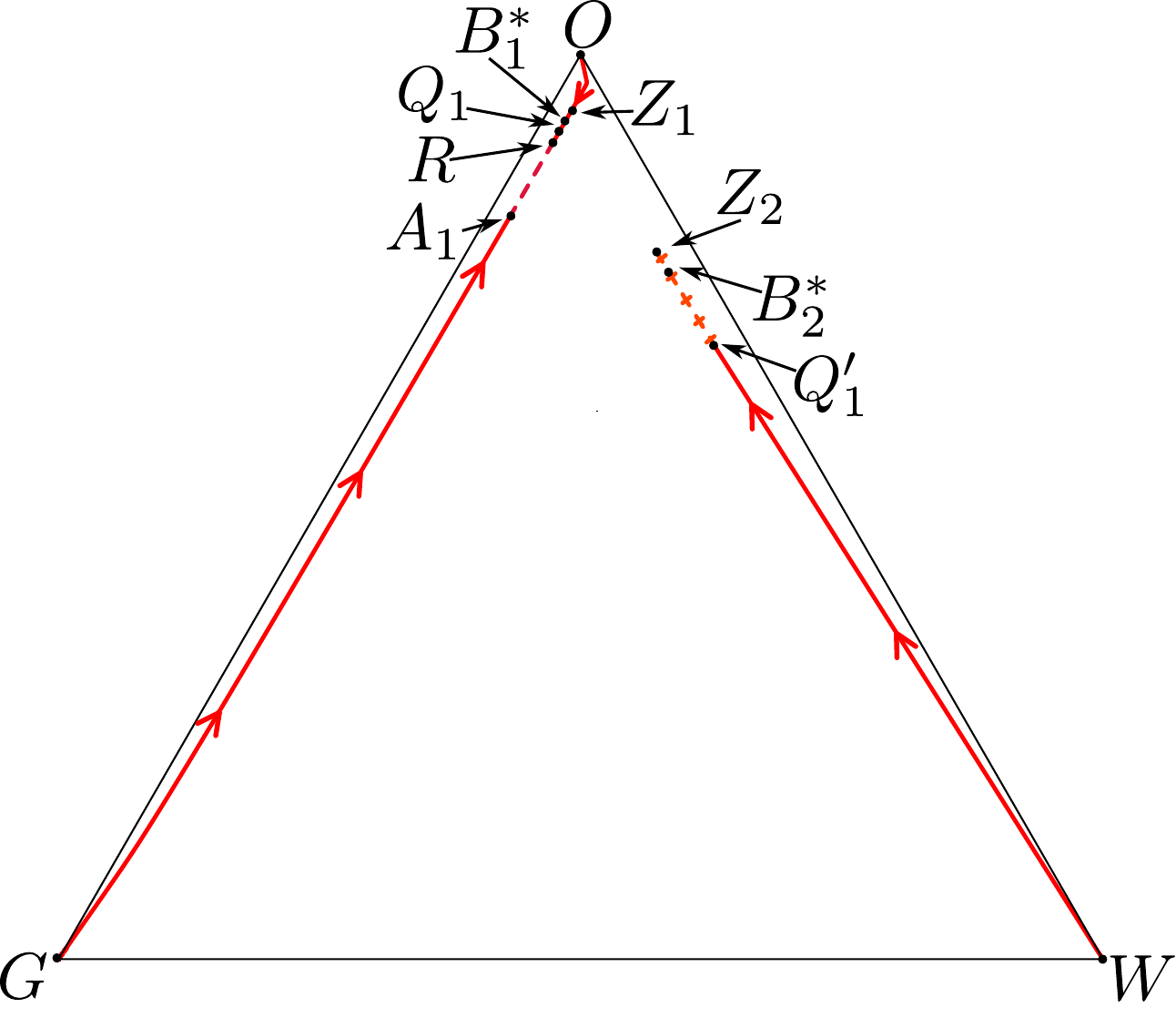}}  
	\hspace{1.75mm}
	\subfigure[Riemann solution with $L$ varying along \GWc.]
	{\includegraphics[width=0.4\linewidth]{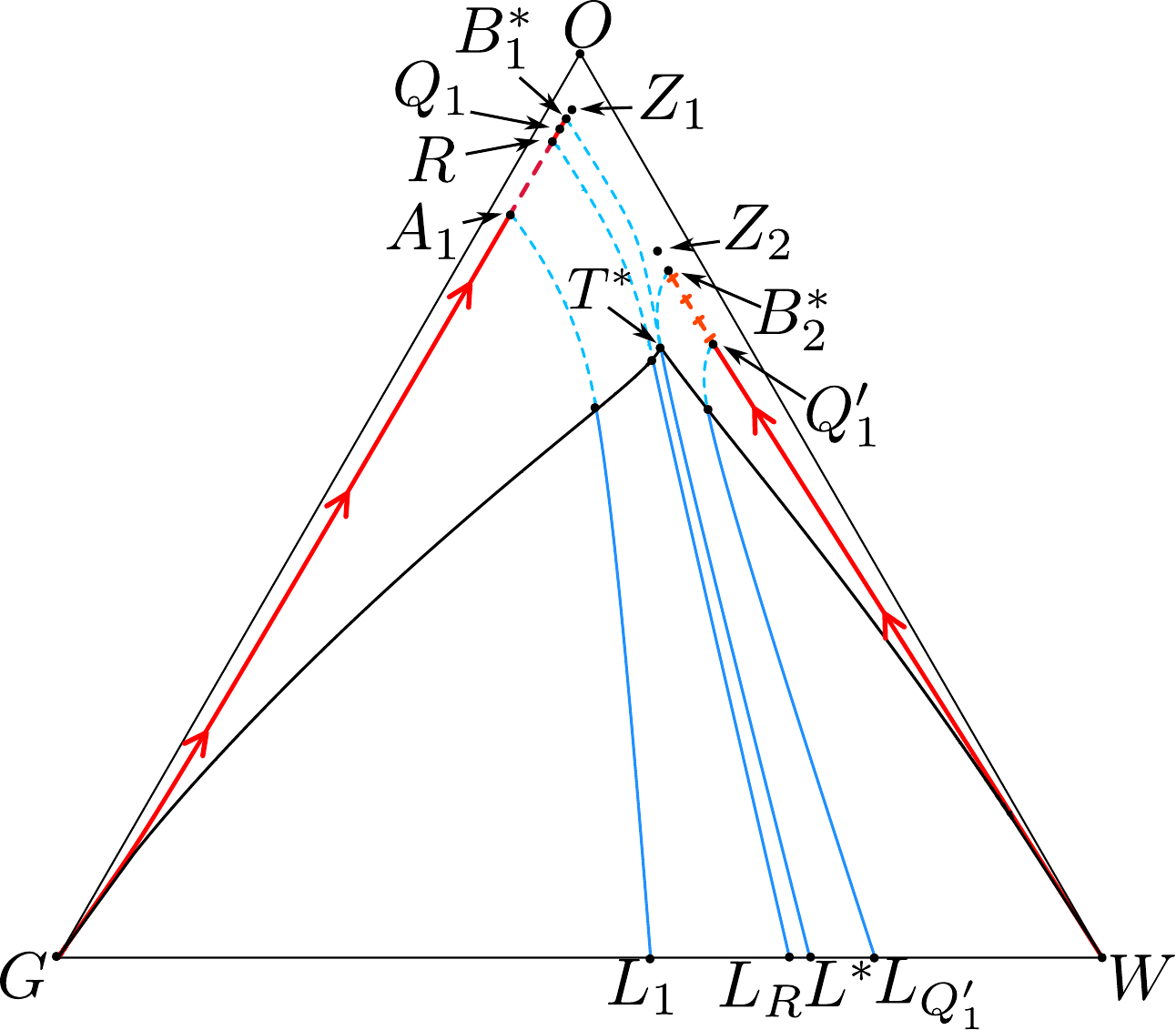}}  
	\caption{{$R$ in subregion of $\Theta_2$ bounded by
	$\Theta_2^a$ of
    Fig.~\ref{fig:RRegions Completed2}(b). We have
    $\sigma(Q_1',Q_1) = \lambdaf(Q_1) = \lambdaf(Q_1')$ and
    $\sigma(B_2^*,B_1^*) = \lambdaf(B_1^*)
    = \lambdas(T^*)$.
Curve $(G,T^*]_{\text{ext}}\cup[T^*,W)_{\text{ext}}$ in (b)
is an extension of the $\wm_f(R)$ segment $(G,B_1^*]\cup[B_2^*,W)$ such that
$\sigma(M', M) = \lambdas(M')$, for all $(M', M) \in [G, T^*]_{\text{ext}}\cup [T^*, W]_{\text{ext}} 
\times (G, B_1^*]\cup [B_2^*, W)$.    }}
	\label{fig:RS-R11}
\end{figure}

Refer to Fig.~\ref{fig:RS-R11}(a). 
The $f$-rarefaction curve for a state $R$ in $\Theta_2^a$ crosses the double contact segment $C_2-\mathcal{I}_2^f$ at a state $Q_1$, then the invariant
segment \EU at a state $Z_1$ before it reaches the vertex $O$.
The segment $[Z_1, Q_1]$ of the $f$-rarefaction curve
defines the admissible composite segment $[Z_2, Q_1']$
satisfying $\sigma(M';M) = \lambdaf(M)$
for all $M'$ in $[Z_2, Q_1']$ and $M$ in $[Z_1, Q_1]$. On the other hand the rarefaction segment $(Q_1, R)$ does not define a composite segment because there is no state $M'$
in $\Hm(M)$ such that $\sigma(M'; M)=\lambdaf(M)$. See the description of the Hugoniot locus for $R$ in $\Theta_2$ above, replacing $R$ with $M$ and $A_3$ with $M'$.


\begin{cla}\label{cla:BWCY3Y4IF2IF1}
	Refer to Fig.~\ref{fig:RS-R11}(a).
	For $R$ in $\Theta_2^a$, $\wm_f(R)$ comprises states
	$M$ along the $f$-rarefaction segments
	$[G, A_1)$, $[O,R)$ and $(Q_1', W]$, states $M$
	along the shock segment $[A_1, R)$,
	and states $M'$ along the composite segment
	$[Z_2, Q_1']$ such that for each state 
	$M'\in [Z_2, Q_1']$ there is a unique state 
	$M \in [Z_1, Q_1]$ with $\sigma(M';M) = \lambdaf(M)$. 
\end{cla}

\medskip

\noindent{\bf The backward fast wave curve $\wm_f(R)$ for $R$ in subregion $\Theta_2^b \cup \Theta_2^c$.}
\medskip

\begin{figure}
	\centering
	\subfigure[Wave Curve for $R$ in subregion $\Theta_2^b$
]{\includegraphics[width=0.4\linewidth]{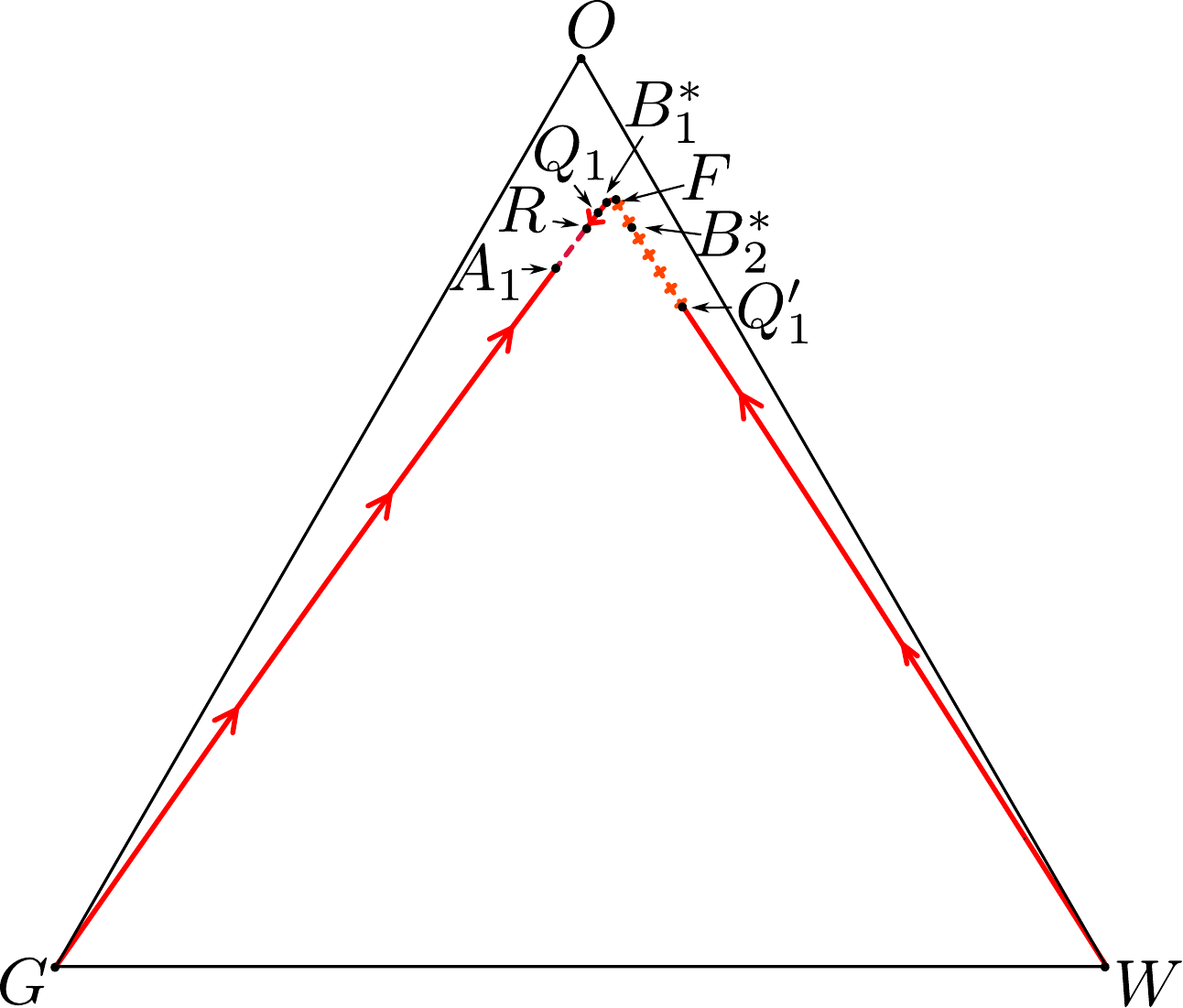}}   \hspace{1.5mm}	
	\subfigure[Wave Curve for $R$ in subregion $\Theta_2^c$
	]{\includegraphics[width=0.4\linewidth]{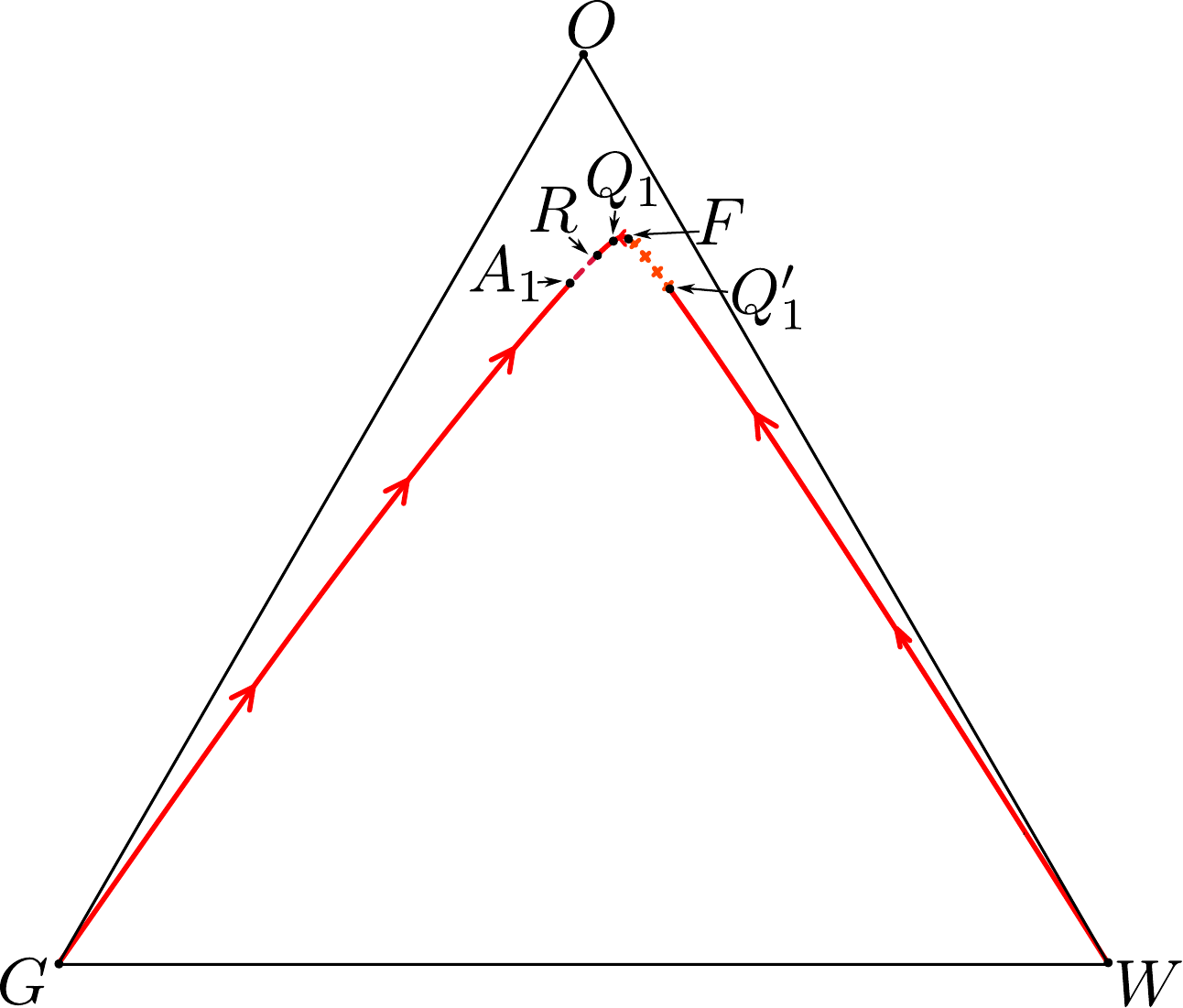}} 
	\caption{Wave Curves}
	\label{fig:Wave_Curves_Theta_2_b_Theta_2_c}
\end{figure}

\begin{figure}
	\centering
	\subfigure[Riemann solution for $R$ in subregion $\Theta_2^b$
]{\includegraphics[width=0.4\linewidth]{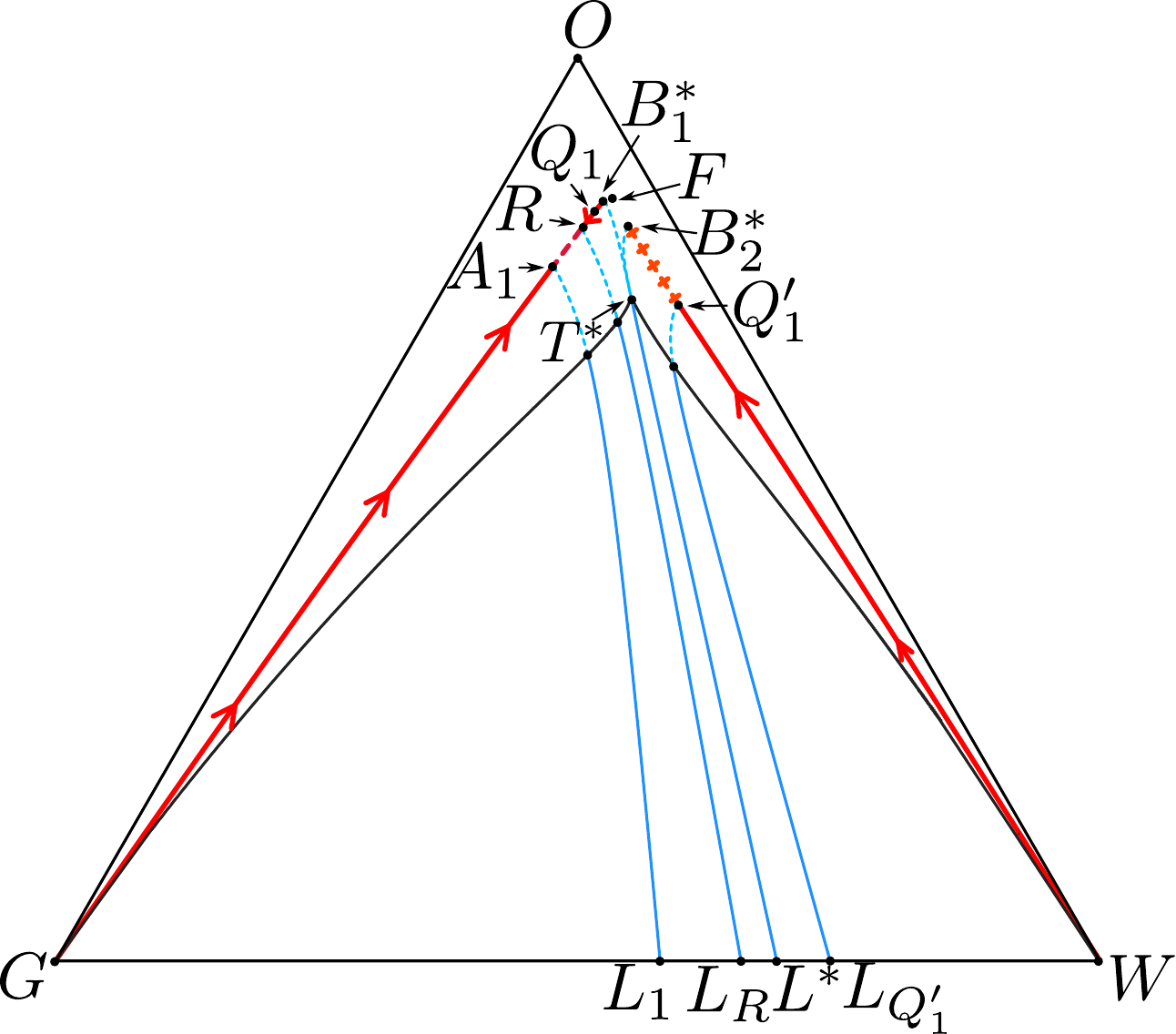}}   \hspace{1.5mm}	
	\subfigure[Riemann solution for $R$ in subregion $\Theta_2^c$
	]{\includegraphics[width=0.4\linewidth]{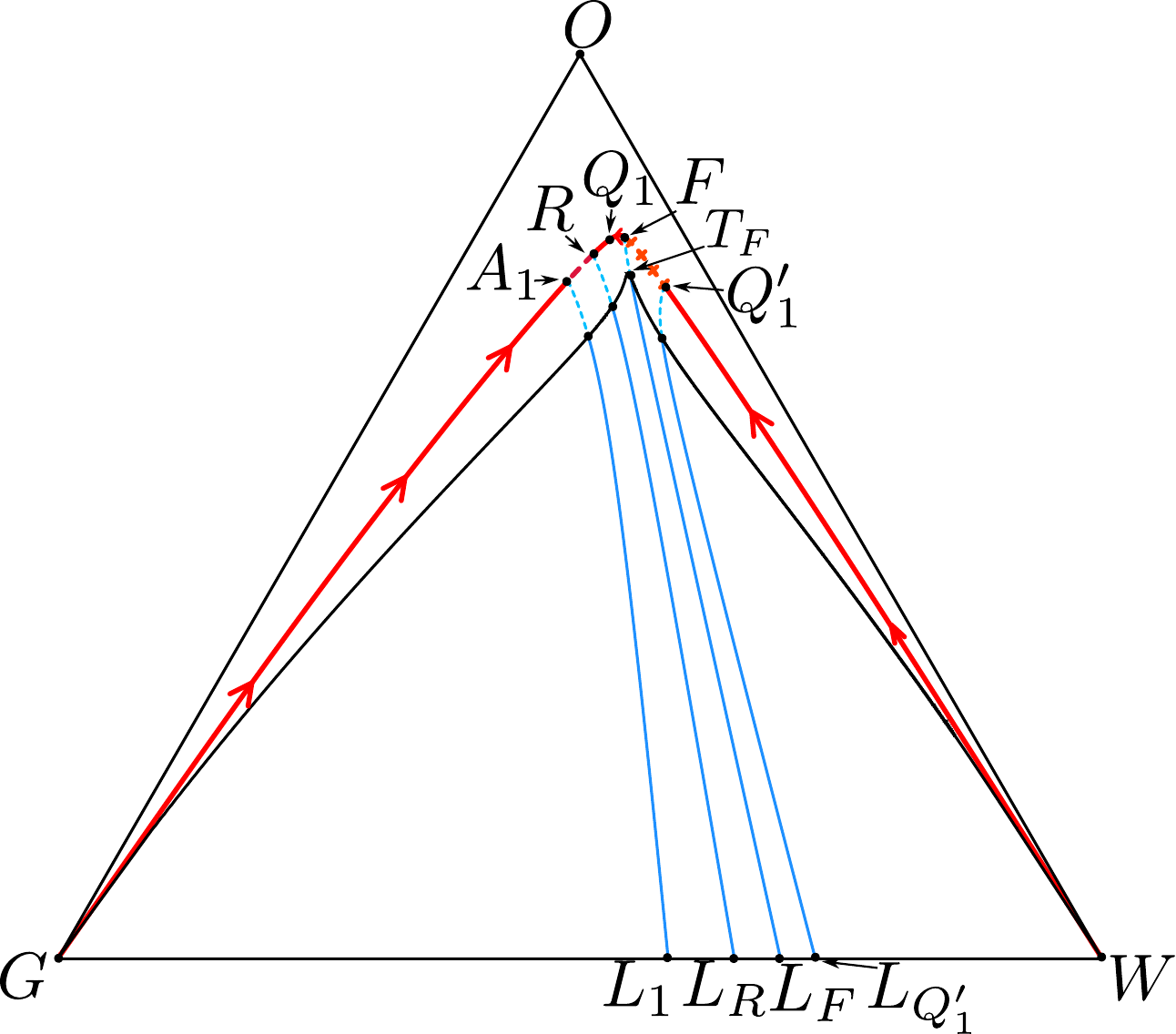}} 
	\caption{Riemann solution for $R$ in $\Theta_2$ with $R$ below the invariant segment $C_2$ - $\mathcal{I}_2^f$ of
    Fig.~\ref{fig:RRegions Completed}(b)
}
	\label{fig:RegionR11BeC}
\end{figure}

For $R$ in subregion $\Theta_2^b \cup \Theta_2^c$, the local $f$-rarefaction curve through $R$ extends from $R$ to a state $F$ on the $f$-inflection segment $\um$-$K$-$P$. Now the state $Q_1$ lies on the intersection of the $f$-rarefaction segment $(R, F]$ with the double contact $C$. 
The segment $[Q_1, F]$ of the $f$-rarefaction curve
defines the admissible composite segment $(F, Q_1']$
satisfying $\sigma(M';M) = \lambdaf(M)$
for all $M'$ in $(F, Q_1']$ and $M$ in $(Q_1, F]$. The state $Q_1'$, which lies in the double contact segment $P-C_4$ shown in  Fig.~\ref{fig:BoundaryExt_MixedDoubleCont}, corresponds to the state $Q_1$.
See Fig.~\ref{fig:Wave_Curves_Theta_2_b_Theta_2_c}. Consequently, we have:

\begin{cla}\label{cla:Y4Y5I3fI2f}
Refer to Fig.~\ref{fig:RegionR11BeC}.
	For $R$ in subregion $\Theta_2^b \cup \Theta_2^c$
	in Fig.\ref{fig:RRegions Completed2}(b),
	$\wm_f(R)$ comprises states
	$M$ along the  rarefaction segments
	$[G, A_1)$, $(R, F]$ and $[Q_1', W]$, states
	$M$
	along the shock segment $[A_1, R)$,
	and states $M'$ along the composite segment
	$(F, Q_1']$ such that for each state 
	$M'\in (F, Q_1']$ there is a unique state 
	$M \in [Q_1, F)$ with $\sigma(M';M) = \lambdaf(M)$. 
\end{cla}

\medskip
\noindent{\bf The Riemann solution for $R$ in subregion $\Theta_2 = \Theta_2^a \cup \Theta_2^b \cup \Theta_2^c$.}

\medskip

Compared with the Riemann solutions previously discussed, the main differences in the
Riemann solutions for $R$ in $\Theta_2$ lie in the role of states $Q_1$ and $Q_1'$ on the double contact segment $C$.

\medskip

\noindent{\bf Riemann solution for $R$ in subregion $\Theta_2^a \cup \Theta_2^b$.}
\medskip

It is worth noticing that for $R$ in $\Theta_2^a \cup \Theta_2^b$ the $f$-rarefaction segment through $R$ intersects the mixed contact segment $MC$ at a state $B_1^*$, which defines the states $T^*$ and $B_2^*$ satisfying
$ \lambdas(T^*) = \sigma(T^*;B_1^*) = \lambdaf(B_1^*) 
= \sigma(B_2^*; B_1^*) = \sigma(T^*; B_2^*)$.
\begin{cla}\label{cla:RSolution-RinR11}
Refer to Fig.~\ref{fig:RS-R11}(b). Let $L$ be a state on the edge \GW of the saturation triangle and  $R$ be a state in subregion $\Theta_2^a\cup \Theta_2^b$ in
Figs.~\ref{fig:RegioesThetas}(a) and \ref{fig:RRegions Completed2}(b).

Let $A_1$, $Q_1$, $Q_1'$, $B_1^*$ and $B_2^*$ be the states on  $\wm_f(R)$ satisfying
$\sigma(A_1;R) = \lambdaf(A_1)$, $\lambdaf(Q_1) = \lambdaf(Q_1') = \sigma(Q_1';Q_1)$, $ \lambdas(T^*) = \sigma(T^*;B_1^*) = \lambdaf(B_1^*) 
= \sigma(B_2^*; B_1^*) = \sigma(T^*; B_2^*)$, with $T^*\in \mathcal{H}( B_1^*$).

Let $L_1$, $L_R$, $L^*$ and $L_{Q_1'}$ be the intersection points of the backward slow wave curves through
$A_1$, $R$,
$B_1^*$ (or $B_2^*$) and $Q_1'$ with the edge \GWc, respectively.
Then, 
\begin{itemize}
\item[(i)] if $L=G$, the Riemann solution is \;
$
    G\testright{R_f} A_1  \xrightarrow{'S_f} R;
$
\item[(ii)] if $L=W$, the Riemann solution is \;
$
   W\testright{R_f} Q_1' \xrightarrow{'S_f\,'} Q_1  \xrightarrow{R_f} R;
$
\item[(iii)] if $L \in(G, L_1)$, the Riemann solution is \;
$
     L\testright{R_s} T_1 \xrightarrow{'S_s} M_1 \testright{R_f} A_1 \xrightarrow{'S_f} R,
$
where $T_1\in(G, T^*)_{\text{ext}}$ and $M_1\in(G, A_1)$; 
\item[(iv)] if $L \in[L_1, L_R]$, the Riemann solution is \;
$
     L\testright{R_s} T_1 \xrightarrow{'S_s} M_1 \testright{S_f}  R,
$
where $T_1\in(G, T^*)_{\text{ext}}$ and $M_1\in[A_1,R)$; 
\item[(v)] if $L \in (L_R, L^*)$, 
the Riemann solution is \;
$
     L\testright{R_s} T_1 \xrightarrow{'S_s} M_1 \testright{R_f}  R,
$
where $T_1\in(G,T^*)_{\text{ext}}$ and $M_1\in(R,B_1^*)$; 
\item[(vi)] if $L \in(L^*, L_{Q_1'})$,
the Riemann solution is \;
$
     L\testright{R_s} T_1 \xrightarrow{'S_s} M_1' \testright{S_f\,'} M_1 \xrightarrow{R_f} R,
$
where $T_1\in(W, T^*)_{\text{ext}}$, $M_1'\in(B_2^*, Q_1')$ and $M_1\in(B_1^*, Q_1)$; 
\item[(vii)] if $L \in(L_{Q_1'},W)$,
the Riemann solution is \;
$
     L\testright{R_s} T_1 \xrightarrow{'S_s} M_1' \testright{R_s}Q_1' \xrightarrow{'S_f\,'} Q_1 \xrightarrow{R_f} R,
$
where $T_1\in(W, T^*)_{\text{ext}}$  and $M_1'\in(W, Q_1')$. 
\end{itemize}
\end{cla}
The case $L = L^*$ is similar to that described in 
Remark~\ref{rem:L*B*}.

\medskip

\noindent{ \bf Riemann solution for $R$ in subregion $\Theta_2^c$.}
\medskip

For $R$ in $\Theta_2^c$ the $f$-rarefaction segment $(R, F]$ does not intersect
the mixed contact segment $MC$. Consequently, the states $B_1^*$, $T^*$, and $B_2^*$, which occur when $R$ lies in region $\Theta_2^a \cup \Theta_2^b$, no longer exist. Hence only the double contact segment $C$ interferes in the construction of Riemann solutions.

\begin{cla}\label{cla:RSolution-RinY5PIf4If3}
Refer to Fig.~\ref{fig:RegionR11BeC}(b). Let $L$ be a state on the edge \GW of the saturation triangle and  $R$ be a state in subregion $\Theta_2^c$ in
Figs.~\ref{fig:RegioesThetas}(a) and \ref{fig:RRegions Completed2}(b).
Let $A_1$, $Q_1$, $Q_1'$ and $F$ be the states on the
$\wm_f(R)$ satisfying
$\sigma(A_1; R) = \lambdaf(A_1)$, $\lambdaf(Q_1)= \lambdaf(Q_1') = \sigma(Q_1';Q_1)$ and $F$ on the intersection of the $f$-rarefaction segment through $R$ with the $f$-inflection segment $K$-$P$.

Let $L_1$, $L_R$, $L_F$ and $L_{Q_1'}$
be the intersection points of
the backward slow wave curves through
$A_1$, $R$,
$F$ and $Q_1'$ with the edge \GWc, respectively.
Then, 
\begin{itemize}
\item[(i)] if $L=G$, the Riemann solution is
$
    G\testright{R_f} A_1  \xrightarrow{'S_f} R;
$
\item[(ii)] if $L=W$, the Riemann solution is
$
   W\testright{R_f} Q_1' \xrightarrow{'S_f\,'} Q_1  \xrightarrow{R_f} R;
$
\item[(iii)] if $L \in[L_1, L_R]$, the Riemann solution is
$
     L\testright{R_s} T_1 \xrightarrow{'S_s} M_1 \testright{S_f}  R,
$
where $T_1\in(G, T_F)_{\text{ext}}$ and $M_1\in[A_1,R)$;
\item[(iv)] if $L \in (L_R, L_F]$, the Riemann solution is
$
     L\testright{R_s} T_1 \xrightarrow{'S_s} M_1 \testright{R_f}  R,
$
where $T_1\in(G, T_F)_{\text{ext}}$ and $M_1\in[F,R)$; 
\item[(v)] if $L \in (L_F, L_{Q_1'}]$, 
the Riemann solution is
$
     L\testright{R_s} T_1 \xrightarrow{'S_s} M_1' \testright{S_f\,'} M_1 \xrightarrow{R_f} R,
$
where $T_1\in(W, T_F)_{\text{ext}}$, $M_1'\in(F, Q_1')$ and $M_1\in(R, Q_1)$; 

\item[(vi)] if $L \in(L_{Q_1'},W)$, the Riemann solution is
$
     L\testright{R_s} T_1 \xrightarrow{'S_s} M_1' \testright{R_s}Q_1' \xrightarrow{'S_f\,'} Q_1 \xrightarrow{R_f} R,
$
where $T_1\in(W, T^*)_{\text{ext}}$  and $M_1'\in(W, Q_1')$. 
\end{itemize}
\end{cla}



%

\subsubsection{Subregion \texorpdfstring{$\Theta_3$}{O}.}
\label{subsec:RinTheta_3}
The boundary of subregion $\Theta_3$
in Figs.~\ref{fig:RegioesThetas}(a) and
\ref{fig:RRegions Completed2}(a) consists of the segment $\mathcal{I}_E$-$R_P$ of the $f$-inflection locus $\mathcal{I}_f$; the $f$-composite segment $R_P$-$P'$ corresponding to the $f$-rarefaction segment $P$-$R_P$ in
Fig.~\ref{fig:RRegions Completed2}(b);
the segment $P'$-$C'_3$, which is the extension of the segment $P$-$C_3$ of the double contact locus $C$ in Fig.~\ref{fig:RRegions Completed2}(b); and the segment $C'_3$-$\mathcal{I}_E$
of the edge \GOc.

The region $\Theta_3$ is subdivided into three subregions: $\Theta_3^i$, with $i\in \{a,b,c\}$. Their shared boundaries are the invariant segment $\mathcal{I}_2^f$-$C_2'$, and the $f$-composite segment $\mathcal{I}_1^f$-$C_1'$ associated with the $f$-rarefaction segment $\mathcal{I}_1^f$-$C_1$.
 
\bigskip

\noindent{\bf The backward Hugoniot curve and the admissible shocks for $R$ in subregion $\Theta_3$.}
\medskip 

\begin{figure}[]
	\centering
	\subfigure[$R = (0.0401428,0.722888)$. $R \in \Theta_3^a$]
	{\includegraphics[width=0.4\linewidth]{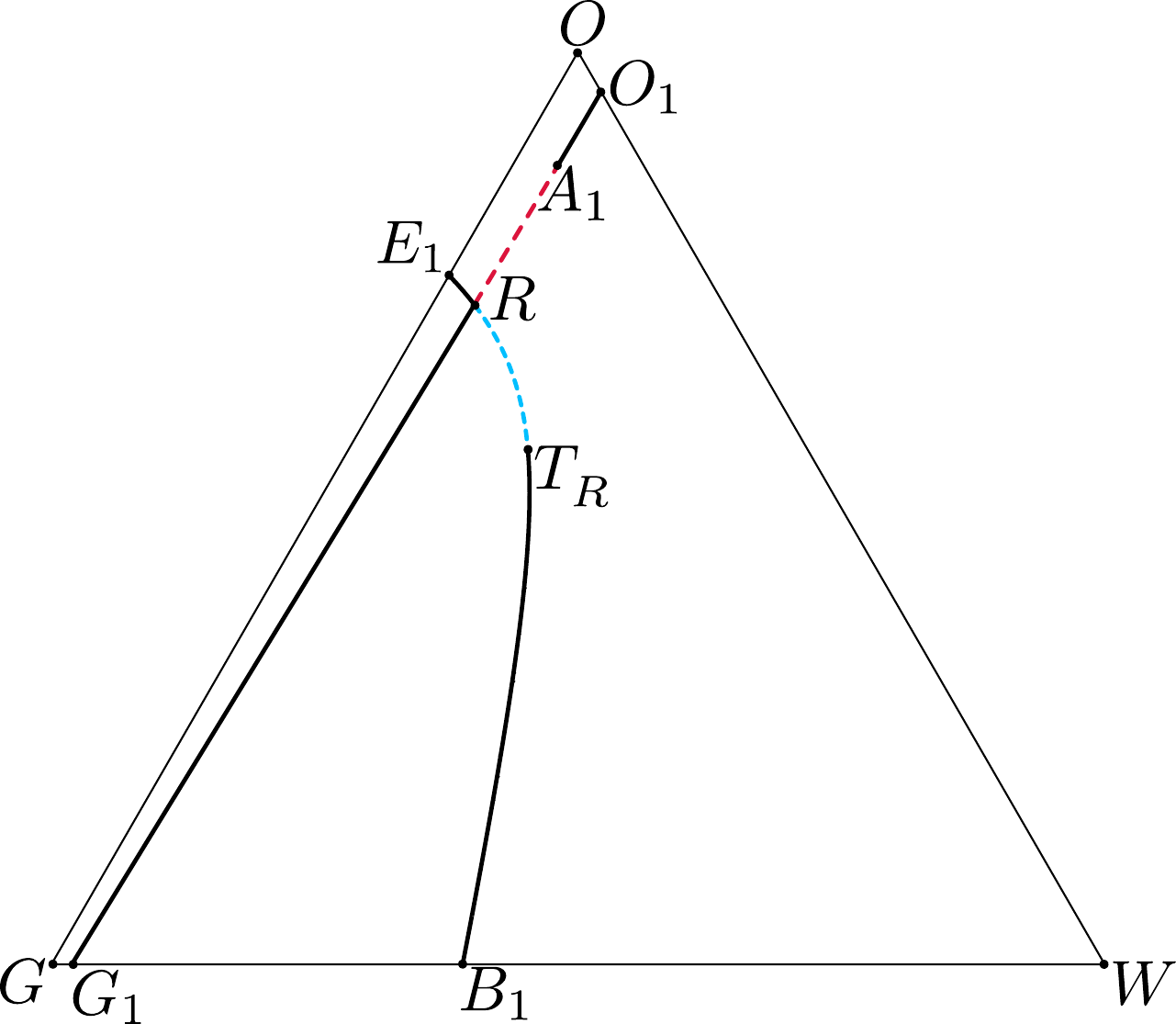}}  
	\hspace{1.75mm}
	\subfigure[$R = (0.0323924,0.578612)$. $R \in \Theta_3^a$]
	{\includegraphics[width=0.4\linewidth]{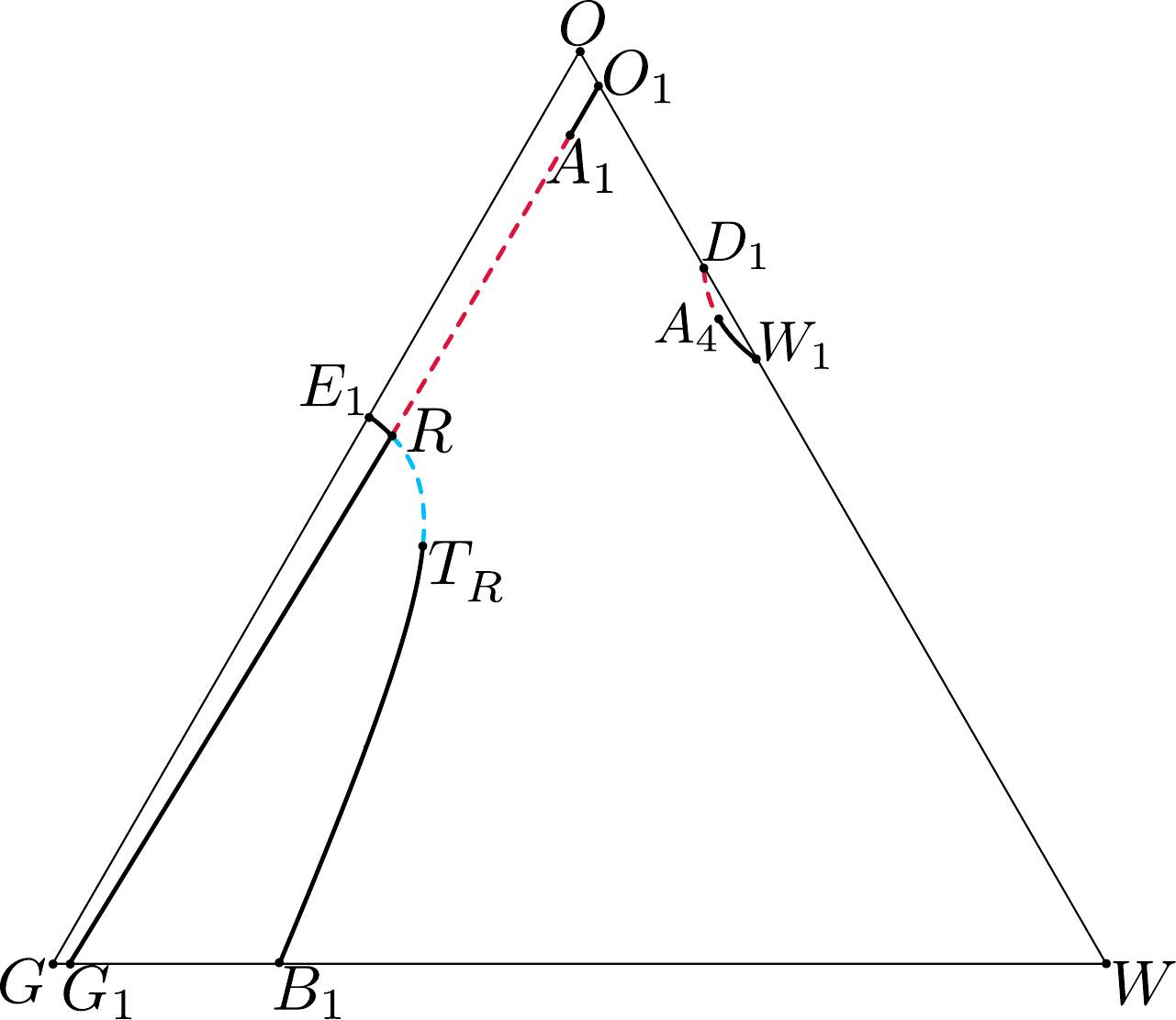}}  
	\caption{Backward Hugoniot curve based on $R$ in subregion $\Theta_3^a$. The detached branch $[D_1, W_1]$ in (b) is not admissible according to the viscous profile criterion.
	We have: $\sigma(A_1; R) = \lambdaf(A_1)$,
	$\sigma(A_4; R) = \lambdaf(A_4)$,
	$\sigma(T_R; R) = \lambdas(T_R)$
	}
	\label{fig:BFWC-R4_Above}
\end{figure}

\begin{figure}[]
	\centering
	\subfigure[$R = (0.0847167,0.709481)$.  $R \in \Theta_3^b$]
	{\includegraphics[width=0.4\linewidth]{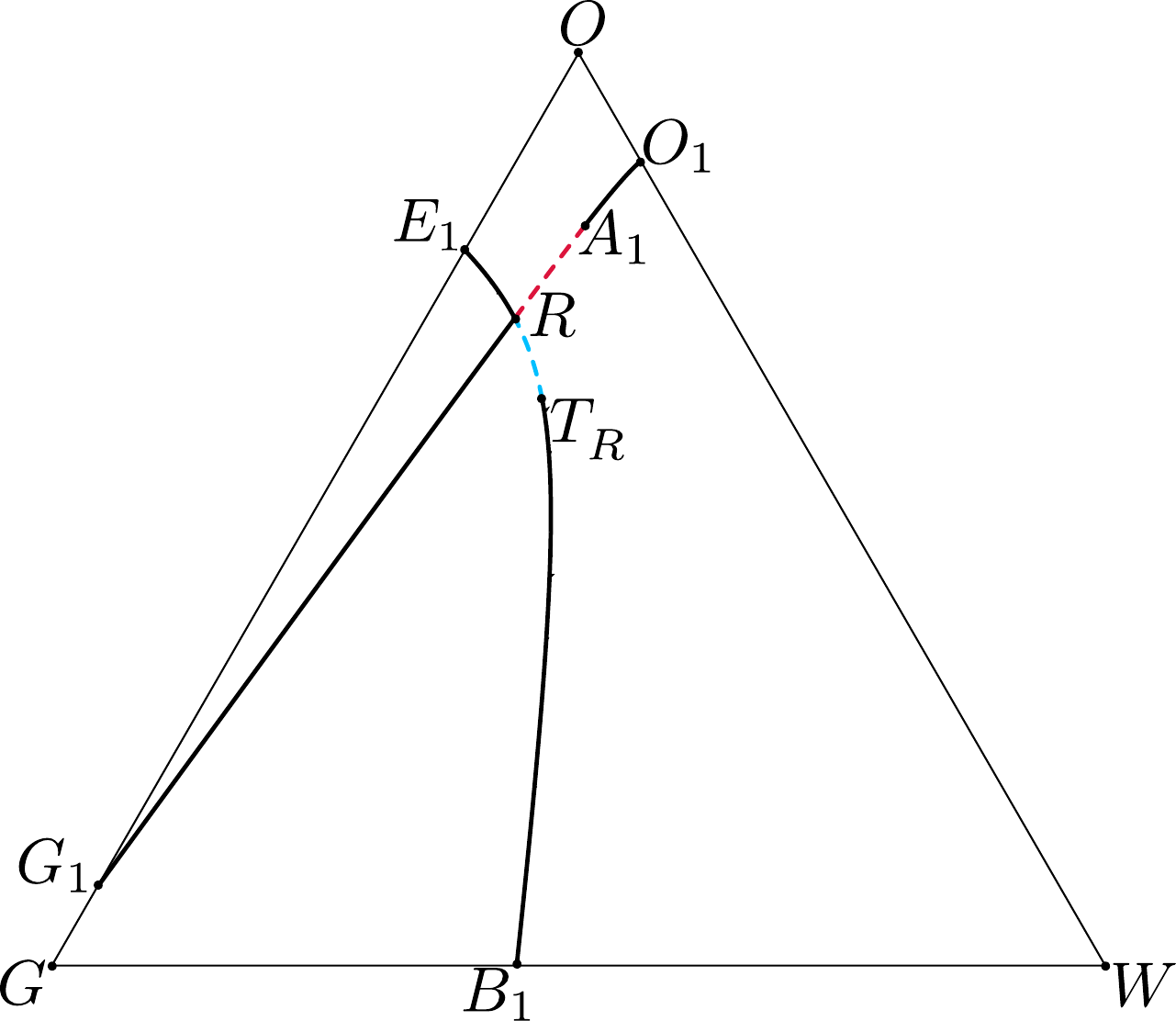}}  
	\hspace{1.75mm}
	\subfigure[$R = (0.0705098,0.630641)$.  $R \in \Theta_3^b$]
	{\includegraphics[width=0.4\linewidth]{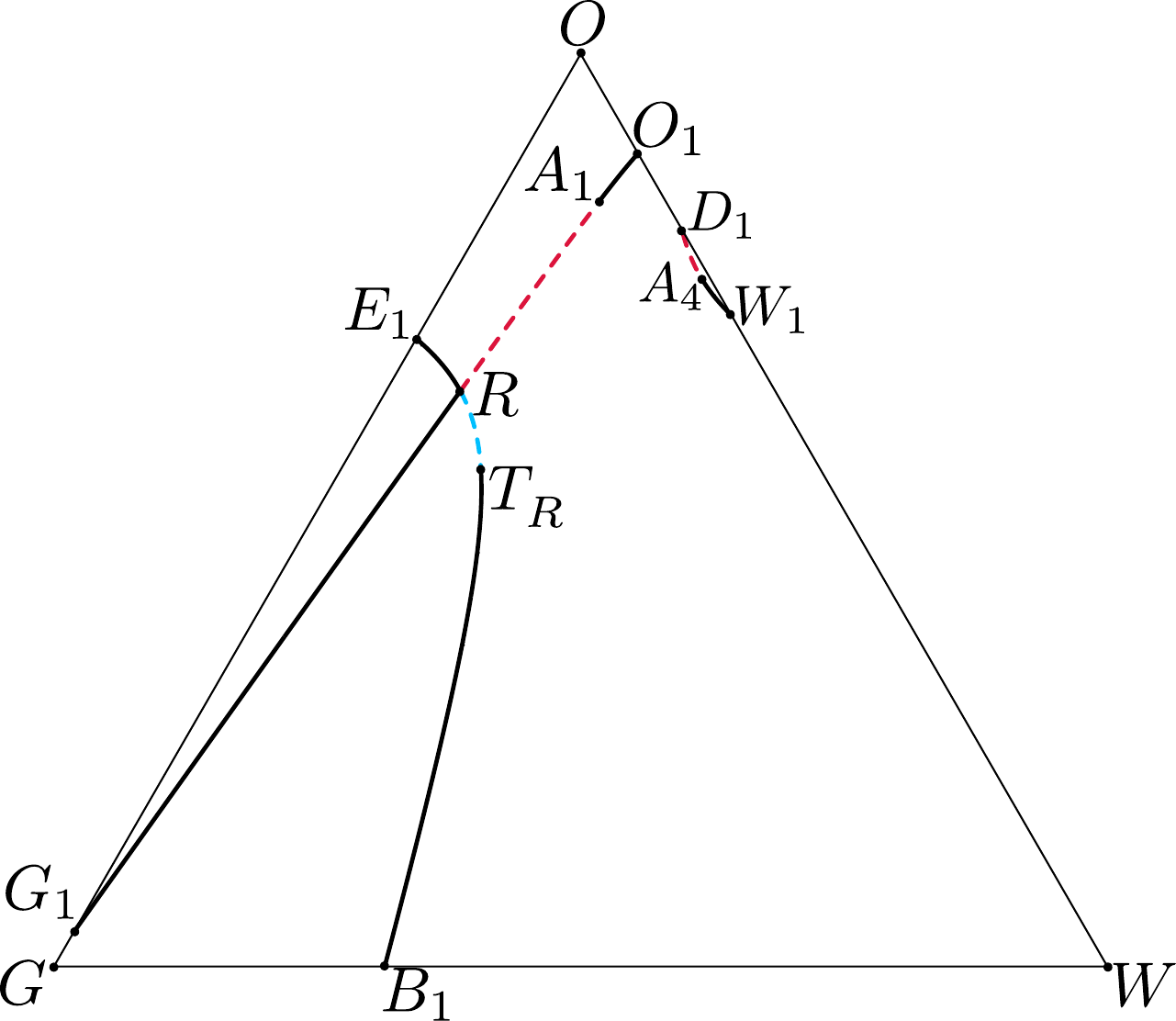}}  
	\caption{Backward Hugoniot curve based on $R$ in subregion $\Theta_3^b$. The detached branch $[D_1, W_1]$ in (b) is not admissible according to the viscous profile criterion.
	We have: $\sigma(A_1; R) = \lambdaf(A_1)$,
	$\sigma(A_4; R) = \lambdaf(A_4)$,
	$\sigma(T_R; R) = \lambdas(T_R)$.
	}
	\label{fig:BFWC-R4_belowB}
\end{figure}

\begin{figure}[]
	\centering
	\subfigure[$R = (0.109168,0.715865)$.  $R$ in $\Theta_3^c$]
	{\includegraphics[width=0.4\linewidth]{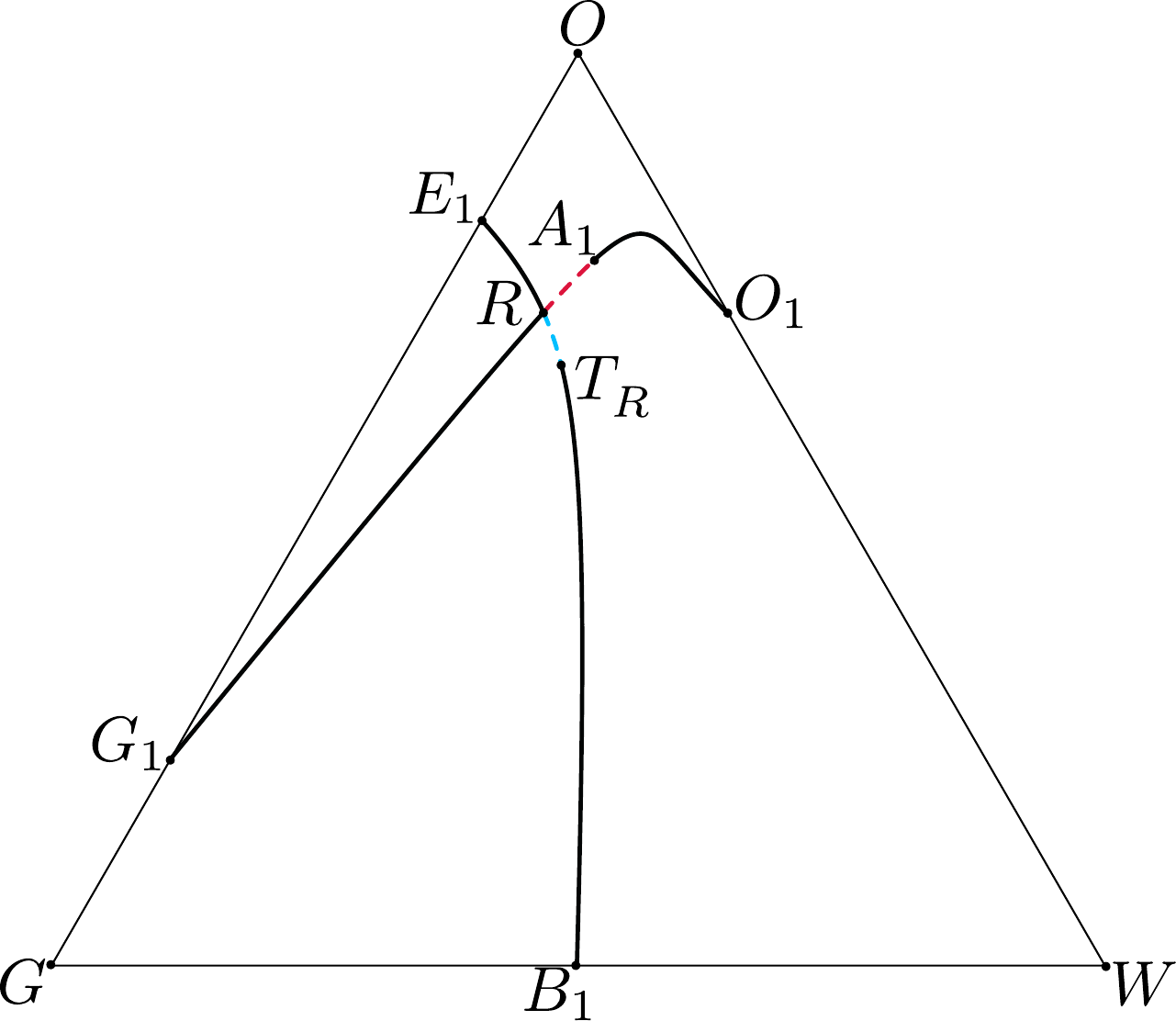}}  
	\hspace{1.75mm}
	\caption{{Backward Hugoniot curve for $R$ in subregion $\Theta_3^c$.
	We have: $\sigma(A_1; R) = \lambdaf(A_1)$ and
	$\sigma(T_R; R) = \lambdas(T_R)$.
	}}
	\label{figBlowFWC-R4_beC}
\end{figure}

Refer to Figs.~\ref{fig:BFWC-R4_Above}, \ref{fig:BFWC-R4_belowB}, and
\ref{figBlowFWC-R4_beC}.
As the state $R$ enters into region $\Theta_3$
the local $f$-shock segment $[A_1, R)$, that
was part of the attached branch $[G_1, R)$ of $\mathcal{H}(R)$
for $R$ in $\Theta_2$, now is part of  
the attached branch $(R, O_1]$ of $\mathcal{H}(R)$.
Due to this inversion of direction, we now denote this local shock segment by $(R, A_1]$.
This is the main difference in the construction of the
Riemann solution when compared with the case for $R$ in
$\Theta_2$.



Referring to Fig.~\ref{figBlowFWC-R4_beC}(b) the Hugoniot branch $[R,O_1]$ has two $f$-shock segments, $[R, A_1]$ and $[A_5,A_4]$. Only the segment $[R, A_1]$ is admissible. 

As for $R$ in $\Theta_2$,  we also have two cases
for the $\wm_f(R)$ to be considered, according to whether $R$
lies above or below the invariant segment \GUc.

\medskip

\noindent{\bf The backward fast wave curve $\wm_f(R)$ for $R$ in subregion $\Theta_3^a$.}
\medskip

\begin{figure}[ht]
	\centering
	\subfigure[Admissible $\wm_f(R)$.]
	{\includegraphics[width=0.4\linewidth]{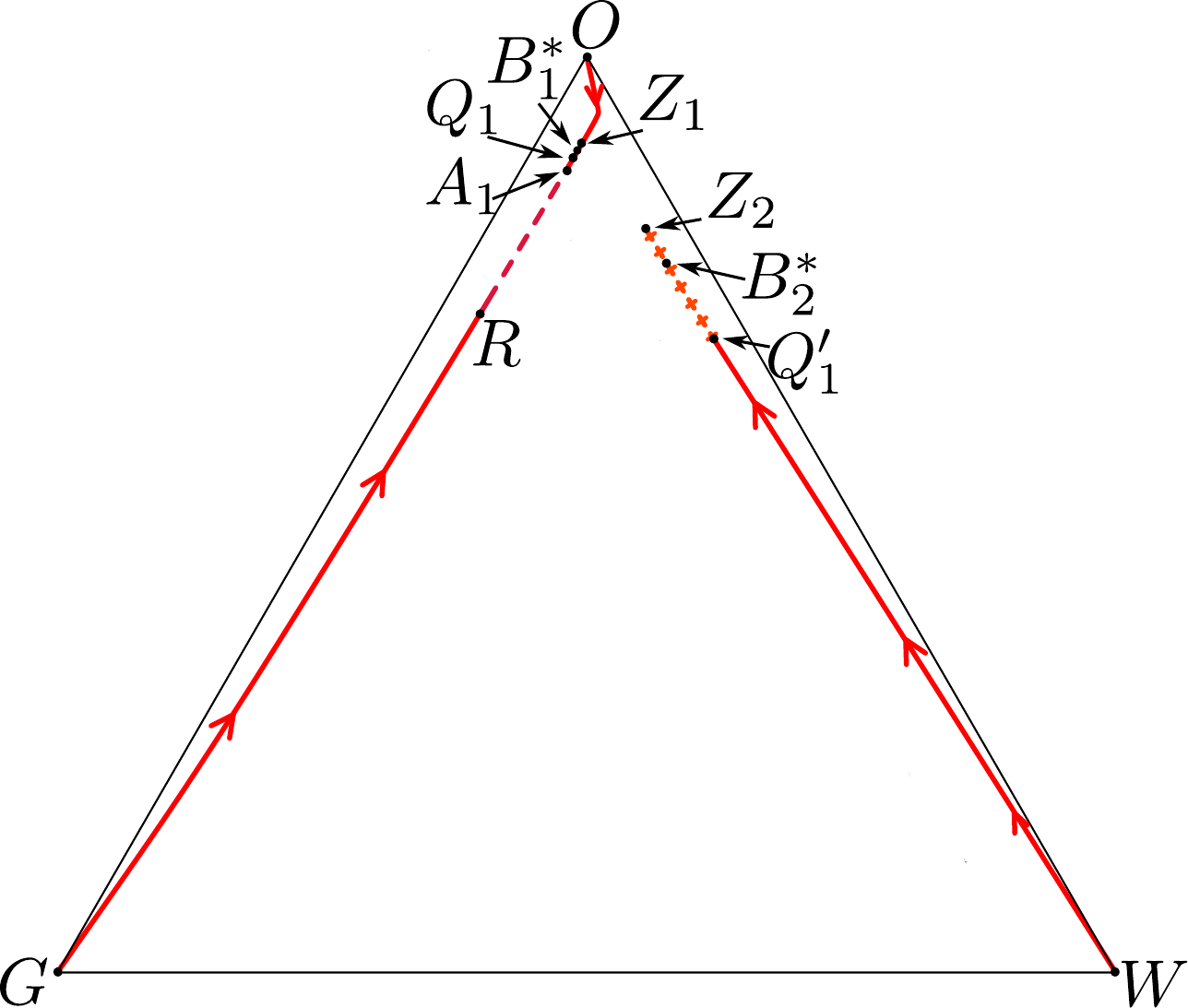}}  
	\hspace{1.75mm}
	\subfigure[Riemann solution. ]
	{\includegraphics[width=0.4\linewidth]{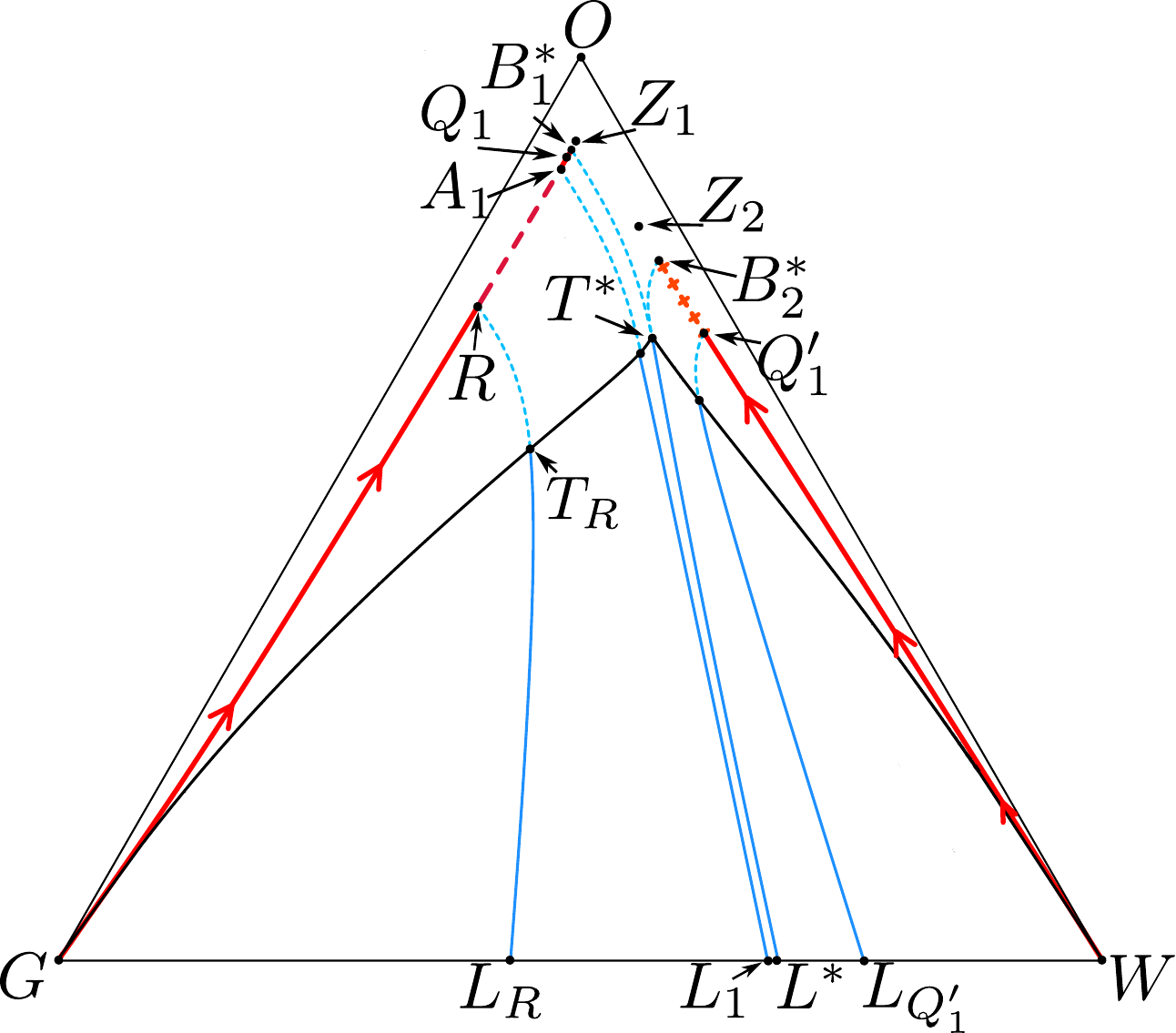}}  
 	\caption{(a) Admissible $\wm_f(R)$;  (b) Riemann solution for $R\in \Theta_3^a$. 
 	 The curve $[G,T^*]_{\text{ext}}\cup[T^*,W]_{\text{ext}}$
is an extension of the $\wm_f(R)$ segment $[G,B_1^*]\cup[B_2^*,W]$ such that
$\sigma(M', M) = \lambdas(M')$ for all $M' \in [G, T^*]_{\text{ext}}\cup [T^*, W]_{\text{ext}}$ and $ 
M \in [G, B_1^*]\cup [B_2^*, W]$. For this case, we consider $R = (0.0401428,0.722888).$ }
 	\label{fig:RS-R12}
\end{figure}

For the state $R$ in the subregion $\Theta_3^a$ the endpoint $A_1$ of the $f$-shock segment $(R, A_1]$ in the attached branch $(R, O_1]$
of $\mathcal{H}(R)$ lies between the $f$-inflection segment $\m{I}_f$ and the 
double contact segment $C$ ( i.e., in $\Theta_2^a$).
Since $\sigma(A_1; R) = \lambdaf(A_1)$, the
$f$-shock segment $(R, A_1]$ of the
$\wm_f(R)$ is followed by
the $f$-rarefaction segment $(A_1, O]$.
In the other direction, $\wm_f(R)$ comprises the
$f$-rarefaction segment $(R, G]$.
The $f$-rarefaction segment $(A_1, O]$ intersects
the double contact segment $C$
at a state $Q_1$ and the invariant segment \EU
at a state $Z_1$. As shown in  Fig.~\ref{fig:RS-R12}(a),
there are states $Z_2$, in the invariant segment 
\WUc, and  $Q_1'$ corresponding to  the states $Z_1$ and $Q_1$ such that
$\sigma(Z_2;Z_1) = \lambdaf(Z_1)$ and $\sigma(Q_1';Q_1) = \lambdaf(Q_1) = \lambdaf(Q_1')$. Notice that $Q_1'$ lies
in a part of the double contact locus not displayed in the figures.
Consequently the $f$-rarefaction segment
$[Z_1, Q_1]$ defines a $f$-composite segment
$[Z_2, Q_1']$ such that $\sigma(M';M) = \lambdaf(M)$
for all states $M'$ along $[Z_2, Q_1']$
and $M$ along $[Z_1, Q_1]$.
Since $\sigma(Q_1';Q_1) = \lambdaf(Q_1')$,
from $Q_1'$ the $\wm_f(R)$ is continued
by the $f$-rarefaction segment $(Q_1', W]$.

In summary, we have:
\begin{cla}\label{cla:BWCIF1IF2Y4Y3hat}
	(Refer to Fig.~\ref{fig:RS-R12}(a).)
	For $R$ in subregion $\Theta_3^a$
	in Figs.~\ref{fig:RegioesThetas}(a), \ref{fig:RRegions Completed2}(a),
	$\wm_f(R)$ comprises states
	$M$ along the $f$-rarefaction segments
	$[G, R)$, $[O,A_1)$ and $[W, Q_1')$, states $M$
	along the $f$-shock segment $(R, A_1]$,
	and states $M'$ along the $f$-composite segment
	$[Q_1', Z_2]$ such that for each state 
	$M'\in [Q_1', Z_2]$ there is a unique state 
	$M \in [Q_1, Z_1]$ with $\sigma(M';M) = \lambdaf(M)$.
\end{cla}

\medskip

\noindent{\bf The backward fast wave curve $\wm_f(R)$ for $R$ in subregion $\Theta_3^b\cup \Theta_3^c$.}
\medskip

\begin{figure}[ht]
	\centering
	\subfigure[Admissible $\wm_f(R)$.]
	{\includegraphics[width=0.4\linewidth]{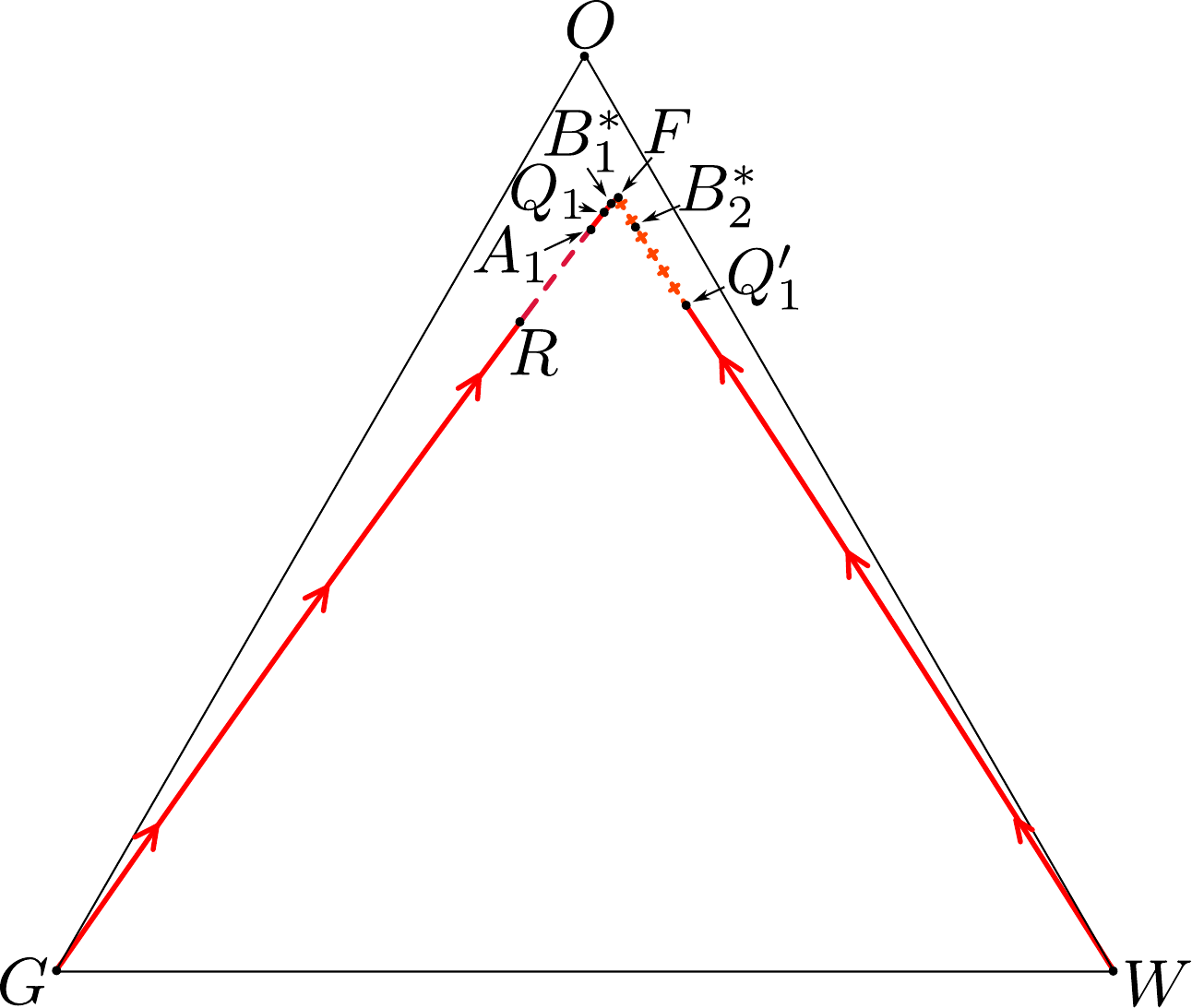}}  
	\hspace{1.75mm}
	\subfigure[Riemann solution. ]
	{\includegraphics[width=0.4\linewidth]{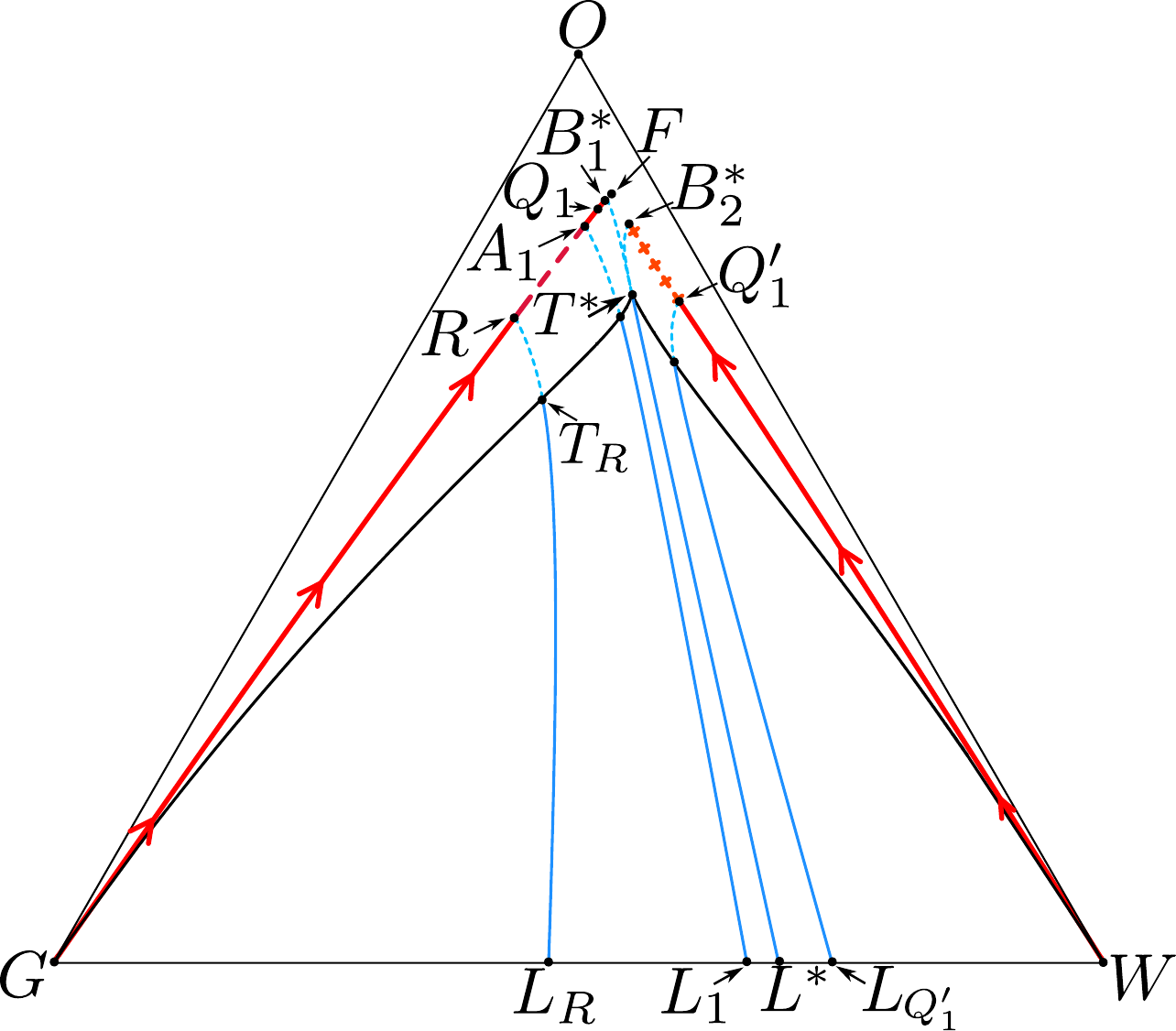}}  
 	\caption{(a) $\wm_f(R)$; (b) Riemann solution for $R\in \Theta_3^b$.
 	 The curve $[G,T^*]_{\text{ext}}\cup[T^*,W]_{\text{ext}}$
is an extension of the $\wm_f(R)$ segment $[G,B_1^*]\cup[B_2^*,W]$ such that
$\sigma(M', M) = \lambdas(M')$ for all $M'\in [G, T^*]_{\text{ext}}\cup [T^*, W]_{\text{ext}}$ and $
M \in [G, B_1^*]\cup [B_2^*, W]$. For this case, we consider $R = (0.0847167,0.709481)$. }
 	\label{fig:RS-R4_A4}
\end{figure}
\begin{figure}[ht]
	\centering
	\subfigure[Admissible $\wm_f(R)$.]
	{\includegraphics[width=0.4\linewidth]{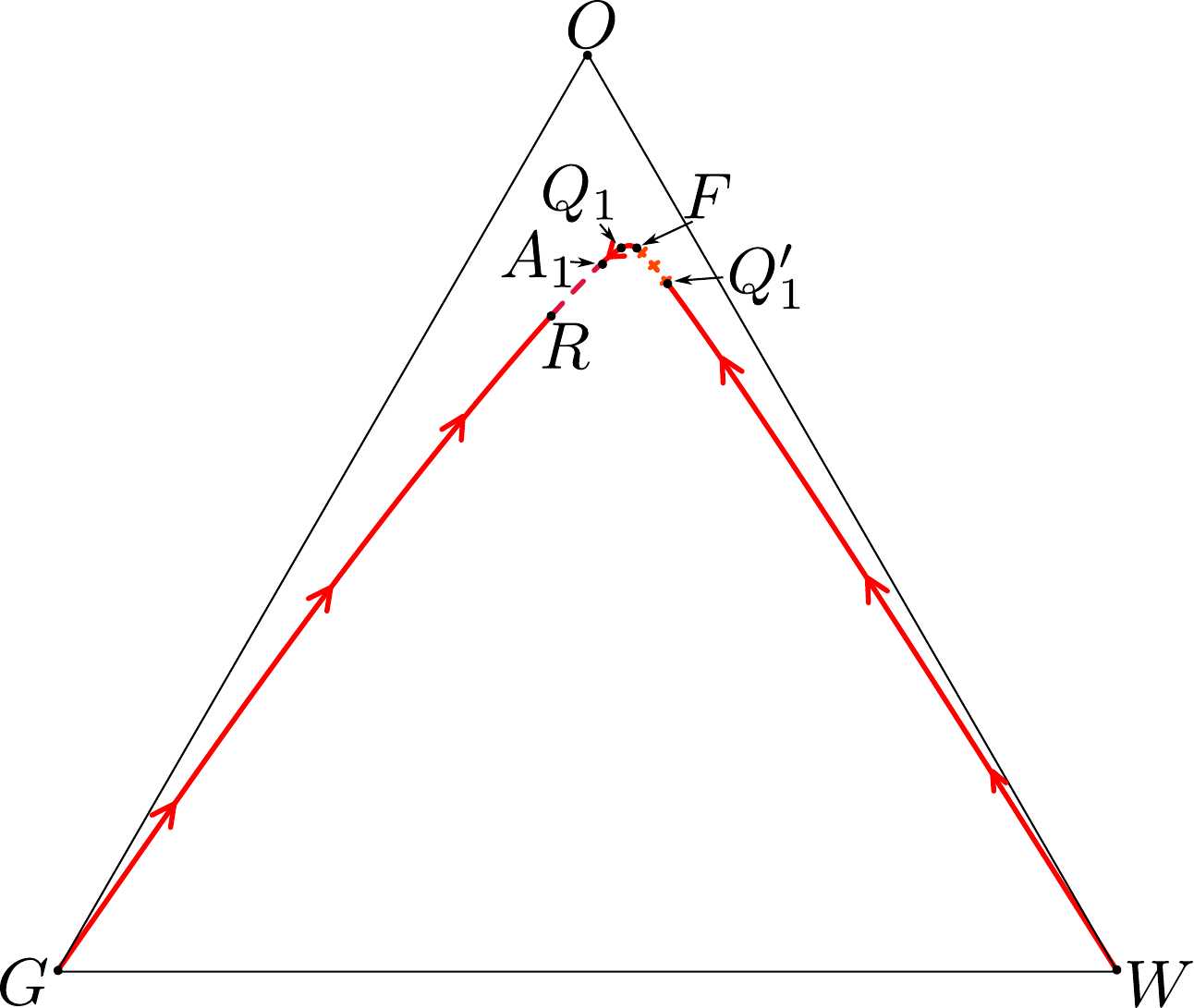}}  
	\hspace{1.75mm}
	\subfigure[Riemann solution.]
	{\includegraphics[width=0.4\linewidth]{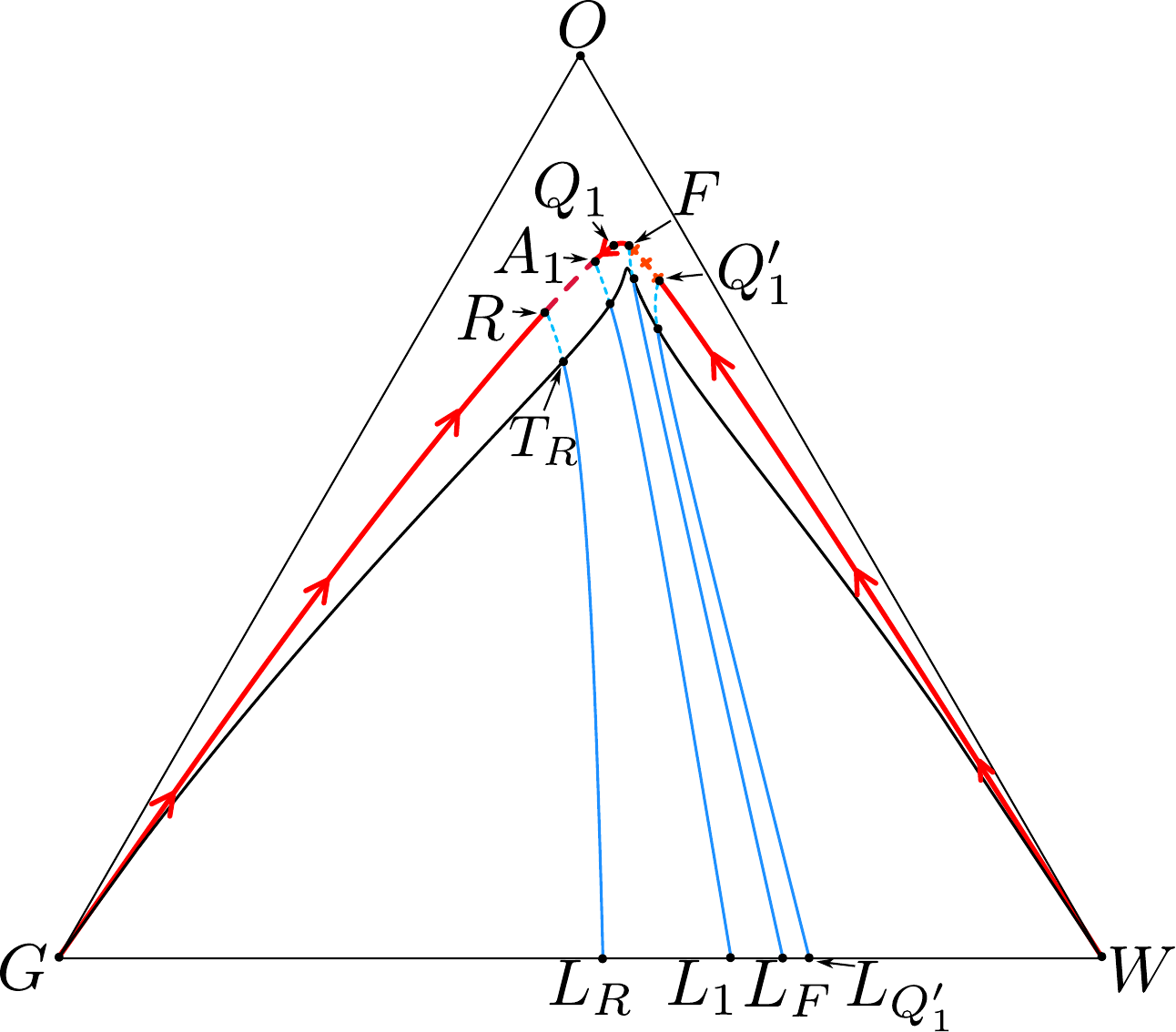}}  
 	\caption{Admissible $\wm_f(R)$ and Riemann solution for $R\in \Theta_3$ in subregion bounded by
 	 $\Theta_3^c$. 
 	 Curve $[G,T_1^*]_{\text{ext}}\cup[T_1^*,W]_{\text{ext}}$
is an extension of the $\wm_f(R)$ segment $[G,B_1^*]\cup[B_2^*,W]$ such that
$\sigma(M', M) = \lambdas(M')$, for all $(M', M) \in [G, T_1^*]_{\text{ext}}\cup [T_1^*, W]_{\text{ext}} 
\times [G, B_1^*]\cup [B_2^*, W]$. For this case, we consider $R = (0.109168,0.715865)$. }
 	\label{fig:RS-R4_A5}
\end{figure}

The behavior of $\wm_f(R)$ for
$R$ in subregion $\Theta_3^b\cup \Theta_3^c$ is similar to that for
$R$ in subregion $\Theta_2^b\cup \Theta_2^c$.
It differs from the case with
$R \in \Theta_3^a$
in that now instead of intersecting the invariant segment
\EUc, the $f$-rarefaction segment
through $A_1$
intersects the $f$-inflection segment $\um$-$K$ or $K$-$P$ depending on whether $R \in \Theta_3^b$ or
$R \in \Theta_3^c$.
We have:

\begin{cla}\label{cla:I2fI4fPhatY4hat}
	Refer to Figs.~\ref{fig:RS-R4_A4}(a) and
	\ref{fig:RS-R4_A5}(a).
	For $R$ in subregion
	$\Theta_3^b\cup \Theta_3^c$	in Figs.~\ref{fig:RegioesThetas}(a), \ref{fig:RRegions Completed2}(a)
	$\wm_f(R)$ comprises states
	$M$ along the $f$-rarefaction segments
	$[G, R)$, $[F, A_1)$ and $[W, Q_1']$, states $M$
	along the $f$-shock segment $[A_1, R)$,
	and states $M'$ along the $f$-composite segment
	$[Q_1', F)$ such that for each state 
	$M'\in [Q_1', F)$ there is a unique state 
	$M \in (F, Q_1]$ with $\sigma(M';M) = \lambdaf(M)$. 
\end{cla}

\medskip
\noindent{\bf The Riemann solution for $R$ in subregion $\Theta_3 = \Theta_3^a \cup \Theta_3^b \cup \Theta_3^c$.}

\medskip

\medskip
\noindent{\bf The Riemann solution for $R$ in subregion $\Theta_3^a \cup \Theta_3^b$.}

\medskip

We recall the intersection points $B_1^*$, and $Q_1$ of the $f$-rarefaction segment 
$[O, A_1)$ in $\wm_f(R)$ with 
the mixed contact segment $MC$, and the double contact segment $C$, respectively, which are displayed in Figs.~\ref{fig:RS-R12} and \ref{fig:RS-R4_A4}.
They play an important role in the construction of Riemann solutions as will be seen in the following.

\begin{cla}\label{cla:RSolution-RinR12-above}
Refer to Figs.~\ref{fig:RS-R12} and \ref{fig:RS-R4_A4}.
Let $L$ be a state on the edge \GW of the
saturation triangle and $R$ be a state in subregion $\Theta_3^a 
\cup \Theta_3^b$ in
Figs.~\ref{fig:RegioesThetas}(a), \ref{fig:RRegions Completed2}(a). 
Let $A_1$, $Q_1$, and $Q_1'$ be the states on
$\wm_f(R)$ such that $\sigma(A_1; R) = \lambdaf(A_1)$, $\lambdaf(Q_1)= \lambdaf(Q_1') = \sigma(Q_1';Q_1)$.
Let $T^*$, $B_1^*$ and $B_2^*$ be states satisfying
$\lambdaf(T^*) = \sigma(T^*; B_1^*) = \sigma(T^*; B_2^*) = \sigma(B_2^*; B_1^*) = \lambdaf(B_1^*)$.

Let  $L_R$, $L_1$, $L^*$ and $L_{Q_1'}$
be the intersection points of
the backward slow wave curves through
$R$, $A_1$,
$B_1^*$ (or $B_2^*$) and ${Q_1'}$ with the edge \GWc, respectively.
Then, 

\begin{itemize}
\item[(i)] if $L=G$, the Riemann solution is $L\testright{R_f} R$; 
\item[(ii)] if $L \in(G, L_R)$, 
the Riemann solution is \;
$
     L\testright{R_s} T_1 \xrightarrow{'S_s} M_1 \testright{R_f}  R
$, \; where $T_1\in(G, T_R)_{\text{ext}}$ and $M_1\in(G,R)$;
 \item[(iii)] if $L=L_R$, the Riemann
solution is \;
$
     L\testright{R_s} T_R \xrightarrow{'S_s} R
$;
 \item[(iv)] if $L \in(L_R, L_1]$ the Riemann solution is \;
$
     L\testright{R_s} T_1 \xrightarrow{'S_s} M_1 \testright{S_f}  R
$, \; where $T_1\in(T_R, T^*)_{\text{ext}}$ and $M_1\in(R, A_1)$;
\item[(v)] if $L \in(L_1, L^*)$ the Riemann solution is \;
$
     L\testright{R_s} T_1 \xrightarrow{'S_s} M_1 \testright{R_f}  A_1 \xrightarrow{'S_f} R
$, \; where $T_1\in(T_R, T^*)_{\text{ext}}$ and $M_1\in(A_1, B_1^*)$;
\item[(vi)] if $L \in (L^*,L_{Q_1'}]$ the Riemann solution is \;
$
     L\testright{R_s} T_1 \xrightarrow{'S_s} M_1 \testright{S_f'} M_1' \xrightarrow{R_f}  A_1 \xrightarrow{'S_f} R
$, \; where $T_1\in(T^*, W)_{\text{ext}}$, $M_1\in(B_2^*, Q_1')$ and $M_1'\in(B_1^*, Q_1)$;
\item[(vii)] if $L \in[L_{Q_1'},W)$ the solution is
$
     L\testright{R_s} T_1 \xrightarrow{'S_s} M_1 \testright{R_f} Q_1' \xrightarrow{'S_f'} Q_1 \xrightarrow{R_f}  A_1 \xrightarrow{'S_f} R
$, \; where $T_1\in(T^*, W)_{\text{ext}}$ and $M_1\in(Q_1', W)$;
\item[(viii)] If  $L=W$, the solution is
$
     L \testright{R_f} Q_1' \xrightarrow{'S_f'} Q_1 \xrightarrow{R_f}  A_1 \xrightarrow{'S_f} R
$.
\end{itemize}
\begin{rem}\label{rem:L*B*3}
If $L = L^*$, see Figs.~\ref{fig:RS-R12} and \ref{fig:RS-R4_A4}, there are two paths to reach the right state
$R$, namely 
$
     L^*\testright{R_s} T^* \xrightarrow{'S_o'} B_1^* \xrightarrow{R_f}  A_1 \xrightarrow{'S_f} R
$ or 
$
     L^*\testright{R_s} T^* \xrightarrow{'S_s} B_2^* \xrightarrow{S_f'} B_1^* \xrightarrow{R_f}  A_1 \xrightarrow{'S_f} R
$
but the triple shock rule guarantees they represent the same solution in $xt$-space.
In other words, the solution consists of a single wave group.
\end{rem}

\end{cla}

\medskip
\noindent{\bf The Riemann solution for $R$ in subregion $\Theta_3^c$.}

\medskip

For $R$ in subregion $\Theta_3^c$, 
the $f$-rarefaction segment $(R, F]$ 
no longer intersects the mixed contact segment $MC$, but continues to intersect
the double contact segment $C$
at a state $Q_1$. 

In summary, we have:
\begin{cla}\label{cla:RSolution-RinI3FI4fPhatY5hat}
Refer to Fig.~\ref{fig:RS-R4_A5}(b).
Let $L$ be a state on the edge \GW of the
saturation triangle and $R$ be a state in subregion $\Theta_3^c$ in
Figs.~\ref{fig:RegioesThetas}(a), \ref{fig:RRegions Completed2}(a).
Let $A_1$, $Q_1$, $Q_1'$, and $F$ be the states in $\wm_f(R)$ with
$\lambdaf(A_1)=\sigma(A_1;R)$,
$\lambdaf(Q_1')=\sigma(Q_1';A_1) = \lambdaf(Q_1)$, and $F$ on the $f$-inflection segment $K$-$P$ in Fig~\ref{fig:RRegions Completed2}(c).

Let  $L_R$, $L_1$, $L_F$ and $L_{Q_1'}$
be the intersection points of
the backward slow wave curves through
$R$, $A_1$,
$F$ and $Q_1'$ with the edge \GWc, respectively.
Then, 

\begin{itemize}
\item[(i)]if $L=G$, the Riemann solution is $L\testright{R_f} R$. 
\item[(ii)] if $L \in(G, L_R)$, the Riemann
solution is $
     L\testright{R_s} T_1 \xrightarrow{'S_s} M_1 \testright{R_f}  R
$, \; where $T_1\in(G, T_R)$ and $M_1\in(G,R)$;
 \item[(iii)] if $L=L_R$, the Riemann
solution is \;
$
     L\testright{R_s} T_R \xrightarrow{'S_s} R
$;
 \item[(iv)] if $L \in(L_R, L_1]$ the Riemann solution is \;
$
     L\testright{R_s} T_1 \xrightarrow{'S_s} M_1 \testright{S_f}  R
$, \; where $T_1$ lies on the slow extension of the $f$-shock segment $(R, A_1]$ and $M_1\in(R, A_1]$;
\item[(v)] if $L \in(L_1, L_F]$ the Riemann solution is \;
$
     L\testright{R_s} T_1 \xrightarrow{'S_s} M_1 \testright{R_f}  A_1 \xrightarrow{'S_f} R
$, \;
where $T_1$ lies on the slow extension of the $f$-rarefaction segment $(A_1, F]$ and $M_1\in(A_1, F]$;
\item[(vi)] if $L \in(L_F,L_{Q_1'})$
the Riemann solution is \;
$
     L\testright{R_s} T_1 \xrightarrow{'S_s} M_1 \testright{S_f'} M_1' \xrightarrow{R_f}  A_1 \xrightarrow{'S_f} R
$, \; where $T_1$ lies on the slow extension of the $f$-composite segment $(F,Q_1']$, $M_1\in (F, Q_1']$ and  $M_1'\in(F, Q_1)$;
\item[(vii)] if $L \in[L_{Q_1'},W)$ the solution is
$
     L\testright{R_s} T_1 \xrightarrow{'S_s} M_1 \testright{R_f} Q_1' \xrightarrow{'S_f'} Q_1 \xrightarrow{R_f}  A_1 \xrightarrow{'S_f} R
$, \; where $T_1$ lies on the slow extension of the $f$-rarefaction segment $(Q_1', W)$ and $M_1\in(Q_1', W)$;

\item[(viii)] if  $L=W$, the solution is
$
     L \testright{R_f} Q_1' \xrightarrow{'S_f'} Q_1 \xrightarrow{R_f}  A_1 \xrightarrow{'S_f} R
$.
\end{itemize}
\end{cla}




\subsection{Right state in region \texorpdfstring{$\Omega$}{O}}
\label{subsec:Omega}
We divide region $\Omega$ into subregions $\Omega_1$  and $\Omega_2$; see Fig~\ref{fig:RegioesThetas}(a). These subregions are defined below. 

 \begin{figure}[]
	\centering
	\subfigure[R-Region $\Omega_1$. ]
	{\includegraphics[scale=0.6]{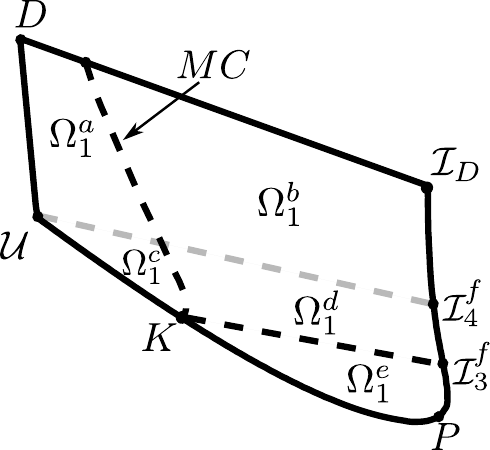}}  
	 \hspace{6mm}
	\subfigure[R-Region $\Omega_2$.]
	{\includegraphics[scale=0.5]{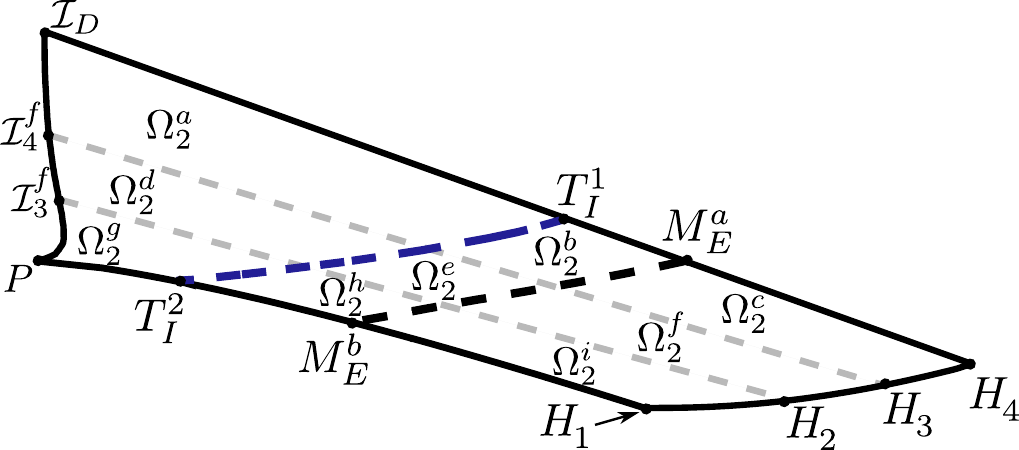}}   
	\caption{Zoom of the $R$-regions $\Omega_1$ and $\Omega_2$
 in Fig.~\ref{fig:RegioesThetas}
	displaying in detail their subdivisions.}
		\label{fig:RRegionsOmegas}
\end{figure}

\subsubsection{Subregion \texorpdfstring{$\Omega_1$}{O}}
The boundary of subregion $\Omega_1$ in Fig~\ref{fig:RegioesThetas}(b)
consists of the following segments: the invariant segment \DUc, the segments $\um$-$P$, $P$-$\mathcal{I}_D$ of the $f$-inflection locus, and the segment $\mathcal{I}_D$-$D$ of the edge \WOc.

The region $\Omega_1$ is subdivided into five subregions: $\Omega_1^i$, with $i\in \{a,b,c,d,e\}$, see Fig.~\ref{fig:RRegionsOmegas}(a). These subregions are defined by
the mixed double contact segment $MC$, the invariant segment $\um$-$\mathcal{I}_4^f$, and
the $f$-rarefaction segment $K$-$\mathcal{I}_3^f$.

The Hugoniot locus $\mathcal{H}(R)$ and the fast-backward wave curve $\wm_f(R)$ for right states $R$ in regions $\Omega_1^i$ are similar in structure to $\mathcal{H}(R)$ and $\wm_f(R)$ for $R$ in $\Theta_1^i$, $i \in \{a, b, c, d,e\}$. As a consequence, the Riemann solutions for $R$ in $\Omega_1^a \cup \Omega_1^c$, $\Omega_1^b \cup \Omega_1^d$,
$\Omega_1^e$ are similar to the Riemann solutions for $R$ in $\Theta_1^a \cup \Theta_1^c$, $\Theta_1^b \cup \Theta_1^d$,
$\Theta_1^e$, respectively.

\subsubsection{Subregion \texorpdfstring{$\Omega_2$}{O}}

The boundary of subregion $\Omega_2$ in Fig~ \ref{fig:RRegionsOmegas}(b)
consists of the following segments: the segment $\mathcal{I}_D$-$P$ of the $f$-inflection locus, the $f$-Hysteresis segment $P$-$H_1$, the segment of $s$-Hysteresis $H_1$-$H_4$, and the segment $H_4$-$\mathcal{I}_D$
on the edge \WOc.
The region $\Omega_2$ is further subdivided into nine subregions: $\Omega_2^\alpha$, with $\alpha\in \{a, b, c, d, e, f, g, h, i \}$ which will be detailed below.

\begin{figure}[ht]
	\centering
	\subfigure[Hugoniot for generic state $R$ at $\Omega_2^a$. $R=(0.271633, 0.711087)$.]
	{\includegraphics[scale=0.27]{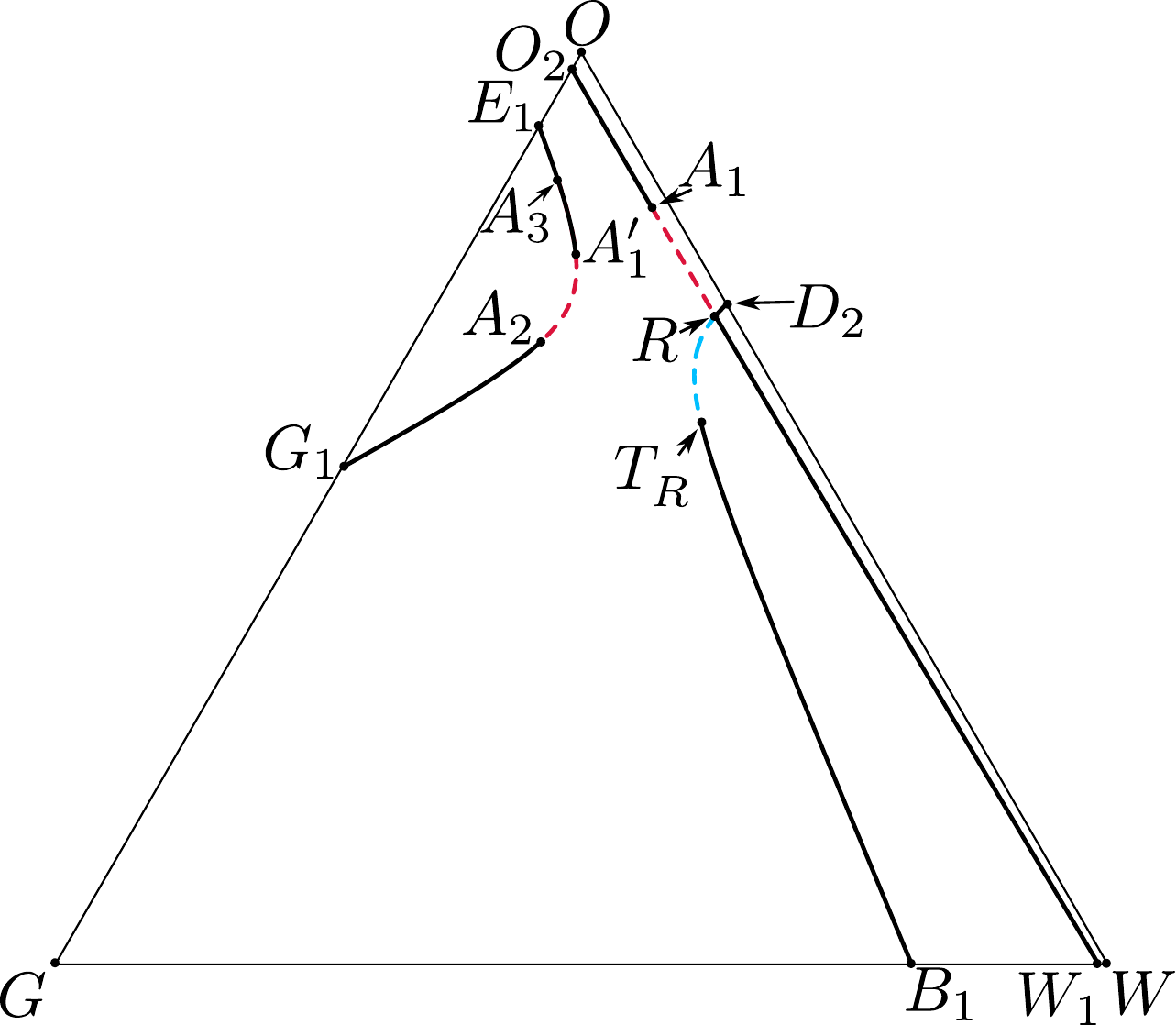}}  
	\hspace{1.75mm}
	\subfigure[Hugoniot for generic state in the region between the curves $EB$ and $M_E$ ,$\Omega_2^b$. $R=(0.29452, 0.684057)$.]
	{\includegraphics[scale=0.27]{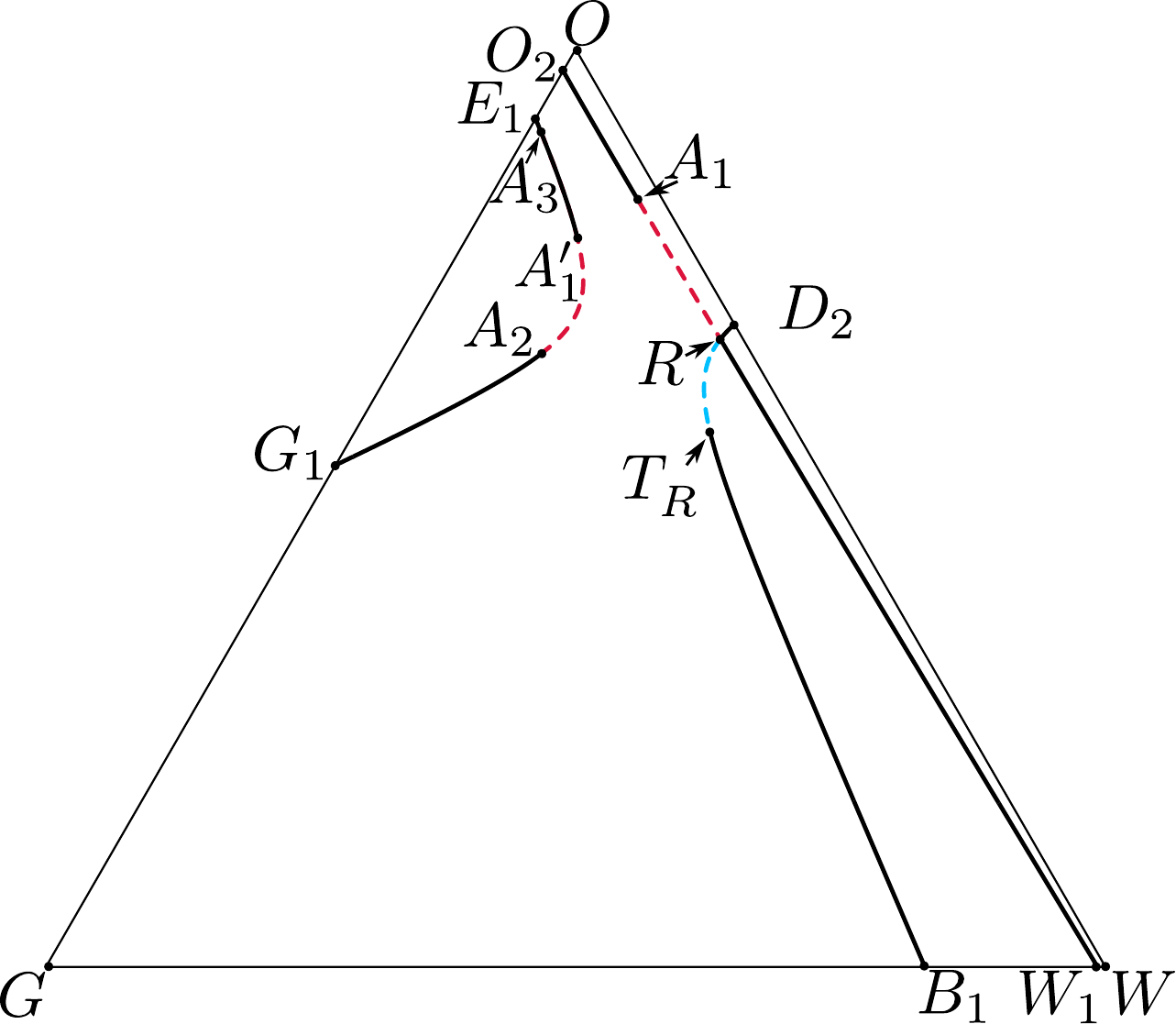}}  
 \caption{Hugoniot for generic states $R$ in $\Omega_2$}
 	\label{fig:Hugoniot_Omega_2A}
\end{figure}

\begin{figure}[ht]
	\centering
     \subfigure[Hugoniot for generic state $R$ at $\Omega_2^d$. $R=(0.230643, 0.731226)$.]
	{\includegraphics[scale=0.27]{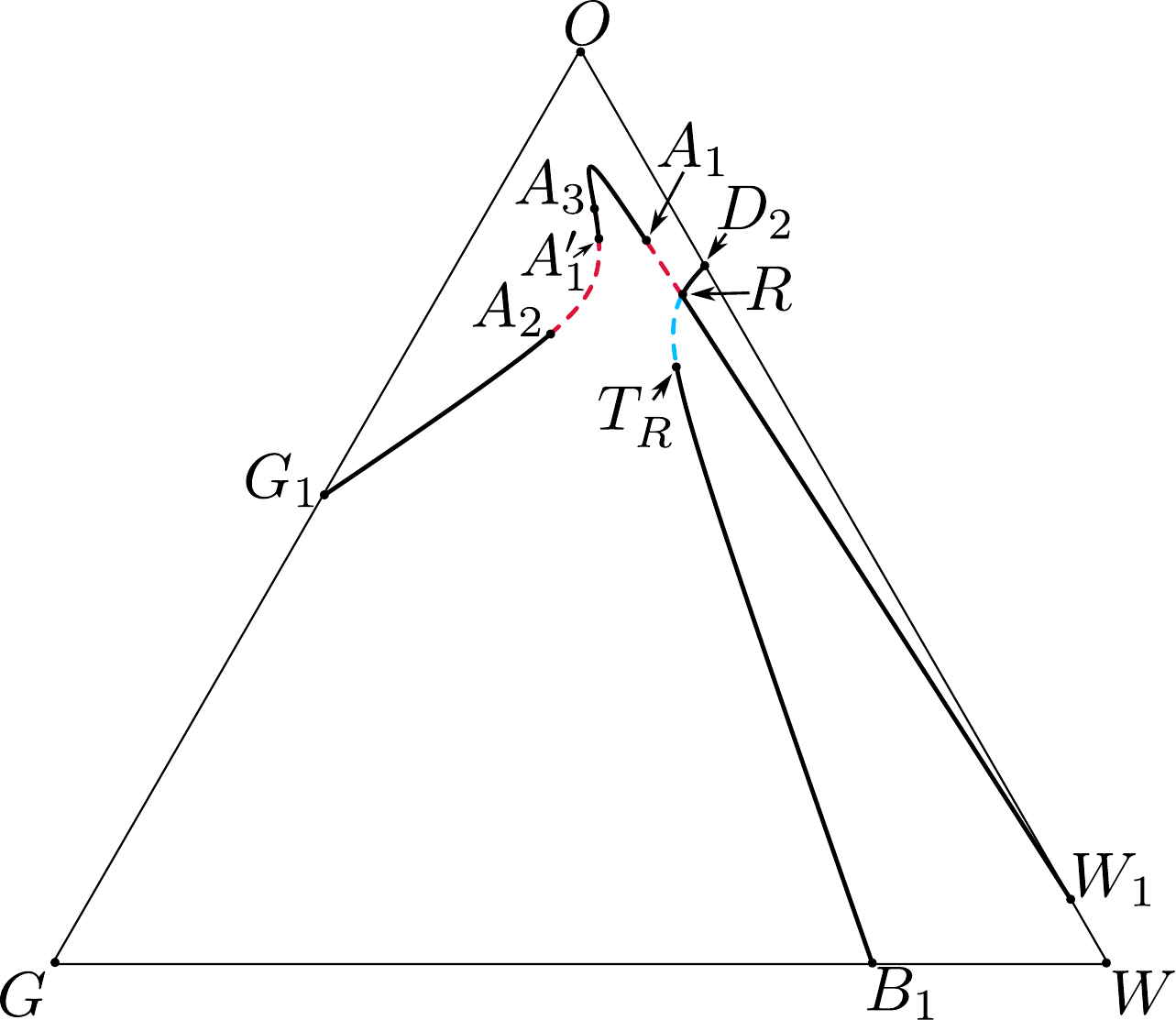}}  
	\hspace{1.75mm}
	\subfigure[Hugoniot for generic state in the region between the curves $EB$ and $M_E$, $\Omega_2^e$. $R=(0.258922, 0.704903)$.]
	{\includegraphics[scale=0.27]{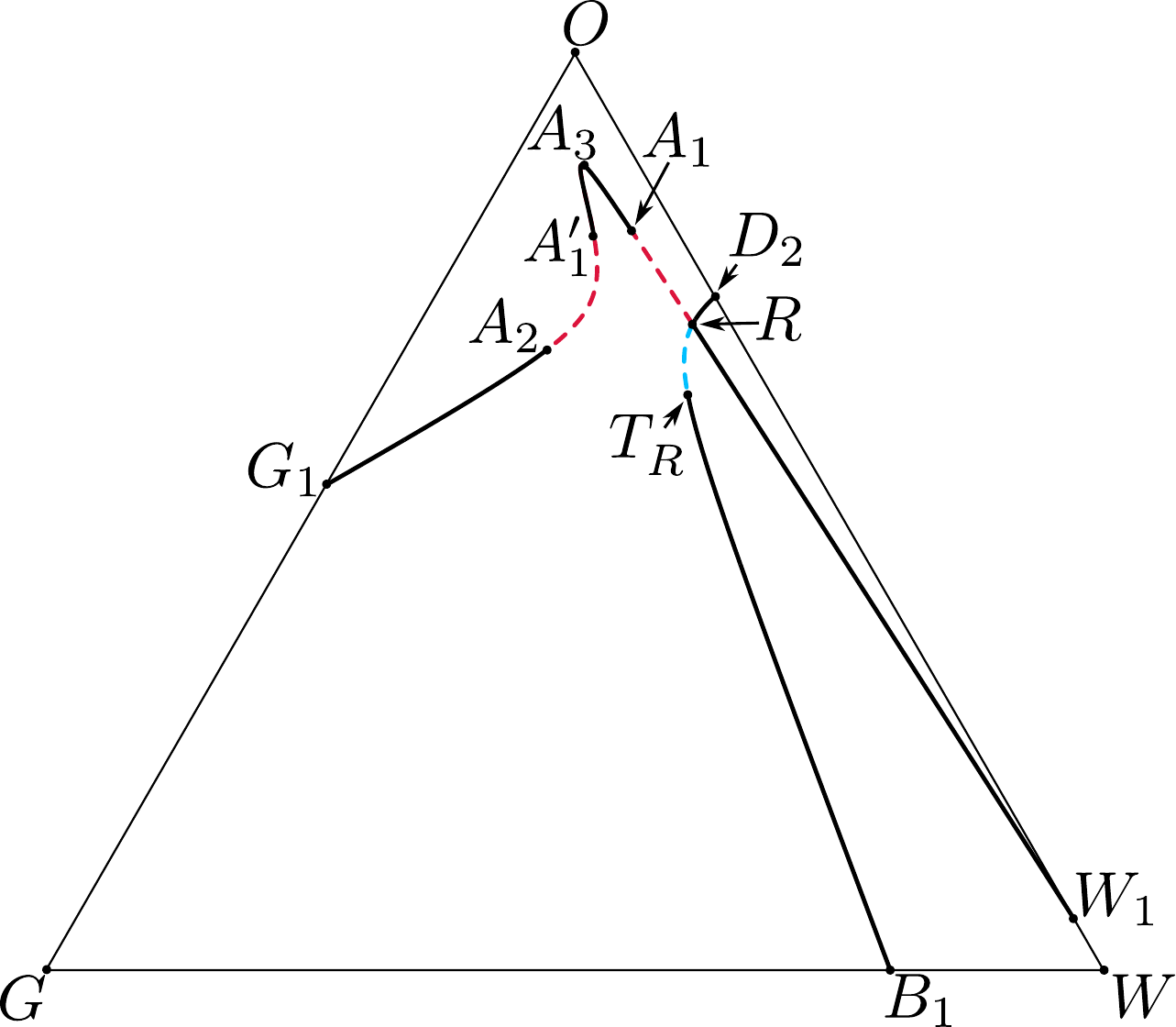}}  
 	\caption{Hugoniot for generic states $R$ in $\Omega_2$.}
 	\label{fig:Hugoniot_Omega_2B}
\end{figure}

\begin{figure}[ht]
	\centering
	\subfigure[Wave curve for a generic state in the region between the curves $EB$ and $M_E$, $\Omega_2^b$. $R=(0.29452, 0.684057)$.]
	{\includegraphics[scale=0.27]{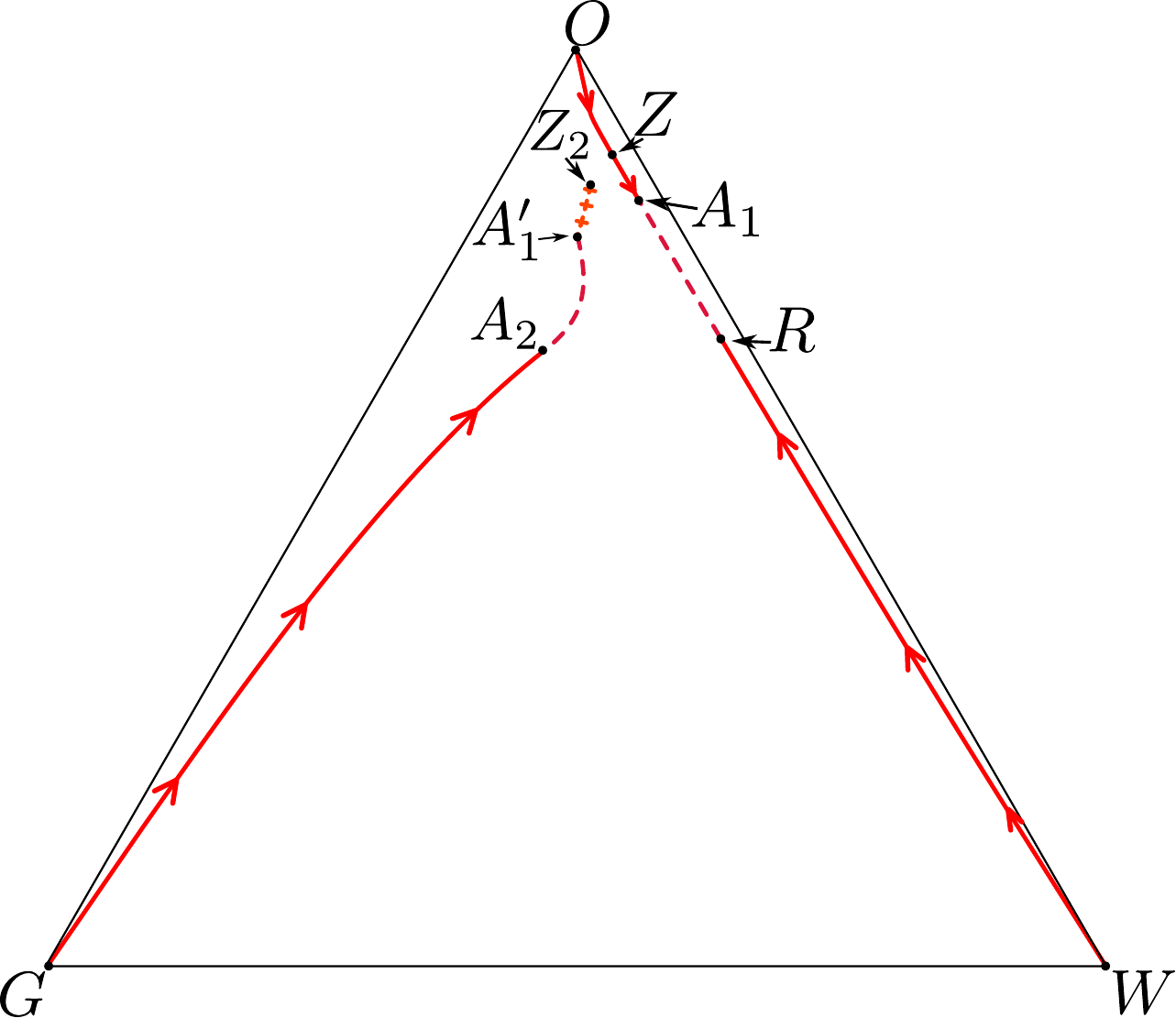}}  
	\hspace{1.75mm}
     \subfigure[Riemann solution for a generic state in the region between the curves $EB$ and $M_E$, $\Omega_2^b$. $R=(0.29452,  0.684057)$.
      The segment ${[A_1', B_2^*]}$ is a composite defined by the rarefaction ${[A_1, B_1^*]}$.]
	{\includegraphics[scale=0.27]{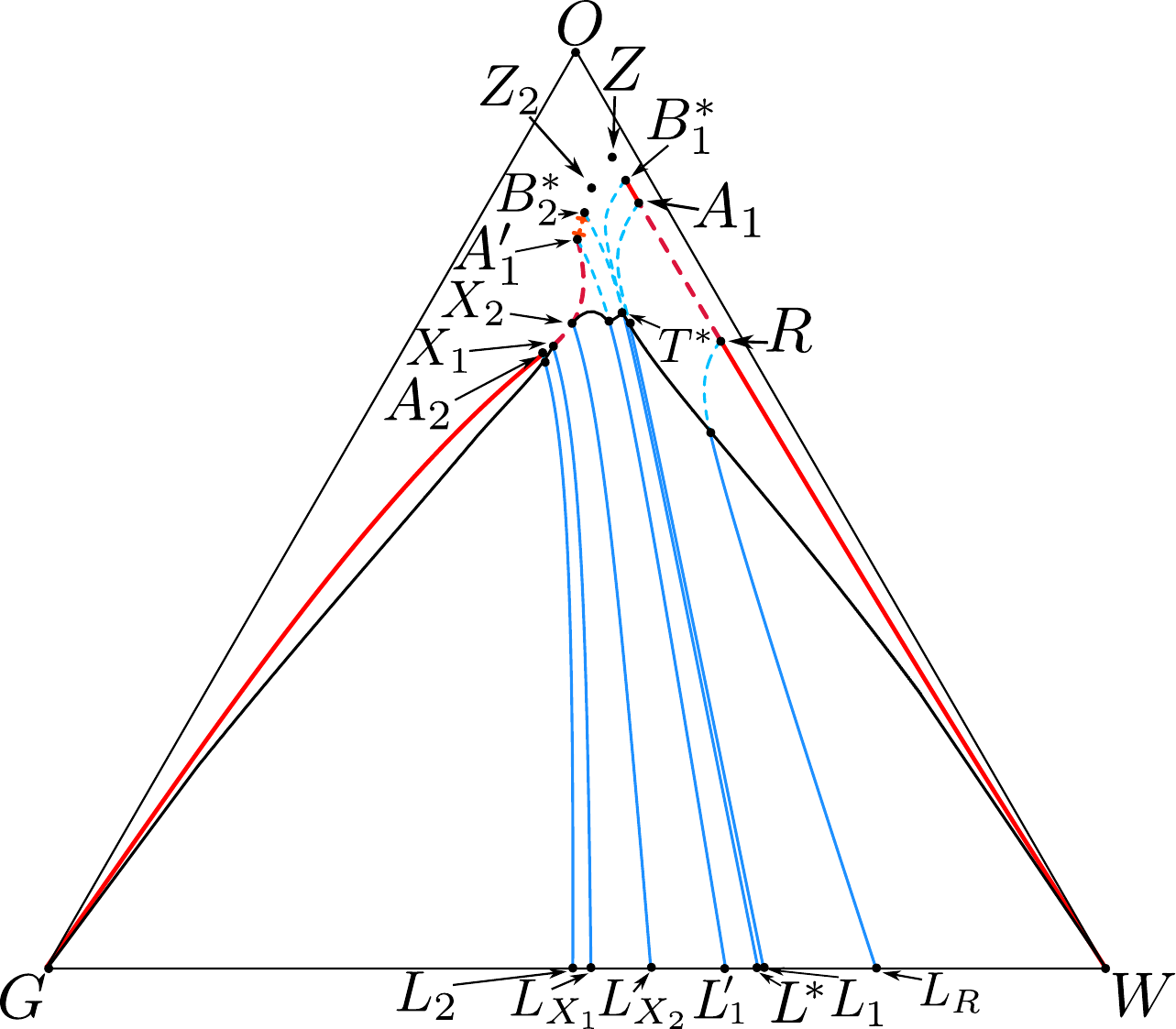}}  
	\hspace{1.75mm}
	 \caption{Riemann solution for a generic states $R$ in $\Omega_2$ }
 	\label{fig:Wave_Curve_2C}
\end{figure}

\begin{figure}[ht]
	\centering
	\subfigure[Wave curve for a generic state in region $\Omega_2^e$, between curves $EB$ and $M_E$. $R=(0.258922,0,704903)$.]
	{\includegraphics[scale=0.27]{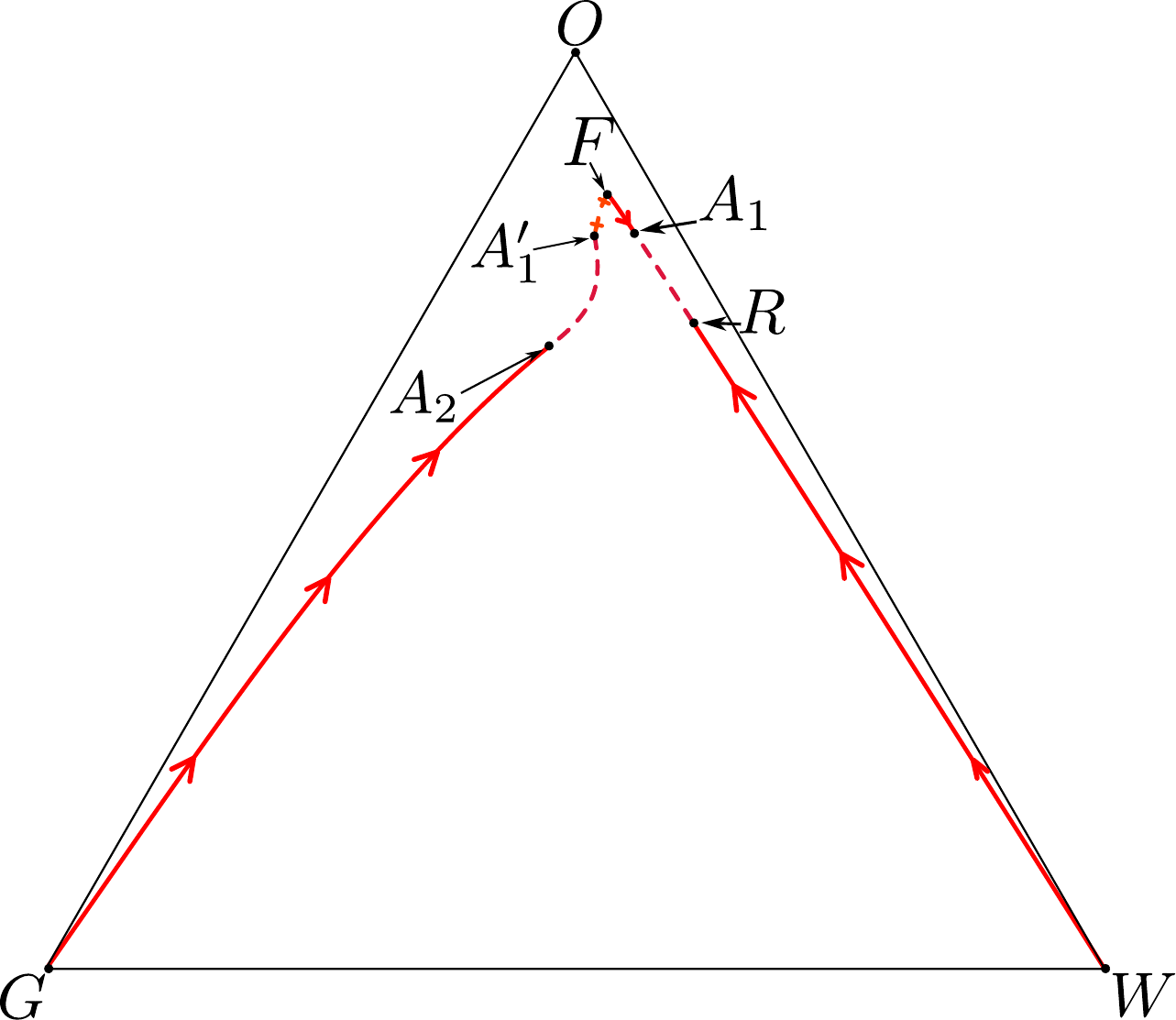}}  
	\hspace{1.75mm}
     \subfigure[Riemann solution for a generic state in region $\Omega_2^e$, between curves $EB$ and $M_E$. $R=(0.258922,0,704903)$.]
	{\includegraphics[scale=0.27]{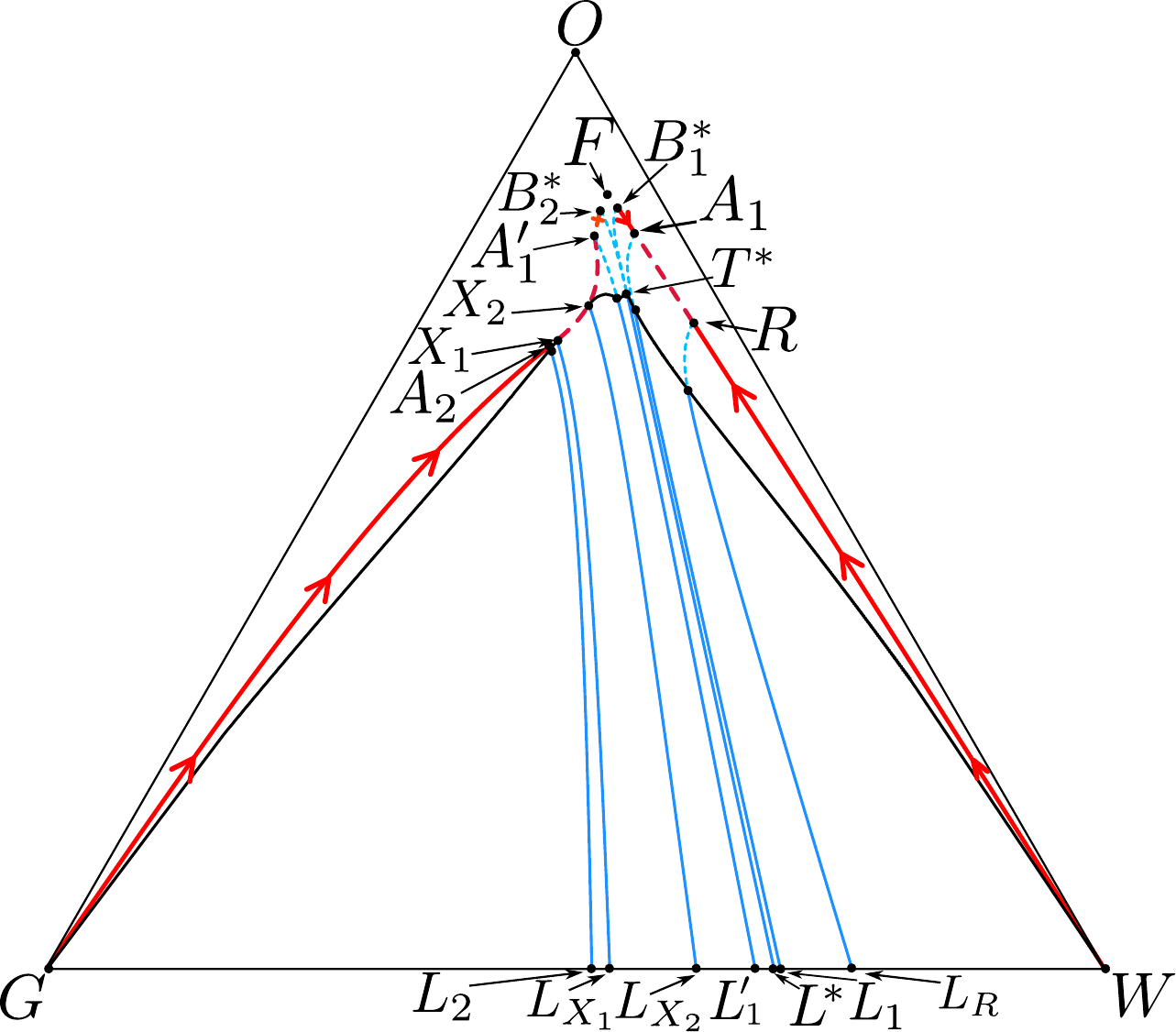}}  
	\hspace{1.75mm}
	 \caption{Riemann solution for a generic states $R$ in $\Omega_2$ }
 	\label{fig:Wave_Curve_2F}
\end{figure}

In region $\Omega_2$, the Hugoniot curves and the $f$-wave curves are of two types distinguished by their topological structure.
These distinct topological structures are caused by the position of state $R$ relative
to the invariant segment $\mathcal{I}_4^f$-$H_3$,
as shown in Figs.~\ref{fig:Hugoniot_Omega_2A}(b), \ref{fig:Hugoniot_Omega_2B}(b) for Hugoniot curves, and in Figs.~\ref{fig:Wave_Curve_2C}(a), \ref{fig:Wave_Curve_2F}(b) for $f$-wave curves.

\begin{figure}[ht]
	\centering
     \subfigure[$M\in (A_2,A_1')$  and $R=(0.271633, 0.711087)$.]
	{\includegraphics[scale=0.315]{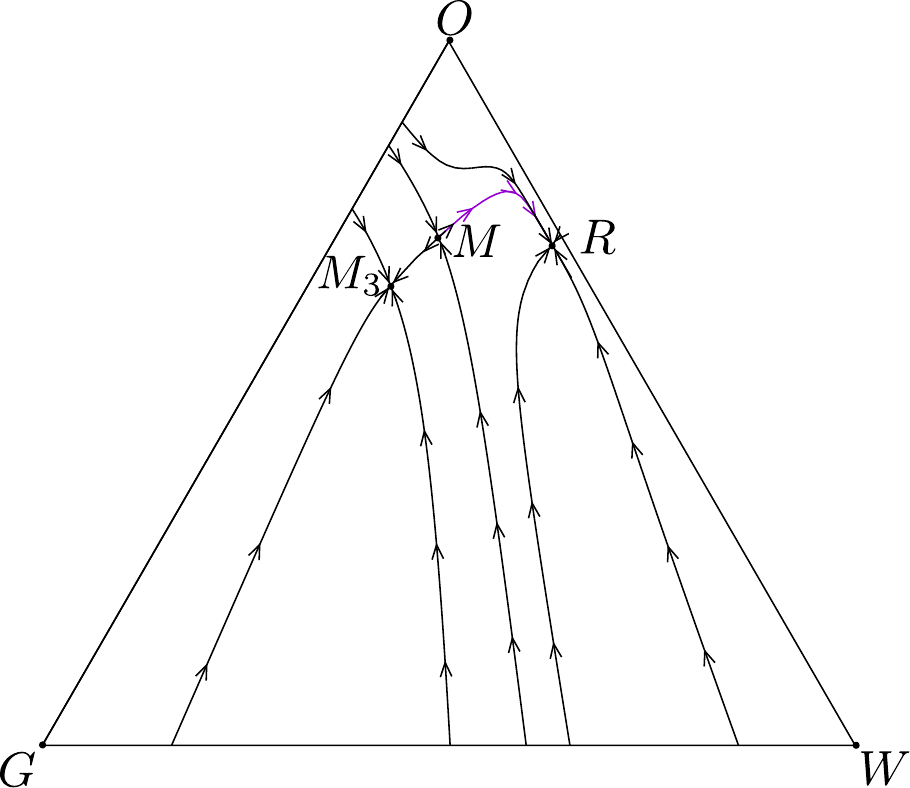}}  
	\hspace{1.75mm}
	\subfigure[$M=A_1'$  and $R=(0.271633, 0.711087)$.]
	{\includegraphics[scale=0.315]{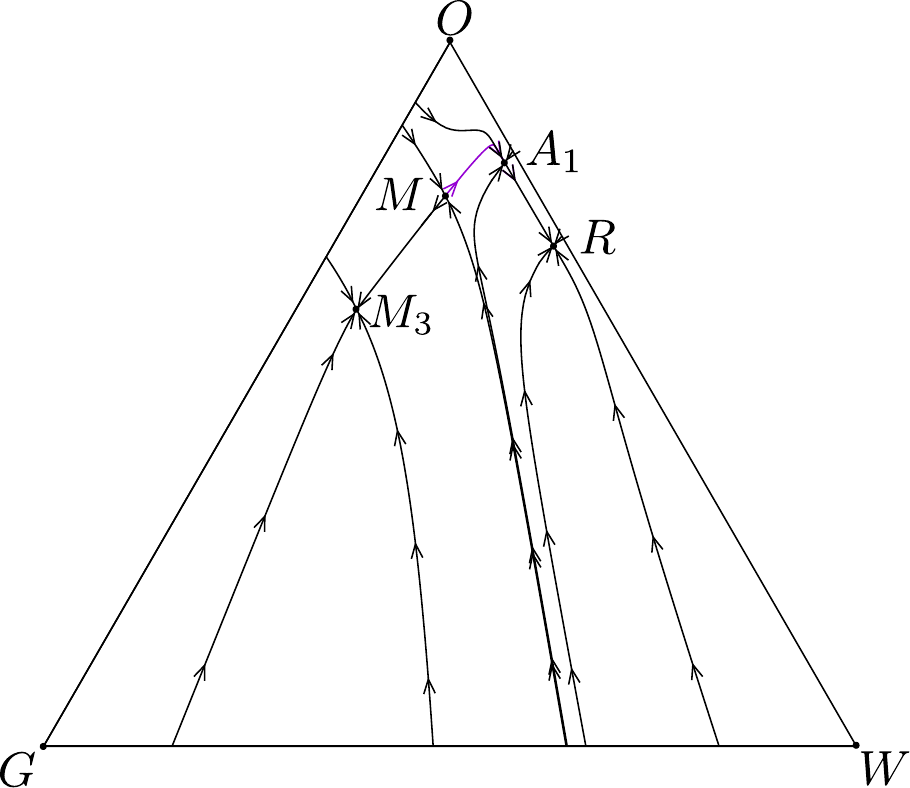}}  
 \subfigure[$M\in (A_1',A_3)$  and $R=(0.271633, 0.711087)$.]
	{\includegraphics[scale=0.315]{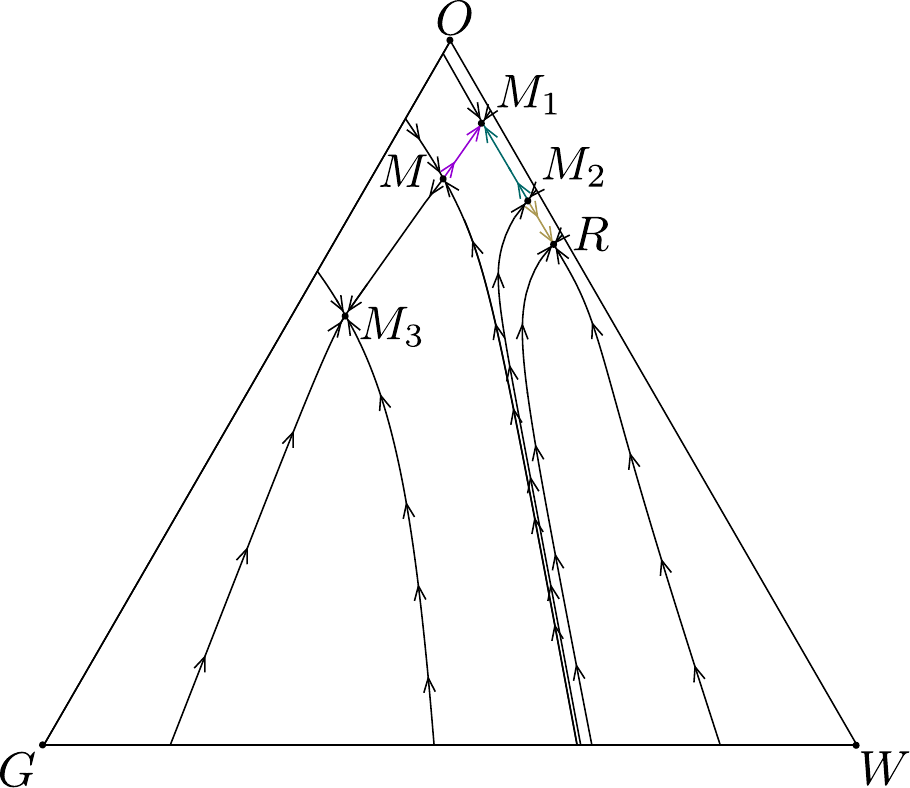}}  
 	\caption{ Some orbits of the ODE system (4.3) for M in the detached segment $(A_2, A_3)$ of $\mathcal{H}(R)$ in Fig. 3(b). The discontinuity joining state $M$ to state $R$ has a viscous profile (i.e. is admissible)
    only if $M$ is in  $[A2,A_1')$.}
 	\label{fig:dynamical_system_A_1p}
\end{figure}

The Hugoniot curves contain the $f$-shock segments $(R,A_1]$ and $[A_2,A_3]$,
    with $\lambda_f(A_1)=\sigma(A_1;R)$, $\sigma(A_2;R) = \lambdaf(A_2)$, and $\sigma(A_3;R) = \lambdaf(R)$. 
    The segment $(R, A_1]$ is admissible. On the other hand only the portion $[A_2,A_1')$ of the segment $[A_2,A_3]$ is admissible, with the state $A_1'$ satisfying $\sigma(A_1;R)=  \sigma(A_1';A_1)=\sigma(A_1';R)$.
    From the viewpoint of dynamical systems, $A_1'$ is associated with a saddle-node bifurcation. 
    If $M \in [A_2, A_1')$ there is an orbit from
    $M$ to $R$. When $M = A_1'$ there is a
    saddle-node at $A_1$. If $M \in (A_1', A_3]$ the saddle-node splits into an attractor $M_1$
    and a saddle $M_2$; there is an orbit from $M$ to $M_1$, an orbit from $M_2$ to $M_1$, and an orbit from $M_2$ to $R$; see Fig.~\ref{fig:dynamical_system_A_1p}. 
A similar analysis of dynamical systems related to the admissibility of
Lax shock waves was performed in \cite{Andrade2018}.
Regarding the two types of $f$-wave curves
shown in Figs.~\ref{fig:Wave_Curve_2C}(a), \ref{fig:Wave_Curve_2F}(b), as in the region $\Theta$ their differences play no role in the Riemann solutions because of the left states $L$ we are considering.

\begin{figure}[ht]
	\centering
	\subfigure[Hugoniot for generic state $R$ at $\Omega_2^g$. $R=(0.21181, 0.743343)$.]
	{\includegraphics[scale=0.27]{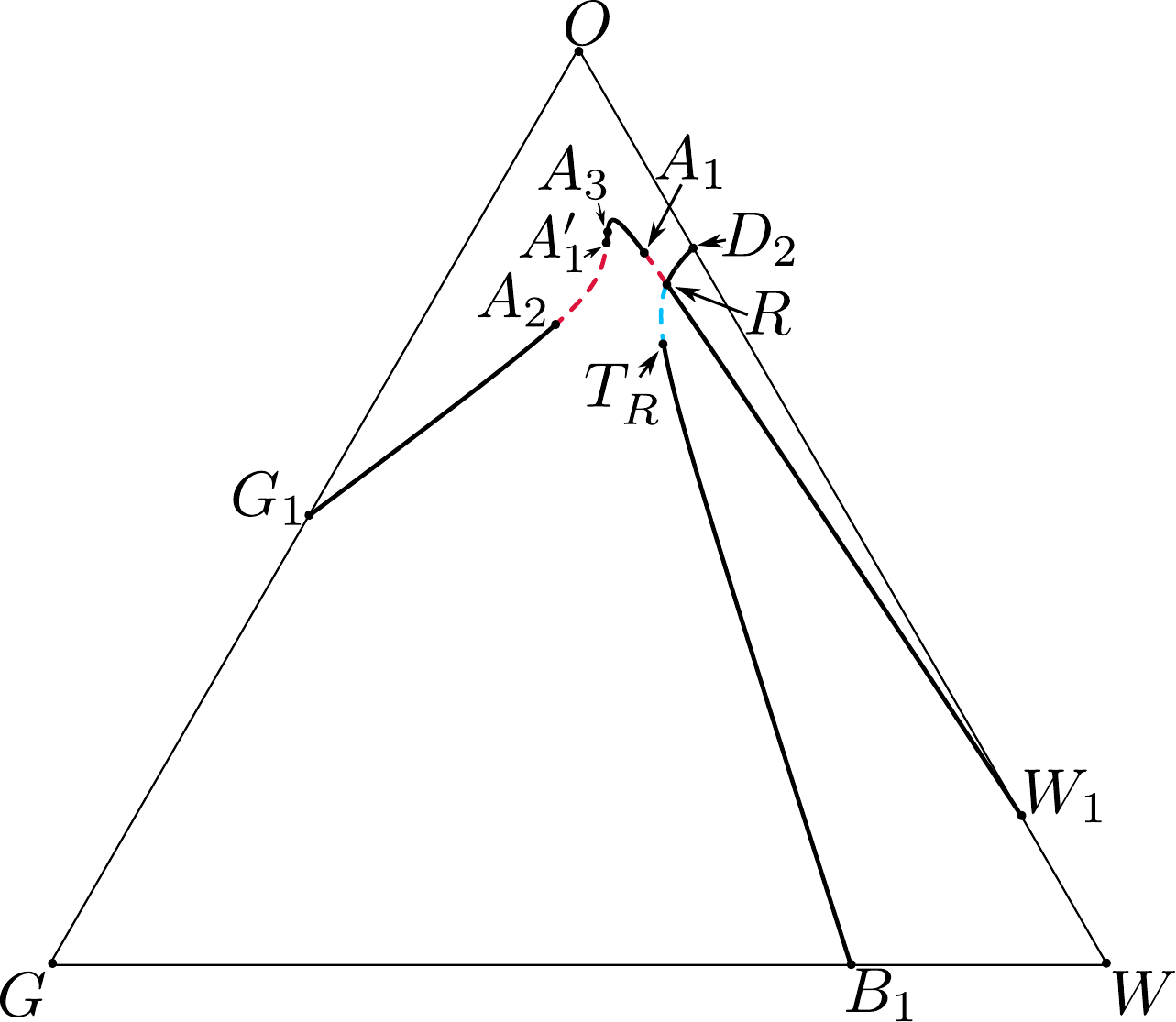}}  
	\hspace{1.75mm}
	\subfigure[Hugoniot for generic state in the region between the curves $EB$ and $M_E$, $\Omega_2^h$. $R=(0.246289 0.712158)$.]
	{\includegraphics[scale=0.27]{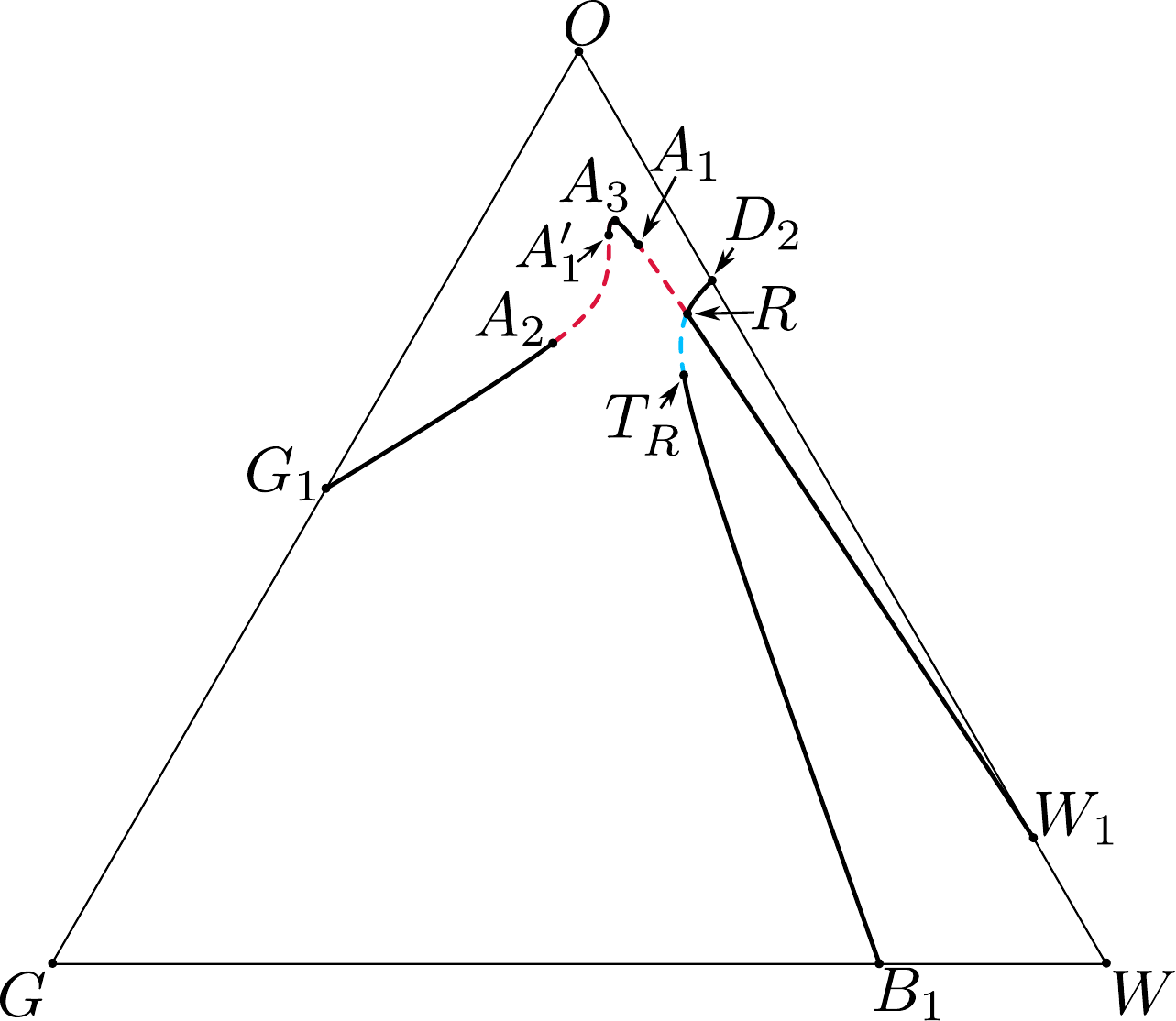}}  
 \caption{Hugoniot for generic states $R$ in $\Omega_2$.}
 	\label{fig:Hugoniot_Omega_2C}
\end{figure}

 \begin{figure}[ht]
	\centering
     \subfigure[Hugoniot for generic state $R$ at $\Omega_2^c$. $R=(0.295244,0.67847)$.]
	{\includegraphics[scale=0.27]{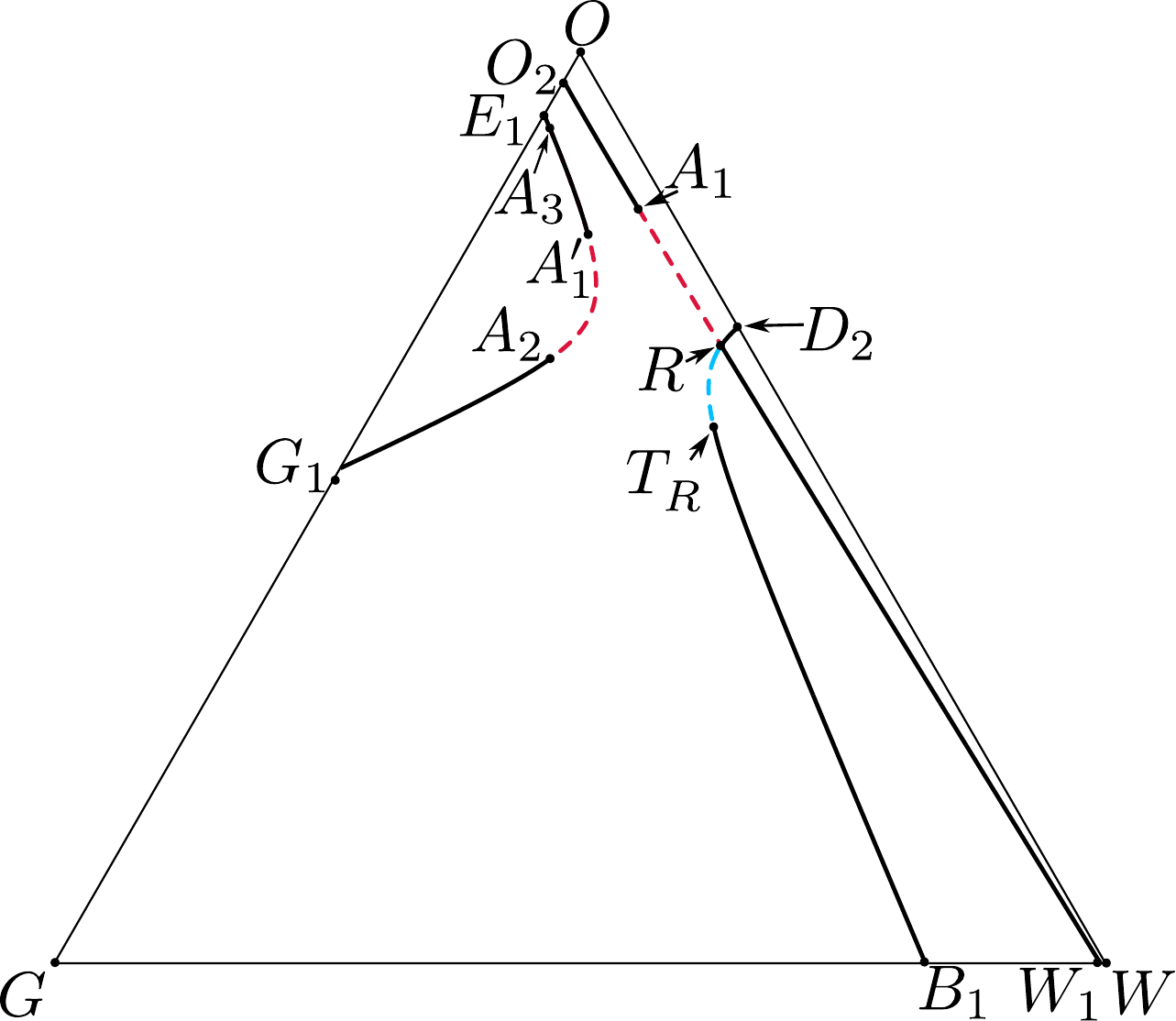}}  
	\hspace{1.75mm}
	\subfigure[Hugoniot for generic state $R$ at $\Omega_2^f$. $R=(0.275999,0.69048)$.]
	{\includegraphics[scale=0.27]{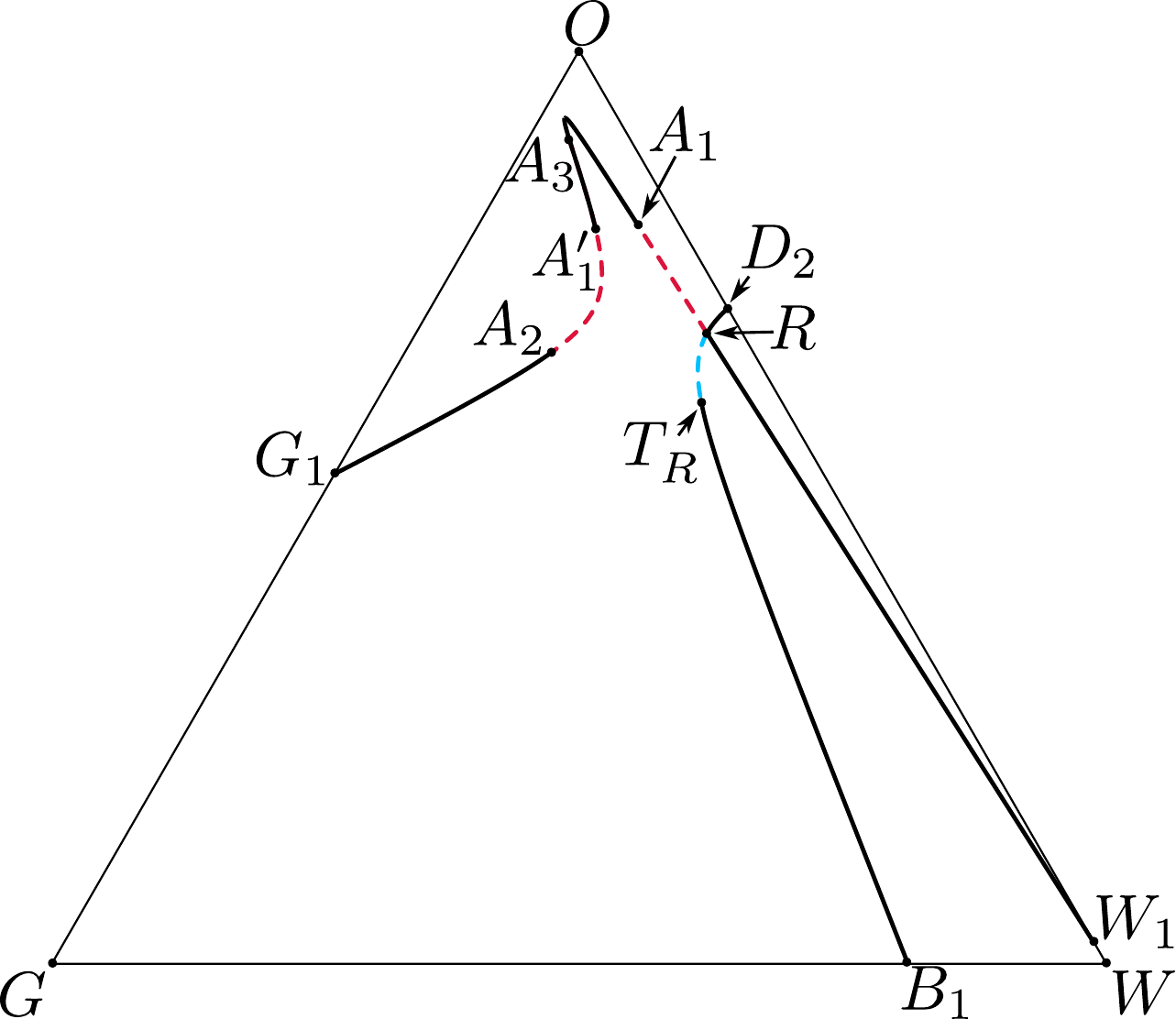}}  
 \subfigure[Hugoniot for generic state $R$ at $\Omega_2^i$. $R=(0.2683, 0.691978)$.]
	{\includegraphics[scale=0.27]{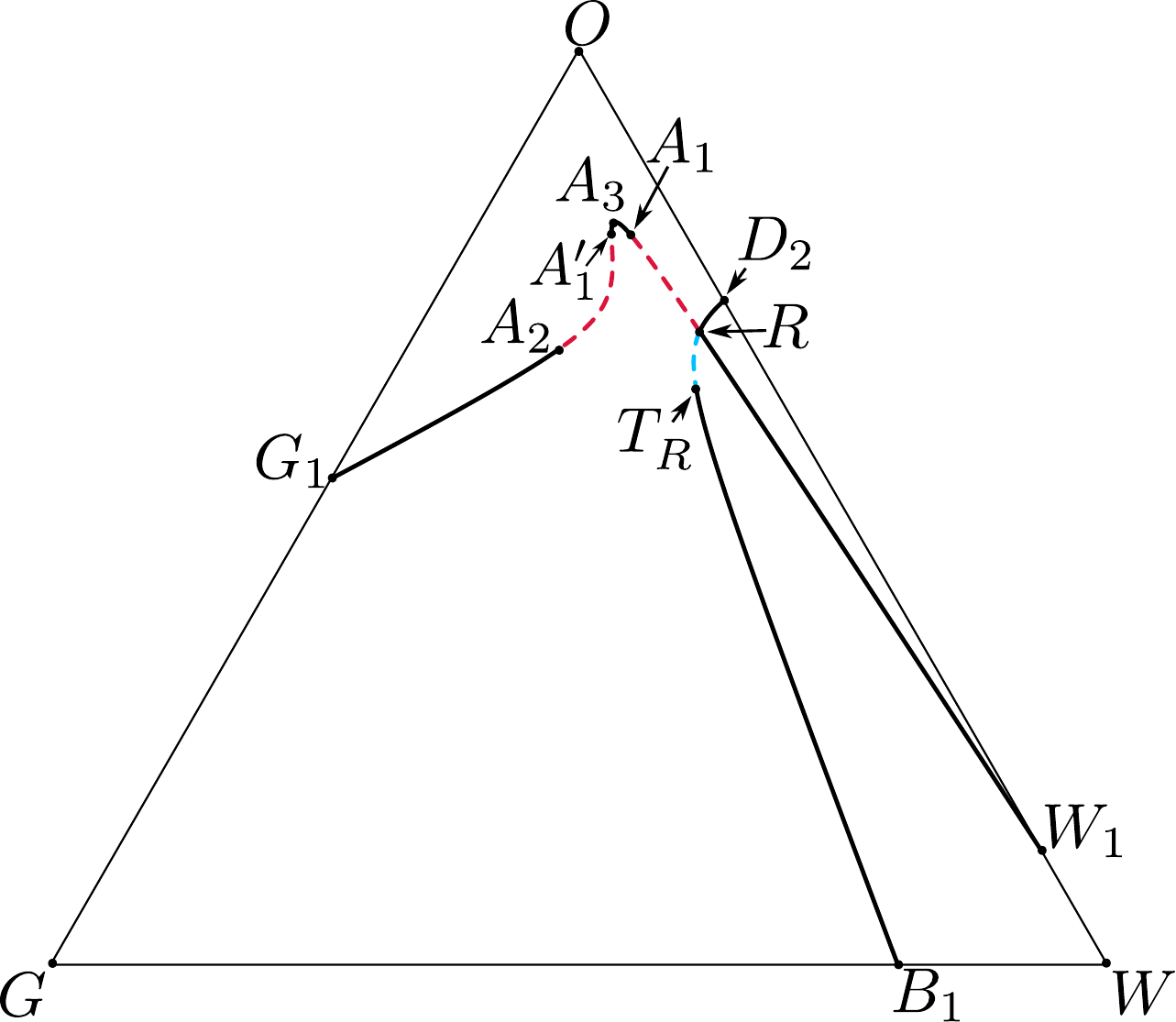}}  
 	\caption{Hugoniot for generic states $R$ in 
  $\Omega_2$. 
  }
 	\label{fig:Hugoniot_Omega_2D}
\end{figure}

\begin{figure}[ht]
	\centering
	\subfigure[Wave Curve for a generic state $R$ at $\Omega_2^a$. $R=(0.271633, 0.711087)$.]
	{\includegraphics[scale=0.27]{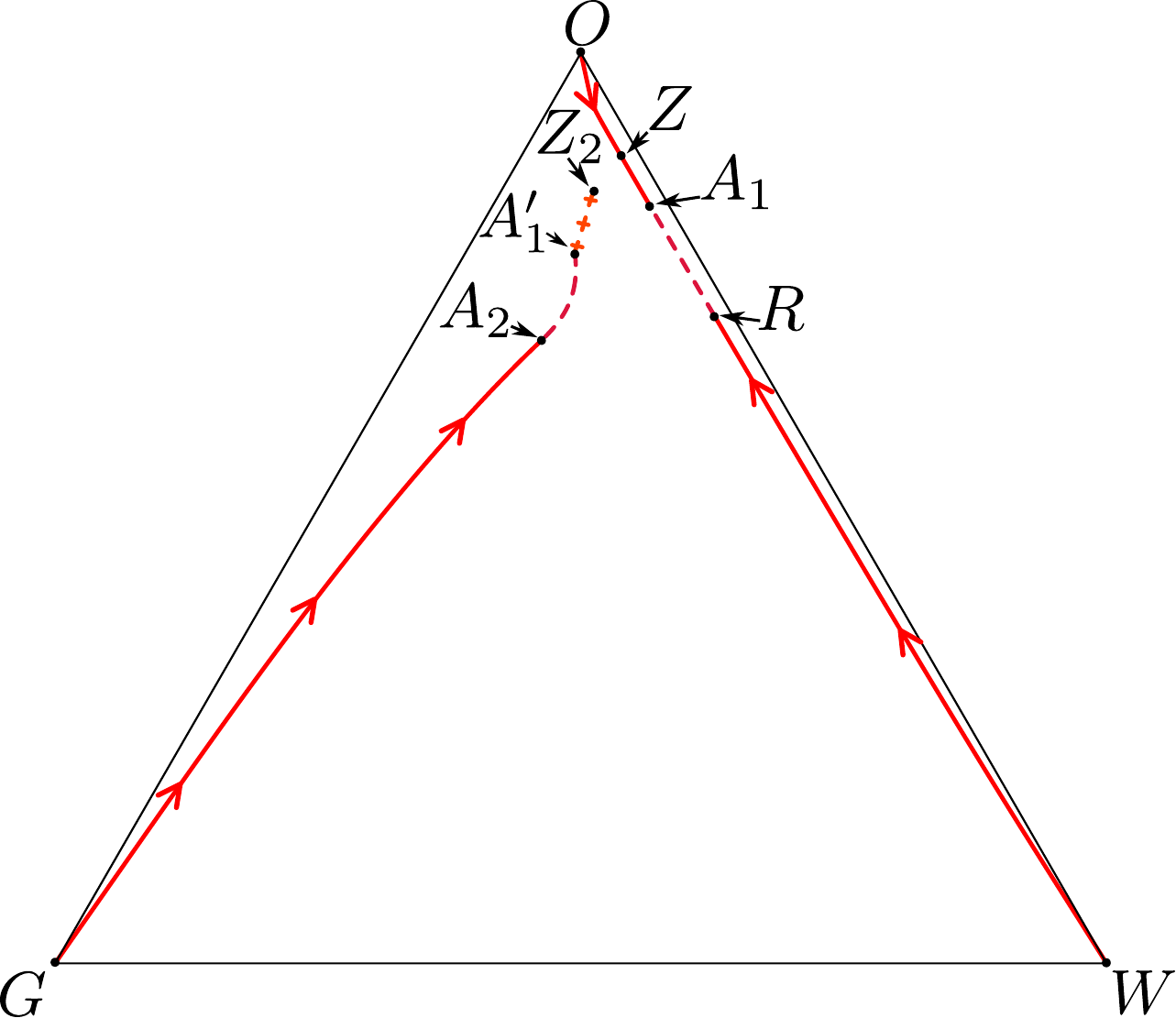}}  
	\hspace{1.75mm}
     \subfigure[Wave Curve for a generic state $R$ at $\Omega_2^d$. $R=(0.230643, 0.731226)$.]
	{\includegraphics[scale=0.27]{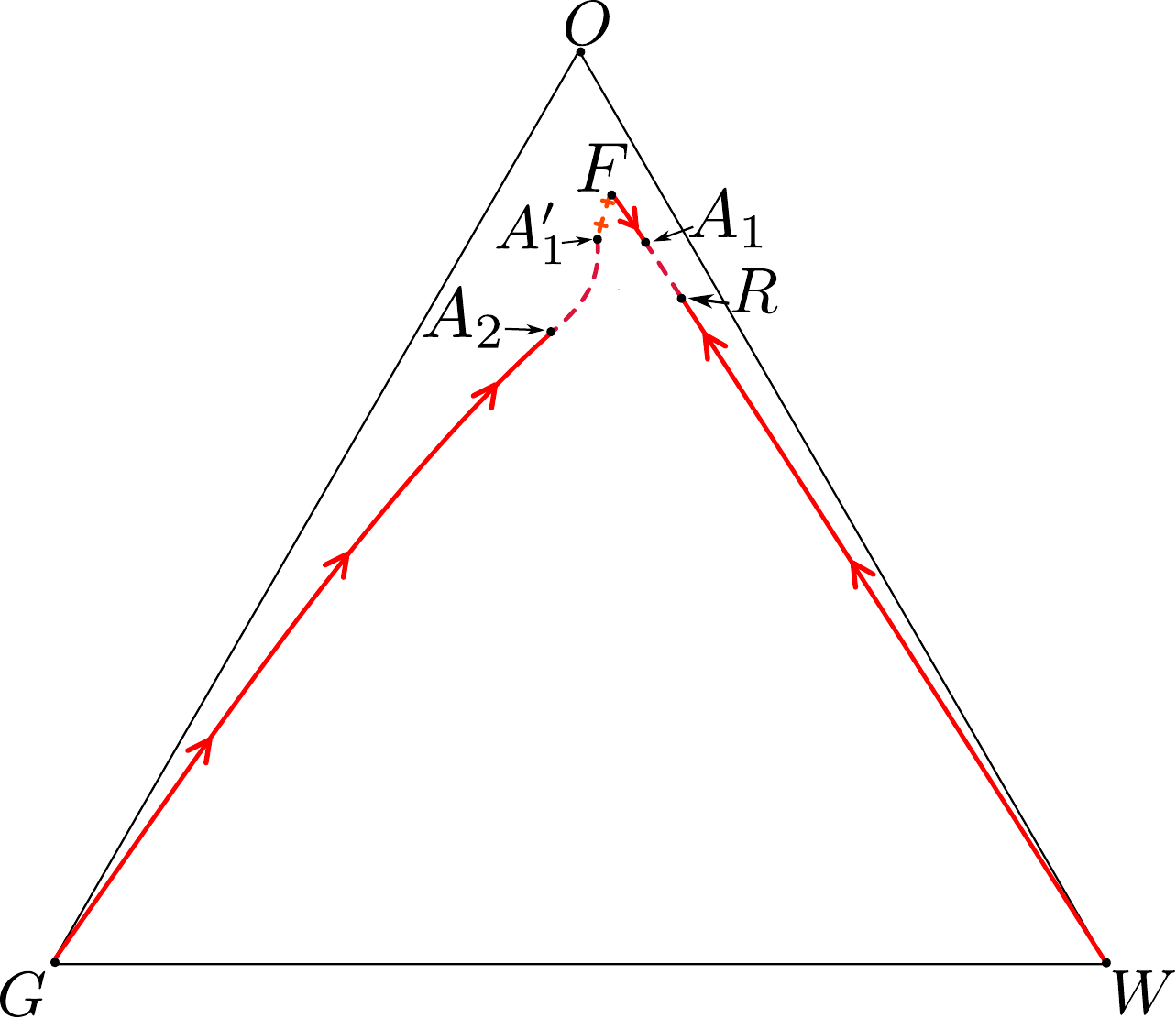}}  
	\hspace{1.75mm}
	 \caption{Wave Curve for a generic states $R$ in $\Omega_2$ }
 	\label{fig:Wave_Curve_2A}
\end{figure}
\begin{figure}[ht]
	\centering
	\subfigure[Riemann solution for a generic state $R$ at $\Omega_2^a$. $R=(0.271633, 0.711087)$.]
	{\includegraphics[scale=0.27]{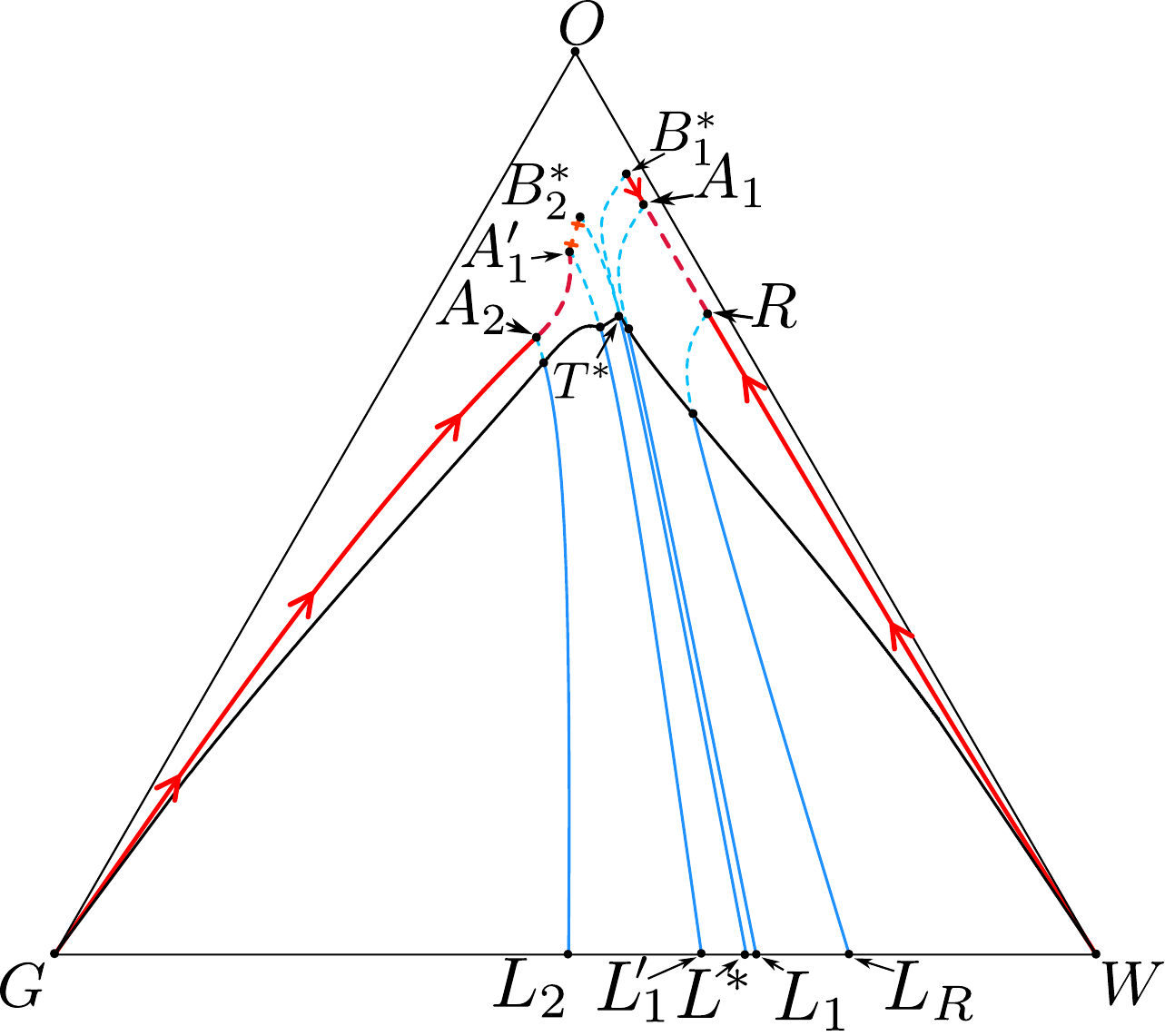}}  
	\hspace{1.75mm}
     \subfigure[Riemann solution for a generic state $R$ at $\Omega_2^d$. $R=(0.230643, 0.731226)$.]
	{\includegraphics[scale=0.27]{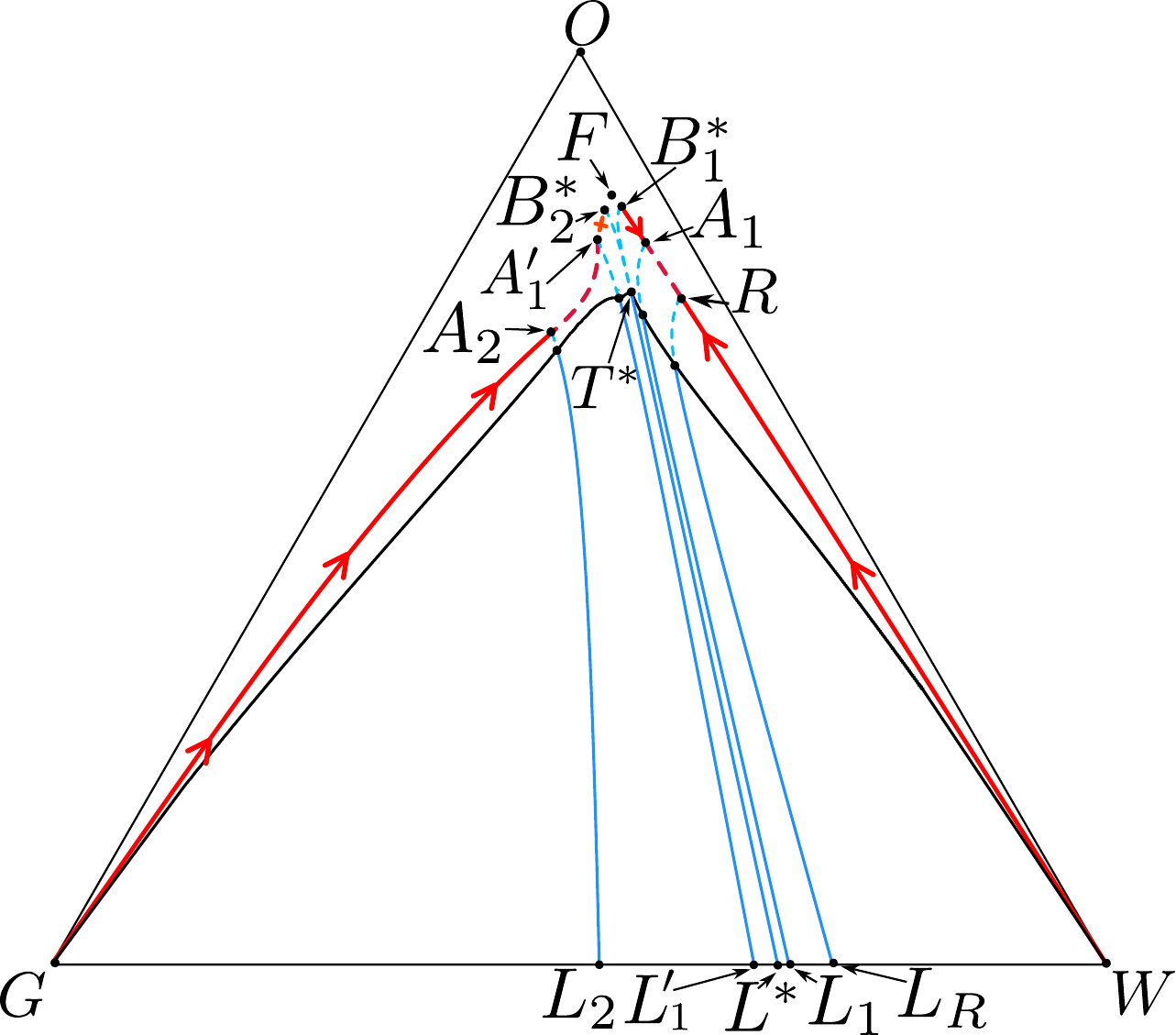}}  
	\hspace{1.75mm}
	 \caption{Riemann solution for a generic states $R$ in $\Omega_2$.}
 	\label{fig:Wave_Curve_2B}
\end{figure}

\begin{figure}[ht]
	\centering
	\subfigure[Wave Curve for a generic state $R$ at $\Omega_2^c$. $R=(0.295244,0.67847)$.]
	{\includegraphics[scale=0.27]{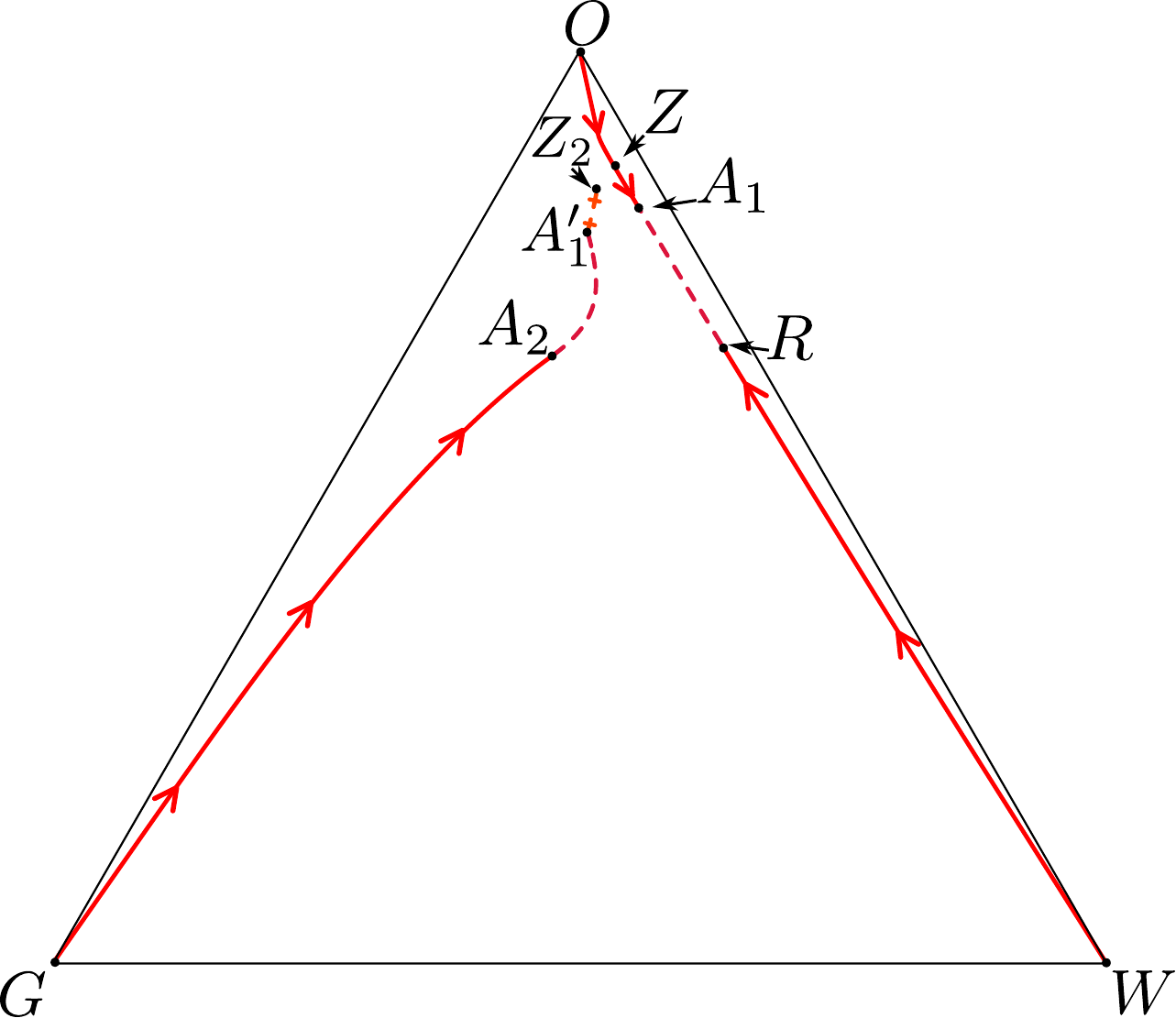}}  
	\hspace{1.75mm}
     \subfigure[Wave Curve for a generic state $R$ at $\Omega_2^f$. $R=(0.275999,0.69048)$.]
	{\includegraphics[scale=0.27]{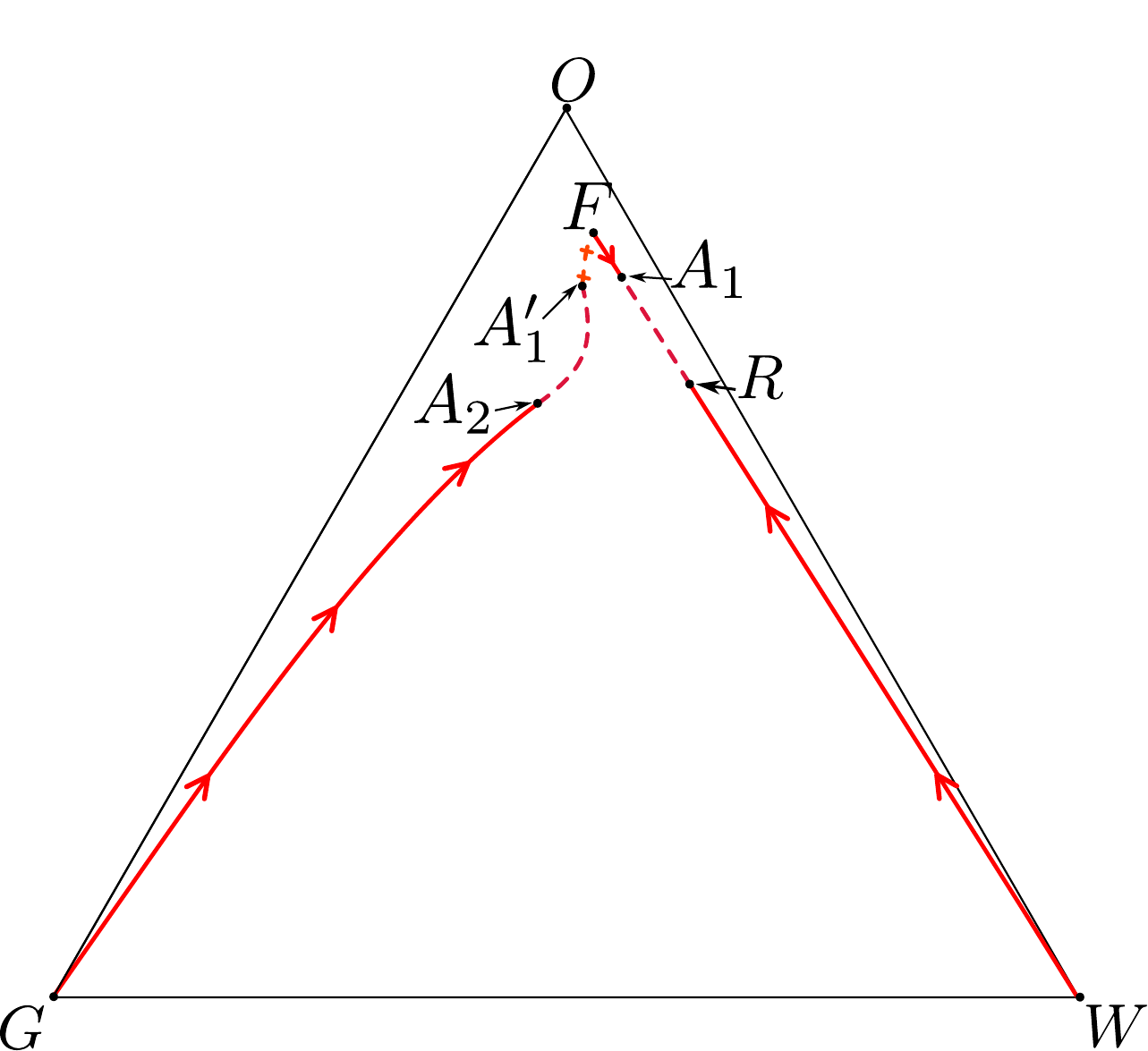}}  
	\hspace{1.75mm}
	 \caption{Wave Curve for a generic states $R$ in $\Omega_2$ }
 	\label{fig:Wave_Curve_2D}
\end{figure}
\begin{figure}[ht]
	\centering
	\subfigure[Riemann solution for a generic state $R$ at $\Omega_2^c$. $R=(0.295244,0.67847)$.]
	{\includegraphics[scale=0.27]{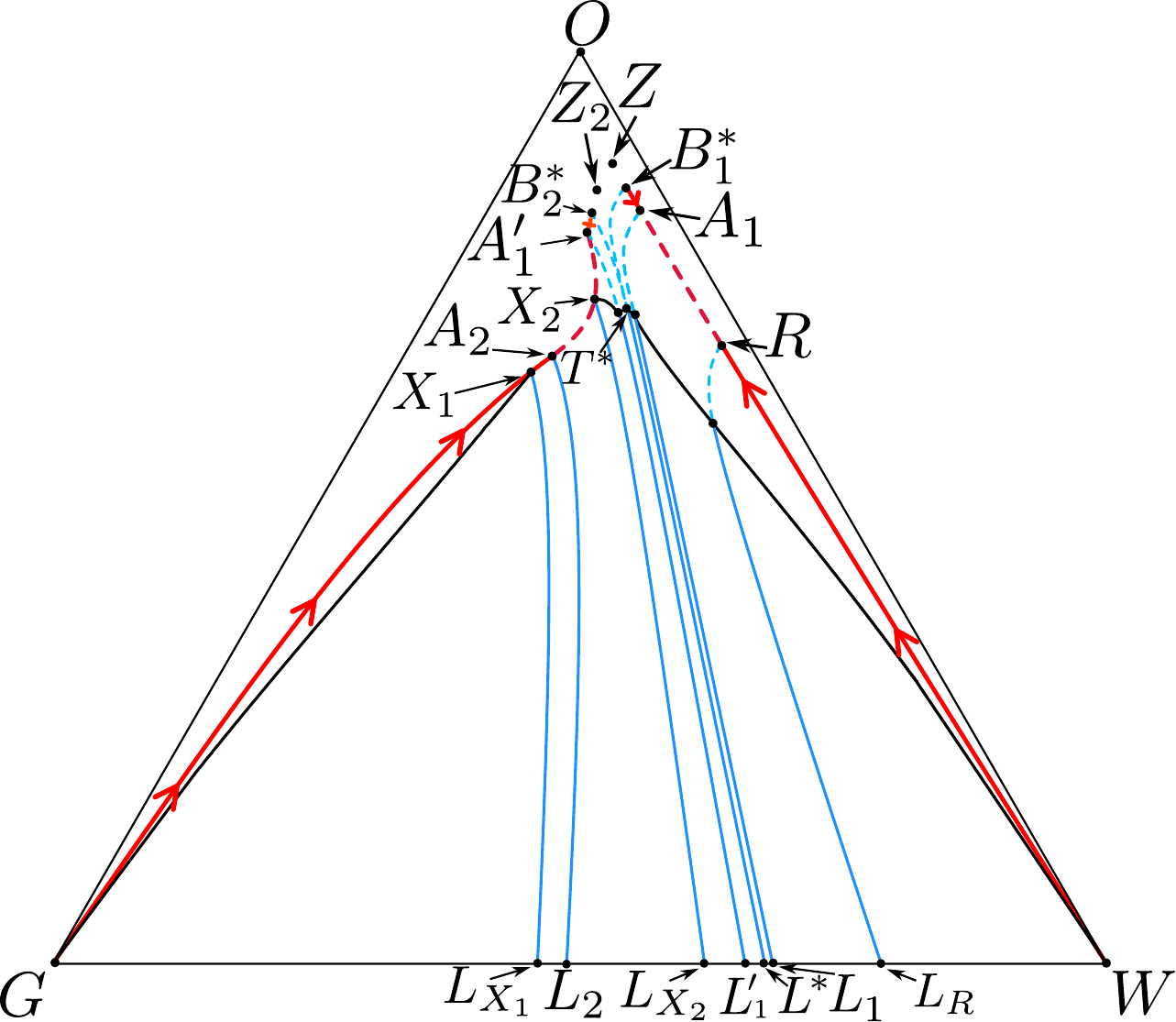}}  
	\hspace{1.75mm}
     \subfigure[Riemann solution for a generic state $R$ at $\Omega_2^f$. $R=(0.275999,0.69048)$.]
	{\includegraphics[scale=0.27]{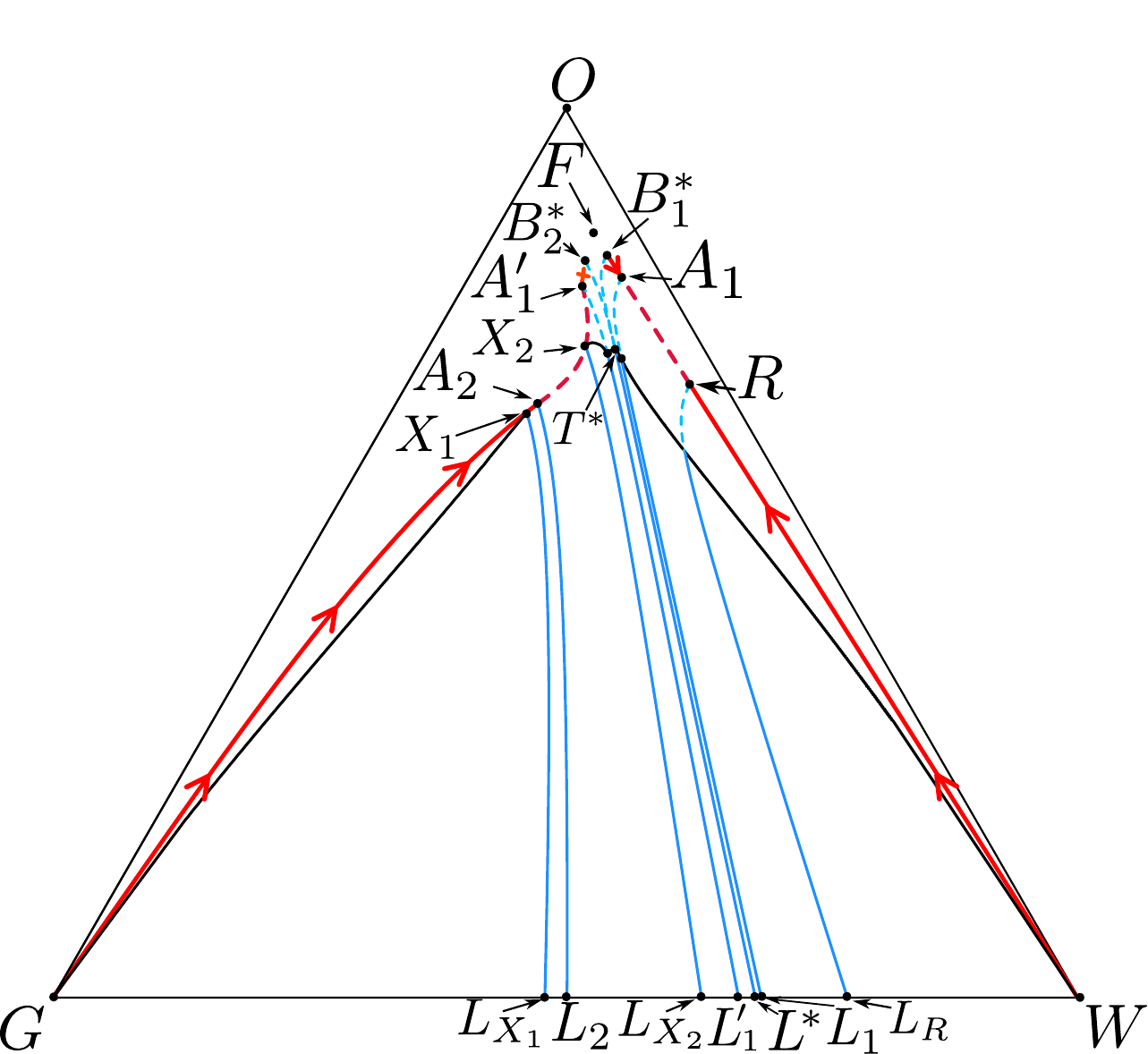}}  
	\hspace{1.75mm}
	 \caption{Riemann solution for a generic states $R$ in $\Omega_2$ }
 	\label{fig:Wave_Curve_2E}
\end{figure}

\begin{figure}[ht]
	\centering
	\subfigure[Wave Curve for a generic state $R$ in $\Omega_2^g$. $R=(0.21181, 0.743343)$.]
	{\includegraphics[scale=0.27]{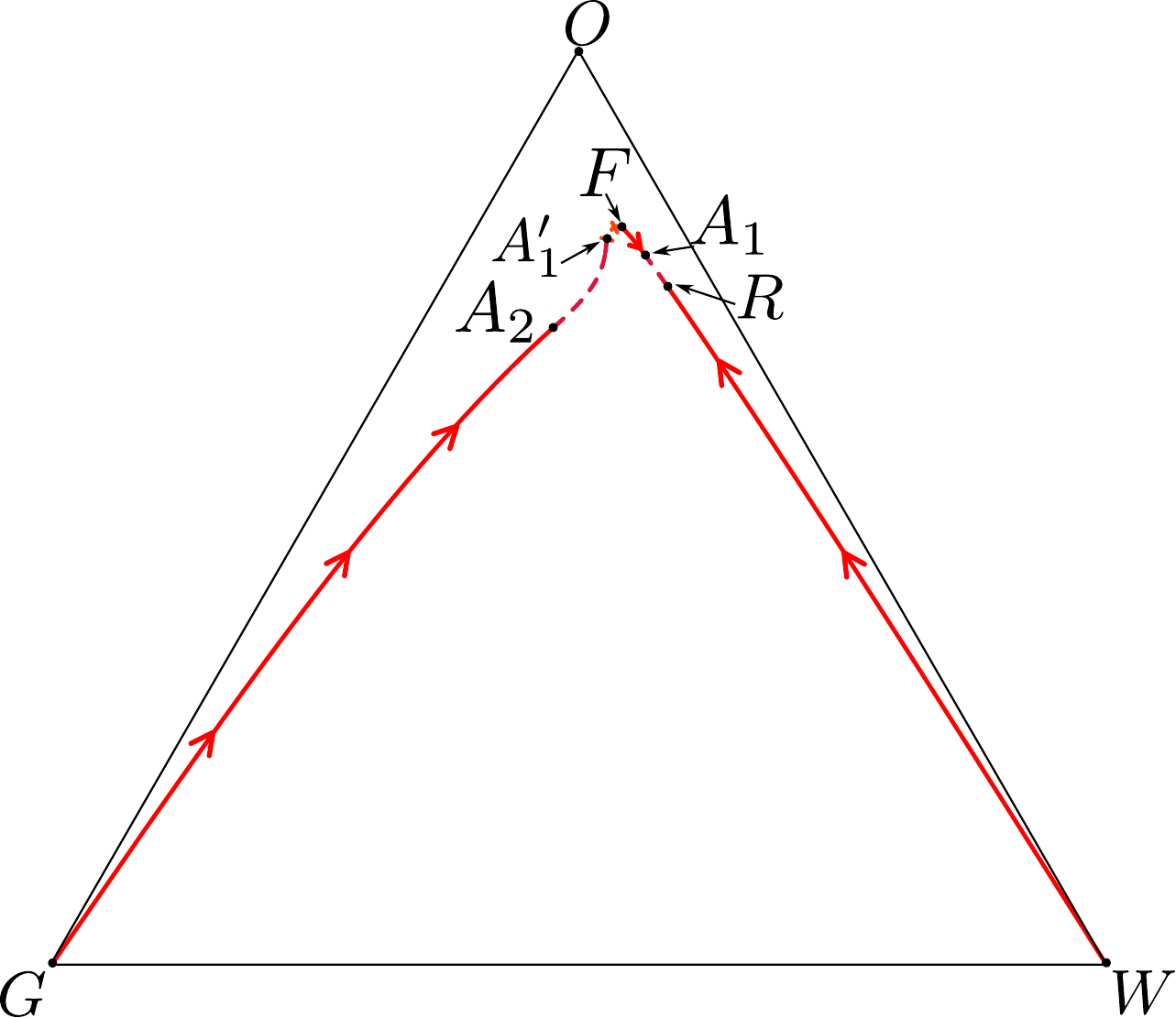}}  
	\hspace{1.75mm}
     \subfigure[Wave Curve for a generic state $R$ in $\Omega_2^i$. $R=(0.2683, 0.691978)$.]
	{\includegraphics[scale=0.27]{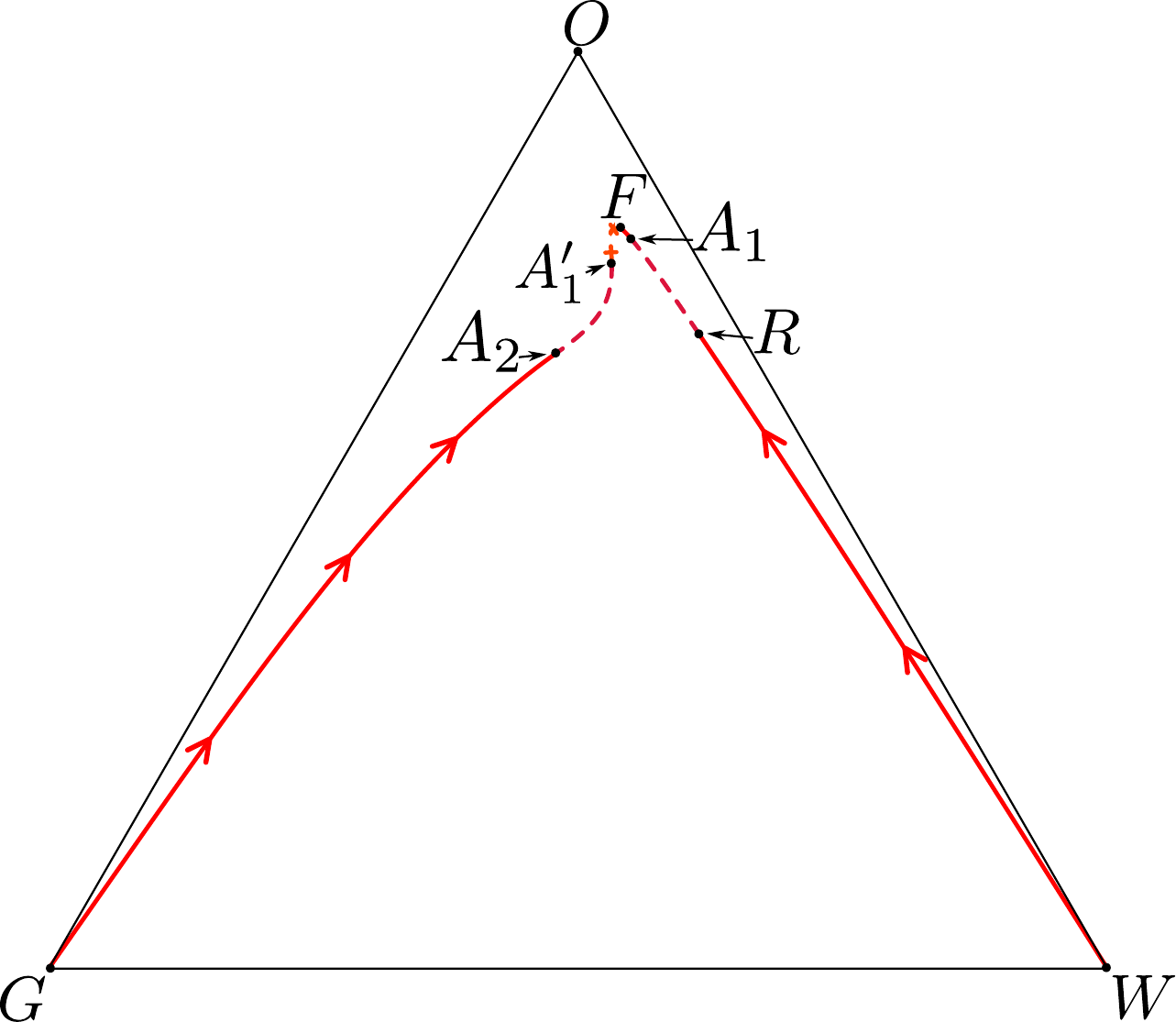}}  
	\hspace{1.75mm}
	 \caption{Wave Curve for a generic states $R$ in $\Omega_2$ }
 	\label{fig:Wave_Curve_2G}
\end{figure}
\begin{figure}[ht]
	\centering
	\subfigure[Riemann solution for a generic state $R$ in $\Omega_2^g$. $R=(0.21181, 0.743343)$.]
	{\includegraphics[scale=0.27]{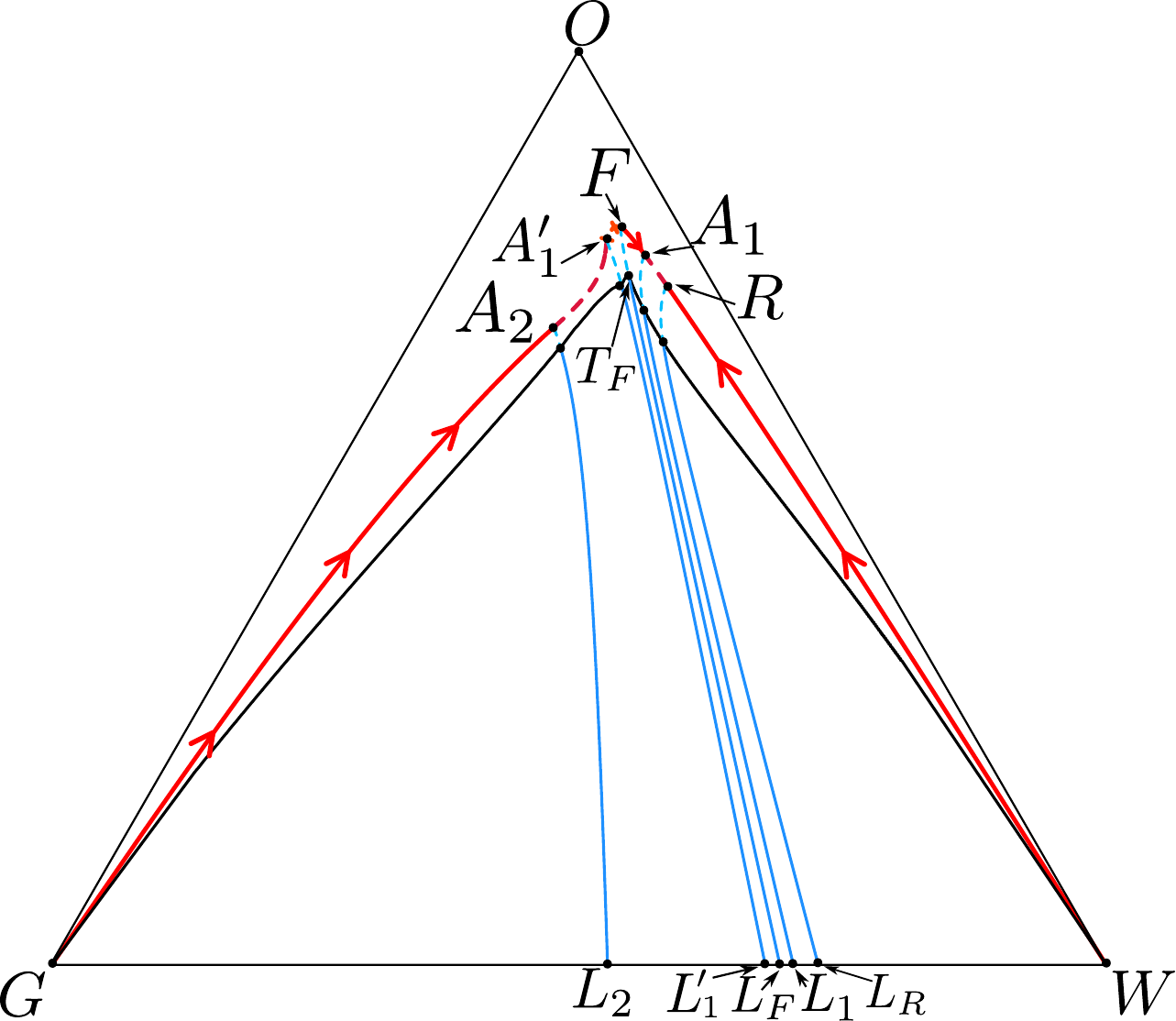}}  
	\hspace{1.75mm}
     \subfigure[Riemann solution for a generic state $R$ in $\Omega_2^i$. $R=(0.2683, 0.691978)$.]
	{\includegraphics[scale=0.27]{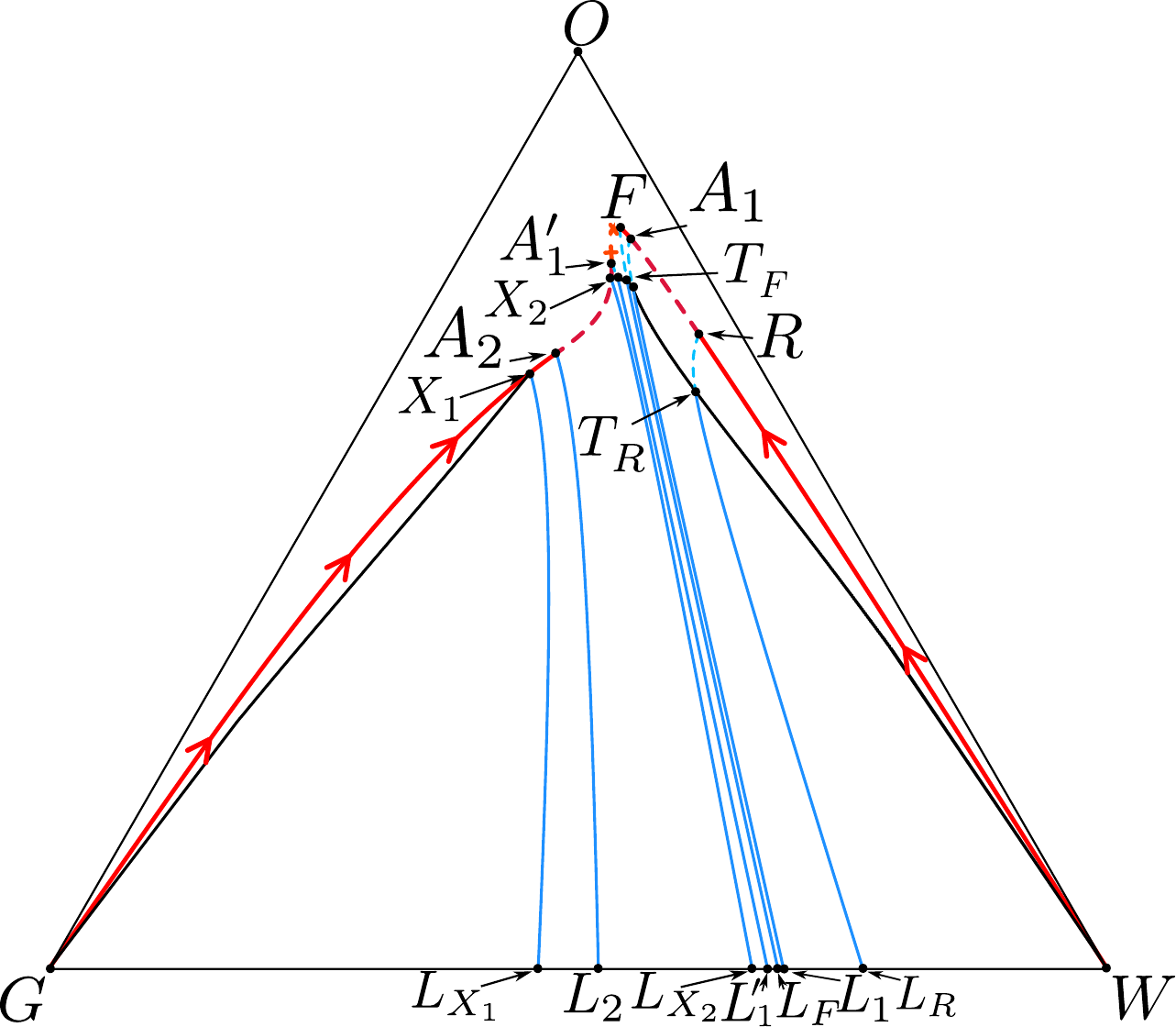}}  
	\hspace{1.75mm}
	 \caption{Riemann solution for a generic state $R$ in $\Omega_2$.}
 	\label{fig:Wave_Curve_2H}
\end{figure}

Next we discuss Riemann solutions.
For this purpose region $\Omega_2$ is subdivided into subregions $\Omega_2^\alpha$, $\alpha \in \{a, b, \dots, i\}$.
This subdivision involves some of the loci introduced in
Section~\ref{sec:Bifurcation}.

The boundaries of the subregions $\Omega_2^\alpha$, $\alpha \in \{a, b, \dots, i\}$ in 
Fig.~\ref{fig:RRegionsOmegas}(b) are defined by
the invariant segment $\mathcal{I}_4^f$-$H_3$;
the $f$-rarefaction segment $\mathcal{I}_3^f$-$H_2$;
the segment $T_I^1$-$T_I^2$ of curve $T_I$; and the segment $M_E^a$-$M_E^b$ of curve $M_E$.
The boundary $T_I^1$-$T_I^2$, which is part of curve $T_I$, separates states $R$ of two types: those for which the fast wave curve does not cross the $s$-inflection locus and those for which it does.
The boundary $M_E^a$-$M_E^b$, which is part of curve $M_E$, determines whether the fast wave curve intersects the slow inflection locus in a shock or a rarefaction segment. 

\medskip
\noindent{ \bf Riemann solution for $R$ in subregion $\Omega_2^a\cup\Omega_2^d$.}
\medskip


\begin{cla}\label{cla:RSolution-Omega_2aU2c}
Refer to Fig.~\ref{fig:Wave_Curve_2B}.
Let $L$ be a state on the edge \GW of the
saturation triangle and $R$ be a state in subregion $\Omega_2^a\cup\Omega_2^d$ in
Figs.~\ref{fig:RegioesThetas}(b), \ref{fig:RRegionsOmegas}(b).
Let $A_1$, $A_2$, $A_1'$, $B_1^*$ and $B_2^*$ be the states in $\wm_f(R)$ with
$\sigma(A_1;R) = \lambdaf(A_1)$, $\lambda_f(A_1)=\sigma(A_1;R)=  \sigma(A_1';A_1)=\sigma(A_1';R) $, $ \lambdas(T^*) = \sigma(T^*;B_1^*) = \lambdaf(B_1^*) 
= \sigma(B_2^*; B_1^*) = \sigma(T^*; B_2^*)$, with $T^*\in\mathcal{H}(B_1^*)$.

Let  $L_2$, $L_1'$, $L^*$, $L_1$ and $L_R$
be the intersection points of
the backward slow wave curves through
$A_2$, $A_1'$, $B_1^*$ ($B_2^*$), $A_1$  and $R$ with the edge \GWc, respectively.
Then, 
 
\begin{itemize}
\item[(i)] if $L=G$, the Riemann solution is 
$
   G\testright{R_f} A_2 \xrightarrow{'S_f\,} R ;
$
\item[(ii)] if $L \in(G, L_2)$, the Riemann solution is 
$
    L\testright{R_s} T_1 \xrightarrow{'S_s} M_1  \testright{R_f}A_2 \xrightarrow{'S_f} R ,
$
where $T_1\in(G, T^*)_{\text{ext}}$ and $M_1\in(G, A_2)$; 
\item[(iii)] if $L \in[L_2, L_1']$ or $L \in[L_1, L_R)$, the Riemann solutions are 
$ L\testright{R_s} T_1 \xrightarrow{'S_s} M_1  \testright{S_f} R$,
   or
   $L\testright{R_s} T_2 \xrightarrow{'S_s} M_2  \testright{S_f} R 
$,
where $T_1\in(G, T^*)_{\text{ext}}$, $T_2\in(T^*,W)_{\text{ext}}$, $M_1\in[A_2,A_1']$ and  $M_2\in[A_1,R)$;
\item[(iv)] if $L \in(L_1', L^*)$, the Riemann solution is 
$
    L\testright{R_s} T_1 \xrightarrow{'S_s} M_1'  \testright{S_f\,'}M_1 \xrightarrow{R_f} A_1\xrightarrow{'S_f} R ,
$
where $T_1\in(G, T^*)_{\text{ext}}$,  $M_1'\in(A_1',B_2^*)$ and $M_1\in(B_1^*, A_1)$;
\item[(v)] if $L \in(L^*, L_1)$, the Riemann solution is
$
    L\testright{R_s} T_1 \xrightarrow{'S_s} M_1  \testright{R_f}A_1 \xrightarrow{'S_f} R,
$
where $T_1\in(T^*,W)_{\text{ext}}$ and $M_1\in(B_1^*,A_1)$;
\item[(vi)] if $L \in (L_R,W)$, the Riemann solution is
$
     L\testright{R_s} T_1 \xrightarrow{'S_s}  M_1  \testright{R_f} R,
$
where $T_1\in(T^*,W)_{\text{ext}}$ and $M_1\in(R,W)$;
\item[(vii)] if $L=W$, the Riemann solution is
$
   W\testright{R_f} R
$.
\end{itemize}
\end{cla}

\begin{rem}\label{rem:L*B*4}
If $L = L^*$, see Fig.~\ref{fig:Wave_Curve_2B}, there are two paths to reach the right state
$R$, namely 
$
     L^*\testright{R_s} T^* \xrightarrow{'S_s} B_2^* \xrightarrow{S_f'} B_1^* \xrightarrow{R_f}  A_1 \xrightarrow{'S_f} R
$
or
$
     L^*\testright{R_s} T^* \xrightarrow{'S_o'} B_1^* \xrightarrow{R_f}  A_1 \xrightarrow{'S_f} R
$ 
but the triple shock rule guarantees they represent the same solution in $xt$-space.
In other words, the solution consists of a single wave group.
\end{rem}

Next, we consider $R$ in region $\Omega_2^g$.
This case differs from the previous one in that
the states $B_1^*$ and $B_2^*$ collapse to
state $F$ and the corresponding state $T^*$ becomes $T_F$. We then have
\begin{cla}\label{cla:RSolution-Omega_2e}
Refer to Fig.~\ref{fig:Wave_Curve_2H}.
Let $L$ be a state on the edge \GW of the
saturation triangle and $R$ be a state in subregion $\Omega_2^g$ in
Figs.~\ref{fig:RegioesThetas}(b), \ref{fig:RRegionsOmegas}(b).
Then the Riemann solution is the one described in Claim~\ref{cla:RSolution-Omega_2aU2c},
replacing $T^*$ with $T_F$ and $L^*$ with $L_F$.
\end{cla}

\medskip
\noindent{ \bf Riemann solution for $R$ in subregion $\Omega_2^b\cup\Omega_2^e$.}
\medskip


Since the state $R$ has now crossed the curve $T_I$,  a new feature arises, namely, the interaction between the wave curve $W_f(R)$ and the $s$-inflection locus. More precisely, 
as shown in Figs.~\ref{fig:Wave_Curve_2C}(b) and \ref{fig:Wave_Curve_2F}(b),  $W_f(R)$ intersects the $s$-inflection locus at two points, denoted as $X_1$ and $X_2$, with $X_1$ and $X_2$
in the admissible $f$-shock segment $[A_2, A_1')$.
This fact leads to changes in the structure of the certain slow wave curves used in the Riemann solutions. 

\begin{rem}\label{rem:X1X2}
We notice that the states $X_1$ and $X_2$ are also intersection
points of the slow extension of $W_f(R)$ with the
$s$-inflection locus.
\end{rem}

\medskip
\begin{cla}\label{cla:RSolution-Omega_2bU2e}
Refer to Fig.~\ref{fig:Wave_Curve_2C}.
Let $L$ be a state on the edge \GW of the
saturation triangle and $R$ be a state in subregion $\Omega_2^b\cup\Omega_2^e$ in
Figs.~\ref{fig:RegioesThetas}(b), \ref{fig:RRegionsOmegas}(b).
Let $A_1$, $A_2$, $A_1'$, $B_1^*$ and $B_2^*$  be the states in $\wm_f(R)$ such that
$\sigma(A_1;R) = \lambdaf(A_1)$, $\sigma(A_1;R)=  \sigma(A_1';A_1)=\sigma(A_1';R) $, $ \lambdas(T^*) = \sigma(T^*;B_1^*) = \lambdaf(B_1^*) 
= \sigma(B_2^*; B_1^*) = \sigma(T^*; B_2^*)$, with $T^*\in\mathcal{H}(B_1^*)$.

Let  $L_2$, $L_{X_1}$, $L_{X_2}$, $L_1'$,  $L^*$, 
$L_1$ and $L_R$
be the intersection points of
the backward slow wave curves through
$X_1$, $A_2$, $X_2$, $A_1'$, $B_1^*$ ($B_2^*$),
$A_1$ and $R$ with the edge \GWc, respectively.
Then, 
\begin{itemize}
\item[(i)] if $L=G$, the Riemann solution is 
$
   G\testright{R_f} A_2 \xrightarrow{'S_f\,} R ;
$
\item[(ii)] if $L \in(G, L_2)$, the Riemann solution is 
$
    L\testright{R_s} T_1 \xrightarrow{'S_s} M_1  \testright{R_f}A_2 \xrightarrow{'S_f} R ,
$
where $T_1\in(G, X_1)_{\text{ext}}$ and $M_1\in(G, A_2)$;

\item[(iii)] if $L = L_2$, the Riemann solution is 
$
     L\testright{R_s} T_1 \xrightarrow{'S_s} A_2 \testright{'S_f} R ,
$
where $T_1\in(G, X_1)_{\text{ext}}$;

\item[(iv)] if $L \in (L_2, L_{X_1})$, the Riemann solution is 
$
     L\testright{R_s} T_1 \xrightarrow{'S_s} M_1 \testright{S_f} R ,
$
where $T_1\in(G, X_1)_{\text{ext}}$ and $M_1\in (A_2, X_1)$;

\item[(v)] if $L \in[L_{X_1},L_{X_2}]$, the Riemann solution is 
$
    L\testright{R_s}  M_1  \testright{S_f} R ,
$
where $M_1\in[X_1,X_2]$; 
\item[(vi)] if $L \in (L_{X_2}, L_1']$ or $L \in[L_1, L_R)$, the Riemann solutions are 

$
L\testright{R_s} T_1 \xrightarrow{'S_s} M_1  \testright{S_f} R ,\quad\mbox{ or }\quad 
   L\testright{R_s} T_2 \xrightarrow{'S_s} M_2  \testright{S_f} R 
$
where $T_1\in(X_2, T^*)_{\text{ext}}$, $M_1\in[X_2,A_1']$, 
\quad or \quad $T_2\in(T^*,W)_{\text{ext}}$  and  $M_2\in[A_1,R)$;
\item[(vii)] if $L \in(L_1', L^*)$, the Riemann solution is 
$
    L\testright{R_s} T_1 \xrightarrow{'S_s} M_1'  \testright{S_f\,'}M_1 \xrightarrow{R_f} A_1\xrightarrow{'S_f} R ,
$
where $T_1\in(X_2, T^*]_{\text{ext}}$,  $M_1'\in(A_1',B_2^*)$ and $M_1\in (B_1^*, A_1)$; 
\item[(viii)] if $L \in(L^*, L_1)$, the Riemann solution is
$
    L\testright{R_s} T_1 \xrightarrow{'S_s} M_1  \testright{R_f}A_1 \xrightarrow{'S_f} R,
$
where $T_1\in(T^*,W)_{\text{ext}}$ and $M_1\in(B_1^*,A_1)$;
\item[(ix)] if $L=W$, the Riemann solution is
$
   W\testright{R_f} R.
$
\end{itemize}
\end{cla}

\begin{rem}\label{rem:L*B*5}
If $L = L^*$, see Fig.~\ref{fig:Wave_Curve_2C}, there are two paths to reach the right state
$R$, namely 
$
     L^*\testright{R_s} T^* \xrightarrow{'S_s} B_2^* \xrightarrow{S_f'} B_1^* \xrightarrow{R_f}  A_1 \xrightarrow{'S_f} R
$
or
$
     L^*\testright{R_s} T^* \xrightarrow{'S_o'} B_1^* \xrightarrow{R_f}  A_1 \xrightarrow{'S_f} R
$ 
but the triple shock rule guarantees they represent the same solution in $xt$-space.
In other words, the solution consists of a single wave group.
\end{rem}

We now consider $R$ in region $\Omega_2^h$.
The changes in the
Riemann solutions when $R$ moves from
$\Omega_2^b\cup \Omega_2^e$ to $\Omega_2^h$ are similar to those 
described in Claims~\ref{cla:RSolution-Omega_2aU2c}, \ref{cla:RSolution-Omega_2e} when $R$ moves from
$\Omega_2^a\cup \Omega_2^d$  to $\Omega_2^g$.
Namely,
\begin{cla}\label{cla:RSolution-Omega_2h}
Refer to Fig.~\ref{fig:Wave_Curve_2F}.
Let $L$ be a state on the edge \GW of the
saturation triangle and $R$ be a state in subregion $\Omega_2^h$ in
Figs.~\ref{fig:RegioesThetas}(b), \ref{fig:RRegionsOmegas}(b).
Then the Riemann solution is the one described in Claim~\ref{cla:RSolution-Omega_2bU2e},
replacing $T^*$ with $T_F$ and $L^*$ with $L_F$.
\end{cla}

\medskip

\noindent{ \bf Riemann solution for $R$ in subregion $\Omega_2^c\cup\Omega_2^f$.}
\medskip

The Riemann solutions that remain to be considered arise when $R$ crosses the boundary $M_E$
in Fig.~\ref{fig:RRegionsOmegas}, where $A_2$ and $X_1$ coalesce,
differ from the previous ones by the relative position of states $A_2$ and $X_1$. 

We start with $R$ in region $\Omega_2^c\cup\Omega_2^f$. 

\medskip
\begin{cla}\label{cla:RSolution-Omega_2bU2d}
Refer to Fig.~\ref{fig:Wave_Curve_2E}.
Let $L$ be a state on the edge \GW of the
saturation triangle and $R$ be a state in subregion $\Omega_2^c\cup\Omega_2^f$ in
Figs.~\ref{fig:RegioesThetas}(b), \ref{fig:RRegionsOmegas}(b).
Let $A_1$, $A_2$, $A_1'$, $B_1^*$ and $B_2^*$  be the states in $\wm_f(R)$ such that
$\sigma(A_1;R) = \lambdaf(A_1)$, $\sigma(A_1;R)=  \sigma(A_1';A_1)=\sigma(A_1';R) $, $ \lambdas(T^*) = \sigma(T^*;B_1^*) = \lambdaf(B_1^*) 
= \sigma(B_2^*; B_1^*) = \sigma(T^*; B_2^*)$, with $T^*\in\mathcal{H}(B_1^*)$.

Let $L_{X_1}$, $L_2$, $L_{X_2}$, $L_1'$,  $L^*$, 
$L_1$ and $L_R$
be the intersection points of
the backward slow wave curves through
$X_1$, $A_2$, $X_2$, $A_1'$, $B_1^*$ ($B_2^*$),
$A_1$ and $R$ with the edge \GWc, respectively.
Then, 
\begin{itemize}
\item[(i)] if $L=G$, the Riemann solution is 
$
   G\testright{R_f} A_2 \xrightarrow{'S_f\,} R ;
$
\item[(ii)] if $L \in(G, L_{X_1})$, the Riemann solution is 
$
    L\testright{R_s} T_1 \xrightarrow{'S_s} M_1  \testright{R_f}A_2 \xrightarrow{'S_f} R ,
$
where $T_1\in(G, X_1)_{\text{ext}}$ and $M_1\in(G, X_1)$; 
\item[(iii)] if $L \in[L_{X_1},L_2)$, the Riemann solution is 
$
    L\testright{R_s}  M_1  \testright{R_f}A_2 \xrightarrow{'S_f} R ,
$
where $M_1\in[X_1, A_2)$;
\item[(iv)] if $L \in[L_2,L_{X_2}]$, the Riemann solution is 
$
    L\testright{R_s}  M_1  \testright{S_f} R ,
$
where $M_1\in[A_2,X_2]$; 
\item[(v)] if $L \in (L_{X_2}, L_1']$ or $L \in[L_1, L_R)$, the Riemann solutions are 
$
L\testright{R_s} T_1 \xrightarrow{'S_s} M_1 \testright{S_f} R$, or $
   L\testright{R_s} T_2 \xrightarrow{'S_s} M_2  \testright{S_f} R 
$
where $T_1\in(X_2, T^*)_{\text{ext}}$, $M_1\in (X_2,A_1']$, $T_2\in(T^*,W)_{\text{ext}}$,  and  $M_2\in[A_1,R)$;
\item[(vi)] if $L \in(L_1', L^*)$, the Riemann solution is 
$
    L\testright{R_s} T_1 \xrightarrow{'S_s} M_1'  \testright{S_f\,'}M_1 \xrightarrow{R_f} A_1\xrightarrow{'S_f} R ,
$
where $T_1\in(X_2, T^*)_{\text{ext}}$,  $M_1'\in(A_1',B_2^*)$ and $M_1\in(B_1^*, A_1)$;
\item[(vii)] if $L \in(L^*, L_1)$, the Riemann solution is
$
    L\testright{R_s} T_1 \xrightarrow{'S_s} M_1  \testright{R_f}A_1 \xrightarrow{'S_f} R,
$
where $T_1\in(T^*,W)_{\text{ext}}$ and $M_1\in(B_1^*,A_1)$;
\item[(viii)] if $L=W$, the Riemann solution is
$
   W\testright{R_f} R.
$
\end{itemize}
\end{cla}


\begin{rem}\label{rem:L*B*6}
If $L = L^*$, see Fig.~\ref{fig:Wave_Curve_2E}, there are two paths to reach the right state
$R$, namely 
$
     L^*\testright{R_s} T^* \xrightarrow{'S_s} B_2^* \xrightarrow{S_f'} B_1^* \xrightarrow{R_f}  A_1 \xrightarrow{'S_f} R
$
or
$
     L^*\testright{R_s} T^* \xrightarrow{'S_o'} B_1^* \xrightarrow{R_f}  A_1 \xrightarrow{'S_f} R
$ 
but the triple shock rule guarantees they represent the same solution in $xt$-space.
In other words, the solution consists of a single wave group.
\end{rem}

\noindent{ \bf Riemann solution for $R$ in subregion $\Omega_2^i$.}
\medskip

We now consider $R$ in region $\Omega_2^i$.
The changes in the
Riemann solution when $R$ moves from
$\Omega_2^c\cup \Omega_2^f$ to $\Omega_2^i$
are similar to those 
described in Claims~\ref{cla:RSolution-Omega_2aU2c}, \ref{cla:RSolution-Omega_2e} when $R$ moves from
$\Omega_2^a\cup \Omega_2^d$  to $\Omega_2^g$, or in Claims~\ref{cla:RSolution-Omega_2bU2e}, \ref{cla:RSolution-Omega_2h}
when $R$ moves from
$\Omega_2^b\cup \Omega_2^e$  to $\Omega_2^h$.
Namely,

\begin{cla}\label{cla:RSolution-Omega_2f}
Refer to Fig.~\ref{fig:Wave_Curve_2H}(b).
For $L$ in the edge \GW of the
saturation triangle and $R$ in subregion $\Omega_2^i$ in
Figs.~\ref{fig:RegioesThetas}(b), \ref{fig:RRegionsOmegas}(b),
the Riemann solution is the one described in Claim~\ref{cla:RSolution-Omega_2bU2d},
replacing $T^*$ with $T_F$ and $L^*$ with $L_F$.
\end{cla}


\subsection{Right state in region \texorpdfstring{$\Gamma$}{O}}
\label{subsec:Gamma}
Region $\Gamma$, see Fig~\ref{fig:RegioesThetas}(c), is subdivided into subregions $\Gamma_1$  and $\Gamma_2$ by the fast-inflection locus $\mathcal{R}_P$-$\mathcal{I}_B$.

\begin{figure}[]
	\centering
  \subfigure[$R = (0.214134,0.71306)$ in $\Gamma_1$, below slow inflection locus. ]
	{\includegraphics[scale=0.27]{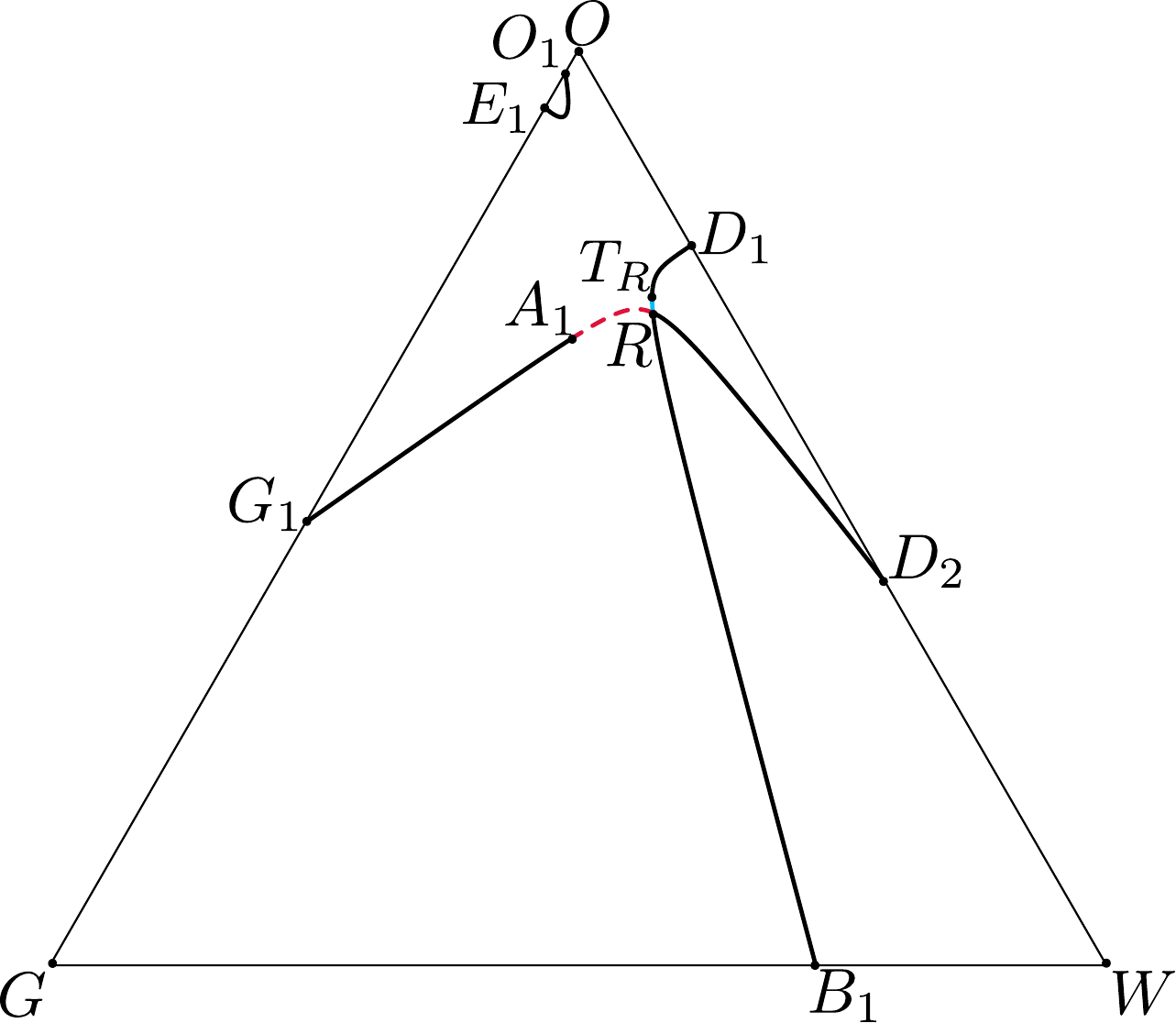}} 
   \hspace{6mm}
 \subfigure[$R = (0.280838,0.67656)$ in $\Gamma_1$, above slow inflection locus.  ]
	{\includegraphics[scale=0.27]{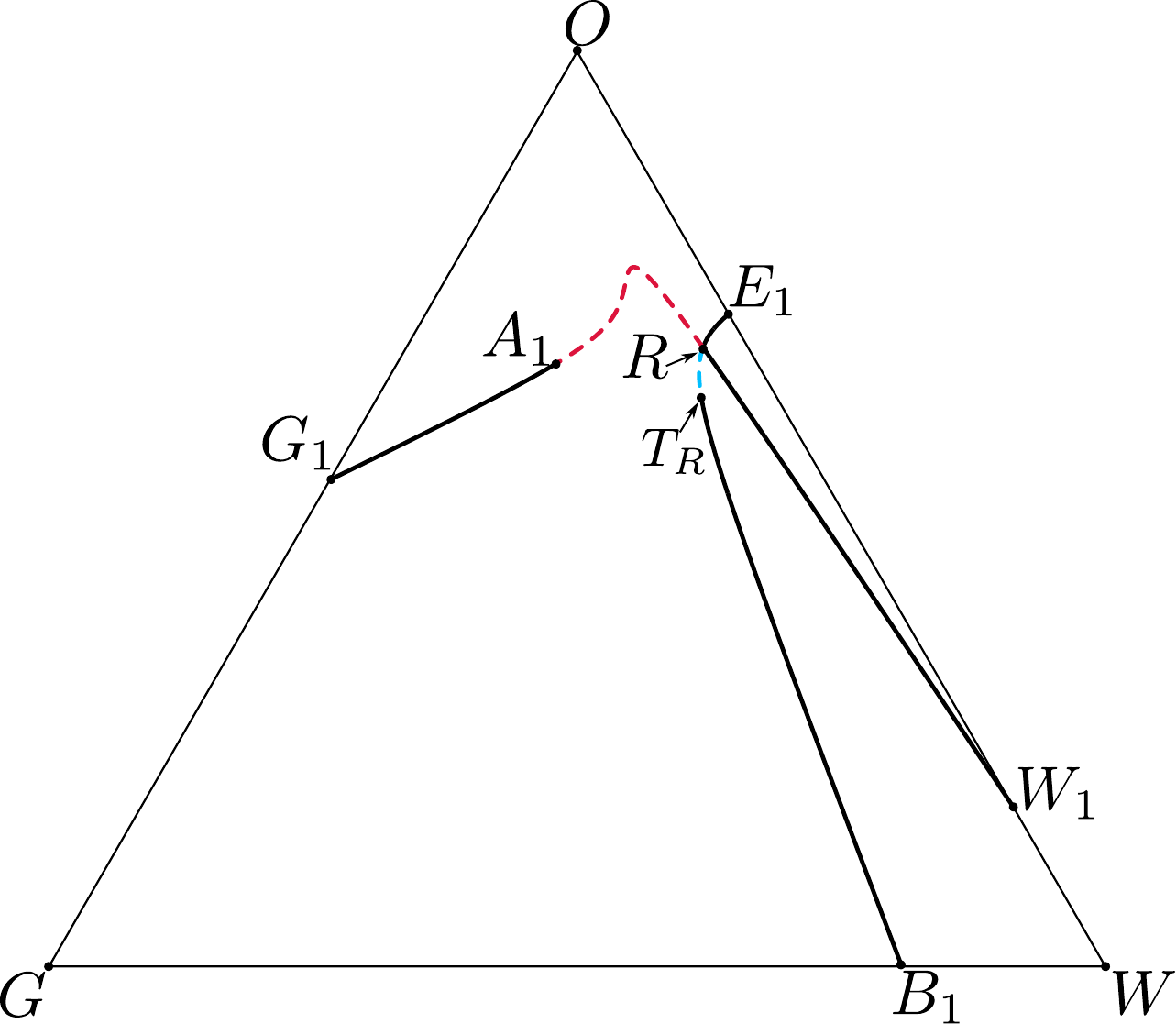}}  
	\subfigure[$R = (0.179281,0.587257)$ in $\Gamma_2$, below slow inflection locus. ]
	{\includegraphics[scale=0.27]{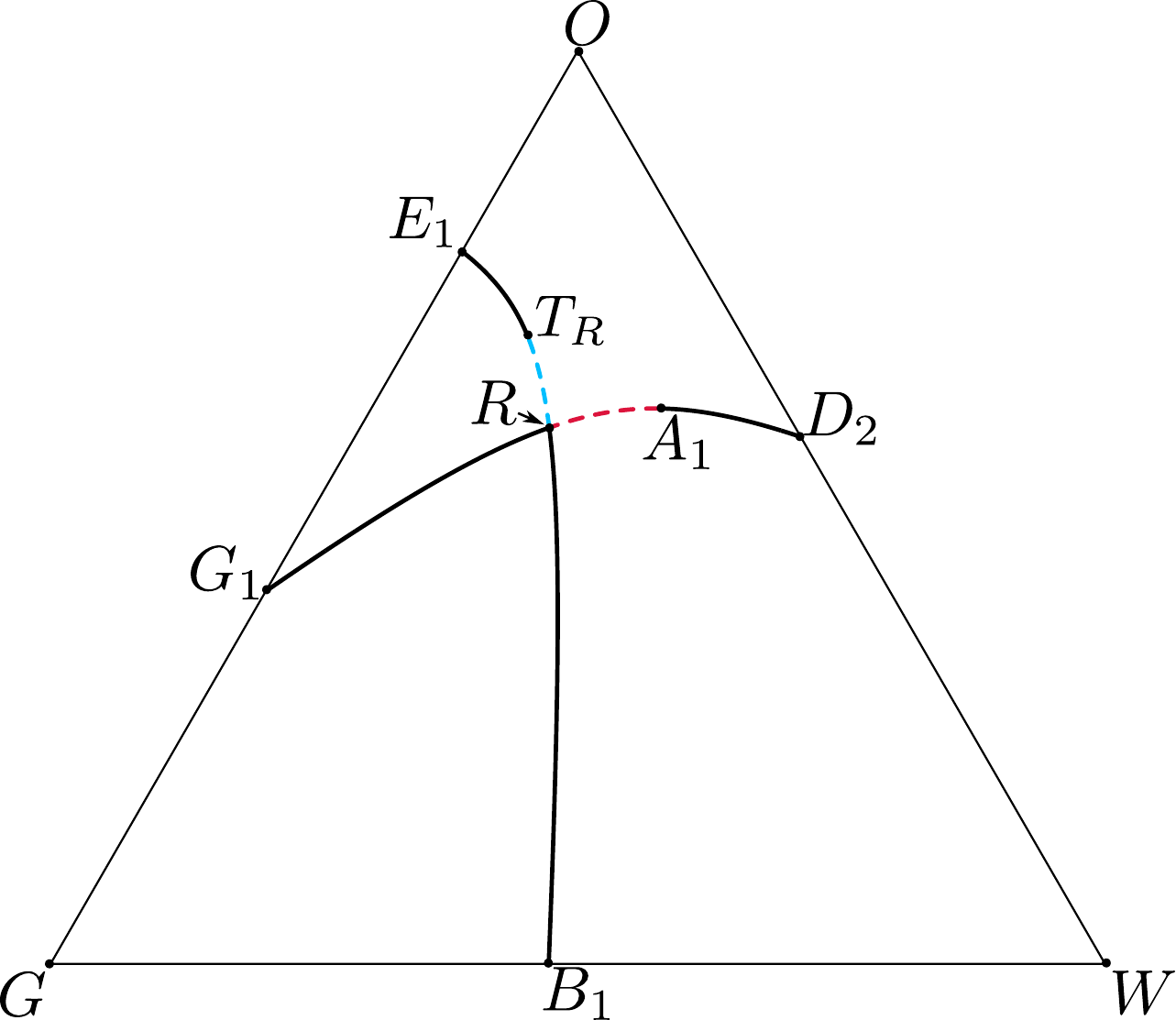}}  
	 \hspace{6mm}
	\subfigure[$R = (0.0827537,0,60369)$ in $\Gamma_2$, above slow inflection locus.]
	{\includegraphics[scale=0.27]{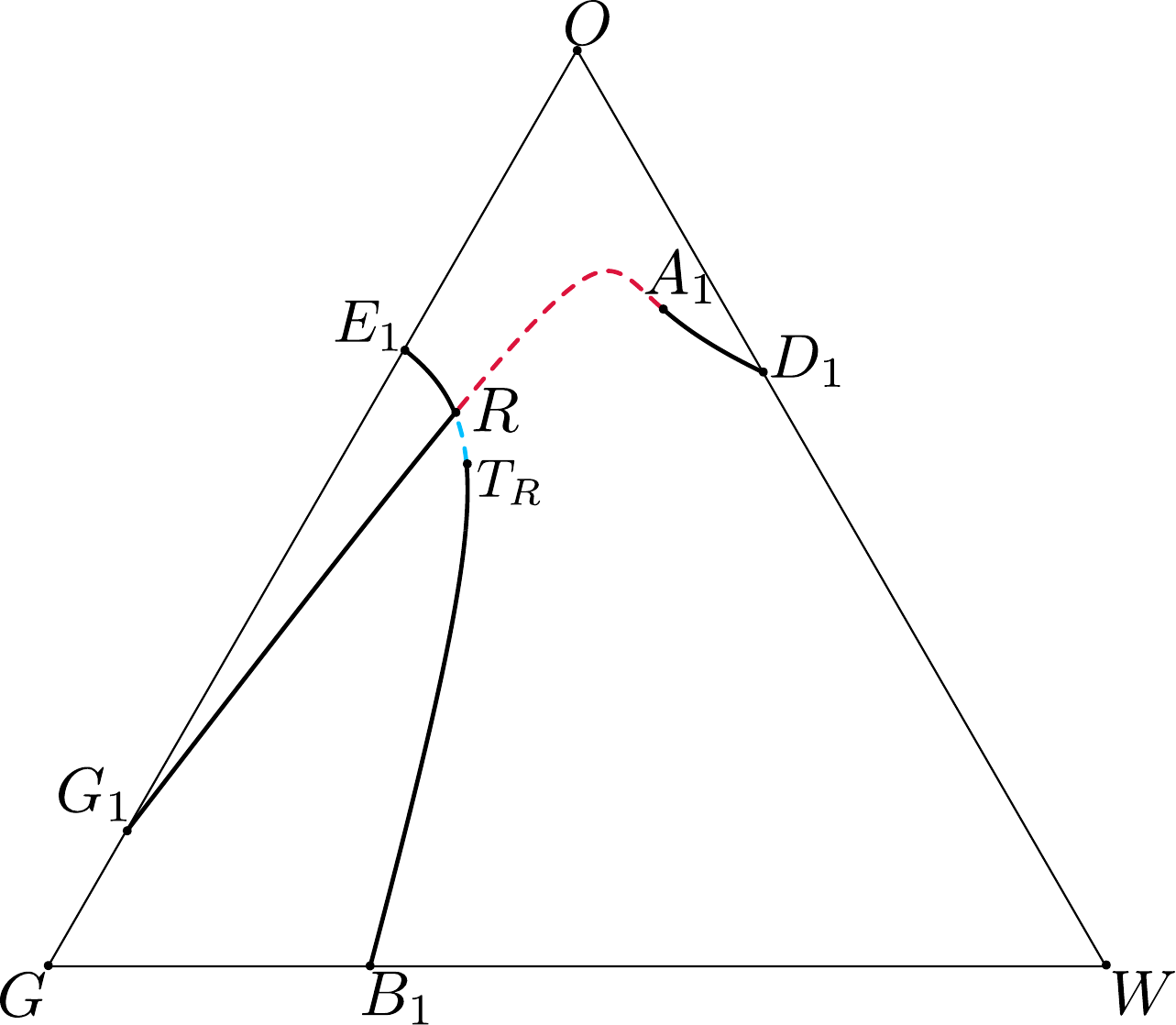}}   
		\caption{Hugoniot curves for states $R$ in $\Gamma_1$ and $\Gamma_2$, far
  from edge \GWc.
  	} 
	\label{fig:Hugoniot_Gamma_regions}
\end{figure}

Typical Hugoniot curves for $R$ in regions $\Gamma_1$ and $\Gamma_2$ far from the edge \GW are depicted in Figs.~\ref{fig:Hugoniot_Gamma_regions}. The admissible shocks are parameterized by the local $s$-shock segment $(R, T_R]$, with $\sigma(T_R; R) = \lambdas(T_R)$, and the local $f$-shock segment $(R, A_1]$, with $\sigma(A_1; R) = \lambdaf(A_1)$.

 \begin{figure}[]
	\centering
	\subfigure[$R$-Region $\Gamma_2$. ]
	{\includegraphics[scale=0.5]{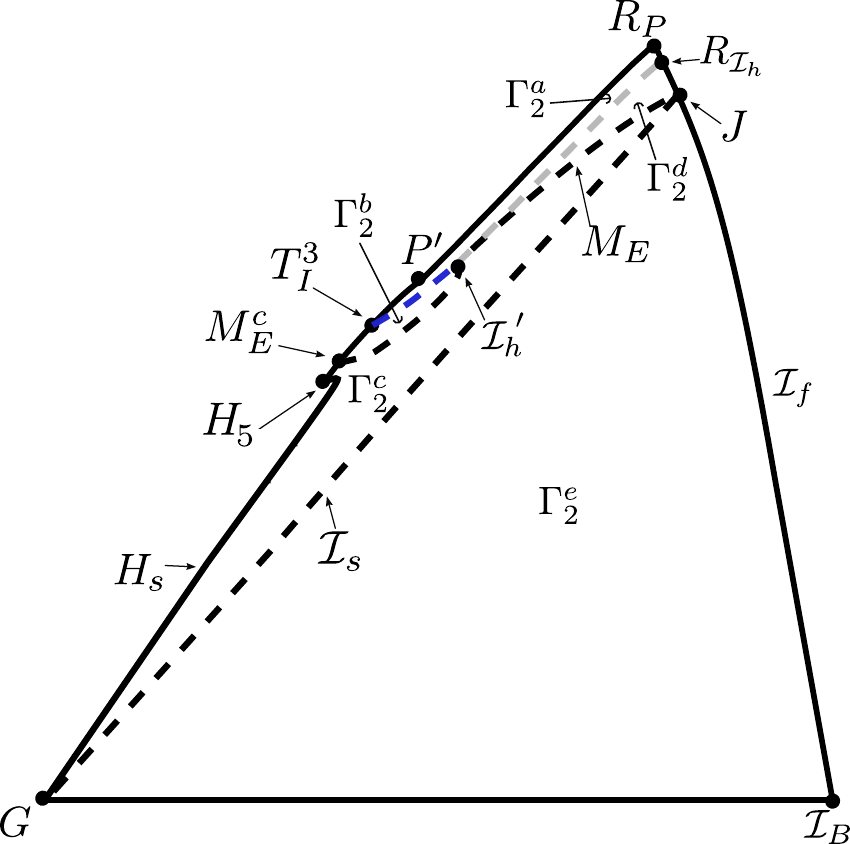}}  
	 \hspace{6mm}
	\subfigure[$R$-Region $\Gamma_1$.]
	{\includegraphics[scale=0.5]{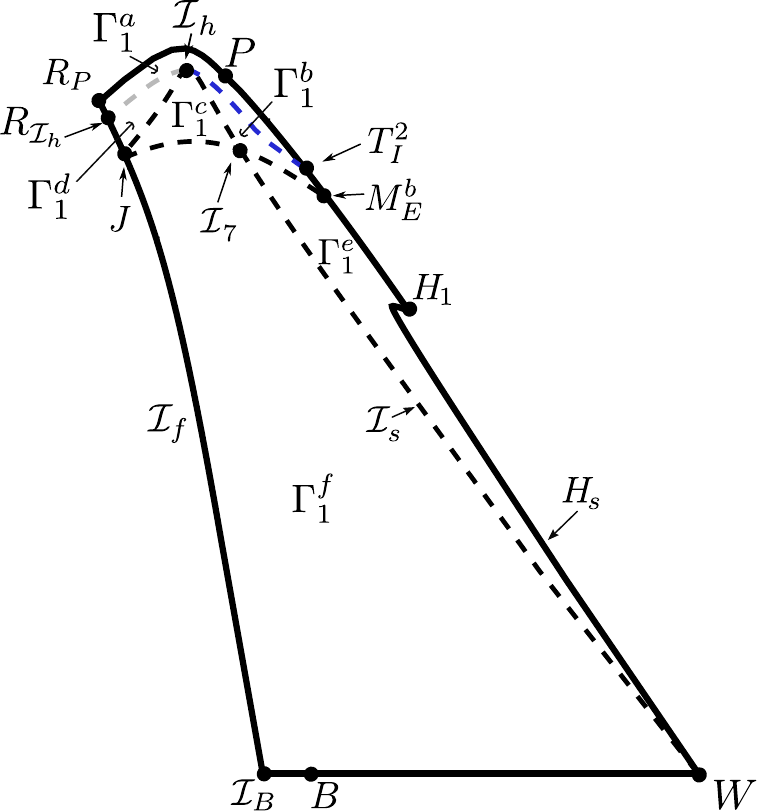}}   
	\caption{ Zoom of the $R$-regions $\Gamma_1$ and $\Gamma_2$
 in Fig.~\ref{fig:RegioesThetas}
	displaying their subdivisions $\Gamma_1^{\alpha}$, $\alpha \in \{a, b, c, d, e, f\}$ and $\Gamma_2^{\beta}$, $\beta \in \{a, b, c, d, e\}$.
 The boundary between $\Gamma_1^a$ and  $\Gamma_1^b$ consists of the segment $\mathcal{I}_h$-$T_I^2$ of curve $T_I$.
 The boundary between $\Gamma_1^b$ and $\Gamma_1^c$ is the segment $\mathcal{I}_{7}$-$\mathcal{I}_h$ of the
 $s$-inflection $\mathcal{I}_s$.
 The boundary between $\Gamma_1^c$ and $\Gamma_1^d$ is the segment $J$-$\mathcal{I}_h$ of the
 $s$-inflection $\mathcal{I}_s$.
Here, the state $J$, also shown in Fig.~\ref{fig:IntegralCurves_extensions}, 
is the intersection of the loci $M_E$,  $\mathcal{I}_s$, and $\mathcal{I}_f$.
 The boundary between $\Gamma_1^b$ and $\Gamma_1^e$ is the segment $\mathcal{I}_{7}$-$M_E^b$ of curve $M_E$.
 The boundary between $\Gamma_1^c$ and $\Gamma_1^f$ is the segment $J$-$\mathcal{I}_{7}$ of curve $M_E$.
 The boundary between $\Gamma_1^e$ and $\Gamma_1^f$ is the segment $\mathcal{I}_{7}$-$W$ of the slow inflection locus.
 The boundary between $\Gamma_2^a$ and $\Gamma_2^d$ is the composite segment $R_{\mathcal{I}_h}$-${\mathcal{I}_h}'$.
 The boundary between $\Gamma_2^a$ and $\Gamma_2^b$ is the segment ${\mathcal{I}_h}'$-$T_I^3$ of curve $T_I$.
 The boundary between $\Gamma_2^b$ and $\Gamma_2^c$ is the segment ${\mathcal{I}_h}'$-$M_E^c$ of $M_E$.
 The boundary between $\Gamma_2^c$ and $\Gamma_2^d$ is the segment ${\mathcal{I}_h}'$-$J$ of $M_E$.
 The boundary between $\Gamma_2^c$ and $\Gamma_2^e$ is the segment
 $J$-$G$ of the $s$-inflection locus.
 }
	\label{fig:RRegionsGamma}
\end{figure}
\begin{figure}[]
	\centering
 \subfigure[$R = (0.0803214, 0.609742)$ in $\Gamma_2^a$.]
	{\includegraphics[scale=0.27]{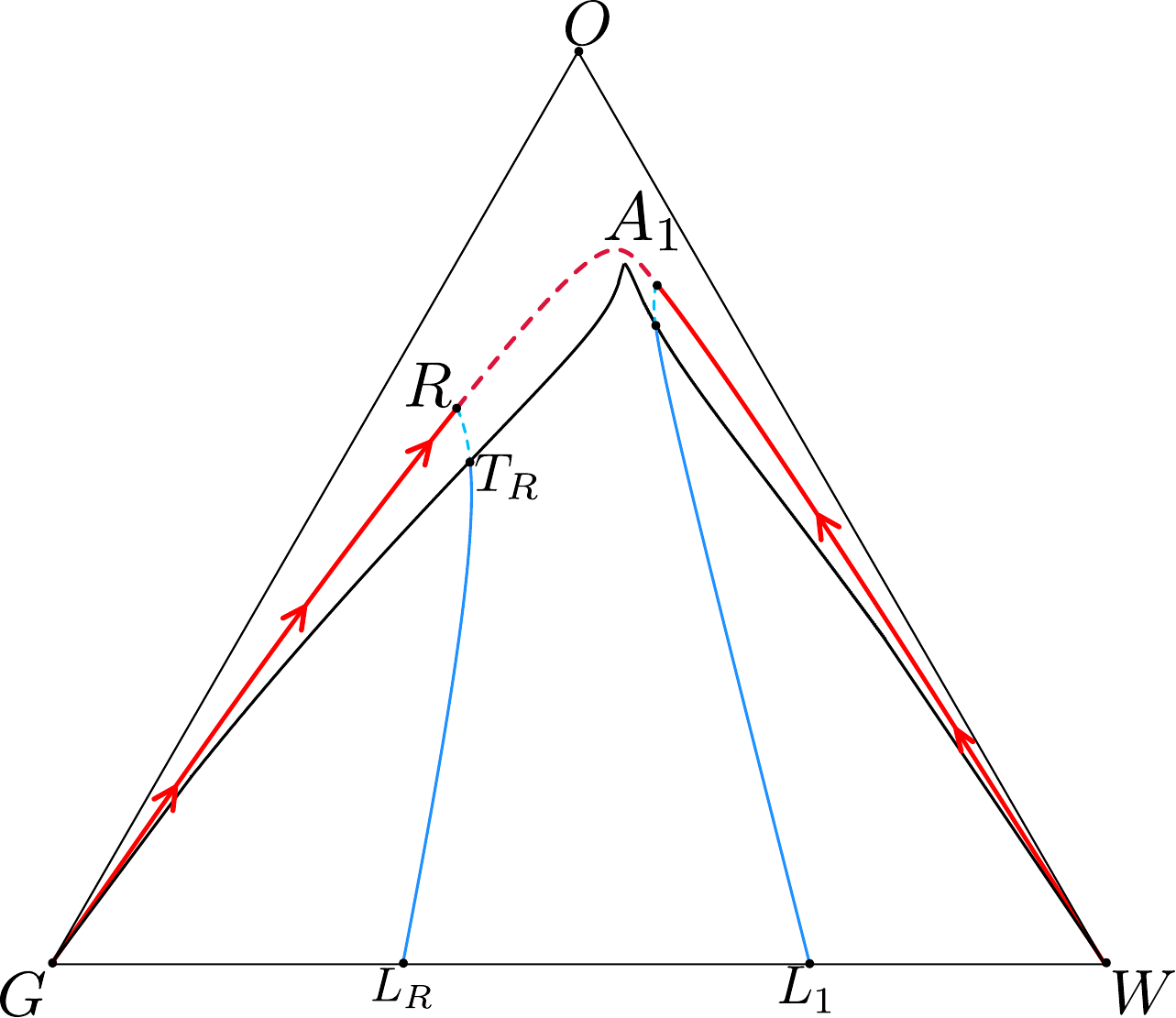}}  
 \hspace{6mm}
 \subfigure[$R = (0.0827537,  0.60369)$ in $\Gamma_2^b$. ]
	{\includegraphics[scale=0.27]{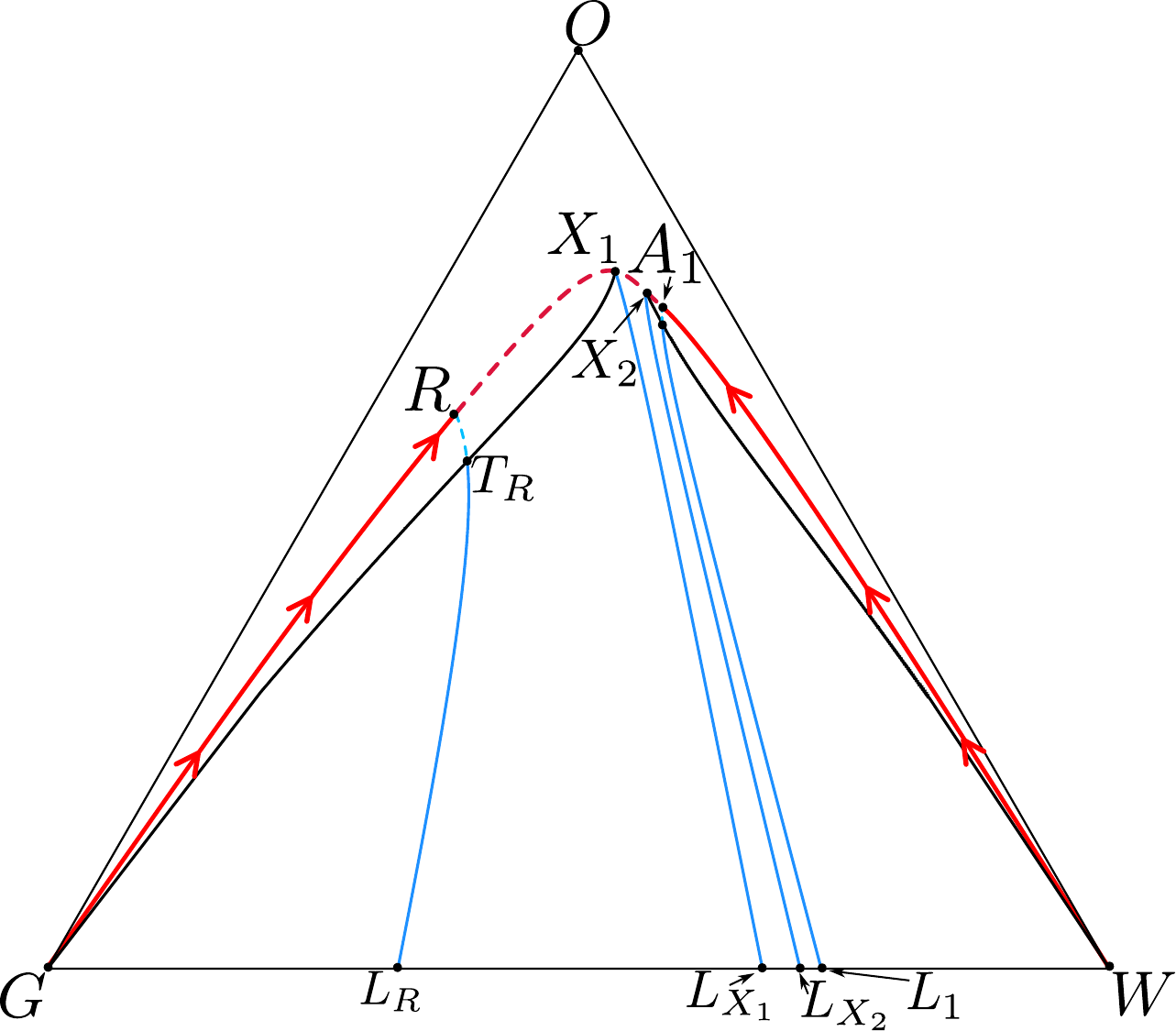}} 
  \hspace{6mm}
	\subfigure[$R = (0.0576305,  0.470357)$ in $\Gamma_2^c$. ]
	{\includegraphics[scale=0.27]{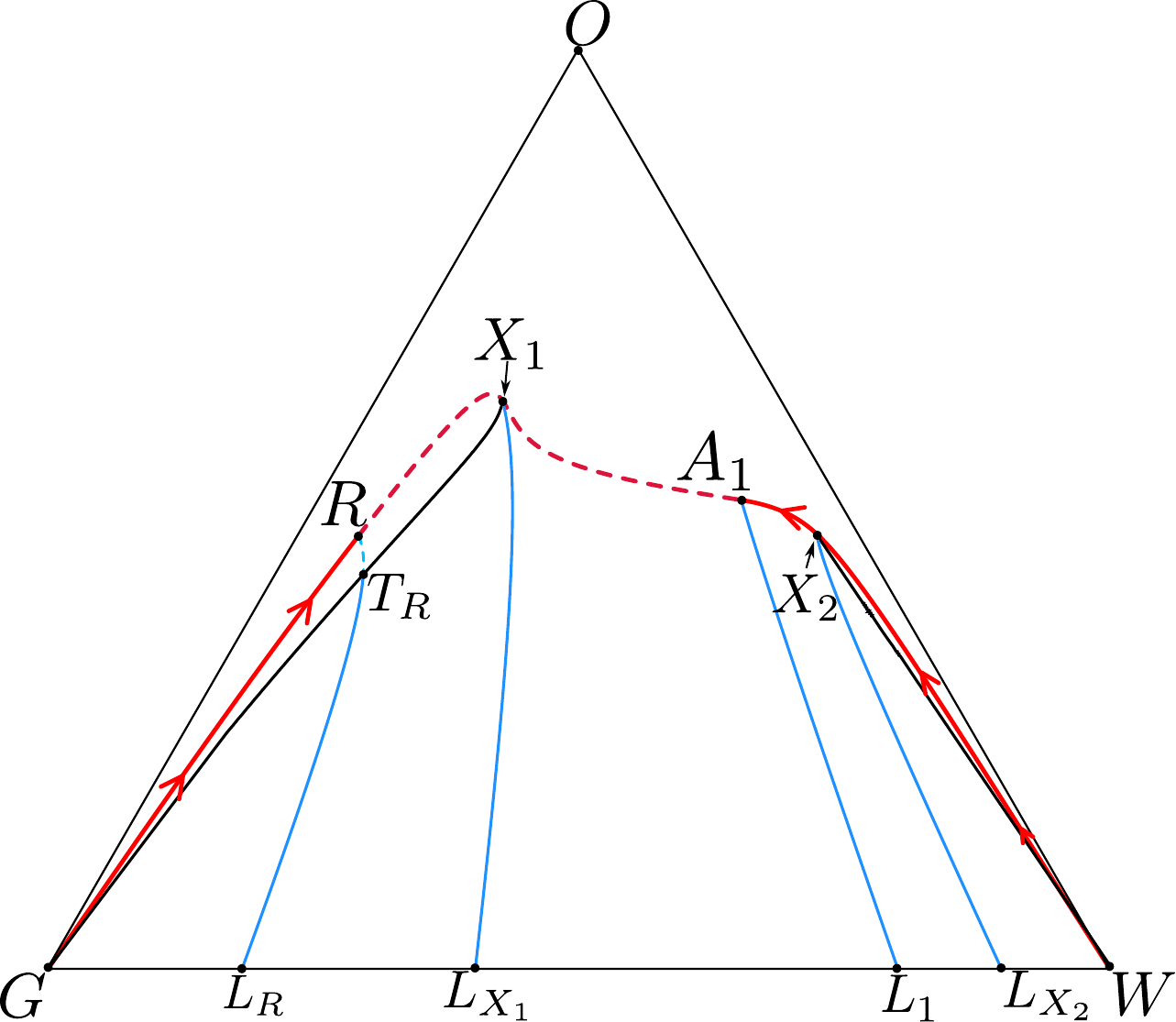}}  
	 \hspace{6mm}
  	\subfigure[$R = (0.125154,  0.716874)$ in $ \Gamma_2^d$.]
	{\includegraphics[scale=0.27]{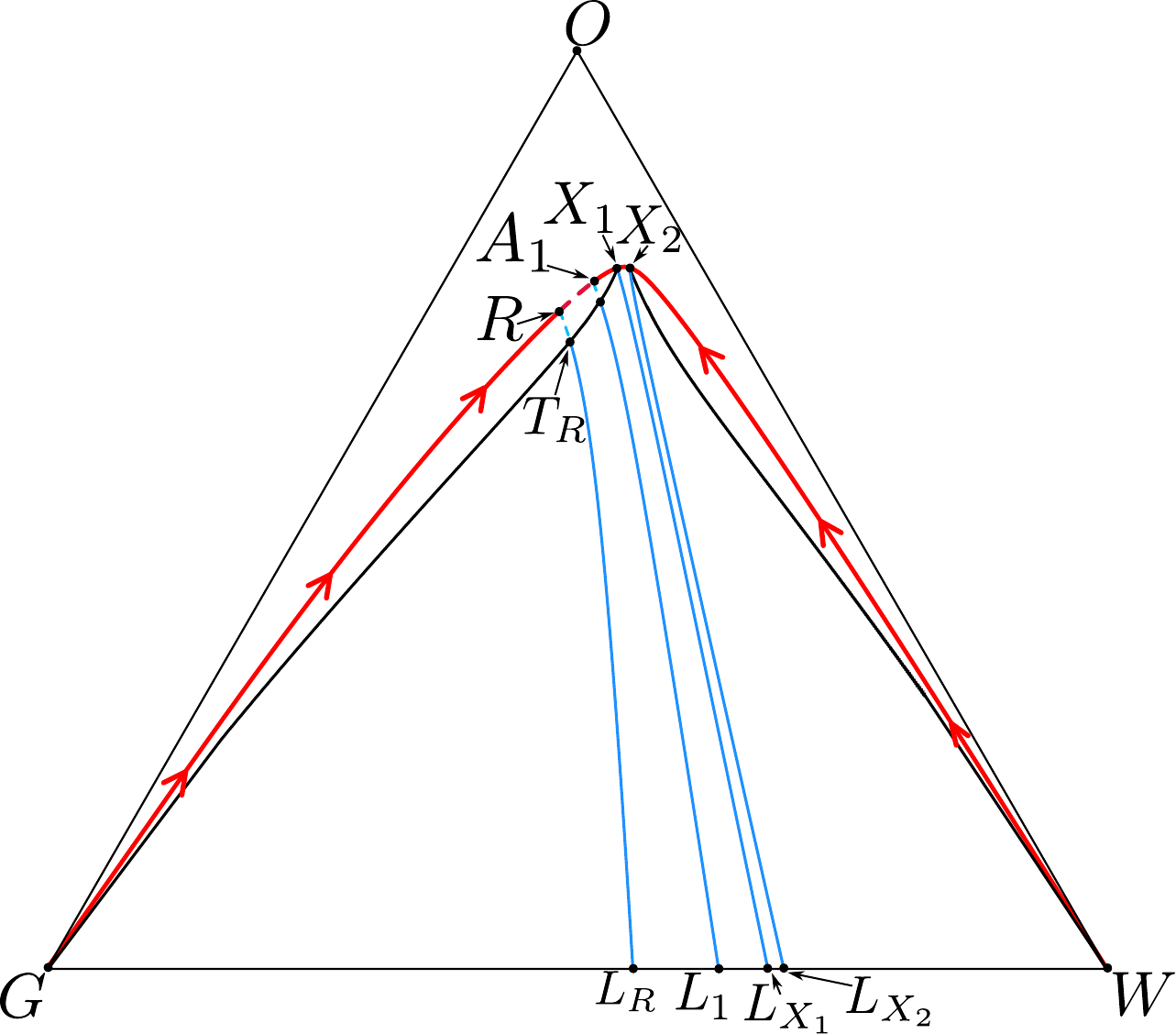}}   
 	 \hspace{6mm}
	\subfigure[$R = (0.179281, 0.587257)$ in $\Gamma_2^e$.]
	{\includegraphics[scale=0.27]{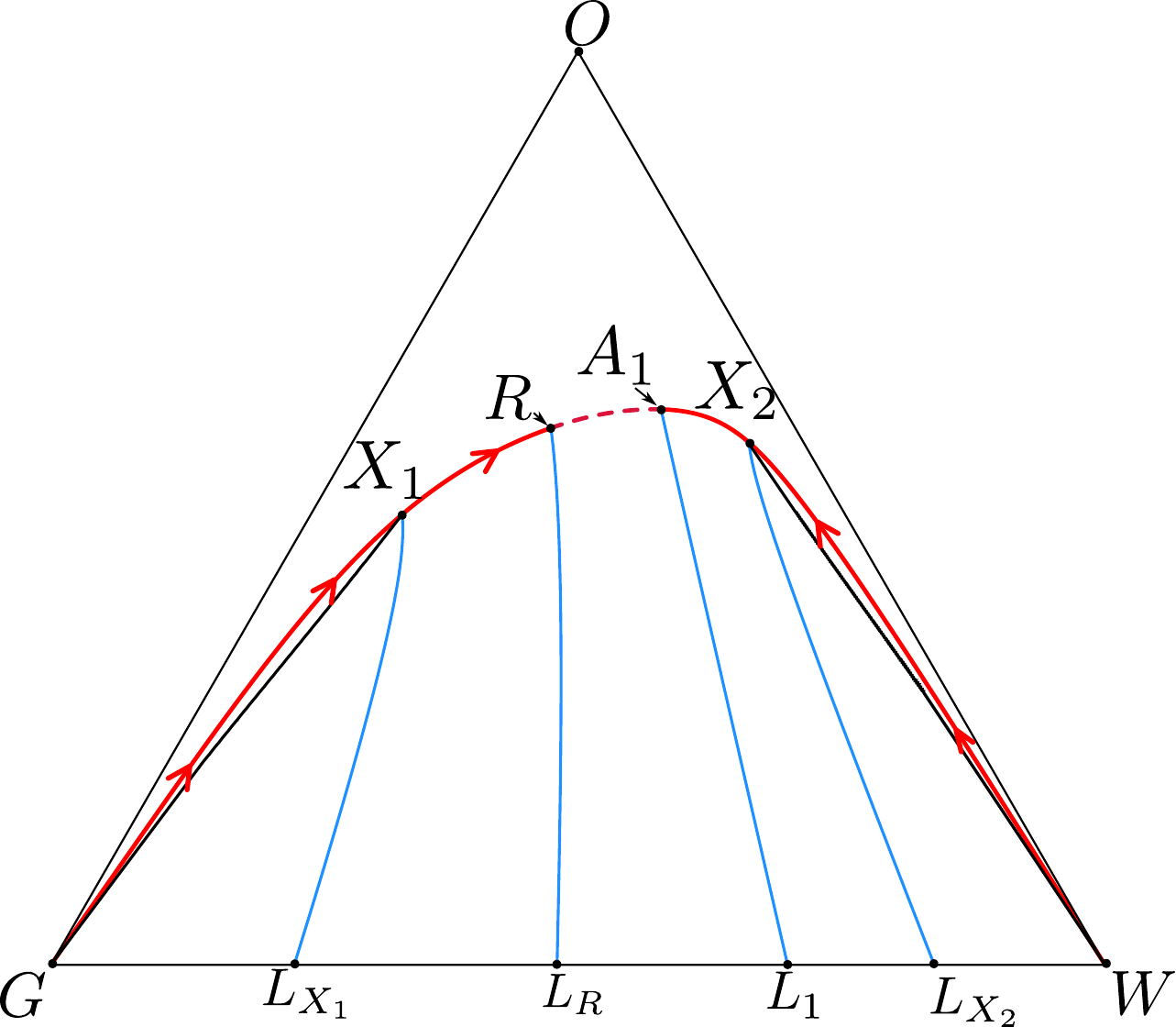}}   
 	\caption{Wave curves $W_f(R)$ and Riemann solutions for states $R$ in $\Gamma_2$.
  	} 
	\label{fig:Solution_Riemann_P_Gamma_2}
\end{figure}

Typical backward fast wave curves $W_f(R)$ for $R$ in $\Gamma_1$ and $\Gamma_2$
far from edge \GW are shown in Figs.~\ref{fig:Solution_Riemann_P_Gamma_1}
and \ref{fig:Solution_Riemann_P_Gamma_2}, where it is also shown geometric representations of Riemann solutions.
Such wave curves consist of two rarefaction segments 
and one shock segment comprising only one local branch.
If $R$ lies in $\Gamma_1$ the rarefaction segments are $[W, R)$ and $(A_1, G]$, whereas if $R$ lies in $\Gamma_2$ the rarefaction segments are $[G, R)$ and $(A_1, W]$. In both cases, the shock segment is $(R, A_1]$.

For $R$ near the edge \GWc, both the Hugoniot curves and the fast wave curves may
have nonlocal branches.
However, such branches play no role in the Riemann problems we consider; hence they are omitted.



\subsubsection{Subregion \texorpdfstring{$\Gamma_1$}{P}}\label{sec:Gamma1}
The region $\Gamma_1$ is further subdivided into subregions $\Gamma_1^\alpha$, $\alpha \in \{a, b, c, d, e,f\}$.
These subregions are primarily defined by two types of boundaries:
the first type, comprising the segment $T_I^2$-$\mathcal{I}_h$ of the curve $T_I$ together with the
$f$-rarefaction segment $\mathcal{I}_h$-$R_{\mathcal{I}_h}$, defines the subregion $\Gamma_1^a$, characterized by states $R$ for which the fast wave curve does not intersect the $s$-inflection locus;
the second type, comprising the segments $M_E^b$-$J$ of the curve $M_E$ and $J$-$\mathcal{I}_h$-$\mathcal{I}_{7}$ of the $s$-inflection locus, determines the type of the fast wave at the intersection points.

 \begin{figure}[]
	\centering
 \subfigure[$R = (0.165936, 0.769909)$ in $\Gamma_1^a$.]
	{\includegraphics[scale=0.27]{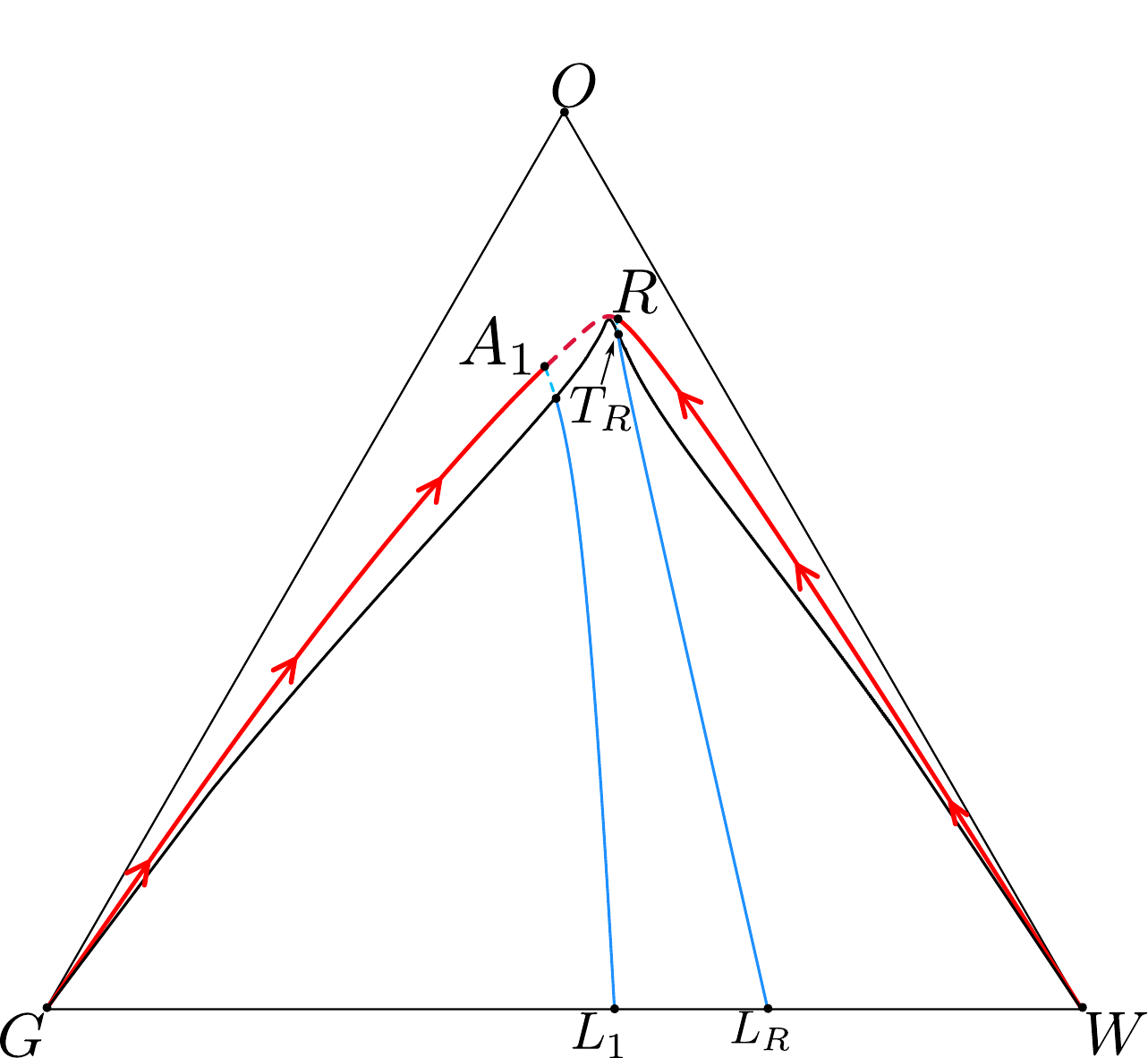}} 
  \hspace{6mm}
  \subfigure[$R = (0.207227, 0.737142)$ in  $\Gamma_1^b$.]
	{\includegraphics[scale=0.27]{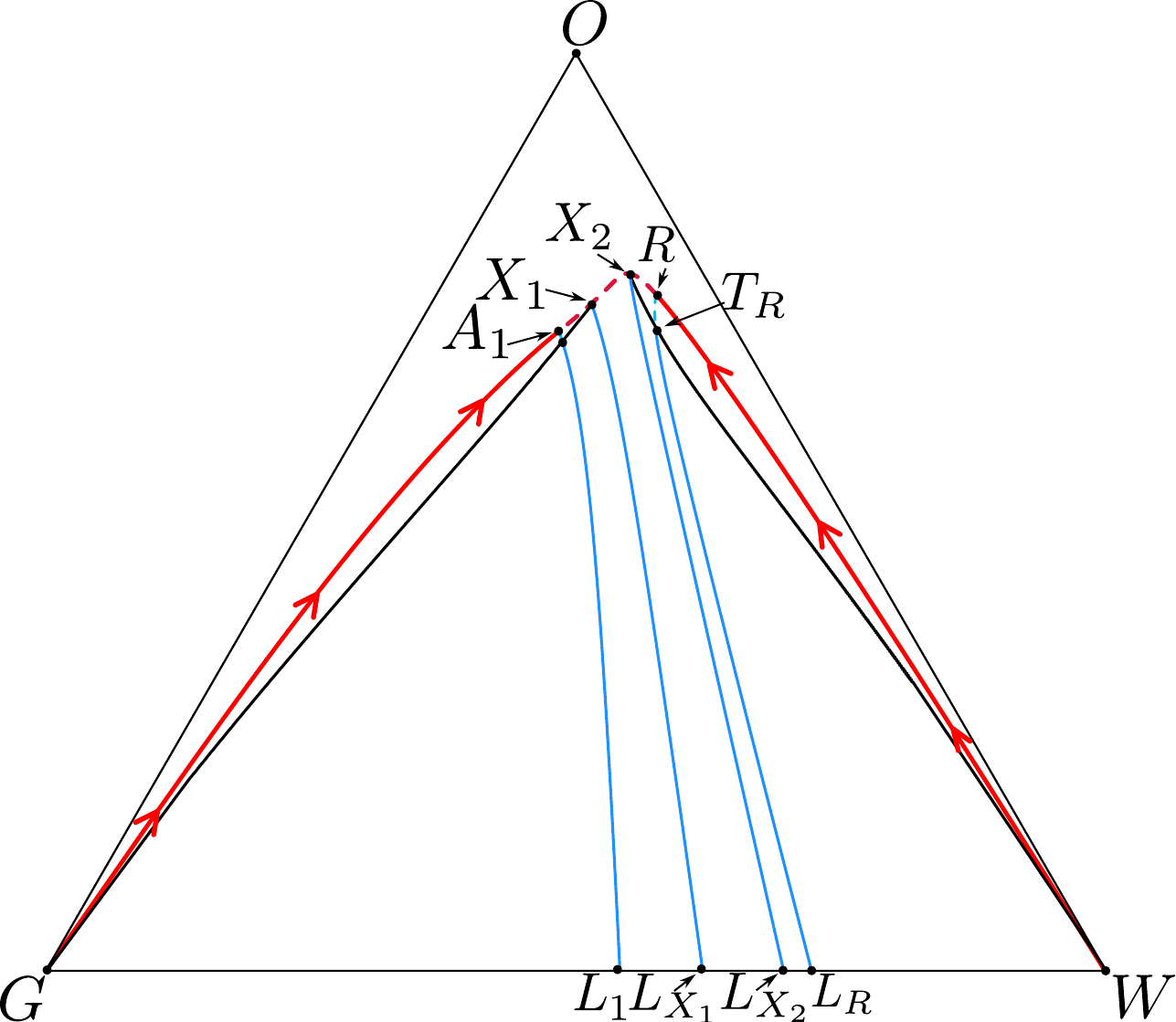}}
 \hspace{6mm}
  \subfigure[ $R = (0.180295, 0.744851)$ in $\Gamma_1^c$. ]
	{\includegraphics[scale=0.27]{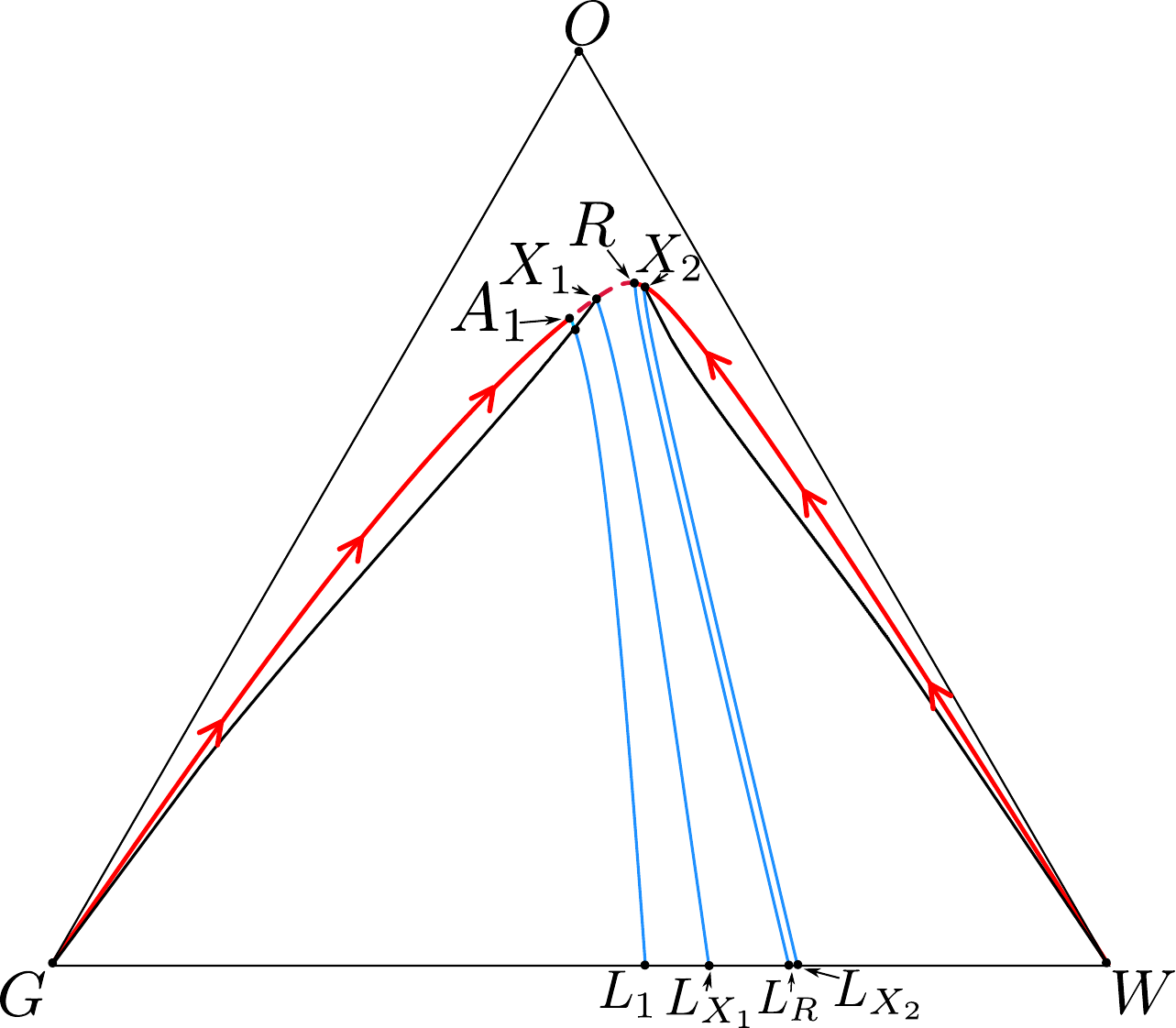}}  
	 \hspace{6mm}
	\subfigure[ $R = (0.14618, 0.752911)$ in  $\Gamma_1^d$.]
	{\includegraphics[scale=0.27]{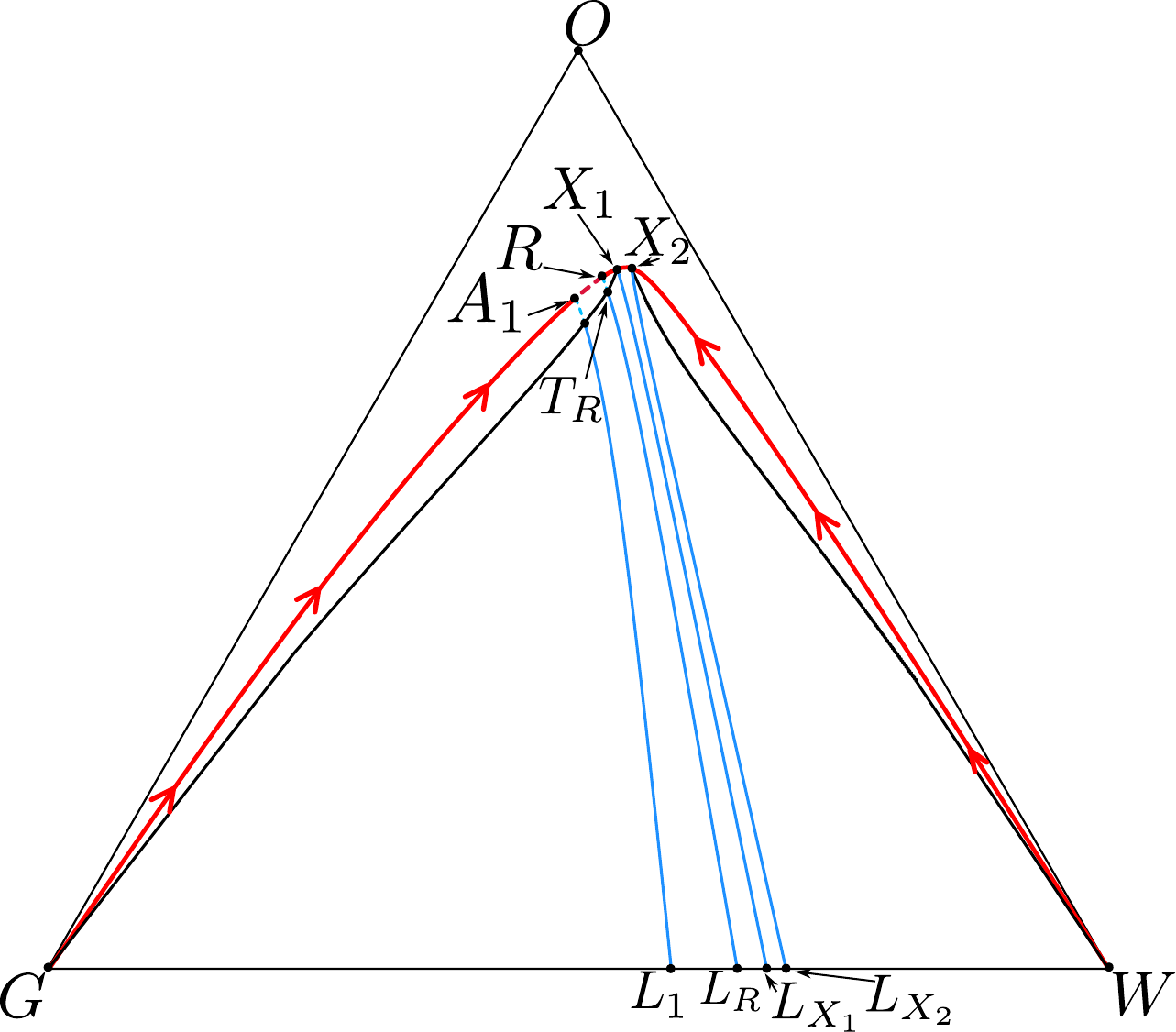}}  
	 \hspace{6mm}
	\subfigure[$R = (0.328456, 0.635357)$ in $\Gamma_1^e$.]
	{\includegraphics[scale=0.27]{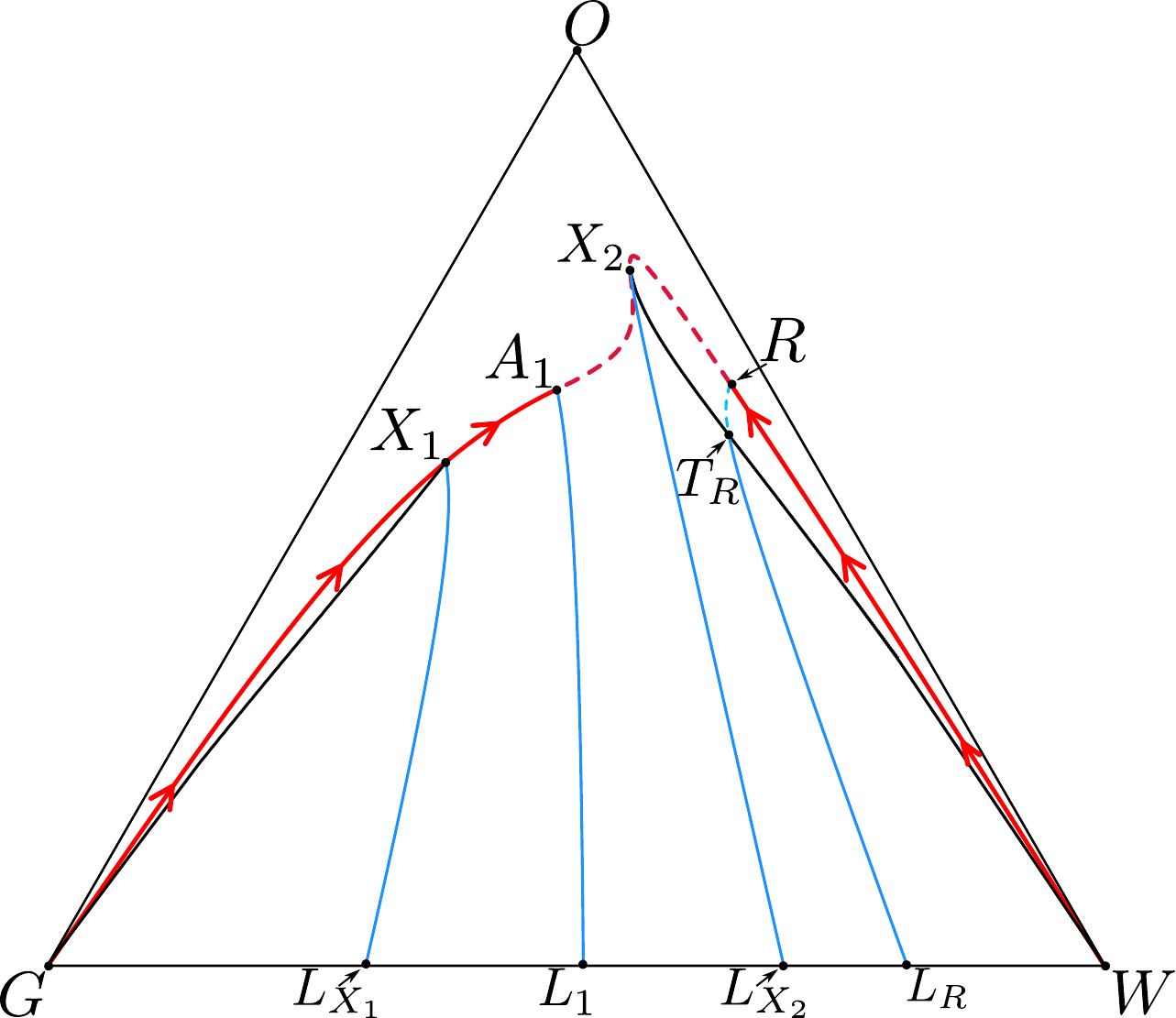}} 
 \hspace{6mm}
	\subfigure[ $R = (0.196822, 0.704416)$ in $\Gamma_1^f$. ]
	{\includegraphics[scale=0.27]{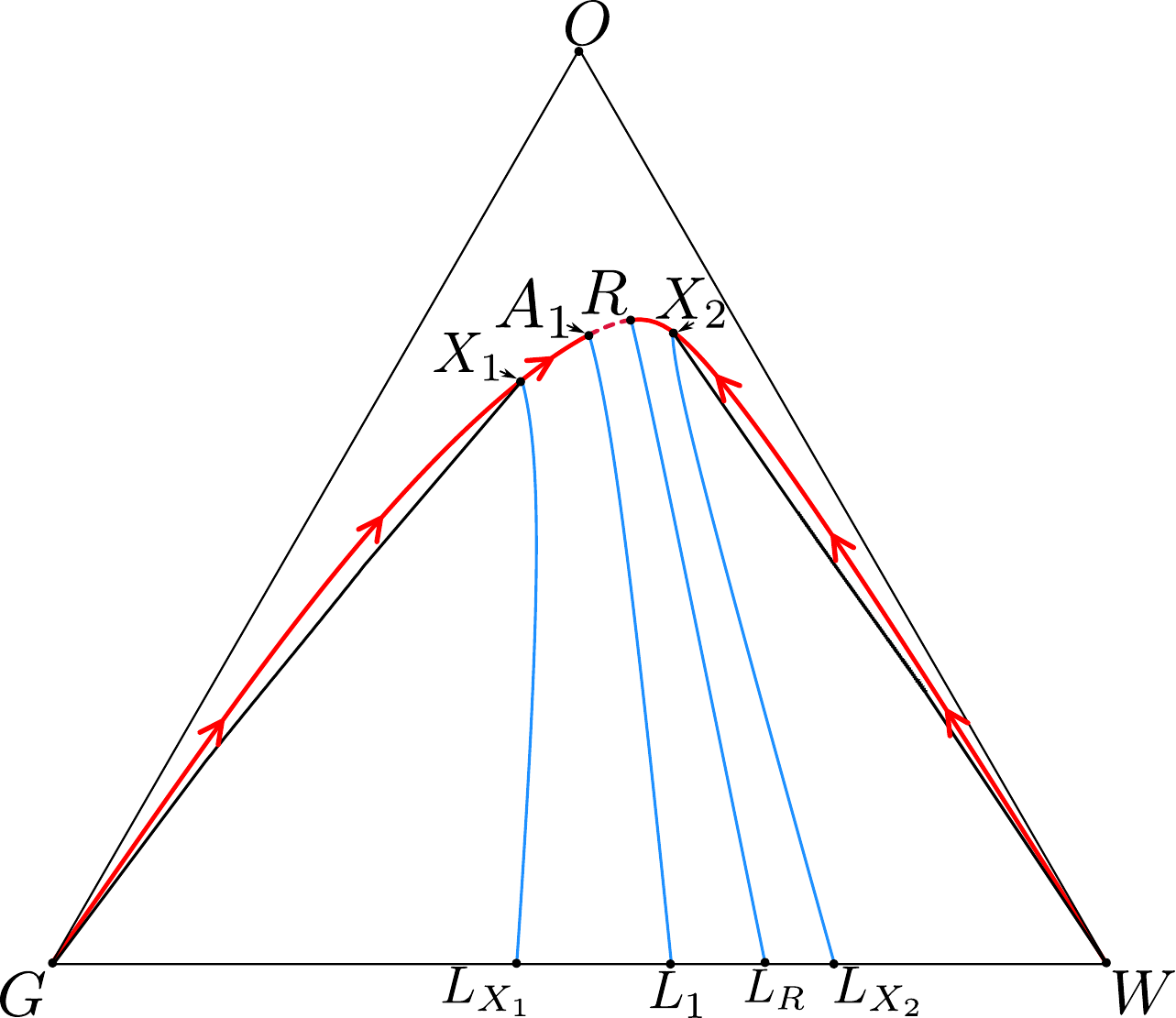}}  
	\caption{Wave curves $W_f(R)$ and Riemann solutions for states $R$ in $\Gamma_1$.
	} 
	\label{fig:Solution_Riemann_P_Gamma_1}
\end{figure}

\medskip


\noindent{ \bf Riemann solution for $R$ in subregion $\Gamma_1^a$.}
\medskip

\begin{cla}\label{cla:RSolution-Gamma_1a-new}
Refer to Fig.~\ref{fig:Solution_Riemann_P_Gamma_1}(a).
Let $L$ be a state in the edge \GW of the saturation triangle and $R$ be a state in subregion $\Gamma_1^a$ in Figs.~\ref{fig:RRegions Completed} and \ref{fig:RRegionsGamma}(b).
Let $A_1$ be the state in $\wm_f(R)$ satisfying
$\sigma(A_1; R) = \lambdaf(A_1)$.
Let $L_1$ and $L_R$ be the intersection points of the backward slow wave curves of $A_1$ and $R$ with the edge \GWc.
Then, 
\begin{itemize}
\item[(i)]if $L=G$, the Riemann solution is
$G\testright{R_f} A_1 \xrightarrow{'S_f\,} R$ ; 
\item[(ii)] if $L \in(G, L_1)$, the Riemann solution is
$
   L\testright{R_s} T_1 \xrightarrow{'S_s} M_1  \testright{R_f}A_1 \xrightarrow{'S_f} R,
$
where $T_1\in(G, T_R)_{\text{ext}}$ and $M_1\in(G, A_1)$; 
\item[(iii)] if $L \in[L_1, L_R)$, the Riemann solution is
$L\testright{R_s} T_1 \xrightarrow{'S_s} M_1 \testright{S_f} R,
$
where $T_1\in(G, T_R)_{\text{ext}}$ and $M_1\in(A_1,R)$;
\item[(iv)] if $L = L_R$, 
the Riemann solution is
$ L\testright{R_s} T_R \xrightarrow{'S_s} R$;
\item[(v)] if $L \in (L_R, W)$, 
the Riemann solution is
$ L\testright{R_s} T_1 \xrightarrow{'S_s} M_1 \testright{R_f}R,
$
where $T_1\in(T_R,W)_{\text{ext}}$ and  $M_1\in(R, W)$;
\item[(vi)]if $L=W$, the Riemann solution is 
$ W\testright{R_f} R$.
\end{itemize}
\end{cla}

\medskip

\noindent{ \bf Riemann solution for $R$ in subregion $\Gamma_1^b$.}

\begin{cla}\label{cla:RSolution-Gamma_1b-new}
Refer to Fig.~\ref{fig:Solution_Riemann_P_Gamma_1}(b).
Let $L$ be a state in the edge \GW of the saturation triangle and $R$ be a state in subregion $\Gamma_1^b$ in Figs.~\ref{fig:RRegions Completed} and \ref{fig:RRegionsGamma}(b).
Let $A_1$ be the state in $\wm_f(R)$ satisfying
$\sigma(A_1; R) = \lambdaf(A_1)$ and
$X_1$, $X_2$ be the intersection states of $\wm_f(R)$
with the $s$-inflection locus.
Let $L_1$, $L_{X_1}$, $L_{X_2}$ and $L_R$ be the intersection points of the backward slow wave curves of $A_1$, $X_1$, $X_2$ and $R$ with the edge \GWc.
Then,
\begin{itemize}
\item[(i)]if $L=G$, the Riemann solution is
$G\testright{R_f} A_1 \xrightarrow{'S_f\,} R$; 
\item[(ii)] if $L \in(G, L_1)$ the Riemann solution is
$
   L\testright{R_s} T_1 \xrightarrow{'S_s} M_1  \testright{R_f}A_1 \xrightarrow{'S_f} R,
$
where $T_1\in(G, X_1)_{\text{ext}}$ and $M_1\in(G, A_1)$; 
\item[(iii)] if $L \in [L_1, L_{X_1})$ or $L \in (L_{X_2}, L_R)$ 
 the Riemann solution is
 $L\testright{R_s} T_1 \xrightarrow{'S_s} M_1  \testright{S_f} R$ or
 $L\testright{R_s} T_2 \xrightarrow{'S_s} M_2  \testright{S_f} R$,
 where $T_1\in(G,X_1)_{\text{ext}}$ and $M_1\in[A_1,X_1)$ or $T_2\in(X_2, W)_{\text{ext}}$ and $M_2\in(X_2, R)$;
 \item[(iv)] if $L \in [L_{X_1}, L_{X_2}]$  
 the Riemann solution is
 $L\testright{R_s} M_1 \testright{S_f} R$,
 where $M_1\in[X_1, X_2]$;
 \item[(v)]if $L=L_R$, the Riemann solution is 
$ L\testright{R_s} T_R \xrightarrow{'S_s} R$;
\item[(vi)] if $L \in(L_R, W)$, the Riemann solution is $
   L\testright{R_s} T_1 \xrightarrow{'S_s} M_1  \testright{R_f} R$,
where $T_1\in( T_R,W)_{\text{ext}}$ and $M_1\in(R, W)$; 
\item[(vii)]if $L=W$, the Riemann solution is 
$ W\testright{R_f} R$.
\end{itemize}
\end{cla}

\noindent{ \bf Riemann solution for $R$ in subregion $\Gamma_1^c$.}
\medskip

\begin{cla}\label{cla:RSolution-Gamma_1c-new}
Refer to Fig.~\ref{fig:Solution_Riemann_P_Gamma_1}(c).
Let $L$ be a state in the edge \GW of the saturation triangle and $R$ be a state in subregion $\Gamma_1^c$ in Figs.~\ref{fig:RRegions Completed} and \ref{fig:RRegionsGamma}(b).
Let $A_1$ be the state in $\wm_f(R)$ satisfying
$\sigma(A_1; R) = \lambdaf(A_1)$ and
$X_1$, $X_2$ be the intersection states of $\wm_f(R)$
with the $s$-inflection locus.
Let $L_1$, $L_{X_1}$, $L_{X_2}$ and $L_R$ be the intersection points of the backward slow wave curves of $A_1$, $X_1$, $X_2$ and $R$ with the edge \GWc.
Then,
\begin{itemize}
\item[(i)]if $L=G$, the Riemann solution is
$G\testright{R_f} A_1 \xrightarrow{'S_f\,} R$ ; 
\item[(ii)] if $L \in(G, L_1)$ the Riemann solution is
$
   L\testright{R_s} T_1 \xrightarrow{'S_s} M_1  \testright{R_f}A_1 \xrightarrow{'S_f} R,
$
where $T_1\in(G, X_1)_{\text{ext}}$ and $M_1\in(G, A_1)$; 
\item[(iii)] if $L \in [L_1, L_{X_1})$ the Riemann solution is
 $L\testright{R_s} T_1 \xrightarrow{'S_s} M_1  \testright{S_f} R$,
 where $T_1\in(G,X_1)_{\text{ext}}$ and $M_1\in[A_1,X_1)$;
 \item[(iv)] if $L \in [L_{X_1}, L_R)$  
 the Riemann solution is
 $L\testright{R_s} M_1 \testright{S_f} R$,
 where $M_1\in[X_1, R)$;
 \item[(v)]if $L=L_R$, the Riemann solution is 
$ L\testright{R_s} R$;
\item[(vi)] if $L \in(L_R, L_{X_2}]$, the Riemann solution is 
$
   L\testright{R_s} M_1  \testright{R_f} R,
$
where  $M_1\in(R, X_2]$; 
\item[(vii)] if $L \in(L_{X_2},W)$, the Riemann solution is 
$
   L\testright{R_s} T_1 \xrightarrow{'S_s} M_1  \testright{R_f} R,
$
where $T_1\in(X_2,W)_{\text{ext}}$ and $M_1\in(R, W)$; 
\item[(viii)]if $L=W$, the Riemann solution is 
$ W\testright{R_f} R$.
\end{itemize}
\end{cla}

\noindent{ \bf Riemann solution for $R$ in subregion $\Gamma_1^d$.}
\medskip

\begin{cla}\label{cla:RSolution-Gamma_1d-new}
Refer to Fig.~\ref{fig:Solution_Riemann_P_Gamma_1}(d).
Let $L$ be a state in the edge \GW of the saturation triangle and $R$ be a state in subregion $\Gamma_1^d$ in Figs.~\ref{fig:RRegions Completed} and \ref{fig:RRegionsGamma}(b).
Let $A_1$ be the state in $\wm_f(R)$ satisfying
$\sigma(A_1; R) = \lambdaf(A_1)$ and
$X_1$, $X_2$ be the intersection states of $\wm_f(R)$
with the $s$-inflection locus.
Let $L_1$, $L_{X_1}$, $L_{X_2}$ and $L_R$ be the intersection points of the backward slow wave curves of $A_1$, $X_1$, $X_2$ and $R$ with the edge \GWc.
Then,
\begin{itemize}
\item[(i)]if $L=G$, the Riemann solution is
$G\testright{R_f} A_1 \xrightarrow{'S_f\,} R$; 
\item[(ii)] if $L \in(G, L_1)$ the Riemann solution is
$
   L\testright{R_s} T_1 \xrightarrow{'S_s} M_1  \testright{R_f}A_1 \xrightarrow{'S_f} R,
$
where $T_1\in(G, T_R)_{\text{ext}}$ and $M_1\in(G, A_1)$; 
\item[(iii)] if $L \in [L_1, L_{R})$ the Riemann solution is
 $L\testright{R_s} T_1 \xrightarrow{'S_s} M_1  \testright{S_f} R$,
 where $T_1\in(G,T_R)_{\text{ext}}$ and $M_1\in[A_1,R)$;
 \item[(iv)]if $L=L_R$, the Riemann solution is 
$ L \testright{R_s} T_R \xrightarrow{'S_s}  R$;
 \item[(v)] if $L \in  (L_R,L_{X_1})$  
 the Riemann solution is
 $L\testright{R_s} T_1 \xrightarrow{'S_s} M_1 \testright{R_f} R$,
 where $T_1\in(T_R, X_1)_{\text{ext}}$ and $M_1\in(R, X_1)$;
\item[(vi)] if $L \in[L_{X_1}, L_{X_2}]$, the Riemann solution is 
$
   L\testright{R_s} M_1  \testright{R_f} R,
$
where  $M_1\in[X_1, X_2]$; 
\item[(vii)] if $L \in(L_{X_2},W)$, the Riemann solution is 
$
   L\testright{R_s} T_1 \xrightarrow{'S_s} M_1  \testright{R_f} R,
$
where $T_1\in(X_2,W)_{\text{ext}}$ and $M_1\in(X_2, W)$; 
\item[(viii)]if $L=W$, the Riemann solution is 
$ W\testright{R_f} R$.
\end{itemize}
\end{cla}

\noindent{ \bf Riemann solution for $R$ in subregion $\Gamma_1^e$.}
\medskip

\begin{cla}\label{cla:RSolution-Gamma_1e-new}
Refer to Fig.~\ref{fig:Solution_Riemann_P_Gamma_1}(e).
Let $L$ be a state in the edge \GW of the saturation triangle and $R$ be a state in subregion $\Gamma_1^e$ in Figs.~\ref{fig:RRegions Completed} and \ref{fig:RRegionsGamma}(b).
Let $A_1$ be the state in $\wm_f(R)$ satisfying
$\sigma(A_1; R) = \lambdaf(A_1)$ and
$X_1$, $X_2$ be the intersection states of $\wm_f(R)$
with the $s$-inflection locus.
Let $L_1$, $L_{X_1}$, $L_{X_2}$ and $L_R$ be the intersection points of the backward slow wave curves of $A_1$, $X_1$, $X_2$ and $R$ with the edge \GWc.
Then,
\begin{itemize}
\item[(i)]if $L=G$, the Riemann solution is
$G\testright{R_f} A_1 \xrightarrow{'S_f\,} R$; 
\item[(ii)] if $L \in(G, L_{X_1})$ the Riemann solution is
$
   L\testright{R_s} T_1 \xrightarrow{'S_s} M_1  \testright{R_f}A_1 \xrightarrow{'S_f} R,
$
where $T_1\in(G, X_1)_{\text{ext}}$ and $M_1\in(G, X_1)$; 
\item[(iii)] if $L \in [L_{X_1}, L_{1})$ the Riemann solution is
 $L\testright{R_s} M_1  \testright{R_f} A_1 \xrightarrow{'S_f} R$,
 where $M_1\in[X_1,A_1)$;
 \item[(iv)] if $L \in  [L_1,L_{X_2}]$  
 the Riemann solution is
 $L\testright{R_s} M_1 \testright{S_f} R$,
 where $M_1\in[A_1, X_2]$  ;
\item[(v)] if $L \in(L_{X_2}, L_{R})$, the Riemann solution is 
$
L\testright{R_s} T_1 \xrightarrow{'S_s} M_1  \testright{S_f} R,
$
where $T_1\in(X_2, T_R)_{\text{ext}}$ and $M_1\in(X_2, R)$;
\item[(vi)]if $L=L_R$, the Riemann solution is 
$ L \testright{R_s} T_R \xrightarrow{'S_s}  R$;
\item[(vii)] if $L \in(L_{R},W)$, the Riemann solution is 
$
   L\testright{R_s} T_1 \xrightarrow{'S_s} M_1  \testright{R_f} R,
$
where $T_1\in(T_R,W)_{\text{ext}}$ and $M_1\in(R, W)$; 
\item[(viii)]if $L=W$, the Riemann solution is 
$ W\testright{R_f} R$.
\end{itemize}
\end{cla}

\noindent{ \bf Riemann solution for $R$ in subregion $\Gamma_1^f$.}
\medskip

\begin{cla}\label{cla:RSolution-Gamma_1f-new}
Refer to Fig.~\ref{fig:Solution_Riemann_P_Gamma_1}(f).
Let $L$ be a state in the edge \GW of the saturation triangle and $R$ be a state in subregion $\Gamma_1^f$ in Figs.~\ref{fig:RRegions Completed} and \ref{fig:RRegionsGamma}(b).
Let $A_1$ be the state in $\wm_f(R)$ satisfying
$\sigma(A_1; R) = \lambdaf(A_1)$ and
$X_1$, $X_2$ be the intersection states of $\wm_f(R)$
with the $s$-inflection locus.
Let $L_1$, $L_{X_1}$, $L_{X_2}$ and $L_R$ be the intersection points of the backward slow wave curves of $A_1$, $X_1$, $X_2$ and $R$ with the edge \GWc.
Then,
\begin{itemize}
\item[(i)]if $L=G$, the Riemann solution is
$G\testright{R_f} A_1 \xrightarrow{'S_f\,} R$; 
\item[(ii)] if $L \in(G, L_{X_1})$ the Riemann solution is
$
   L\testright{R_s} T_1 \xrightarrow{'S_s} M_1  \testright{R_f}A_1 \xrightarrow{'S_f} R,
$
where $T_1\in(G, X_1)_{\text{ext}}$ and $M_1\in(G, X_1)$; 
\item[(iii)] if $L \in [L_{X_1}, L_{1})$ the Riemann solution is
 $L\testright{R_s} M_1  \testright{R_f} A_1 \xrightarrow{'S_f} R$,
 where $M_1\in[X_1,A_1)$;
 \item[(iv)] if $L \in  [L_1,L_{R})$  
 the Riemann solution is
 $L\testright{R_s} M_1 \testright{S_f} R$,
 where $M_1\in[A_1, R)$;
 \item[(v)]if $L=L_R$, the Riemann solution is 
$ L \testright{R_s}  R$;
\item[(vi)] if $L \in(L_{R},L_{X_2}]$ the Riemann solution is 
$
L\testright{R_s} M_1  \testright{R_f} R,
$
where $M_1\in(R,X_2]$; 
\item[(vii)] if $L \in(L_{X_2},W)$ the Riemann solution is 
$
   L\testright{R_s} T_1 \xrightarrow{'S_s} M_1  \testright{R_f} R,
$
where $T_1\in(X_2,W)_{\text{ext}}$ and $M_1\in(X_2, W)$; 
\item[(viii)]if $L=W$, the Riemann solution is 
$ W\testright{R_f} R$.
\end{itemize}
\end{cla}

\medskip

\subsubsection{Subregion \texorpdfstring{$\Gamma_2$}{P}}\label{sec:Gamma2}

The region $\Gamma_2$ is further subdivided into subregions $\Gamma_2^\alpha$, $\alpha \in \{a, b, c, d, e\}$.

Similar to region $\Gamma_1$, these subregions are primarily defined by two types of boundaries.
The first type comprises
the composite segment $R_{\mathcal{I}_h}$-${\mathcal{I}_h}'$, which separates
$\Gamma_2^a$ and $\Gamma_2^d$, together with the segment ${\mathcal{I}_h}'$-$T_I^3$ of curve $T_I$, which separates subregions $\Gamma_2^a$ and $\Gamma_2^b$. The subregion $\Gamma_2^a$ consists of states $R$ for which the fast wave curve does not intersect the $s$-inflection locus.
The second type comprises two segments of the $M_E$ curve; namely, $J$-${\mathcal{I}_h}'$
separating subregions $\Gamma_2^d$ and $\Gamma_2^c$ and
${\mathcal{I}_h}'$-$M_E^c$ separating
$\Gamma_2^c$ and $\Gamma_2^b$.
For $R$ in region $\Gamma_2^d$ the fast wave curve intersects the $s$-inflection at two points, both in the rarefaction segment $(A_1, W]$ in Fig.~\ref{fig:Solution_Riemann_P_Gamma_2}(d); for $R$ in region $\Gamma_2^c$ one of the intersection points occurs in the shock segment $(R, A_1]$ and the other in the rarefaction segment $(A_1, W]$ in Fig.~\ref{fig:Solution_Riemann_P_Gamma_2}(c); for $R$ in region $\Gamma_2^b$ both intersection points occur in the shock segment $(R, A_1]$ in Fig.~\ref{fig:Solution_Riemann_P_Gamma_2}(b).
Finally the $s$-inflection segment $G$-$J$ separates regions $\Gamma_2^c$ and $\Gamma_2^e$; for $R$ in region $\Gamma_2^e$ the fast wave curve intersects the $s$-inflection at two points, one in the rarefaction segment $[G, R)$ and the other in the rarefaction segment $[W, A_1)$ in Fig.~\ref{fig:Solution_Riemann_P_Gamma_2}(e).

\medskip

\noindent{ \bf Riemann solution for $R$ in subregion $\Gamma_2^a$.}
\medskip

\begin{cla}\label{cla:RSolution-Gamma_2a-new}
Refer to Fig.~\ref{fig:Solution_Riemann_P_Gamma_2}(a).
Let $L$ be a state in the edge \GW of the saturation triangle and $R$ be a state in subregion $\Gamma_2^a$ in Figs.~\ref{fig:RRegions Completed} and \ref{fig:RRegionsGamma}(a).
Let $A_1$ be the state in $\wm_f(R)$ satisfying
$\sigma(A_1; R) = \lambdaf(A_1)$.
Let $L_1$ and $L_R$ be the intersection points of the backward slow wave curves of $A_1$ and $R$ with the edge \GWc.
Then, 
\begin{itemize}
\item[(i)]if $L=G$, the Riemann solution is
$G\testright{R_f} R$; 
\item[(ii)] if $L \in(G, L_R)$, the Riemann solution is
$
   L\testright{R_s} T_1 \xrightarrow{'S_s} M_1  \testright{R_f} R,
$
where $T_1\in(G, T_R)_{\text{ext}}$ and $M_1\in(G, R)$; 
\item[(iii)] if $L = L_R$, 
the Riemann solution is
$ L\testright{R_s} T_R \xrightarrow{'S_s} R$;
\item[(iv)] if $L \in(L_R, L_1]$, the Riemann solution is
$L\testright{R_s} T_1 \xrightarrow{'S_s} M_1 \testright{S_f} R,
$
where $T_1\in(T_R,W)_{\text{ext}}$ and $M_1\in(R,A_1]$;

\item[(v)] if $L \in (L_1, W)$, 
the Riemann solution is
$ L\testright{R_s} T_1 \xrightarrow{'S_s} M_1 \testright{R_f} A_1 \xrightarrow{'S_f} R,
$
where $T_1\in(T_R,W)_{\text{ext}}$ and  $M_1\in(A_1, W)$;
\item[(vi)]if $L=W$, the Riemann solution is 
$ W\testright{R_f} A_1 \xrightarrow{'S_f} R$.
\end{itemize}
\end{cla}

\medskip

\noindent{ \bf Riemann solution for $R$ in subregion $\Gamma_2^b$.}

\begin{cla}\label{cla:RSolution-Gamma_2b-new}
Refer to Fig.~\ref{fig:Solution_Riemann_P_Gamma_2}(b).
Let $L$ be a state in the edge \GW of the saturation triangle and $R$ be a state in subregion $\Gamma_2^b$ in Figs.~\ref{fig:RRegions Completed} and \ref{fig:RRegionsGamma}(a).
Let $A_1$ be the state in $\wm_f(R)$ satisfying
$\sigma(A_1; R) = \lambdaf(A_1)$ and
$X_1$, $X_2$ be the intersection states of $\wm_f(R)$
with the $s$-inflection locus.
Let $L_1$, $L_{X_1}$, $L_{X_2}$ and $L_R$ be the intersection points of the backward slow wave curves of $A_1$, $X_1$, $X_2$ and $R$ with the edge \GWc.
Then,
\begin{itemize}
\item[(i)]if $L=G$, the Riemann solution is
$G\testright{R_f} R$; 
\item[(ii)] if $L \in(G, L_R)$, the Riemann solution is
$
   L\testright{R_s} T_1 \xrightarrow{'S_s} M_1  \testright{R_f} R,
$
where $T_1\in(G, T_R)_{\text{ext}}$ and $M_1\in(G, R)$;
\item[(iii)]if $L=L_R$, the Riemann solution is 
$ L\testright{R_s} T_R \xrightarrow{'S_s} R$;
\item[(iv)] if $L \in (L_R, L_{X_1})$ or $L \in (L_{X_2}, L_1]$ 
 the Riemann solution is
 $L\testright{R_s} T_1 \xrightarrow{'S_s} M_1  \testright{S_f} R$ or
 $L\testright{R_s} T_2 \xrightarrow{'S_s} M_2  \testright{S_f} R$,
 where $T_1\in(T_R,X_1)_{\text{ext}}$ and $M_1\in(R,X_1)$ or $T_2\in(X_2, W)_{\text{ext}}$ and $M_2\in(X_2,A_1]$;
 \item[(v)] if $L \in [L_{X_1}, L_{X_2}]$  
 the Riemann solution is
 $L\testright{R_s} M_1 \testright{S_f} R$,
 where $M_1\in[X_1, X_2]$;
\item[(vi)] if $L \in(L_1, W)$ the Riemann solution is $
   L\testright{R_s} T_1 \xrightarrow{'S_s} M_1 \testright{R_f} A_1 \xrightarrow{'S_f} R,
$
where $T_1\in( X_2,W)_{\text{ext}}$ and $M_1\in(A_1, W)$; 
\item[(vii)]if $L=W$, the Riemann solution is 
$ W \testright{R_f} A_1 \xrightarrow{'S_f} R$.
\end{itemize}
\end{cla}

\noindent{ \bf Riemann solution for $R$ in subregion $\Gamma_2^c$.}
\medskip

\begin{cla}\label{cla:RSolution-Gamma_2c-new}
Refer to Fig.~\ref{fig:Solution_Riemann_P_Gamma_2}(c).
Let $L$ be a state in the edge \GW of the saturation triangle and $R$ be a state in subregion $\Gamma_2^c$ in Figs.~\ref{fig:RRegions Completed} and \ref{fig:RRegionsGamma}(a).
Let $A_1$ be the state in $\wm_f(R)$ satisfying
$\sigma(A_1; R) = \lambdaf(A_1)$ and
$X_1$, $X_2$ be the intersection states of $\wm_f(R)$
with the $s$-inflection locus.
Let $L_1$, $L_{X_1}$, $L_{X_2}$ and $L_R$ be the intersection points of the backward slow wave curves of $A_1$, $X_1$, $X_2$ and $R$ with the edge \GWc.
Then,
\begin{itemize}
\item[(i)]if $L=G$, the Riemann solution is
$G\testright{R_f} R$; 
\item[(ii)] if $L \in(G, L_R)$, the Riemann solution is
$
   L\testright{R_s} T_1 \xrightarrow{'S_s} M_1  \testright{R_f} R,
$
where $T_1\in(G, T_R)_{\text{ext}}$ and $M_1\in(G, R)$;
\item[(iii)]if $L=L_R$, the Riemann solution is 
$ L\testright{R_s} T_R \xrightarrow{'S_s} R$;
\item[(iv)] if $L \in (L_R, L_{X_1})$ the Riemann solution is
 $L\testright{R_s} T_1 \xrightarrow{'S_s} M_1  \testright{S_f} R$,
 where $T_1\in(T_R,X_1)_{\text{ext}}$ and $M_1\in(R,X_1)$;
\item[(v)] if $L \in [L_{X_1},L_1]$ the Riemann solution is
 $L\testright{R_s} M_1  \testright{S_f} R$,
 where $M_1\in[X_1,A_1]$;
 \item[(vi)] if $L \in (L_{1}, L_{X_2}]$  
 the Riemann solution is
 $L\testright{R_s} M_1 \testright{R_f} A_1 \xrightarrow{'S_f} R$,
 where $M_1\in(A_1, X_2]$;
\item[(vii)] if $L \in(L_{X_2},W)$, the Riemann solution is 
$
   L\testright{R_s} T_1 \xrightarrow{'S_s} M_1 \testright{R_f} A_1 \xrightarrow{'S_f} R,
$
where $T_1\in(X_2,W)_{\text{ext}}$ and $M_1\in(X_2, W)$; 
\item[(viii)]if $L=W$, the Riemann solution is 
$ W \testright{R_f} A_1 \xrightarrow{'S_f} R$.
\end{itemize}
\end{cla}

\noindent{ \bf Riemann solution for $R$ in subregion $\Gamma_2^d$.}
\medskip

\begin{cla}\label{cla:RSolution-Gamma_2d-new}
Refer to Fig.~\ref{fig:Solution_Riemann_P_Gamma_2}(d).
Let $L$ be a state in the edge \GW of the saturation triangle and $R$ be a state in subregion $\Gamma_2^d$ in Figs.~\ref{fig:RRegions Completed} and \ref{fig:RRegionsGamma}(a).
Let $A_1$ be the state in $\wm_f(R)$ satisfying
$\sigma(A_1; R) = \lambdaf(A_1)$ and
$X_1$, $X_2$ be the intersection states of $\wm_f(R)$
with the $s$-inflection locus.
Let $L_1$, $L_{X_1}$, $L_{X_2}$ and $L_R$ be the intersection points of the backward slow wave curves of $A_1$, $X_1$, $X_2$ and $R$ with the edge \GWc.
Then,
\begin{itemize}
\item[(i)]if $L=G$, the Riemann solution is
$G\testright{R_f} R$;
\item[(ii)] if $L \in(G, L_R)$, the Riemann solution is
$
   L\testright{R_s} T_1 \xrightarrow{'S_s} M_1  \testright{R_f} R,
$
where $T_1\in(G, T_R)_{\text{ext}}$ and $M_1\in(G, R)$;
\item[(iii)]if $L=L_R$, the Riemann solution is 
$ L\testright{R_s} T_R \xrightarrow{'S_s} R$;
\item[(iv)] if $L \in (L_R, L_1]$ the Riemann solution is
 $L\testright{R_s} T_1 \xrightarrow{'S_s} M_1  \testright{S_f} R$,
 where $T_1\in(T_R,X_1)_{\text{ext}}$ and $M_1\in(R,A_1]$;
\item[(v)] if $L \in (L_{1},L_{X_1})$ or $L \in(L_{X_2},W)$,
the solution is
 $L\testright{R_s}  T_1 \xrightarrow{'S_s} M_1  \testright{R_f} A_1 \xrightarrow{'S_f} R$ or
 $L\testright{R_s}  T_2 \xrightarrow{'S_s} M_2  \testright{R_f} A_1 \xrightarrow{'S_f} R$,
 where $T_1\in(T_R,X_1)_{\text{ext}}$ and $M_1\in(A_1,X_1)$ or
 $T_2\in(X_2,W)$ and $M_2\in(X_2, W)$;
 \item[(vi)] if $L \in [L_{X_1}, L_{X_2}]$  
 the Riemann solution is
 $L\testright{R_s} M_1 \testright{R_f} A_1 \xrightarrow{'S_f} R$,
 where $M_1\in[X_1, X_2]$; 
\item[(vii)]if $L=W$, the Riemann solution is 
$ W \testright{R_f} A_1 \xrightarrow{'S_f} R$.
\end{itemize}
\end{cla}

\noindent{ \bf Riemann solution for $R$ in subregion $\Gamma_2^e$.}
\medskip

\begin{cla}\label{cla:RSolution-Gamma_2e-new}
Refer to Fig.~\ref{fig:Solution_Riemann_P_Gamma_2}(e).
Let $L$ be a state in the edge \GW of the saturation triangle and $R$ be a state in subregion $\Gamma_2^e$ in Figs.~\ref{fig:RRegions Completed} and \ref{fig:RRegionsGamma}(a).
Let $A_1$ be the state in $\wm_f(R)$ satisfying
$\sigma(A_1; R) = \lambdaf(A_1)$ and
$X_1$, $X_2$ be the intersection states of $\wm_f(R)$
with the $s$-inflection locus.
Let $L_1$, $L_{X_1}$, $L_{X_2}$ and $L_R$ be the intersection points of the backward slow wave curves of $A_1$, $X_1$, $X_2$ and $R$ with the edge \GWc.
Then,
\begin{itemize}
\item[(i)] if $L=G$, the Riemann solution is
$G\testright{R_f} R$; 
\item[(ii)] if $L \in(G, L_{X_1})$ the Riemann solution is
$
   L\testright{R_s} T_1 \xrightarrow{'S_s} M_1  \testright{R_f} R,
$
where $T_1\in(G, X_1)_{\text{ext}}$ and $M_1\in(G, X_1)$; 
\item[(iii)] if $L \in [L_{X_1}, L_{R})$ the Riemann solution is
 $L\testright{R_s} M_1  \testright{R_f} R$,
 where $M_1\in[X_1,R)$;
 \item[(iv)]if $L=L_R$, the Riemann solution is 
$ L\testright{R_s} R$;
 \item[(v)] if $L \in  (L_R,L_{1}]$  
 the Riemann solution is
 $L\testright{R_s} M_1 \testright{S_f} R$,
 where $M_1\in(R,A_1]$;
\item[(vi)] if $L \in(L_{1}, L_{X_2}]$, the Riemann solution is 
$
L\testright{R_s} M_1  \testright{R_f} A_1 \xrightarrow{'S_f} R,
$
where $M_1\in(A_1,X_2]$; 
\item[(vii)] if $L \in(L_{X_2},W)$, the Riemann solution is 
$
   L\testright{R_s} T_1 \xrightarrow{'S_s} M_1 \testright{R_f} A_1 \xrightarrow{'S_f}  R,
$
where $T_1\in(X_2,W)_{\text{ext}}$ and $M_1\in(X_2, W)$; 
\item[(viii)]if $L=W$, the Riemann solution is 
$ W \testright{R_f} A_1 \xrightarrow{'S_f} R$.
\end{itemize}
\end{cla}

\section{Simulations}
\label{sec:Simulations}
This section presents the solution to the Riemann problem for some cases, as presented in previous sections. We compare the analytical and numerical solutions for each case at $t_D = 1$. Direct numerical simulations of System \eqref{eq:system1}-\eqref{eq:flowfunct} were performed using an implicit finite-difference scheme (FDS), which provides second-order accuracy in both space and time. This was combined with Newton's method; further details on the implementation can be found in \cite{lambert2020mathematics}.

\noindent{ \bf Case 1: $L = (0.25, 0.0001)$ and $R = (0.106, 0.888)$.}
\medskip

\begin{figure}[ht]
	\centering
	\subfigure[Riemann problem solution for $ L =(0.25 , 0.0001) $ and $ R = (0.106,0.888) $.]
	{\includegraphics[width=0.45\linewidth]{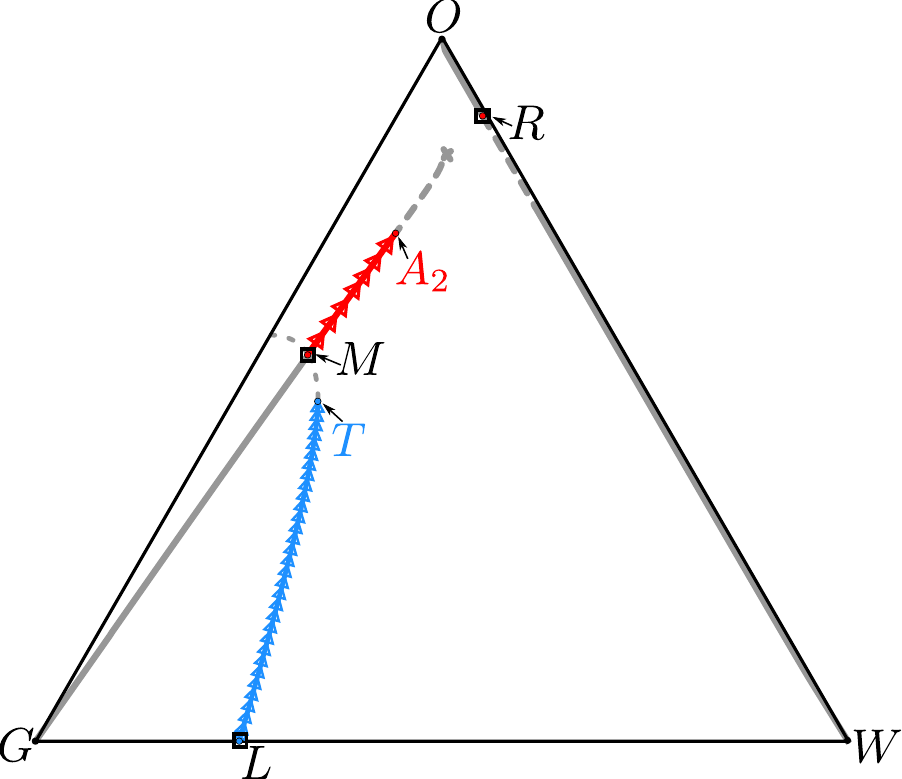}} 
	\hspace{1.5mm}	
 \subfigure[Saturation profiles for water, oil, and gas at $ t_D=1 $.]
	{\includegraphics[width=0.5\linewidth]{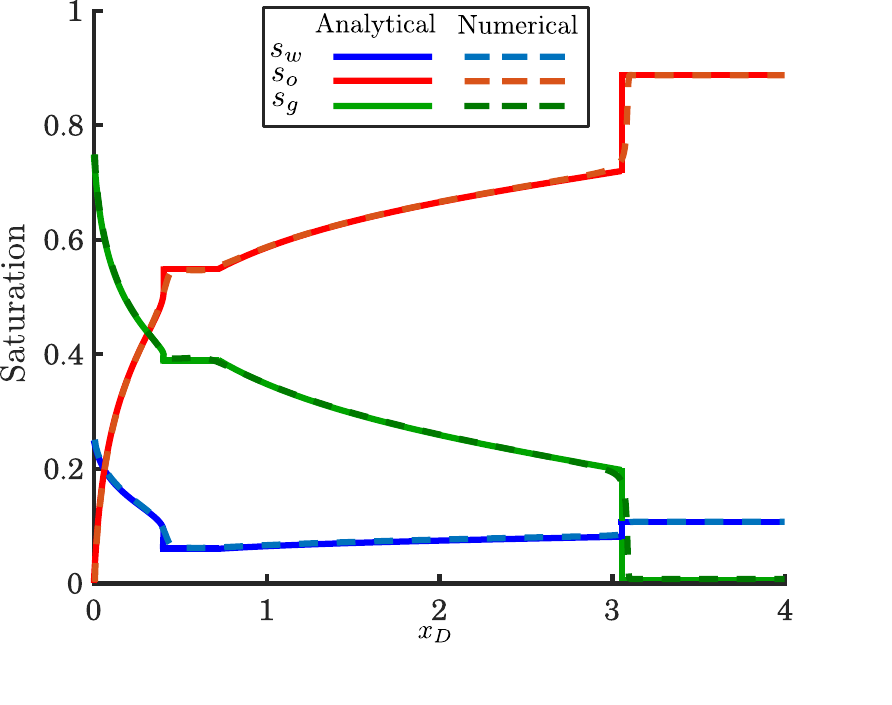}} 
 	\caption{Analytical solution for $R\in \Omega_1^a$ and $L\in{(G, L_2)}$. (a) The composition path in the saturation triangle that consists of $L\testright{R_s} T \xrightarrow{'S_s} M  \testright{R_f}A_2 \xrightarrow{'S_f} R$. The gray curve segments are the parts of the slow and fast wave curves that are not involved in the solution. (b) Analytical profiles (solid curves) compared to numerical simulations (dashed curves).
	}
	\label{fig:Case1}
\end{figure}

Figure~\ref{fig:Case1}(a) illustrates the analytical solution for the Riemann problem with initial states $L \in (G, L_2)$ and $R \in \Omega_1^a$. The composition path (Fig.~\ref{fig:Case1}(a)) consists of an $s$-rarefaction from $L$ to $T$, followed by an $s$-shock from $T$ to $M$, followed  $f$-rarefaction from $M$ to $A_2$, and finally an $f$-shock connecting $A_2$ to $R$. The gray curves represent inactive segments of the slow and fast wave families. Figure~~\ref{fig:Case1}(b) compares the analytical saturation profiles (solid lines) with numerical results (dashed lines) at $t_D = 1$, showing excellent agreement for all phases.

\noindent{ \bf Case 2: $L = (0.635, 0.0001)$ and $R = (0.177, 0.816)$.}
\medskip

\begin{figure}[ht]
	\centering
    \subfigure[Riemann problem solution for $ L =(0.635 , 0.0001) $ and $ R = (0.177,0.816) $.]
	{\includegraphics[width=0.45\linewidth]{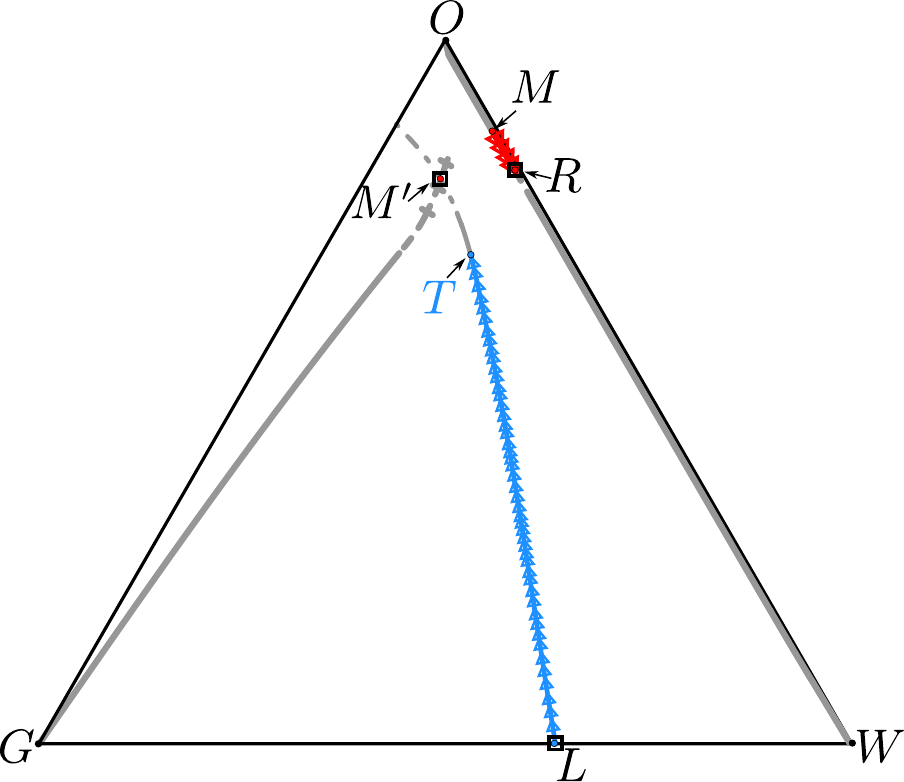}} 
	\hspace{1.5mm}	
 \subfigure[Saturation profiles for water, oil, and gas at $ t_D=1 $.]
	{\includegraphics[width=0.5\linewidth]{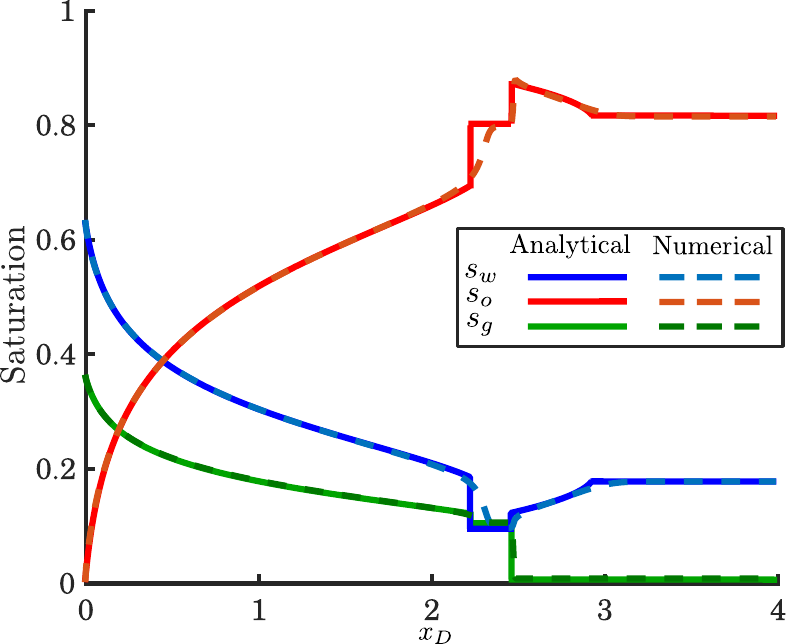}} 
 	\caption{Analytical solution for $R\in \Omega_1^b$ and $L\in{(L^*, L_3)}$. (a) The composition path in the saturation triangle that consists of $L\testright{R_s} T \xrightarrow{'S_s} M'  \testright{{S_f} '}M \xrightarrow{R_f} R$. The gray curve segments are the parts of the slow and fast wave curves that are not involved in the solution. (b) Analytical profiles (solid curves) compared to numerical simulations (dashed curves).
	}
	\label{fig:Case2}
\end{figure}

Figure~\ref{fig:Case2}(a) illustrates the analytical solution for the Riemann problem with initial states $L \in (L^*, L_3)$ and $R \in \Omega_1^b$. The composition path (Fig.~\ref{fig:Case2}(a)) consists of an $s$-rarefaction from $L$ to $T$, followed by an $s$-shock from $T$ to $M'$, followed  $f$-shock from $M'$ to $M$, and finally an $f$-rarefaction connecting $M$ to $R$. The gray curves represent inactive segments of the slow and fast wave families. Figure~~\ref{fig:Case2}(b) compares the analytical saturation profiles (solid lines) with numerical results (dashed lines) at $t_D = 1$.The numerical simulation match the analytical predictions closely, with minor discrepancies near of the $s$-composite front.

\noindent{ \bf Case 3: $L = (0.682, 0.0001)$ and $R = (0.271, 0.711)$.}
\medskip

\begin{figure}[ht]
	\centering
	  \subfigure[Riemann problem solution for $ L =(0.682 , 0.0001) $ and $ R = (0.271,0.711) $.]
	{\includegraphics[width=0.45\linewidth]{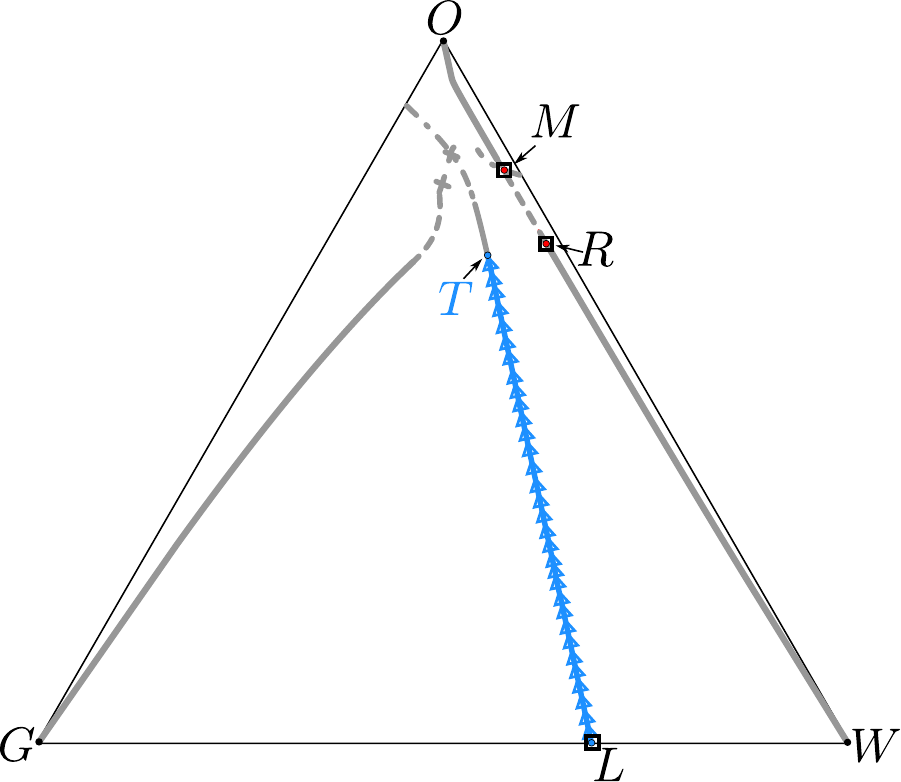}} 
	\hspace{1.5mm}	
 \subfigure[Saturation profiles for water, oil, and gas at $ t_D=1 $.]
	{\includegraphics[width=0.5\linewidth]{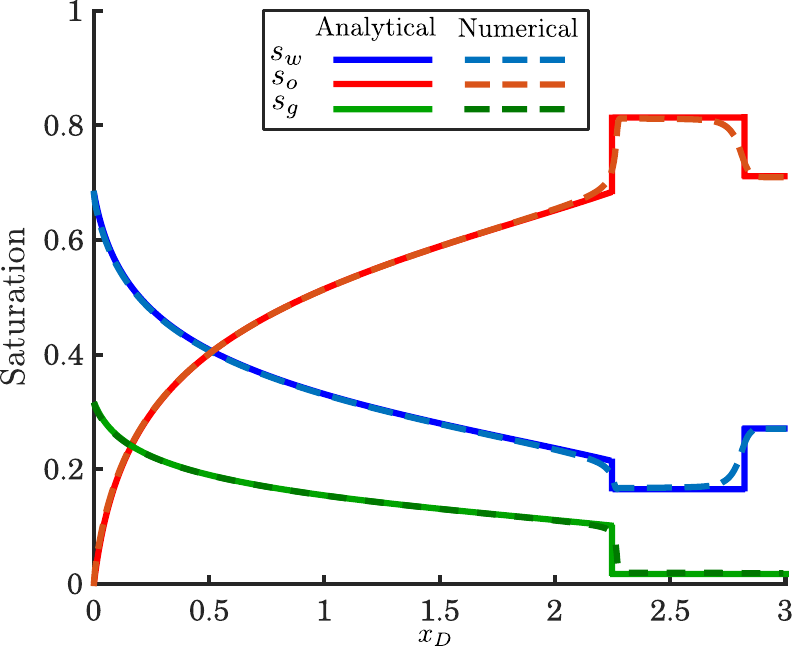}} 
 	\caption{Analytical solution for $R\in \Omega_2^a$ and $L\in{[L_1, L_R)}$. (a) The composition path in the saturation triangle, consisting of an $s$-rarefaction from $L$ to $T$, follows an $s$-shock from $T$ to $M$, and an $f$-shock from $M$ to $R$. The gray curve segments are the parts of the slow and fast wave curves that are not involved in the solution. (b) Analytical profiles (solid curves) compared to numerical simulations (dashed curves).
	}
	\label{fig:Case3}
\end{figure}

Figure~\ref{fig:Case3}(a) addresses the case $R \in \Omega_2^a$ and $L \in [L_1, L_R)$. The solution path (Fig.~\ref{fig:Case3}(a)) involves an $s$-rarefaction ($L \testright{}  T$), an $s$-shock ($T \rightarrow M$), and an $f$-shock ($M \testright{}  R$). The numerical results  (Fig.~\ref{fig:Case3}(b)) match the analytical predictions closely, with minor discrepancies near the $f$-shock front. 


\section{Conclusion}
\label{sec:Conclusion}

In this work, we present the Riemann solution for three-phase flow in porous media under the condition that oil viscosity exceeds that of water and gas. We classify all Riemann solution problems for scenarios where the left states $L$ lie along the edge $G$-$W$, and the right states $R$ span nearly the entire saturation triangle, excluding small regions near the boundaries $G$-$O$ and $W$-$O$. Using the wave curve method, we determine the Riemann solution for initial and injection data within this class. Our solutions were validated with numerical simulations. This study extends previous analytical solutions, which were limited to right states near the corner $O$ or within the quadrilateral $O$-$E$-$\mathcal{U}$-$D$. Notably, this classification remains valid for all viscosity variations satisfying the inequalities \eqref{eq:classical}, corresponding to viscosity regimes where the umbilic point is close to vertex $O$, see Fig.~\ref{fig:ratioviscosities}. 

We verify the $L^1_{loc}$-stability of the Riemann solution with respect to variations in the data. While we do not establish the uniqueness of the Riemann solution, extensive numerical experiments confirm its validity. Since the results of this work were produced for the Corey quadratic exponents model, this analysis serves as a foundation for extending such classifications to more general heavy oil Corey permeability exponents. These extensions have applications in calibrating numerical simulators, uncertainty quantification, and analyzing the effects of physical parameter changes. Our findings provide a comprehensive framework for understanding three-phase flow dynamics in porous media under a wide range of conditions.


\section*{Acknowledgement}

The authors thank Prof. Bradley Plohr for assistance with ELI interactive solver \cite{ELI_web}.

L.L. gratefully acknowledge support from Shell Brasil through the projects ``Avançando na modelagem matemática e computacional para apoiar a implementação da tecnologia `Foam-assisted WAG' em reservatórios do Pré-sal'' (ANP 23518-4) at UFJF, and the strategic importance of the support given by ANP through the R\&D levy regulation.

D.M. was partly supported by CAPES grant 88881.156518/2017-01, by CNPq under grants 405366/2021-3, 306566/2019-2, and by FAPERJ under grants E-26/210.738/2014, E-26/202.764/2017, E-26/201.159/2021.

\bibliographystyle{plainnat}
\bibliography{porousmedianew}

\begin{thebibliography}{32}
\providecommand{\natexlab}[1]{#1}
\providecommand{\url}[1]{\texttt{#1}}
\expandafter\ifx\csname urlstyle\endcsname\relax
  \providecommand{\doi}[1]{doi: #1}\else
  \providecommand{\doi}{doi: \begingroup \urlstyle{rm}\Url}\fi

\bibitem[Andrade et~al.(2016)Andrade, de~Souza, Furtado, and Marchesin]{L.2016}
P.~Andrade, A.~de~Souza, F.~Furtado, and D.~Marchesin.
\newblock {Oil displacement by water and gas in a porous medium: the Riemann
  problem}.
\newblock \emph{Bull. Braz. Math. Soc. (N.S.)}, 47\penalty0 (1):\penalty0
  1--14, 2016.

\bibitem[Andrade et~al.(2018)Andrade, de~Souza, Furtado, and
  Marchesin]{Andrade2018}
P.~Andrade, A.~de~Souza, F.~Furtado, and D.~Marchesin.
\newblock {Three-phase fluid displacement in a porous medium}.
\newblock \emph{Journal of Hyperbolic Differential Equations}, 15\penalty0
  (4):\penalty0 731--753, 2018.

\bibitem[Azevedo et~al.(2010)Azevedo, de~Souza, Furtado, Marchesin, and
  Plohr]{V.2010}
A.~Azevedo, A.~de~Souza, F.~Furtado, D.~Marchesin, and B.~Plohr.
\newblock {The solution by the wave curve method of three-phase flow in virgin
  reservoirs}.
\newblock \emph{Transport in Porous Media}, 83\penalty0 (1):\penalty0 99--125,
  2010.

\bibitem[Azevedo et~al.(2014)Azevedo, de~Souza, Furtado, and
  Marchesin]{Azevedo2014}
A.~Azevedo, A.~de~Souza, F.~Furtado, and D.~Marchesin.
\newblock {Uniqueness of the Riemann solution for three-phase flow in a porous
  medium}.
\newblock \emph{SIAM J. Appl. Math.}, 74\penalty0 (6):\penalty0 1967--1997,
  2014.

\bibitem[Azevedo et~al.(2002)Azevedo, Marchesin, Plohr, and
  Zumbrun]{azevedo2002capillary}
A.~V. Azevedo, D.~Marchesin, B.~Plohr, and K.~Zumbrun.
\newblock Capillary instability in models for three-phase flow.
\newblock \emph{Zeitschrift f{\"u}r angewandte Mathematik und Physik ZAMP},
  53:\penalty0 713--746, 2002.

\bibitem[Bagchi et~al.(2025)Bagchi, Patwardhan, Iglauer, Ben~Mahmud, and
  Ali]{bagchi2025critical}
Chiradip Bagchi, Samarth~D Patwardhan, Stefan Iglauer, Hisham Ben~Mahmud, and
  Muhammad Fazil~Jaffar Ali.
\newblock A critical review on parameters affecting the feasibility of
  underground hydrogen storage.
\newblock \emph{ACS Omega}, 2025.

\bibitem[Barros et~al.(2021)Barros, Pires, and Peres]{barros2021analytical}
Wagner~Q Barros, Adolfo~P Pires, and Alvaro~MM Peres.
\newblock Analytical solution for one-dimensional three-phase incompressible
  flow in porous media for concave relative permeability curves.
\newblock \emph{International Journal of Non-Linear Mechanics}, 137:\penalty0
  103792, 2021.

\bibitem[Davies et~al.(2024)Davies, Ehrmann, and
  Schwenzfeier-Hellkamp]{davies2024safety}
Emma Davies, Andrea Ehrmann, and Eva Schwenzfeier-Hellkamp.
\newblock Safety of hydrogen storage technologies.
\newblock \emph{Processes}, 12\penalty0 (10):\penalty0 2182, 2024.

\bibitem[de~Souza(1992)]{DeSouza1992}
A.~de~Souza.
\newblock {Stability of singular fundamental solutions under perturbations for
  flow in porous media}.
\newblock \emph{Mat. Apl. Comput.}, 11\penalty0 (2):\penalty0 43, 1992.

\bibitem[ELI(2025)]{ELI_web}
ELI.
\newblock {ELI}, {I}nteractive {G}raphical {R}iemann {P}roblem {S}olver.
\newblock \emph{https://eli.fluid.impa.br/}, pages {A}ccessed : 2025--06--08,
  2025.

\bibitem[F.(2018)]{Lozano}
Lozano~L. F.
\newblock \emph{{Diffusive effects in Riemann solutions for the three phase
  flow in porous media, D.Sc. thesis }}.
\newblock PhD thesis, {Instituto de Matem{\'{a}}tica Pura e Aplicada (IMPA),
  Rio de Janeiro, Brazil}, 2018.

\bibitem[Isaacson et~al.(1990)Isaacson, Marchesin, and Plohr]{L.1990}
E.~Isaacson, D.~Marchesin, and B.~Plohr.
\newblock {Transitional waves for conservation laws}.
\newblock \emph{SIAM J. Math. Anal.}, 21\penalty0 (4):\penalty0 837--866, 1990.

\bibitem[Isaacson et~al.(1992)Isaacson, Marchesin, Plohr, and Temple]{L.1992a}
E.~Isaacson, D.~Marchesin, B.~Plohr, and B.~Temple.
\newblock {Multiphase flow models with singular Riemann problems}.
\newblock \emph{Comput. Appl. Math.}, 11\penalty0 (2):\penalty0 147--166, 1992.

\bibitem[Iskandarov et~al.(2022)Iskandarov, Fanourgakis, Ahmed, Alameri,
  Froudakis, and Karanikolos]{iskandarov2022data}
Javad Iskandarov, George~S Fanourgakis, Shehzad Ahmed, Waleed Alameri, George~E
  Froudakis, and Georgios~N Karanikolos.
\newblock Data-driven prediction of in situ co2 foam strength for enhanced oil
  recovery and carbon sequestration.
\newblock \emph{RSC advances}, 12\penalty0 (55):\penalty0 35703--35711, 2022.

\bibitem[Juanes and Patzek(2004)]{juanes2004analytical}
Ruben Juanes and Tadeusz~W Patzek.
\newblock Analytical solution to the riemann problem of three-phase flow in
  porous media.
\newblock \emph{Transport in Porous Media}, 55:\penalty0 47--70, 2004.

\bibitem[Lambert et~al.(2020)Lambert, Alvarez, Ledoino, Tadeu, Marchesin, and
  Bruining]{lambert2020mathematics}
Wanderson Lambert, Amaury Alvarez, Ismael Ledoino, Duilio Tadeu, Dan Marchesin,
  and Johannes Bruining.
\newblock Mathematics and numerics for balance partial differential-algebraic
  equations (pdaes).
\newblock \emph{J. Sci. Comput.}, 84\penalty0 (2):\penalty0 29, 2020.
\newblock \doi{10.1007/s10915-020-01279-w}.

\bibitem[Lax(1952)]{Lax1957}
P.~Lax.
\newblock {Hyperbolic systems of conservation laws II}.
\newblock \emph{Comm. Pure Appl. Math.}, X, 1952.

\bibitem[Liu(1974)]{Liu1974}
T.~P. Liu.
\newblock {The Riemann problem for general 2x2 conservations laws}.
\newblock \emph{Trans. Amer. Math. Soc.}, 199:\penalty0 89--112, 1974.

\bibitem[Liu(1975)]{Liu1975}
T.~P. Liu.
\newblock {The Riemann problem for general systems of conservations laws}.
\newblock \emph{J. Differential Equations}, 18:\penalty0 218--234, 1975.

\bibitem[Lozano et~al.(2024{\natexlab{a}})Lozano, Chapiro, and
  Marchesin]{lozano2024a}
L~Lozano, G~Chapiro, and D~Marchesin.
\newblock Analytical investigation of the three-phase foam flow in porous
  media.
\newblock In \emph{ECMOR 2024}, volume 2024, pages 1--9. European Association
  of Geoscientists \& Engineers, 2024{\natexlab{a}}.

\bibitem[Lozano et~al.(2024{\natexlab{b}})Lozano, Ledoino, Plohr, and
  Marchesin]{lozano2024c}
L.~Lozano, I.~Ledoino, B.~J~. Plohr, and D.~Marchesin.
\newblock Structure of undercompressive shock waves in three-phase flow in
  porous media, 2024{\natexlab{b}}.
\newblock URL \url{https://arxiv.org/abs/2412.04439}.

\bibitem[Matos et~al.(2015)Matos, Azevedo, Mota, and Marchesin]{V.2015}
V.~Matos, A.~Azevedo, J.~Mota, and D.~Marchesin.
\newblock {Bifurcation under parameter change of Riemann solutions for
  nonstrictly hyperbolic systems}.
\newblock \emph{Z. Angew. Math. Phys.}, 66\penalty0 (4):\penalty0 1413--1452,
  2015.

\bibitem[Matos et~al.(2016)Matos, Silva, and Marchesin]{Matos2016}
V.~Matos, J.~D. Silva, and D.~Marchesin.
\newblock {Loss of hyperbolicity changes the number of wave groups in Riemann
  problems}.
\newblock \emph{Bulletin of the Brazilian Mathematical Society}, 47\penalty0
  (2):\penalty0 545--559, 2016.

\bibitem[Mehrabi et~al.(2020)Mehrabi, Sepehrnoori, and
  Delshad]{mehrabi2020solution}
Mehran Mehrabi, Kamy Sepehrnoori, and Mojdeh Delshad.
\newblock Solution construction to a class of {R}iemann problems of multiphase
  flow in porous media.
\newblock \emph{Transport in porous media}, 132:\penalty0 241--266, 2020.

\bibitem[Miao et~al.(2024)Miao, Zhao, Huang, Zhao, Zhao, Guo, and
  Wang]{miao2024core}
Yan Miao, Qiuyang Zhao, Zujie Huang, Keyu Zhao, Hao Zhao, Liejin Guo, and
  Yechun Wang.
\newblock Core flooding experimental study on enhanced oil recovery of heavy
  oil reservoirs with high water cut by sub-and supercritical water.
\newblock \emph{Geoenergy Science and Engineering}, 242:\penalty0 213208, 2024.

\bibitem[Pires et~al.(2024)Pires, Barros, and Peres]{pires2024approximate}
Adolfo~P Pires, Wagner~Q Barros, and Alvaro~MM Peres.
\newblock Approximate analytical solutions for 1-d immiscible water alternated
  gas.
\newblock \emph{Transport in Porous Media}, 151\penalty0 (1):\penalty0
  171--191, 2024.

\bibitem[Salimi et~al.(2012)Salimi, Wolf, and Bruining]{salimi2012influence}
Hamidreza Salimi, Karl-Heinz Wolf, and Johannes Bruining.
\newblock The influence of capillary pressure on the phase equilibrium of the
  co2--water system: Application to carbon sequestration combined with
  geothermal energy.
\newblock \emph{International Journal of Greenhouse Gas Control}, 11:\penalty0
  S47--S66, 2012.

\bibitem[Schaeffer and Shearer(1987)]{schaeffer1987classification}
David~G Schaeffer and Michael Shearer.
\newblock The classification of 2$\times$ 2 systems of non-strictly hyperbolic
  conservation laws, with application to oil recovery.
\newblock \emph{Communications on pure and applied mathematics}, 40\penalty0
  (2):\penalty0 141--178, 1987.

\bibitem[Silva and Marchesin(2014)]{silva2014riemann}
Julio~Daniel Silva and Dan Marchesin.
\newblock Riemann solutions without an intermediate constant state for a system
  of two conservation laws.
\newblock \emph{Journal of Differential Equations}, 256\penalty0 (4):\penalty0
  1295--1316, 2014.

\bibitem[Tang et~al.(2022)Tang, Castaneda, Marchesin, and Rossen]{tang2022foam}
Jinyu Tang, Pablo Castaneda, Dan Marchesin, and William~R Rossen.
\newblock Foam-oil displacements in porous media: Insights from three-phase
  fractional-flow theory.
\newblock In \emph{Abu Dhabi International Petroleum Exhibition and
  Conference}, page D042S195R003. SPE, 2022.

\bibitem[Wendroff(1972)]{Wendroff1972}
B.~Wendroff.
\newblock {The Riemann problem for materials with Non Convex Equations of
  state: I Isentropic flow; II General flow}.
\newblock \emph{J. Math. Anal Appl.}, 38\penalty0 (2--3):\penalty0 454 -- 466;
  640 -- 658, 1972.

\bibitem[Zhang et~al.(2006)Zhang, Sayegh, and Huang]{zhang2006enhanced}
YP~Zhang, S~Sayegh, and S~Huang.
\newblock Enhanced heavy oil recovery by immiscible {WAG} injection.
\newblock In \emph{PETSOC Canadian International Petroleum Conference}, pages
  PETSOC--2006. PETSOC, 2006.

\end{thebibliography}
\end{document}